# Linear Map of $D$-Algebra

Aleks Kleyn

Aleks_Kleyn@MailAPS.org
http://AleksKleyn.dyndns-home.com:4080/
http://sites.google.com/site/AleksKleyn/
http://arxiv.org/a/kleyn_a_1
http://AleksKleyn.blogspot.com/

Abstract. Module is effective representation of ring in Abelian group. Linear map of module over commutative ring is morphism of corresponding representation. This definition is the main subject of the book.

To consider this definition from more general point of view I started the book from consideration of Cartesian product of representations. Polymorphism of representations is a map of Cartesian product of representations which is a morphism of representations with respect to each separate independent variable. Reduced morphism of representations allows us to simplify the study of morphisms of representations. However a representation has to satisfy specific requirements for existence of reduced polymomorphism of representations. It is possible that Abelian group is only $\Omega$-algebra, such that representation in this algebra admits polymorphism of representations. However, today, this statement has not been proved.

Multiplicative $\Omega$-group is $\Omega$-algebra in which product is defined. The definition of tensor product of representations of Abelian multiplicative $\Omega$-group is based on properties of reduced polymorphism of representations of Abelian multiplicative $\Omega$-group.

Since algebra is a module in which the product is defined, then we can use this theory to study linear map of algebra. For instance, we can study the set of linear transformations of $D$-algebra $A$ as representation of algebra $A\otimes A$ in algebra $A$.

# Contents

CHAPTER 1

# Preface

## 1.1. Preface to Version 1

There exist few equivalent definitions of module. For me, the definition of module as effective representation of commutative ring in Abelian group is the most interesting definition. This definition allows us to consider linear algebra from more general point of view and to understand which construction can be properly defined in other algebraic theories. For instance, we may consider a linear map as morphism of representation, i.e. a map which preserves the structure of the representation.

The module in which the product is defined is called algebra. Depending on the problem to be solved, we consider different algebraic structures on algebra. Accordingly, the map, which preserves the structure of algebra, changes.

If we do not consider the representation of ring $D$ in $D$-algebra $A$, then the algebra $A$ is a ring. A map, preserving the structure of algebra as a ring, is called homomorphism of algebra. If we do not consider the product $D$-algebra $A$, then we consider $D$-algebra as module. A map, preserving the structure of algebra as module, is called a linear map of the $D$-algebra. Linear map is important tool to study calculus. A map, preserving the structure of representation, is called a linear homomorphism. In the section 5.3, I proved the existence of a nontrivial linear homomorphism.

February, 2015

## 1.2. Preface to Version 2

The main idea of this book is an attempt to apply to the representation of $\Omega_1$-algebra (which is called $\Omega_1$-algebra of transformations of $\Omega_2$-algebra or just $\Omega_1$-algebra of transformations) some of concepts with which we are familiar in linear algebra. Since linear map is reduced morphism of module when we consider module as representation of ring in Abelian group, then the definition of polymorphism of representations is natural generalization of polylinear map.

Studying the linear map over non commutative division ring $D$, we conclude that in the equation

$$f(ax) = af(x)$$

$a$ cannot be any element of the ring $D$, but should belong to the center of the ring $D$. This statement makes theory of linear maps interesting and rich theory.

So, like in case of polylinear map, we request that transformation which $\Omega_1$-algebra generates in one factor, can be transferred to any other factor. In general, this requirement is not restrictive; but compliance with this requirement allows consideration of $\Omega_1$-algebra of transformations as Abelian multiplicative $\Omega_1$-group.

Tensor product of vector spaces is associative. Consideration of associativity of tensor product of representation generates new constraint on polymorphism of representations. Any two operations of $\Omega_2$-algebra must satisfy the equation (3.1.19). It is possible that Abelian group is only $\Omega_2$-algebra, such that representation in this algebra admits polymorphism of representations. However, today, this statement has not been proved.

If we believe in a statement which we cannot prove or disprove, then discovery of opposite statement may lead us to interesting discoveries. For thousands of years mathematicians tried to prove Euclid's fifth postulate. Lobachevsky and Bolyai proved that there exists geometry which does not hold this postulate. Mathematicians tried to find solution of algebraic equation using radicals. However Galois proved this is imposible to do when the power of equation is greater than 4.

For long time, I believed that existence of reduced morphism of representation imply existence of reduced polymorphism of representation. The theorem 3.1.13 was a surprise to me.

I believe that text about polymorphism and tensor product of representations is important.

- After writing of few papers dedicated to calculus over Banach algebra, I believed that one day I will be able to write similar paper about calculus in a representation. I was confused by the need to refuse addition as estimate how small is distance between maps. Of course, I may state that morphism of representation is in neighborhood of a considered map. However addition is an essential component of the construction. Absence of reduced polymorphism deprive me of the opportunity to define derivatives of second order.
- We can consider linear and polylinear maps of effective representation in Abelian $\Omega_2$-group like we consider linear and polylinear maps of $D$-algebra. The structure of linear map of effective representation in Abelian $\Omega_2$-group is more complicated than the structure of linear map of $D$-algebra. However this problem is interesting for me and I hope to return to this problem in the future.

May, 2015

## 1.3. Conventions

Convention 1.3.1. *We will use Einstein summation convention in which repeated index (one above and one below) implies summation with respect to repeated index. In this case we assume that we know the set of summation index and do not use summation symbol*

$$c^i v_i = \sum_{i\in I} c^i v_i$$

□

Convention 1.3.2. *I assume sum over index $i$ in expression like*

$$a_{i\cdot 0} x a_{i\cdot 1}$$

□

Convention 1.3.3. *Let $A$ be free algebra with finite or countable basis. Considering expansion of element of algebra $A$ relative basis $\overline{\overline{e}}$ we use the same root* *letter to denote this element and its coordinates. In expression* $a^2$*, it is not clear whether this is component of expansion of element a relative basis, or this is operation* $a^2 = aa$*. To make text clearer we use separate color for index of element of algebra. For instance,*

$$a = a^{i} e_{i}$$

□

CONVENTION 1.3.4. *If free finite dimensional algebra has unit, then we identify the vector of basis* $e_{0}$ *with unit of algebra.* □

CHAPTER 2

# Product of Representations

## 2.1. Cartesian Product of Universal Algebras

Definition 2.1.1. *Let $\mathcal{A}$ be a category. Let $\{B_i, i \in I\}$ be the set of objects of $\mathcal{A}$. Object*

$$P = \prod_{i \in I} B_i$$

*and set of morphisms*

$$\{\ f_i : P \longrightarrow B_i\ , i \in I\}$$

*is called a* **product of objects** $\{B_i, i \in I\}$ **in category** $\mathcal{A}$[2.1] *if for any object $R$ and set of morphisms*

$$\{\ g_i : R \longrightarrow B_i\ , i \in I\}$$

*there exists a unique morphism*

$$h : R \longrightarrow P$$

*such that diagram*

$$\begin{array}{ccc} P & \xrightarrow{f_i} & B_i \\ {\scriptstyle h}\uparrow & \nearrow{\scriptstyle g_i} & \\ R & & \end{array} \qquad f_i \circ h = g_i$$

*is commutative for all $i \in I$.*

*If $|I| = n$, then we also will use notation*

$$P = \prod_{i=1}^{n} B_i = B_1 \times ... \times B_n$$

*for product of objects $\{B_i, i \in I\}$ in $\mathcal{A}$.* □

Example 2.1.2. *Let $\mathcal{S}$ be the category of sets.*[2.2] *According to the definition 2.1.1, Cartesian product*

$$A = \prod_{i \in I} A_i$$

*of family of sets $(A_i, i \in I)$ and family of projections on the $i$-th factor*

$$p_i : A \to A_i$$

*are product in the category $\mathcal{S}$.* □

[2.1] I made definition according to [1], page 58.
[2.2] See also the example in [1], page 59.

Theorem 2.1.3. *The product exists in the category $\mathcal{A}$ of $\Omega$-algebras. Let $\Omega$-algebra $A$ and family of morphisms*

$$p_i : A \to A_i \quad i \in I$$

*be product in the category $\mathcal{A}$. Then*

2.1.3.1: *The set $A$ is Cartesian product of family of sets $(A_i, i \in I)$*

2.1.3.2: *The homomorphism of $\Omega$-algebra*

$$p_i : A \to A_i$$

*is projection on $i$-th factor.*

2.1.3.3: *We can represent any $A$-number $a$ as tuple $(p_i(a), i \in I)$ of $A_i$-numbers.*

2.1.3.4: *Let $\omega \in \Omega$ be $n$-ary operation. Then operation $\omega$ is defined componentwise*

$$a_1 ... a_n \omega = (a_{1i} ... a_{ni} \omega, i \in I) \tag{2.1.1}$$

*where $a_1 = (a_{1i}, i \in I)$, ..., $a_n = (a_{ni}, i \in I)$ .*

Proof. Let

$$A = \prod_{i \in I} A_i$$

be Cartesian product of family of sets $(A_i, i \in I)$ and, for each $i \in I$, the map

$$p_i : A \to A_i$$

be projection on the $i$-th factor. Consider the diagram of morphisms in category of sets $\mathcal{S}$

(2.1.2)

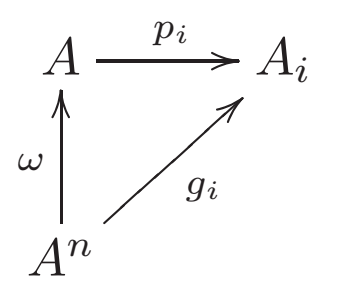

$$p_i \circ \omega = g_i$$

where the map $g_i$ is defined by the equation

$$g_i(a_1, ..., a_n) = p_i(a_1) ... p_i(a_n) \omega$$

According to the definition 2.1.1, the map $\omega$ is defined uniquely from the set of diagrams (2.1.2)

$$a_1 ... a_n \omega = (p_i(a_1) ... p_i(a_n) \omega, i \in I) \tag{2.1.3}$$

The equation (2.1.1) follows from the equation (2.1.3). □

Definition 2.1.4. *If $\Omega$-algebra $A$ and family of morphisms*

$$p_i : A \to A_i \quad i \in I$$

*is product in the category $\mathcal{A}$, then $\Omega$-algebra $A$ is called* **direct** *or* **Cartesian product of $\Omega$-algebras** $(A_i, i \in I)$ . □

Theorem 2.1.5. *Let set $A$ be Cartesian product of sets $(A_i, i \in I)$ and set $B$ be Cartesian product of sets $(B_i, i \in I)$ . For each $i \in I$, let*

$$f_i : A_i \to B_i$$

*be the map from the set $A_i$ into the set $B_i$. For each $i \in I$, consider commutative diagram*

(2.1.4)

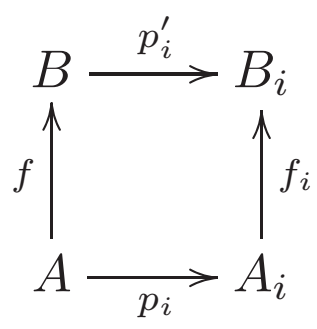

*where maps $p_i$, $p'_i$ are projection on the $i$-th factor. The set of commutative diagrams (2.1.4) uniquely defines map*

$$f : A_1 \to A_2$$

$$f(a_i, i \in I) = (f_i(a_i), i \in I)$$

Proof. For each $i \in I$, consider commutative diagram

(2.1.5)

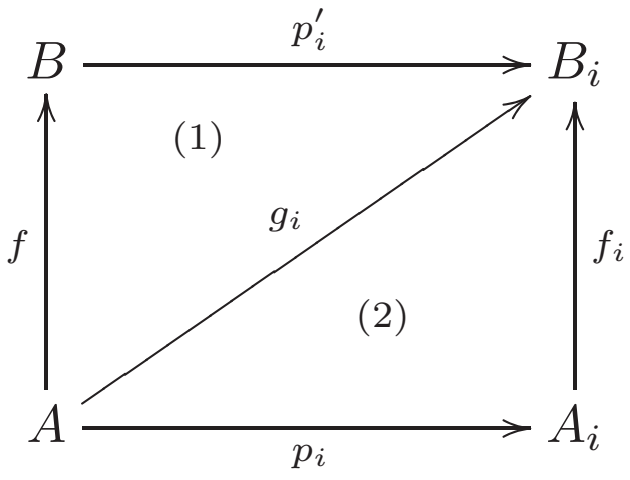

Let $a \in A$. According to the statement 2.1.3.3, we can represent $A$-number $a$ as tuple of $A_i$-numbers

$$a = (a_i, i \in I) \quad a_i = p_i(a) \in A_i \tag{2.1.6}$$

Let

$$b = f(a) \in B \tag{2.1.7}$$

According to the statement 2.1.3.3, we can represent $B$-number $b$ as tuple of $B_i$-numbers

$$b = (b_i, i \in I) \quad b_i = p'_i(b) \in B_i \tag{2.1.8}$$

From commutativity of diagram (1) and from equations (2.1.7), (2.1.8), it follows that

$$b_i = g_i(b) \tag{2.1.9}$$

From commutativity of diagram (2) and from the equation (2.1.6), it follows that

$$b_i = f_i(a_i)$$

□

Theorem 2.1.6. *Let $\Omega$-algebra $A$ be Cartesian product of $\Omega$-algebras $(A_i, i \in I)$ and $\Omega$-algebra $B$ be Cartesian product of $\Omega$-algebras $(B_i, i \in I)$ . For each $i \in I$, let the map*

$$f_i : A_i \to B_i$$

*be homomorphism of $\Omega$-algebra. Then the map*

$$f : A \to B$$

*defined by the equation*

$$f(a_i, i \in I) = (f_i(a_i), i \in I) \tag{2.1.10}$$

*is homomorphism of $\Omega$-algebra.*

Proof. Let $\omega \in \Omega$ be n-ary operation. Let $a_1 = (a_{1i}, i \in I)$, ..., $a_n = (a_{ni}, i \in I)$ , $b_1 = (b_{1i}, i \in I)$, ..., $b_n = (b_{ni}, i \in I)$ . From equations (2.1.1), (2.1.10), it follows that

$$\begin{aligned} f(a_1...a_n\omega) &= f(a_{1i}...a_{ni}\omega, i \in I) \\ &= (f_i(a_{1i}...a_{ni}\omega), i \in I) \\ &= ((f_i(a_{1i}))...(f_i(a_{ni})), i \in I) \\ &= (b_{1i}...b_{ni}\omega, i \in I) \end{aligned}$$

$$f(a_1)...f(a_n)\omega = b_1...b_n\omega = (b_{1i}...b_{ni}\omega, i \in I)$$

□

## 2.2. Cartesian Product of Representations

Lemma 2.2.1. *Let*

$$A = \prod_{i \in I} A_i$$

*be Cartesian product of family of $\Omega_2$-algebras $(A_i, i \in I)$ . For each $i \in I$, let the set ${}^*A_i$ be $\Omega_2$-algebra. Then the set*

$$^\circ A = \{f \in {}^*A : f(a_i, i \in I) = (f_i(a_i), i \in I)\} \tag{2.2.1}$$

*is Cartesian product of $\Omega_2$-algebras ${}^*A_i$.*

Proof. According to the definition (2.2.1), we can represent a map $f \in {}^\circ A$ as tuple

$$f = (f_i, i \in I)$$

of maps $f_i \in {}^*A_i$. According to the definition (2.2.1),

$$(f_i, i \in I)(a_i, i \in I) = (f_i(a_i), i \in I)$$

Let $\omega \in \Omega_2$ be n-ary operation. We define operation $\omega$ on the set ${}^\circ A$ using equation

$$((f_{1i}, i \in I)...(f_{ni}, i \in I)\omega)(a_i, i \in I) = ((f_{1i}(a_i))...(f_{ni}(a_i))\omega, i \in I)$$

□

Definition 2.2.2. *Let $\mathcal{A}_1$ be category of $\Omega_1$-algebras. Let $\mathcal{A}_2$ be category of $\Omega_2$-algebrasz. We define* **category $(\mathcal{A}_1*)\mathcal{A}_2$ of left-side representations**. *Left-side representations of $\Omega_1$-algebra in $\Omega_2$-algebra are objects of this category. Morphisms of corresponding representations are morphisms of this category.* □

Theorem 2.2.3. *In category $(\mathcal{A}_1*)\mathcal{A}_2$ there exists product of single transitive left-side representations of $\Omega_1$-algebra in $\Omega_2$-algebra.*

Proof. For $j=1,2$, let

$$P_j=\prod_{i\in I}B_{ji}$$

be product of family of $\Omega_j$-algebras $\{B_{ji},i\in I\}$ and for any $i\in I$ the map

$$t_{ji}:P_j\longrightarrow B_{ji}$$

be projection onto factor $i$. For each $i\in I$, let

$$h_i:B_{1i}\overset{*}{\longrightarrow}B_{2i}$$

be single transitive $B_{1i}*$-representation in $\Omega_2$-algebra $B_{2i}$.

Let $b_1\in P_1$. According to the statement 2.1.3.3, we can represent $P_1$-number $b_1$ as tuple of $B_{1i}$-numbers

$$b_1=(b_{1i},i\in I)\quad b_{1i}=t_{1i}(b_1)\in B_{1i} \tag{2.2.2}$$

Let $b_2\in P_2$. According to the statement 2.1.3.3, we can represent $P_2$-number $b_2$ as tuple of $B_{2i}$-numbers

$$b_2=(b_{2i},i\in I)\quad b_{2i}=t_{2i}(b_2)\in B_{2i} \tag{2.2.3}$$

Lemma 2.2.4. *For each $i\in I$, consider diagram of maps*

(2.2.4)

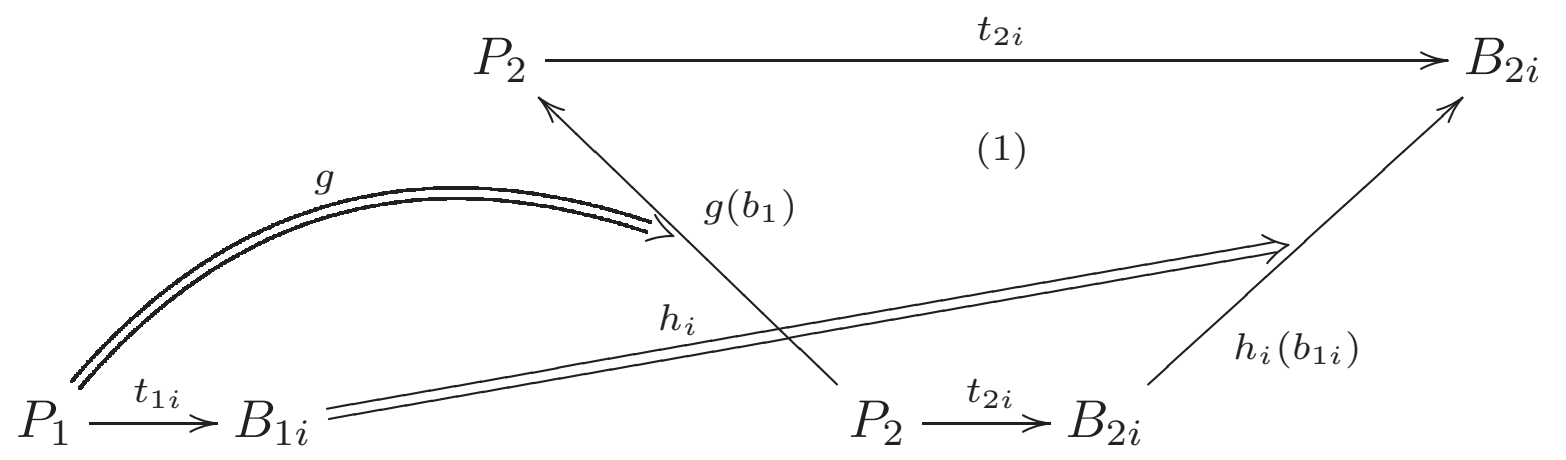

*Let map*

$$g:P_1\rightarrow {}^*P_2$$

*be defined by the equation*

$$g(b_1)\circ b_2=(h_i(b_{1i})\circ b_{2i},i\in I) \tag{2.2.5}$$

*Then the map $g$ is single transitive $P_1*$-representation in $\Omega_2$-algebra $P_2$*

$$g:P_1\overset{*}{\longrightarrow}P_2$$

*The map $(t_{1i},t_{2i})$ is morphism of representation $g$ into representation $h_i$.*

Proof.

2.2.4.1: According to definitions [4]-2.1.1, [4]-2.1.2, the map $h_i(b_{1i})$ is homomorphism of $\Omega_2$-algebra $B_{2i}$. According to the theorem 2.1.6, from commutativity of the diagram (1) for each $i\in I$, it follows that the map

$$g(b_1):P_2\rightarrow P_2$$

defined by the equation (2.2.5) is homomorphism of $\Omega_2$-algebra $P_2$.

2.2.4.2: According to the definition [4]-2.1.2, the set ${}^*B_{2i}$ is $\Omega_1$-algebra. According to the lemma 2.2.1, the set ${}^\circ P_2\subseteq {}^*P_2$ is $\Omega_1$-algebra.

2.2.4.3: According to the definition [4]-2.1.2, the map

$$h_i : B_{1i} \to {}^*B_{2i}$$

is homomorphism of $\Omega_1$-algebra. According to the theorem 2.1.6, the map

$$g : P_1 \to {}^*P_2$$

defined by the equation

$$g(b_1) = (h_i(b_{1i}), i \in I)$$

is homomorphism of $\Omega_1$-algebra.

According to statements 2.2.4.1, 2.2.4.3 and to the definition [4]-2.1.2, the map $g$ is $P_1*$-representation in $\Omega_2$-algebra $P_2$.

Let $b_{21}, b_{22} \in P_2$. According to the statement 2.1.3.3, we can represent $P_2$-numbers $b_{21}, b_{22}$ as tuples of $B_{2i}$-numbers

$$\begin{aligned} b_{21} &= (b_{21i}, i \in I) \quad b_{21i} = t_{2i}(b_{21}) \in B_{2i} \\ b_{22} &= (b_{22i}, i \in I) \quad b_{22i} = t_{2i}(b_{22}) \in B_{2i} \end{aligned} \tag{2.2.6}$$

According to the theorem [4]-2.1.9, since the representation $h_i$ is single transitive, then there exists unique $B_{1i}$-number $b_{1i}$ such that

$$b_{22i} = h_i(b_{1i}) \circ b_{21i}$$

According to definitions (2.2.2), (2.2.5), (2.2.6), there exists unique $P_1$-number $b_1$ such that

$$b_{22} = g(b_1) \circ b_{21}$$

According to the theorems [4]-2.1.9, the representation $g$ is single transitive.

From commutativity of diagram (1) and from the definition [4]-2.2.2, it follows that map $(t_{1i}, t_{2i})$ is morphism of representation $g$ into representation $h_i$. ⊙

Let

$$d_2 = g(b_1) \circ b_2 \quad d_2 = (d_{2i}, i \in I) \tag{2.2.7}$$

From equations (2.2.5), (2.2.7), it follows that

$$d_{2i} = h_i(b_{1i}) \circ b_{2i} \tag{2.2.8}$$

For $j = 1, 2$, let $R_j$ be other object of category $\mathcal{A}_j$. For any $i \in I$, let the map

$$r_{1i} : R_1 \longrightarrow B_{1i}$$

be morphism from $\Omega_1$-algebra $R_1$ into $\Omega_1$-algebra $B_{1i}$. According to the definition 2.1.1, there exists a unique morphism of $\Omega_1$-algebra

$$s_1 : R_1 \longrightarrow P_1$$

such that following diagram is commutative

(2.2.9)

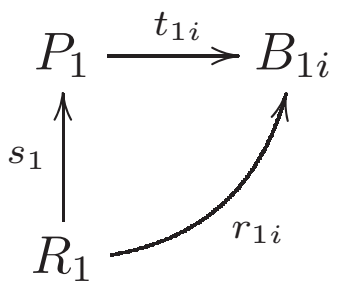

$$t_{1i} \circ s_1 = r_{1i}$$

Let $a_1 \in R_1$. Let

$$b_1 = s_1(a_1) \in P_1 \tag{2.2.10}$$

From commutativity of the diagram (2.2.9) and statements (2.2.10), (2.2.2), it follows that

(2.2.11) $$b_{1i}=r_{1i}(a_1)$$

Let

$$f:R_1\overset{*}{\longrightarrow}R_2$$

be single transitive $R_1*$-representation in $\Omega_2$-algebra $R_2$. According to the theorem [4]-2.2.10, a morphism of $\Omega_2$-algebra

$$r_{2i}:R_2\longrightarrow B_{2i}$$

such that map $(r_{1i},r_{2i})$ is morphism of representations from $f$ into $h_i$ is unique up to choice of image of $R_2$-number $a_2$. According to the remark [4]-2.2.6, in diagram of maps

(2.2.12)

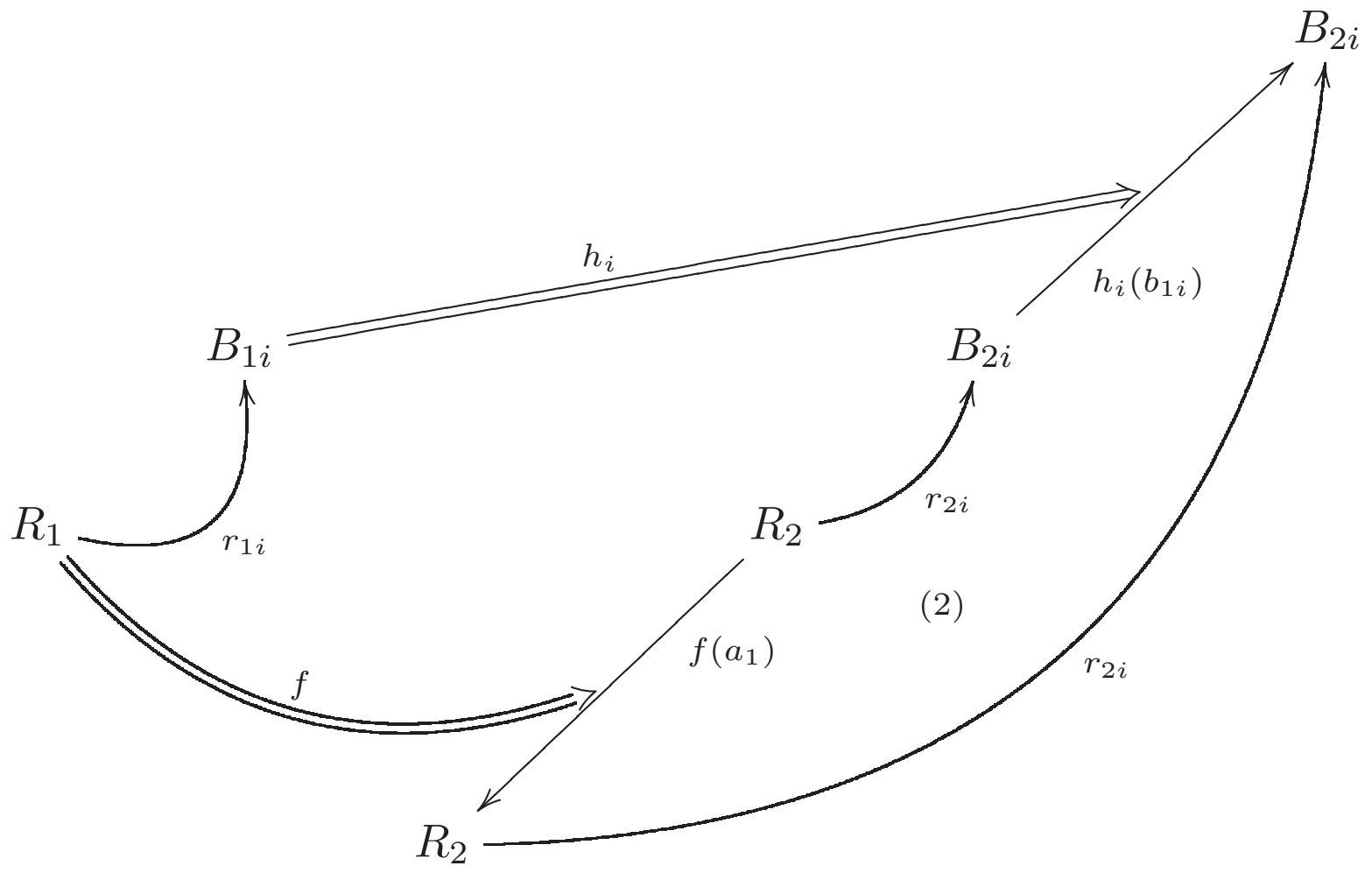

diagram (2) is commutative. According to the definition 2.1.1, there exists a unique morphism of $\Omega_2$-algebra

$$s_2:R_2\longrightarrow P_2$$

such that following diagram is commutative

(2.2.13)

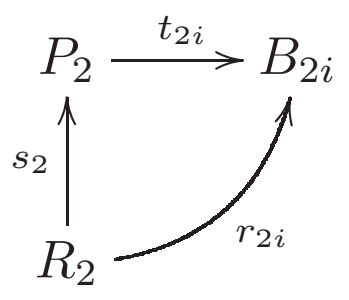

$$t_{2i}\circ s_2=r_{2i}$$

Let $a_2\in R_2$. Let

(2.2.14) $$b_2=s_2(a_2)\in P_2$$

From commutativity of the diagram (2.2.13) and statements (2.2.14), (2.2.3), it follows that

(2.2.15) $$b_{2i}=r_{2i}(a_2)$$

Let

(2.2.16) $$c_2=f(a_1)\circ a_2$$

From commutativity of the diagram (2) and equations (2.2.8), (2.2.15), (2.2.16), it follows that

$$d_{2i}=r_{2i}(c_2) \tag{2.2.17}$$

From equations (2.2.8), (2.2.17), it follows that

$$d_2=s_2(c_2) \tag{2.2.18}$$

and this is consistent with commutativity of the diagram (2.2.13).

For each $i\in I$,we join diagrams of maps (2.2.4), (2.2.9), (2.2.13), (2.2.12)

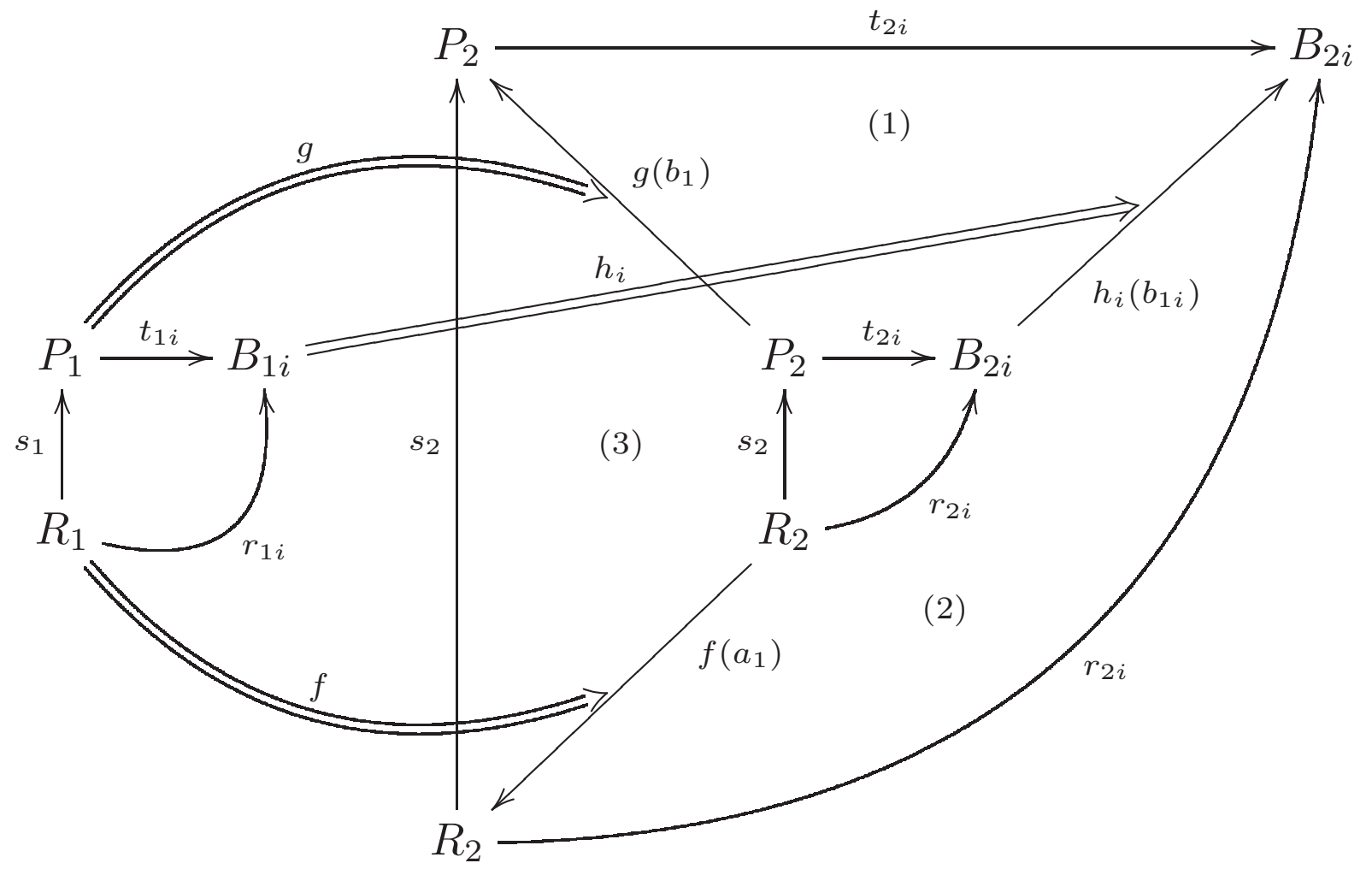

From equations (2.2.7), (2.2.14) and from equations (2.2.16), (2.2.18), commutativity of the diagram (3) follows. Therefore, the map $(s_1,s_2)$ is morphism of representations from $f$ into $g$. According to the theorem [4]-2.2.10, the morphism $(s_1,s_2)$ is defined unambiguously, since we require (2.2.18).

According to the definition 2.1.1, the representation $g$ and family of morphisms of representation $((t_{1i},t_{2i}),i\in I)$ is product in the category $(\mathcal{A}_1*)\mathcal{A}_2$ . □

## 2.3. Reduced Cartesian Product of Representations

Definition 2.3.1. *Let $A_1$ be $\Omega_1$-algebra. Let $\mathcal{A}_2$ be category of $\Omega_2$-algebras. We define* **category $(A_1*)\mathcal{A}_2$ of left-side representations**. *Left-side representations of $\Omega_1$-algebra $A_1$ in $\Omega_2$-algebra are objects of this category. Reduced morphisms of corresponding representations are morphisms of this category.* □

Theorem 2.3.2. *In category $(A_1*)\mathcal{A}_2$ there exists product of effective left-side representations of $\Omega_1$-algebra $A_1$ in $\Omega_2$-algebra and the product is effective left-side representations of $\Omega_1$-algebra $A_1$.*

Proof. Let

$$A_2=\prod_{i\in I}A_{2i}$$

be product of family of $\Omega_2$-algebras $\{A_{2i},i\in I\}$ and for any $i\in I$ the map

$$t_i:A_2\longrightarrow A_{2i}$$

be projection onto factor $i$. For each $i\in I$, let

$$h_i:A_1 \xrightarrow{*} A_{2i}$$

be effective $A_1*$-representation in $\Omega_2$-algebra $A_{2i}$.

Let $b_1\in A_1$. Let $b_2\in A_2$. According to the statement 2.1.3.3, we can represent $A_2$-number $b_2$ as tuple of $A_{2i}$-numbers

$$b_2=(b_{2i},i\in I)\quad b_{2i}=t_i(b_2)\in A_{2i} \tag{2.3.1}$$

Lemma 2.3.3. *For each $i\in I$, consider diagram of maps*

(2.3.2)

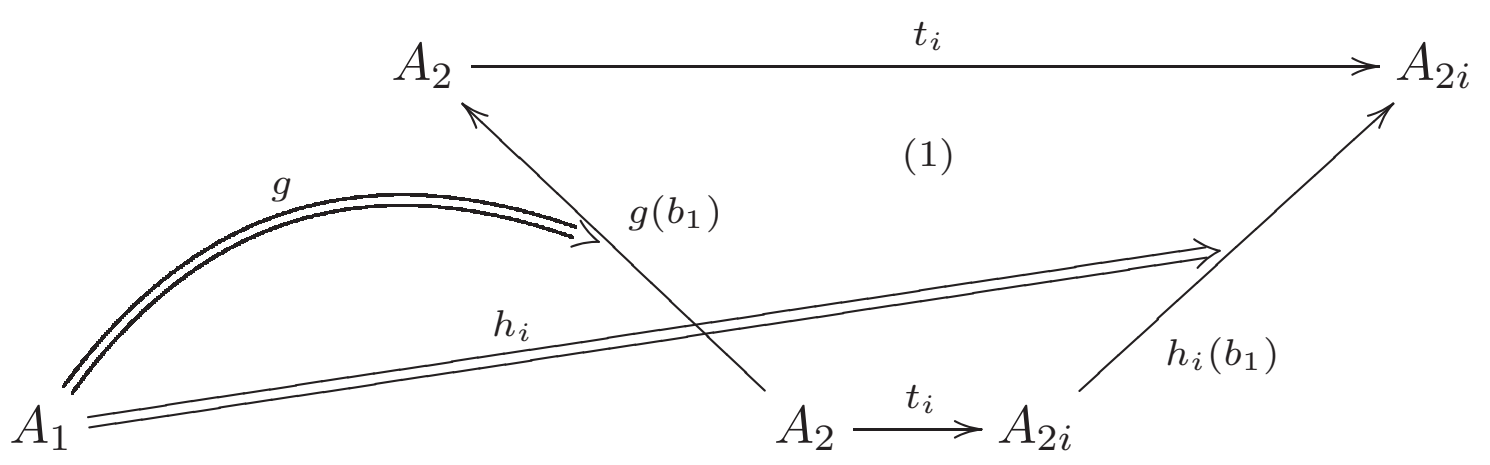

*Let map*

$$g:A_1\rightarrow {}^*A_2$$

*be defined by the equation*

$$g(b_1)\circ b_2=(h_i(b_1)\circ b_{2i},i\in I) \tag{2.3.3}$$

*Then the map $g$ is effective $A_1*$-representation in $\Omega_2$-algebra $A_2$*

$$g:A_1 \xrightarrow{*} A_2$$

*The map $t_i$ is reduced morphism of representation $g$ into representation $h_i$.*

Proof.

2.3.3.1: According to definitions [4]-2.1.1, [4]-2.1.2, the map $h_i(b_1)$ is homomorphism of $\Omega_2$-algebra $A_{2i}$. According to the theorem 2.1.6, from commutativity of the diagram (1) for each $i\in I$, it follows that the map

$$g(b_1):A_2\rightarrow A_2$$

defined by the equation (2.3.3) is homomorphism of $\Omega_2$-algebra $A_2$.

2.3.3.2: According to the definition [4]-2.1.2, the set ${}^*A_{2i}$ is $\Omega_1$-algebra. According to the lemma 2.2.1, the set ${}^\circ A_2\subseteq {}^*A_2$ is $\Omega_1$-algebra.

2.3.3.3: According to the definition [4]-2.1.2, the map

$$h_i:A_1\rightarrow {}^*A_{2i}$$

is homomorphism of $\Omega_1$-algebra. According to the theorem 2.1.6, the map

$$g:A_1\rightarrow {}^*A_2$$

defined by the equation

$$g(b_1)=(h_i(b_1),i\in I)$$

is homomorphism of $\Omega_1$-algebra.

According to statements 2.3.3.1, 2.3.3.3 and to the definition [4]-2.1.2, the map $g$ is $A_1*$-representation in $\Omega_2$-algebra $A_2$.

For any $i \in I$, according to the definition [4]-2.1.6, $A_1$-number $a_1$ generates unique transformation

$$b_{22i} = h_i(b_1) \circ b_{21i} \tag{2.3.4}$$

Let $b_{21}, b_{22} \in A_2$. According to the statement 2.1.3.3, we can represent $A_2$-numbers $b_{21}, b_{22}$ as tuples of $A_{2i}$-numbers

$$\begin{aligned} b_{21} &= (b_{21i}, i \in I) \quad b_{21i} = t_i(b_{21}) \in A_{2i} \\ b_{22} &= (b_{22i}, i \in I) \quad b_{22i} = t_i(b_{22}) \in A_{2i} \end{aligned} \tag{2.3.5}$$

According to the definition (2.3.3) of the representation $g$, from equations (2.3.4), (2.3.5), it follows that $A_1$-number $a_1$ generates unique transformation

$$b_{22} = (h_i(b_1) \circ b_{21i}, i \in I) = g(b_1) \circ b_{21} \tag{2.3.6}$$

According to the definition [4]-2.1.6, the representation $g$ is effective.

From commutativity of diagram (1) and from the definition [4]-2.2.2, it follows that map $t_i$ is reduced morphism of representation $g$ into representation $h_i$. ⊙

Let

$$d_2 = g(b_1) \circ b_2 \quad d_2 = (d_{2i}, i \in I) \tag{2.3.7}$$

From equations (2.3.3), (2.3.7), it follows that

$$d_{2i} = h_i(b_1) \circ b_{2i} \tag{2.3.8}$$

Let $R_2$ be other object of category $\mathcal{A}_2$. Let

$$f : A_1 \xrightarrow{*} R_2$$

be effective $A_1*$-representation in $\Omega_2$-algebra $R_2$. For any $i \in I$, let there exist morphism

$$r_i : R_2 \longrightarrow A_{2i}$$

of representations from $f$ into $h_i$. According to the remark [4]-2.2.6, in diagram of maps

(2.3.9)

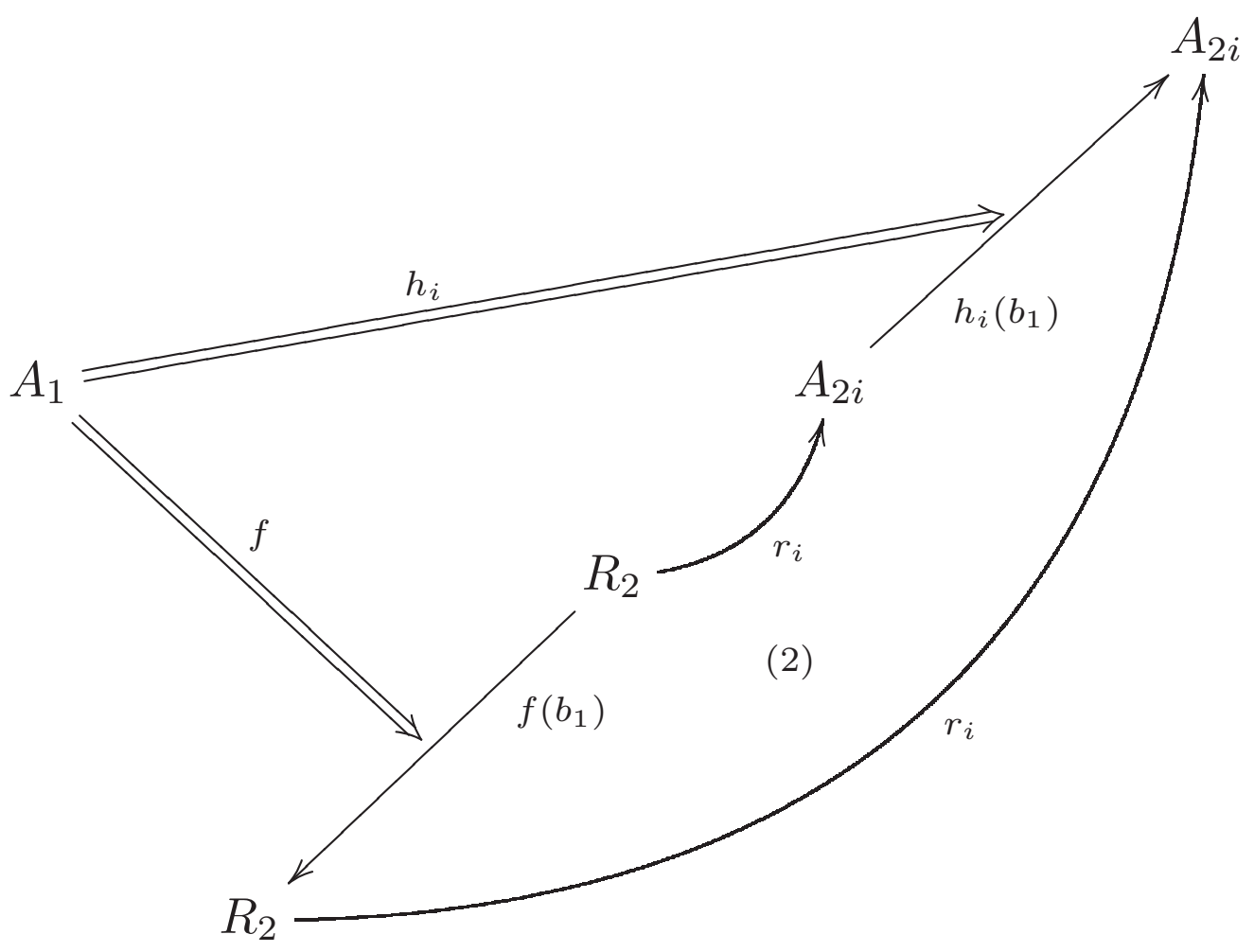

diagram (2) is commutative. According to the definition 2.1.1, there exists a unique morphism of $\Omega_2$-algebra

$$s : R_2 \longrightarrow A_2$$

such that following diagram is commutative

(2.3.10)

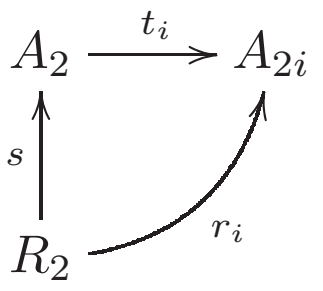

$$t_i \circ s = r_i$$

Let $a_2 \in R_2$. Let

$$b_2 = s(a_2) \in A_2 \tag{2.3.11}$$

From commutativity of the diagram (2.3.10) and statements (2.3.11), (2.3.1), it follows that

$$b_{2i} = r_i(a_2) \tag{2.3.12}$$

Let

$$c_2 = f(a_1) \circ a_2 \tag{2.3.13}$$

From commutativity of the diagram (2) and equations (2.3.8), (2.3.12), (2.3.13), it follows that

$$d_{2i} = r_i(c_2) \tag{2.3.14}$$

From equations (2.3.8), (2.3.14), it follows that

$$d_2 = s(c_2) \tag{2.3.15}$$

and this is consistent with commutativity of the diagram (2.3.10).

For each $i \in I$,we join diagrams of maps (2.3.2), (2.3.10), (2.3.9)

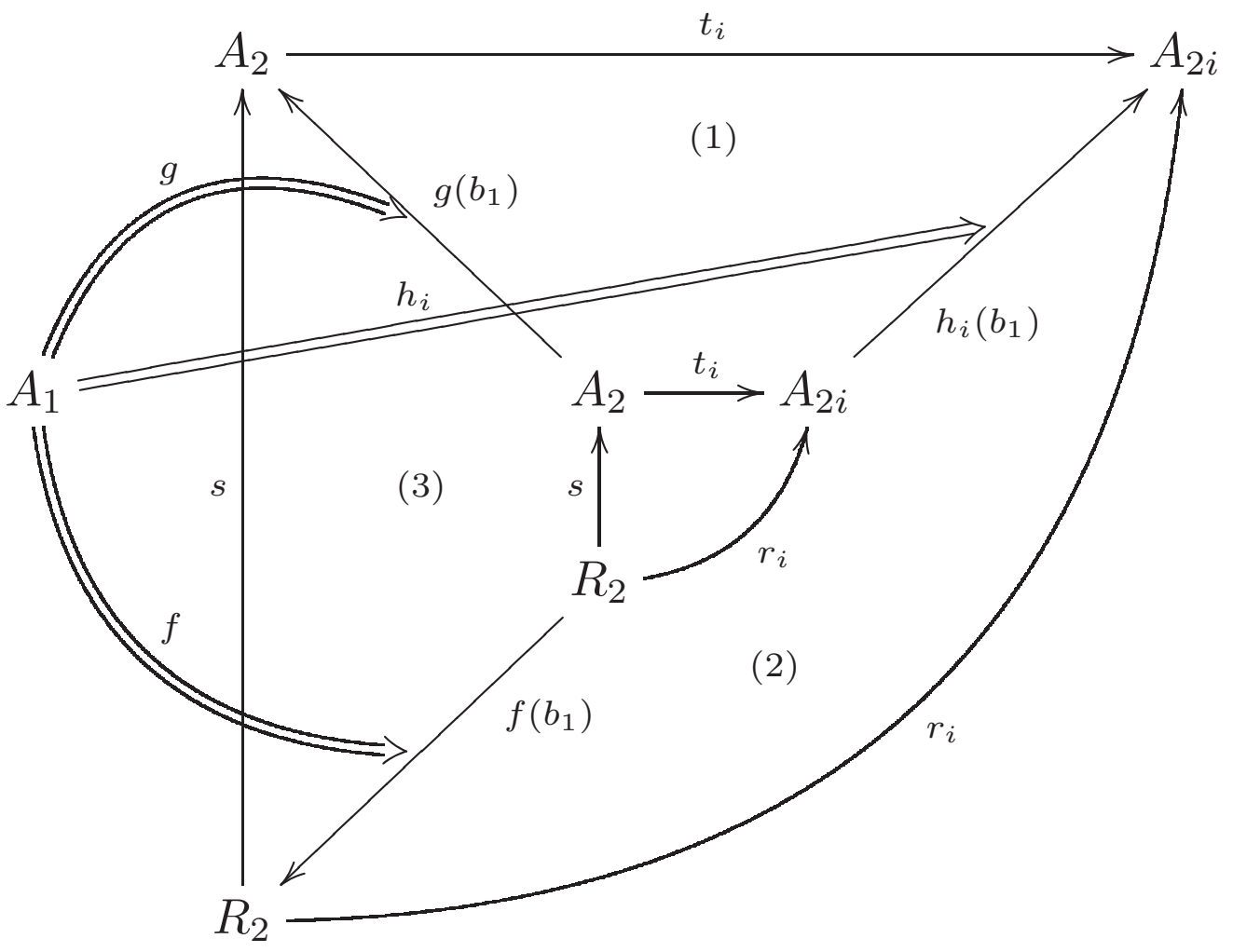

From equations (2.3.7), (2.3.11) and from equations (2.3.13), (2.3.15), commutativity of the diagram (3) follows. Therefore, the map $s$ is reduced morphism of representations from $f$ into $g$. According to the remark [4]-2.3.2, the map $s$ is homomorphism of $\Omega_2$ algebra. According to the theorem 2.1.3 and to the definition 2.1.1, the reduced morphism $s$ is defined unambiguously.

According to the definition 2.1.1, the representation $g$ and family of morphisms of representation $(t_i, i \in I)$ is product in the category $(A_1*)\mathcal{A}_2$ . □

CHAPTER 3

# Tensor Product of Representations

## 3.1. Polymorphism of Representations

Definition 3.1.1. *Let* $A_1$, ..., $A_n$ *be* $\Omega_1$*-algebras. Let* $B_1$, ..., $B_n$, $B$ *be* $\Omega_2$*-algebras. Let, for any* $k$, $k = 1$, ..., $n$,

$$f_k : A_k \dashrightarrow* B_k$$

*be representation of* $\Omega_1$*-algebra* $A_k$ *in* $\Omega_2$*-algebra* $B_k$. *Let*

$$f : A \dashrightarrow* B$$

*be representation of* $\Omega_1$*-algebra* $A$ *in* $\Omega_2$*-algebra* $B$. *The map*

$$r : A_1 \times ... \times A_n \to A \quad R : B_1 \times ... \times B_n \to B$$

*is called* **polymorphism of representations** $f_1$, ..., $f_n$ *into representation* $f$, *if, for any* $k$, $k = 1$, ..., $n$,, *provided that all variables except variables* $a_k \in A_k$, $b_k \in B_k$ *have given value, the map* $(r, R)$ *is a morphism of representation* $f_k$ *into representation* $f$.

*If* $f_1 = ... = f_n$, *then we say that the map* $(r, R)$ *is polymorphism of representation* $f_1$ *into representation* $f$.

*If* $f_1 = ... = f_n = f$, *then we say that the map* $(r, R)$ *is polymorphism of representation* $f$. □

Theorem 3.1.2. *Let the map* $(r, R)$ *be polymorphism of representations* $f_1$, ..., $f_n$ *into representation* $f$. *The map* $(r, R)$ *satisfies to the equality*

$$R(f_1(a_1)(m_1), ..., f_n(a_n)(m_n)) = f(r(a_1, ..., a_n))(R(m_1, ..., m_n)) \tag{3.1.1}$$

*Let* $\omega_1 \in \Omega_1(p)$. *For any* $k$, $k = 1$, ..., $n$,*the map* $r$ *satisfies to the equality*

$$\begin{aligned} & r(a_1, ..., a_{k\cdot 1}...a_{k\cdot p}\omega_1, ..., a_n) \\ = & r(a_1, ..., a_{k\cdot 1}, ..., a_n)...r(a_1, ..., a_{k\cdot p}, ..., a_n)\omega_1 \end{aligned} \tag{3.1.2}$$

*Let* $\omega_2 \in \Omega_2(p)$. *For any* $k$, $k = 1$, ..., $n$,*the map* $R$ *satisfies to the equality*

$$\begin{aligned} & R(m_1, ..., m_{k\cdot 1}...m_{k\cdot p}\omega_2, ..., m_n) \\ = & R(m_1, ..., m_{k\cdot 1}, ..., m_n)...R(m_1, ..., m_{k\cdot p}, ..., m_n)\omega_2 \end{aligned} \tag{3.1.3}$$

Proof. The equality (3.1.1) follows from the definition 3.1.1 and the equality [4]-(2.2.4). The equality (3.1.2) follows from the statement that for any $k$, $k = 1$, ..., $n$,provided that all variables except the variable $x_k \in A_k$ have given value, the map $r$ is homomorphism of $\Omega_1$-algebra $A_k$ into $\Omega_1$-algebra $A$. The equality (3.1.3) follows from the statement that for any $k$, $k = 1$, ..., $n$,provided that all variables except the variable $m_k \in B_k$ have given value, the map $R$ is homomorphism of $\Omega_2$-algebra $B_k$ into $\Omega_2$-algebra $B$. □

Definition 3.1.3. *Let $A$, $B_1$, ..., $B_n$, $B$ be universal algebras. Let, for any $k$, $k = 1$, ..., $n$,*

$$f_k : A \dashrightarrow B_k$$

*be effective representation of $\Omega_1$-algebra $A_k$ in $\Omega_2$-algebra $B_k$. Let*

$$f : A \dashrightarrow B$$

*be effective representation of $\Omega_1$-algebra $A$ in $\Omega_2$-algebra $B$. The map*

$$r_2 : B_1 \times ... \times B_n \to B$$

*is called* **reduced polymorphism of representations** *$f_1$, ..., $f_n$ into representation $f$, if, for any $k$, $k = 1$, ..., $n$,provided that all variables except the variable $x_k \in B_k$ have given value, the map $R$ is a reduced morphism of representation $f_k$ into representation $f$.*

*If $f_1 = ... = f_n$, then we say that the map $R$ is reduced polymorphism of representation $f_1$ into representation $f$.*

*If $f_1 = ... = f_n = f$, then we say that the map $R$ is reduced polymorphism of representation $f$.* □

Theorem 3.1.4. *Let the map $R$ be reduced polymorphism of effective representations $f_1$, ..., $f_n$ into effective representation $f$. For any $k$, $k = 1$, ..., $n$,the map $R$ satisfies to the equality*

$$R(m_1, ..., f_k(a) \circ m_k, ..., m_n) = f(a) \circ R(m_1, ..., m_n) \tag{3.1.4}$$

*Let $\omega_2 \in \Omega_2(p)$. For any $k$, $k = 1$, ..., $n$,the map $R$ satisfies to the equality*

$$\begin{aligned} &R(m_1, ..., m_{k \cdot 1}...m_{k \cdot p}\omega_2, ..., m_n) \\ &= R(m_1, ..., m_{k \cdot 1}, ..., m_n)...R(m_1, ..., m_{k \cdot p}, ..., m_n)\omega_2 \end{aligned} \tag{3.1.5}$$

Proof. The equality (3.1.4) follows from the definition 3.1.3 and the equality [4]-(2.3.3). The equality (3.1.5) follows from the statement that for any $k$, $k = 1$, ..., $n$, provided that all variables except the variable $m_k \in B_k$ have given value, the map $R$ is homomorphism of $\Omega_2$-algebra $B_k$ into $\Omega_2$-algebra $B$. □

We also say that the map $(r, R)$ is polymorphism of representations in $\Omega_2$-algebras $B_1$, ..., $B_n$ into representation in $\Omega_2$-algebra $B$. Similarly, we say that the map $R$ is reduced polymorphism of representations in $\Omega_2$-algebras $B_1$, ..., $B_n$ into representation in $\Omega_2$-algebra $B$.

Comparison of definitions 3.1.1 and 3.1.3 shows that there is a difference between these two forms of polymorphism. This is particularly evident when comparing the difference between equalities (3.1.1) and (3.1.4). If we want to be able to express the reduced polymorphism of representations using polymorphism of representations, then we must require two conditions:

(1) The representation $f$ of universal algebra contains the identity transformation $\delta$. Therefore, there exists $e \in A$ such that $f(e) = \delta$. Without loss of generality, we assume that the choice of $e \in A$ does not depend on whether we consider the representation $f_1$, ..., or $f_n$.

(2) For any $k$, $k=1, ..., n$,

$$r(a_1,...,a_n)=a_k \quad a_i=e \quad i\ne k \tag{3.1.6}$$

Then, provided that $a_i=e$, $i\ne k$, the equality (3.1.1) has form

$$R(m_1,...,f_k(a_k)\circ m_k,...,m_n)=f(r(e,...,a_k,...,e))\circ R(m_1,...,m_n) \tag{3.1.7}$$

It is evident that the equality (3.1.7) coincides with the equality (3.1.4).

A similar problem appears in the analysis of reduced polymorphism of representations. Using the equality (3.1.4), we can write an expression

$$R(m_1,...,f_k(a_k)\circ m_k,...,f_l(a_l)\circ m_l,...,m_n) \tag{3.1.8}$$

either in the following form

$$\begin{aligned}&R(m_1,...,f_k(a_k)\circ m_k,...,f_l(a_l)\circ m_l,...,m_n)\\=&f(a_k)\circ R(m_1,...,m_k,...,f_l(a_l)\circ m_l,...,m_n)\\=&f(a_k)\circ(f(a_l)\circ R(m_1,...,m_k,...,m_l,...,m_n))\\=&(f(a_k)\circ f(a_l))\circ R(m_1,...,m_k,...,m_l,...,m_n)\end{aligned} \tag{3.1.9}$$

or in the following form

$$\begin{aligned}&R(m_1,...,f_k(a_k)\circ m_k,...,f_l(a_l)\circ m_l,...,m_n)\\=&f(a_l)\circ R(m_1,...,f_k(a_k)\circ m_k,...,m_l,...,m_n)\\=&f(a_l)\circ(f(a_k)\circ R(m_1,...,m_k,...,m_l,...,m_n))\\=&(f(a_l)\circ f(a_k))\circ R(m_1,...,m_k,...,m_l,...,m_n)\end{aligned} \tag{3.1.10}$$

Maps $f(a_k)$, $f(a_l)$ are homomorphisms of $\Omega_2$-algebra $B$. Therefore, the map $f(a_k)\circ f(a_l)$ is homomorphism of $\Omega_2$-algebra $B$. However, not every $\Omega_1$-algebra $A$ has such $a$ which depends on $a_k$ and $a_l$ and satisfies to equality

$$f(a)=f(a_k)\circ f(a_l)$$

Since the representation $f$ is single transitive and for any $A$-numbers $a$, $b$ there exists $A$-number $c$ such that

$$f(c)=f(a)\circ f(b) \tag{3.1.11}$$

then the equality (3.1.11) uniquely determines $A$-number $c$. Therefore, we can introduce the product

$$c_1=a_1*b_1$$

such way that

$$f(a*b)=f(a)\circ f(b) \tag{3.1.12}$$

Definition 3.1.5. *Let product*

$$c_1=a_1*b_1$$

*be operation of $\Omega_1$-algebra $A$. Let $\Omega=\Omega_1\setminus\{*\}$. For any operation $\omega\in\Omega(p)$, the product is distributive over the operation $\omega$*

$$a*(b_1...b_n\omega)=(a*b_1)...(a*b_n)\omega \tag{3.1.13}$$

$$(b_1...b_n\omega)*a=(b_1*a)...(b_n*a)\omega \tag{3.1.14}$$

*$\Omega_1$-algebra $A$ is called* **multiplicative $\Omega$-group**.[3.1] □

DEFINITION 3.1.6. *Let $A$, $B$ be multiplicative $\Omega$-groups. The map*

$$f : A \to B$$

*is called* **multiplicative**, *if*

$$f(a * b) = f(a) \circ f(b)$$

□

THEOREM 3.1.7. *Single transitive representation of multiplicative $\Omega$-group is multiplicative map.*

PROOF. Theorem follows from the equality (3.1.12) and from the definition 3.1.6. □

However the statement of the theorem 3.1.7 is not enough to prove equality of expressions (3.1.9) and (3.1.10).

DEFINITION 3.1.8. *If*

$$a * b = b * a \tag{3.1.15}$$

*then multiplicative $\Omega$-group is called* **Abelian**. □

THEOREM 3.1.9. *Let*

$$f : A \xrightarrow{*} M$$

*be effective representation of Abelian multiplicative $\Omega$-group $A$. Then*

$$\begin{aligned} &f(a_k) \circ (f(a_l) \circ R(m_1, ..., m_k, ..., m_l, ..., m_n)) \\ =&f(a_l) \circ (f(a_k) \circ R(m_1, ..., m_k, ..., m_l, ..., m_n)) \end{aligned} \tag{3.1.16}$$

PROOF. From equalities (3.1.9), (3.1.10), (3.1.12), (3.1.15), it follows that

$$\begin{aligned} &f(a_k) \circ (f(a_l) \circ R(m_1, ..., m_k, ..., m_l, ..., m_n)) \\ =&(f(a_k) \circ f(a_l)) \circ R(m_1, ..., m_k, ..., m_l, ..., m_n) \\ =&f(a_k * a_l) \circ R(m_1, ..., m_k, ..., m_l, ..., m_n) \\ =&f(a_l * a_k) \circ R(m_1, ..., m_k, ..., m_l, ..., m_n) \\ =&(f(a_l) \circ f(a_k)) \circ R(m_1, ..., m_k, ..., m_l, ..., m_n) \\ =&f(a_l) \circ (f(a_k) \circ R(m_1, ..., m_k, ..., m_l, ..., m_n)) \end{aligned} \tag{3.1.17}$$

The equality (3.1.16) follows from the equality (3.1.17). □

THEOREM 3.1.10. *Let $A$ be Abelian multiplicative $\Omega$-group. Let $R$ be reduced polymorphism of effective representations $f_1$, ..., $f_n$ into effective representation $f$. Then for any $k$, $l$, $k = 1$, ..., $n$, $l = 1$, ..., $n$,*

$$\begin{aligned} &R(m_1, ..., f_k(a) \circ m_k, ..., m_l, ..., m_n) \\ =&R(m_1, ..., m_k, ..., f_l(a) \circ m_l, ..., m_n) \end{aligned} \tag{3.1.18}$$

PROOF. The equality (3.1.18) directly follows from the equality (3.1.4). □

[3.1] The definition of multiplicative $\Omega$-group is similar to the definition [6]-2.1.3 of $\Omega$-group. However, $\Omega$-group assumes addition as group operation. It is important for us that group operation of multiplicative $\Omega$-group is product. Moreover, operation $\omega$ of $\Omega$-group is distributive over sum. In multiplicative $\Omega$-group, product is distributive over operation $\omega$.

THEOREM 3.1.11. *Let*

$$A \xrightarrow{*} B_1 \qquad A \xrightarrow{*} B_2 \qquad A \xrightarrow{*} B$$

*be effective representations of Abelian multiplicative $\Omega_1$-group $A$ in $\Omega_2$-algebras $B_1$, $B_2$, $B$. Let $\Omega_2$-algebra have 2 operations, namely $\omega_1 \in \Omega_2(p)$, $\omega_2 \in \Omega_2(q)$. The equality*

$$(a_{1\cdot 1}...a_{1\cdot q}\omega_2)...(a_{p\cdot 1}...a_{p\cdot q}\omega_2)\omega_1 = (a_{1\cdot 1}...a_{p\cdot 1}\omega_1)...(a_{1\cdot q}...a_{p\cdot q}\omega_1)\omega_2 \tag{3.1.19}$$

*is necessary condition of existence of reduced polymorphism*

$$R : B_1 \times B_2 \to B$$

PROOF. Let $a_1$, ..., $a_p \in B_1$, $b_1$, ..., $b_q \in B_2$. According to the equality (3.1.5), the expression

$$R(a_1...a_p\omega_1, b_1...b_q\omega_2) \tag{3.1.20}$$

can have 2 values

$$\begin{aligned} &R(a_1...a_p\omega_1, b_1...b_q\omega_2)\\ =&R(a_1, b_1...b_q\omega_2)...R(a_p, b_1...b_q\omega_2)\omega_1\\ =&(R(a_1, b_1)...R(a_1, b_q)\omega_2)...(R(a_p, b_1)...R(a_p, b_q)\omega_2)\omega_1 \end{aligned} \tag{3.1.21}$$

$$\begin{aligned} &R(a_1...a_p\omega_1, b_1...b_q\omega_2)\\ =&R(a_1...a_p\omega_1, b_1)...R(a_1...a_p\omega_1, b_q)\omega_2\\ =&(R(a_1, b_1)...R(a_p, b_1)\omega_1)...(R(a_1, b_q)...R(a_p, b_q)\omega_1)\omega_2 \end{aligned} \tag{3.1.22}$$

From equalities (3.1.21), (3.1.22), it follows that

$$\begin{aligned} &(R(a_1, b_1)...R(a_1, b_q)\omega_2)...(R(a_p, b_1)...R(a_p, b_q)\omega_2)\omega_1\\ =&(R(a_1, b_1)...R(a_p, b_1)\omega_1)...(R(a_1, b_q)...R(a_p, b_q)\omega_1)\omega_2 \end{aligned} \tag{3.1.23}$$

Therefore, the expression (3.1.20) is properly defined iff the equality (3.1.23) is true. Let

$$a_{i\cdot j} = R(a_i, b_j) \in A \tag{3.1.24}$$

The equality (3.1.19) follows from equalities (3.1.23), (3.1.24). □

THEOREM 3.1.12. *There exists reduced polymorphism of effective representations of Abelian multiplicative $\Omega$-group in Abelian group.*

PROOF. Since sum in Abelian group is commutative and associative, then the theorem follows from the theorem 3.1.11. □

THEOREM 3.1.13. *There is no reduced polymorphism of effective representations of Abelian multiplicative $\Omega$-group in ring.*

PROOF. There are two operations in the ring: sum which is commutative and associative and product which is distributive over sum. According to the theorem 3.1.11, the existence of polymorphism of effective representation in the ring implies that sum and product must satisfy the equality

$$a_{1\cdot 1}a_{2\cdot 1} + a_{1\cdot 2}a_{2\cdot 2} = (a_{1\cdot 1} + a_{1\cdot 2})(a_{2\cdot 1} + a_{2\cdot 2}) \tag{3.1.25}$$

However right hand side of the equality (3.1.25) has form

$$\begin{aligned}(a_{1\cdot 1}+a_{1\cdot 2})(a_{2\cdot 1}+a_{2\cdot 2}) &= (a_{1\cdot 1}+a_{1\cdot 2})a_{2\cdot 1}+(a_{1\cdot 1}+a_{1\cdot 2})a_{2\cdot 2}\\ &= a_{1\cdot 1}a_{2\cdot 1}+a_{1\cdot 2}a_{2\cdot 1}+a_{1\cdot 1}a_{2\cdot 2}+a_{1\cdot 2}a_{2\cdot 2}\end{aligned}$$

Therefore, the equality (3.1.25) is not true. □

QUESTION 3.1.14. *It is possible that polymorphism of representations exists only for effective representation in Abelian group. However, this statement has not been proved.* □

## 3.2. Congruence

THEOREM 3.2.1. *Let $N$ be equivalence on the set $A$. Consider category $\mathcal{A}$ whose objects are maps*[3.2]

$$f_1:A\to S_1\quad \ker f_1\supseteq N$$

$$f_2:A\to S_2\quad \ker f_2\supseteq N$$

*We define morphism $f_1\to f_2$ to be map $h:S_1\to S_2$ making following diagram commutative*

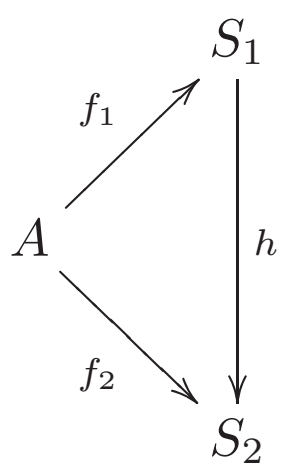

*The map*

$$\mathrm{nat}\,N:A\to A/N$$

*is universally repelling in the category $\mathcal{A}$.*[3.3]

PROOF. Consider diagram

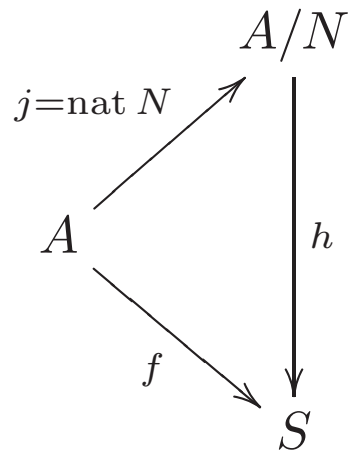

$$\ker f\supseteq N \tag{3.2.1}$$

From the statement (3.2.1) and the equality

$$j(a_1)=j(a_2)$$

it follows that

$$f(a_1)=f(a_2)$$

[3.2]The statement of lemma is similar to the statement on p. [1]-119.
[3.3]See definition of universal object of category in definition on p. [1]-57.

Therefore, we can uniquely define the map $h$ using the equality

$$h(j(b)) = f(b)$$

$\square$

THEOREM 3.2.2. *Let*

$$f : A \twoheadrightarrow^{*} B$$

*be representation of* $\Omega_1$*-algebra* $A$ *in* $\Omega_2$*-algebra* $B$*. Let* $N$ *be such congruence*[3.4] *on* $\Omega_2$*-algebra* $B$ *that any transformation* $h \in {}^*B$ *is coordinated with congruence* $N$*. There exists representation*

$$f_1 : A \twoheadrightarrow^{*} B/N$$

*of* $\Omega_1$*-algebra* $A$ *in* $\Omega_2$*-algebra* $B/N$ *and the map*

$$\mathrm{nat}\, N : B \rightarrow B/N$$

*is reduced morphism of representation* $f$ *into the representation* $f_1$

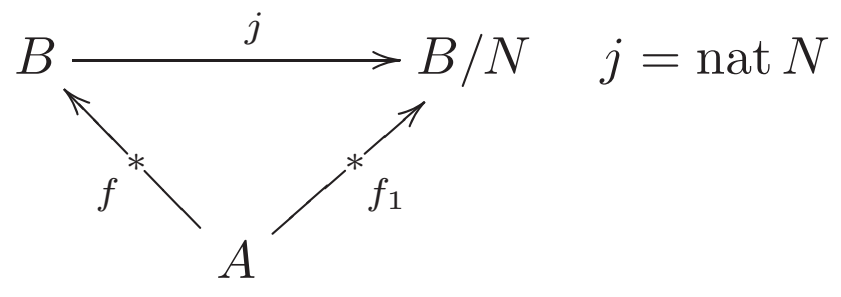

PROOF. We can represent any element of the set $B/N$ as $j(a)$, $a \in B$.

According to the theorem [9]-II.3.5, there exists a unique $\Omega_2$-algebra structure on the set $B/N$. If $\omega \in \Omega_2(p)$, then we define operation $\omega$ on the set $B/N$ according to the equality (3) on page [9]-59

$$j(b_1)...j(b_p)\omega = j(b_1...b_p\omega) \tag{3.2.2}$$

As well as in the proof of the theorem [4]-2.2.16, we can define the representation

$$f_1 : A \twoheadrightarrow^{*} B/N$$

using equality

$$f_1(a) \circ j(b) = j(f(a) \circ b) \tag{3.2.3}$$

We can represent the equality (3.2.3) using diagram

(3.2.4)

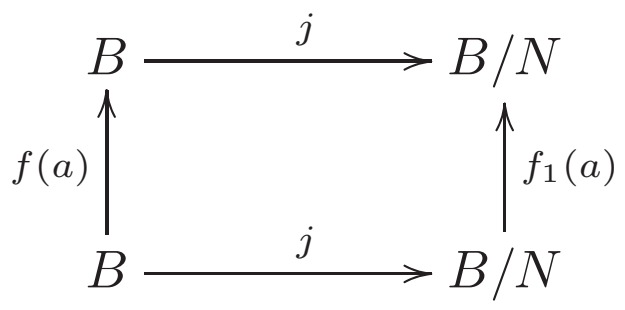

Let $\omega \in \Omega_2(p)$. Since the maps $f(a)$ and $j$ are homomorphisms of $\Omega_2$-algebra, then

$$\begin{aligned} f_1(a) \circ (j(b_1)...j(b_p)\omega) &= f_1(a) \circ j(b_1...b_p\omega) \\ &= j(f(a) \circ (b_1...b_p\omega)) \\ &= j((f(a) \circ b_1)...(f(a) \circ b_p)\omega) \\ &= j(f(a) \circ b_1)...j(f(a) \circ b_p)\omega \\ &= (f_1(a) \circ j(b_1))...(f_1(a) \circ j(b_p))\omega \end{aligned} \tag{3.2.5}$$

[3.4] See the definition of congruence on p. [9]-57.

From the equality (3.2.5), it follows that the map $f_1(a)$ is homomorphism of $\Omega_2$-algebra. From the equality (3.2.3) and from the [4]-2.3.2, it follows that the map $j$ is reduced morphism of the representation $f$ into the representation $f_1$. □

THEOREM 3.2.3. *Let*

$$f : A \xrightarrow{*} B$$

*be representation of $\Omega_1$-algebra $A$ in $\Omega_2$-algebra $B$. Let $N$ be such congruence on $\Omega_2$-algebra $B$ that any transformation $h \in {}^*B$ is coordinated with congruence $N$. Consider category $\mathcal{A}$ whose objects are reduced morphisms of representations*[3.5]

$$R_1 : B \to S_1 \quad \ker R_1 \supseteq N$$

$$R_2 : B \to S_2 \quad \ker R_2 \supseteq N$$

*where $S_1$, $S_2$ are $\Omega_2$-algebras and*

$$g_1 : A \xrightarrow{*} S_1 \qquad g_2 : A \xrightarrow{*} S_2$$

*are representations of $\Omega_1$-algebra $A$. We define morphism $R_1 \to R_2$ to be reduced morphism of representations $h : S_1 \to S_2$ making following diagram commutative*

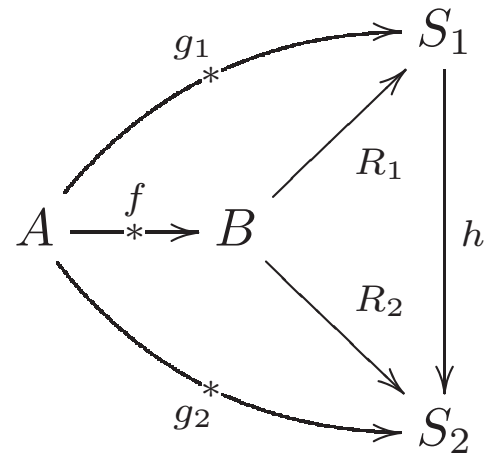

*The reduced morphism* nat $N$ *of representation $f$ into representation $f_1$ (the theorem 3.2.2) is universally repelling in the category $\mathcal{A}$.*[3.6]

PROOF. From the theorem 3.2.1, it follows that there exists and unique the map $h$ for which the following diagram is commutative

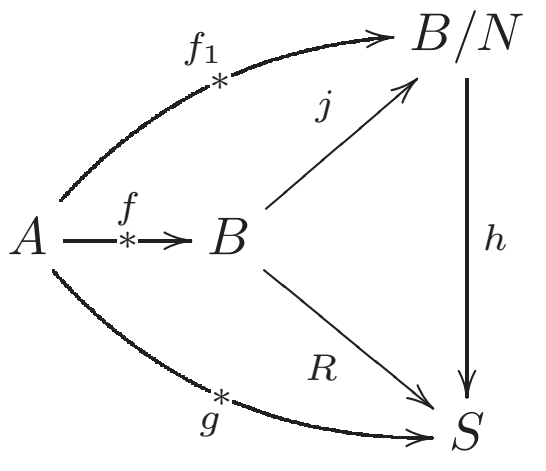

$j = \mathrm{nat}\, N \quad \ker R \supseteq N$

Therefore, we can uniquely define the map $h$ using equality

(3.2.6) $$h(j(b)) = R(b)$$

[3.5] The statement of lemma is similar to the statement on p. [1]-119.

[3.6] See definition of universal object of category in definition on p. [1]-57.

Let $\omega \in \Omega_2(p)$. Since maps $R$ and $j$ are homomorphisms of $\Omega_2$-algebra, then

$$\begin{aligned} h(j(b_1)...j(b_p)\omega) &= h(j(b_1...b_p\omega)) \\ &= R(b_1...b_p\omega) \\ &= R(b_1)...R(b_p)\omega \\ &= h(j(b_1))...h(j(b_p))\omega \end{aligned} \tag{3.2.7}$$

From the equality (3.2.7), it follows that the map $h$ is homomorphism of $\Omega_2$-algebra.

Since the map $R$ is reduced morphism of the representation $f$ into the representation $g$, then the following equality is satisfied

$$g(a)(R(b)) = R(f(a)(b)) \tag{3.2.8}$$

From the equality (3.2.6) it follows that

$$g(a)(h(j(b))) = g(a)(R(b)) \tag{3.2.9}$$

From the equalities (3.2.8), (3.2.9) it follows that

$$g(a)(h(j(b))) = R(f(a)(b)) \tag{3.2.10}$$

From the equalities (3.2.6), (3.2.10) it follows that

$$g(a)(h(j(b))) = h(j(f(a)(b))) \tag{3.2.11}$$

From the equalities (3.2.3), (3.2.11) it follows that

$$g(a)(h(j(b))) = h(f_1(a)(j(b))) \tag{3.2.12}$$

From the equality (3.2.12) it follows that the map $h$ is reduced morphism of representation $f_1$ into the representation $g$. □

## 3.3. Tensor Product of Representations

Definition 3.3.1. *Let $A$ be Abelian multiplicative $\Omega_1$-group. Let $B_1$, ..., $B_n$ be $\Omega_2$-algebras.*[3.7] *Let, for any $k$, $k = 1$, ..., $n$,*

$$f_k : A \xrightarrow{*} B_k$$

*be effective representation of multiplicative $\Omega_1$-group $A$ in $\Omega_2$-algebra $B_k$. Consider category $\mathcal{A}$ whose objects are reduced polymorphisms of representations $f_1$, ..., $f_n$*

$$r_1 : B_1 \times ... \times B_n \longrightarrow S_1 \qquad r_2 : B_1 \times ... \times B_n \longrightarrow S_2$$

*where $S_1$, $S_2$ are $\Omega_2$-algebras and*

$$g_1 : A \xrightarrow{*} S_1 \qquad g_2 : A \xrightarrow{*} S_2$$

*are effective representations of multiplicative $\Omega_1$-group $A$. We define morphism $R_1 \to R_2$ to be reduced morphism of representations $h : S_1 \to S_2$ making following*

[3.7] I give definition of tensor product of representations of universal algebra following to definition in [1], p. 601 - 603.

*diagram commutative*

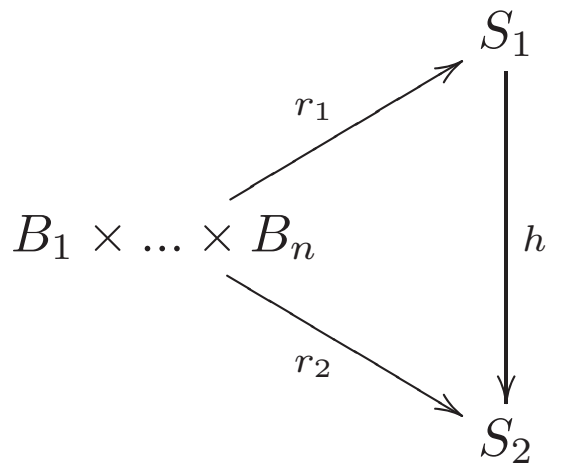

*Universal object* $B_1 \otimes ... \otimes B_n$ *of category* $\mathcal{A}$ *is called* **tensor product** *of representations* $B_1, ..., B_n$. □

Theorem 3.3.2. *Since there exists tensor product of effective representations, then tensor product is unique up to isomorphism of representations.*

Proof. Let $A$ be Abelian multiplicative $\Omega_1$-group. Let $B_1, ..., B_n$ be $\Omega_2$-algebras. Let, for any $k$, $k = 1, ..., n$,

$$f_k : A \xrightarrow{*} B_k$$

be effective representation of multiplicative $\Omega_1$-group $A$ in $\Omega_2$-algebra $B_k$. Let effective representations

$$g_1 : A \xrightarrow{*} S_1 \qquad g_2 : A \xrightarrow{*} S_2$$

be tensor product of representations $B_1, ..., B_n$. From commutativity of the diagram

(3.3.1)

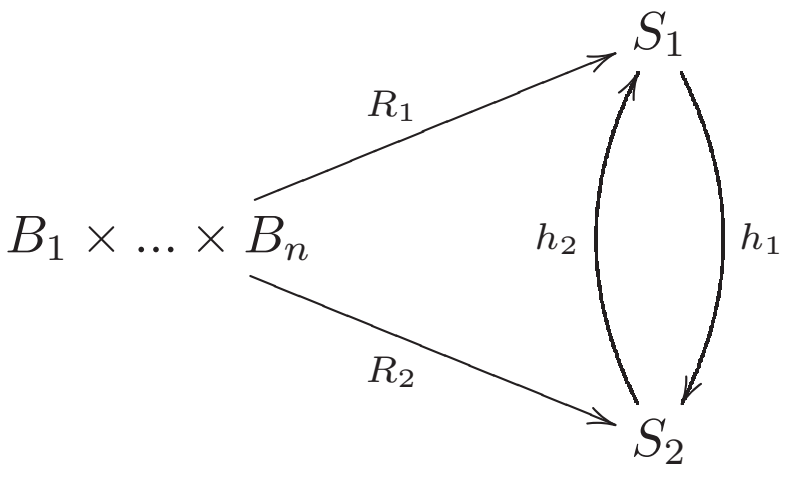

it follows that

$$R_1 = h_2 \circ h_1 \circ R_1$$
$$R_2 = h_1 \circ h_2 \circ R_2 \tag{3.3.2}$$

From equalities (3.3.2), it follows that morphisms of representation $h_1 \circ h_2$, $h_2 \circ h_1$ are identities. Therefore, morphisms of representation $h_1$, $h_2$ are isomorphisms. □

Convention 3.3.3. *Algebras* $S_1$, $S_2$ *may be different sets. However they are indistinguishable for us when we consider them as isomorphic representations. In such case, we write the statement* $S_1 = S_2$. □

Definition 3.3.4. *Tensor product*

$$B^{\otimes n} = B_1 \otimes ... \otimes B_n \quad B_1 = ... = B_n = B$$

*is called* **tensor power** *of representation* $B$. □

Theorem 3.3.5. *Since there exists polymorphism of representations, then there exists tensor product of representations.*

PROOF. Let

$$f : A \xrightarrow{*} M$$

be representation of $\Omega_1$-algebra $A$ generated by Cartesian product $B_1 \times ... \times B_n$ of sets $B_1$, ..., $B_n$.[3.8] Injection

$$i : B_1 \times ... \times B_n \longrightarrow M$$

is defined according to rule[3.9]

$$i \circ (b_1, ..., b_n) = (b_1, ..., b_n) \tag{3.3.3}$$

Let $N$ be equivalence generated by following equalities[3.10]

$$(b_1, ..., b_{i\cdot 1}...b_{i\cdot p}\omega, ..., b_n) = (b_1, ..., b_{i\cdot 1}, ..., b_n)...(b_1, ..., b_{i\cdot p}, ..., b_n)\omega \tag{3.3.4}$$

$$(b_1, ..., f_i(a) \circ b_i, ..., b_n) = f(a) \circ (b_1, ..., b_i, ..., b_n) \tag{3.3.5}$$

$$b_k \in B_k \quad k = 1, ..., n \quad b_{i\cdot 1}, ..., b_{i\cdot p} \in B_i \quad \omega \in \Omega_2(p) \quad a \in A$$

LEMMA 3.3.6. *Let $\omega \in \Omega_2(p)$. Then*

$$\begin{aligned} &f(c) \circ (b_1, ..., b_{i\cdot 1}...b_{i\cdot p}\omega, ..., b_n)\\ =&f(c) \circ ((b_1, ..., b_{i\cdot 1}, ..., b_n)...(b_1, ..., b_{i\cdot p}, ..., b_n)\omega) \end{aligned} \tag{3.3.6}$$

PROOF. From the equality (3.3.5), it follows that

$$f(c) \circ (b_1, ..., b_{i\cdot 1}...b_{i\cdot p}\omega, ..., b_n) = (b_1, ..., f_i(c) \circ (b_{i\cdot 1}...b_{i\cdot p}\omega), ..., b_n) \tag{3.3.7}$$

Since $f_i(c)$ is endomorphism of $\Omega_2$-algebra $B_i$, then from the equality (3.3.7), it follows that

$$f(c) \circ (b_1, ..., b_{i\cdot 1}...b_{i\cdot p}\omega, ..., b_n) = (b_1, ..., (f_i(c) \circ b_{i\cdot 1})...(f_i(c) \circ b_{i\cdot p})\omega, ..., b_n) \tag{3.3.8}$$

From equalities (3.3.8), (3.3.4), it follows that

$$\begin{aligned} &f(c) \circ (b_1, ..., b_{i\cdot 1}...b_{i\cdot p}\omega, ..., b_n)\\ =&(b_1, ..., f_i(c) \circ b_{i\cdot 1}, ..., b_n)...(b_1, ..., f_i(c) \circ b_{i\cdot p}, ..., b_n)\omega \end{aligned} \tag{3.3.9}$$

From equalities (3.3.9), (3.3.5), it follows that

$$\begin{aligned} &f(c) \circ (b_1, ..., b_{i\cdot 1}...b_{i\cdot p}\omega, ..., b_n)\\ =&(f(c) \circ (b_1, ..., b_{i\cdot 1}, ..., b_n))...(f(c) \circ (b_1, ..., b_{i\cdot p}, ..., b_n))\omega \end{aligned} \tag{3.3.10}$$

Since $f(c)$ is endomorphism of $\Omega_2$-algebra $B$, then the equality (3.3.6) follows from the equality (3.3.10). ⊙

LEMMA 3.3.7.

$$f(c) \circ (b_1, ..., f_i(a) \circ b_i, ..., b_n) = f(c) \circ (f(a) \circ (b_1, ..., b_i, ..., b_n)) \tag{3.3.11}$$

---

[3.8] According to theorems 2.1.3, 2.3.2, the set generated by reduced Cartesian product of representations $B_1$, ..., $B_n$ coincides with Cartesian product $B_1 \times ... \times B_n$ of sets $B_1$, ..., $B_n$. At this point of the proof, we do not consider any algebra structure on the set $B_1 \times ... \times B_n$.

[3.9] The equality (3.3.3) states that we identify the basis of the representation $M$ with the set $B_1 \times ... \times B_n$.

[3.10] I considered generating of elements of representation according to the theorem [4]-2.6.4. The theorem 3.3.11 requires the fulfillment of conditions (3.3.4), (3.3.5).

Proof. From the equality (3.3.5), it follows that

$$
\begin{aligned}
f(c) \circ (b_1, ..., f_i(a) \circ b_i, ..., b_n) &= (b_1, ..., f_i(c) \circ (f_i(a) \circ b_i), ..., b_n) \\
&= (b_1, ..., (f_i(c) \circ f_i(a)) \circ b_i, ..., b_n) \\
&= (f(c) \circ f(a)) \circ (b_1, ..., b_i, ..., b_n) \\
&= f(c) \circ (f(a) \circ (b_1, ..., b_i, ..., b_n))
\end{aligned}
\tag{3.3.12}
$$

The equality (3.3.11) follows from the equality (3.3.12). ⊙

Lemma 3.3.8. *For any $c \in A$, endomorphism $f(c)$ of $\Omega_2$-algebra $M$ is coordinated with equivalence $N$.*

Proof. The lemma follows from lemmas 3.3.6, 3.3.7 and from the definition [4]-2.2.14. ⊙

From the lemma 3.3.8 and the theorem [7]-3.2.15, it follows that $\Omega_1$-algebra is defined on the set ${}^*M/N$. Consider diagram

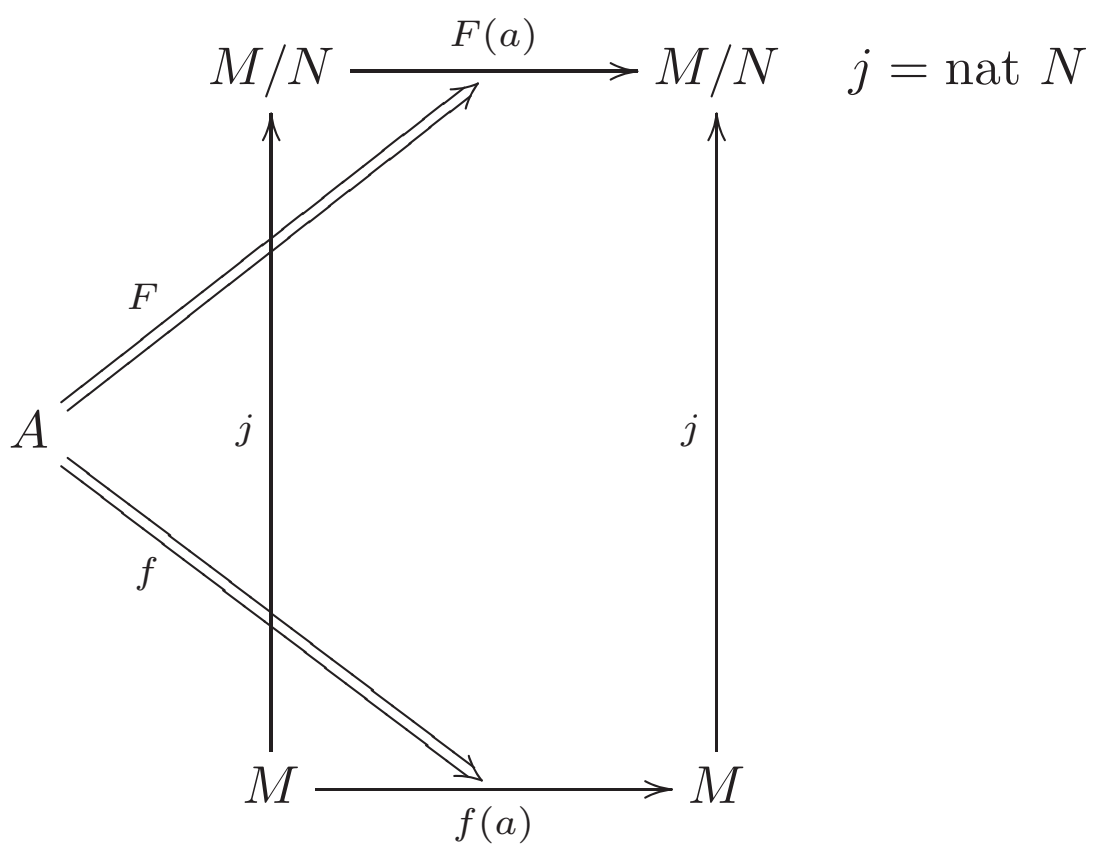

According to lemma 3.3.8, from the condition

$$j \circ b_1 = j \circ b_2$$

it follows that

$$j \circ (f(a) \circ b_1) = j \circ (f(a) \circ b_2)$$

Therefore, transformation $F(a)$ is well defined and

$$F(a) \circ j = j \circ f(a) \tag{3.3.13}$$

If $\omega \in \Omega_1(p)$, then we assume

$$(F(a_1)...F(a_p)\omega) \circ (J \circ b) = J \circ ((f(a_1)...f(a_p)\omega) \circ b)$$

Therefore, map $F$ is representations of $\Omega_1$-algebra $A$. From (3.3.13) it follows that $j$ is reduced morphism of representations $f$ and $F$.

Consider commutative diagram

(3.3.14)

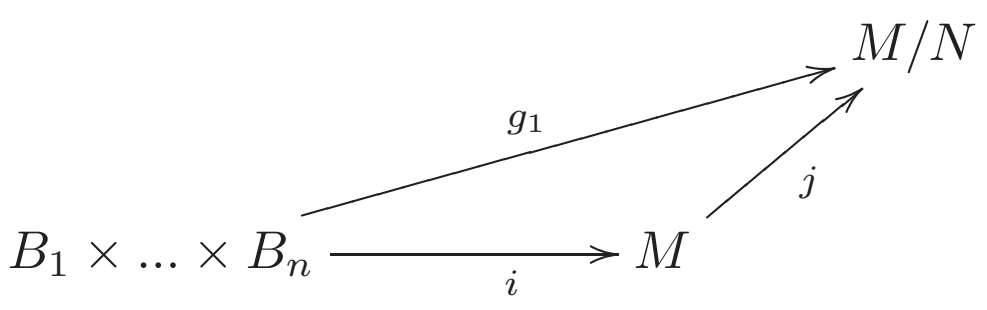

From commutativity of the diagram (3.3.14) and from the equality (3.3.3), it follows that

(3.3.15) $$g_1 \circ (b_1, ..., b_n) = j \circ (b_1, ..., b_n)$$

From equalities (3.3.3), (3.3.4), (3.3.5), it follows that

(3.3.16) $$\begin{aligned} & g_1 \circ (b_1, ..., b_{i\cdot 1} ... b_{i\cdot p}\omega, ..., b_n) \\ = & (g_1 \circ (b_1, ..., b_{i\cdot 1}, ..., b_n)) ... (g_1 \circ (b_1, ..., b_{i\cdot p}, ..., b_n))\omega \end{aligned}$$

(3.3.17) $$g_1 \circ (b_1, ..., f_i(a) \circ b_i, ..., b_n) = f(a) \circ (g_1 \circ (b_1, ..., b_i, ..., b_n))$$

From equalities (3.3.16) and (3.3.17) it follows that map $g_1$ is reduced polymorphism of representations $f_1$, ..., $f_n$.

Since $B_1 \times ... \times B_n$ is the basis of representation $M$ of $\Omega_1$ algebra $A$, then, according to the theorem [4]-2.7.7, for any representation

$$A \xrightarrow{*} V$$

and any reduced polymorphism

$$g_2 : B_1 \times ... \times B_n \longrightarrow V$$

there exists a unique morphism of representations $k : M \to V$, for which following diagram is commutative

(3.3.18)

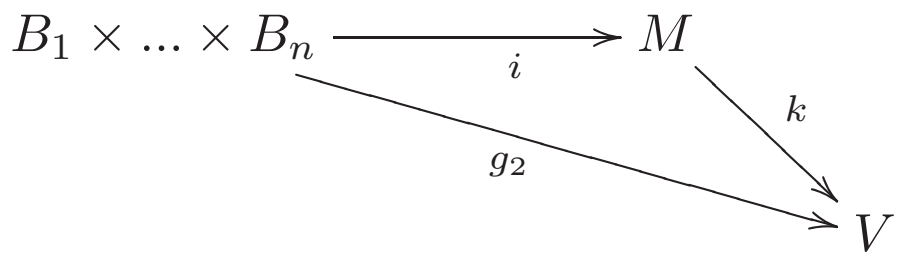

Since $g_2$ is reduced polymorphism, then $\ker k \supseteq N$.

According to the theorem 3.2.3, map $j$ is universal in the category of morphisms of representation $f$ whose kernel contains $N$. Therefore, we have morphism of representations

$$h : M/N \to V$$

which makes the following diagram commutative

(3.3.19)

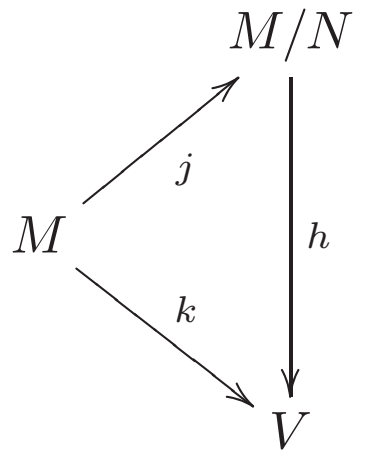

We join diagrams (3.3.14), (3.3.18), (3.3.19), and get commutative diagram

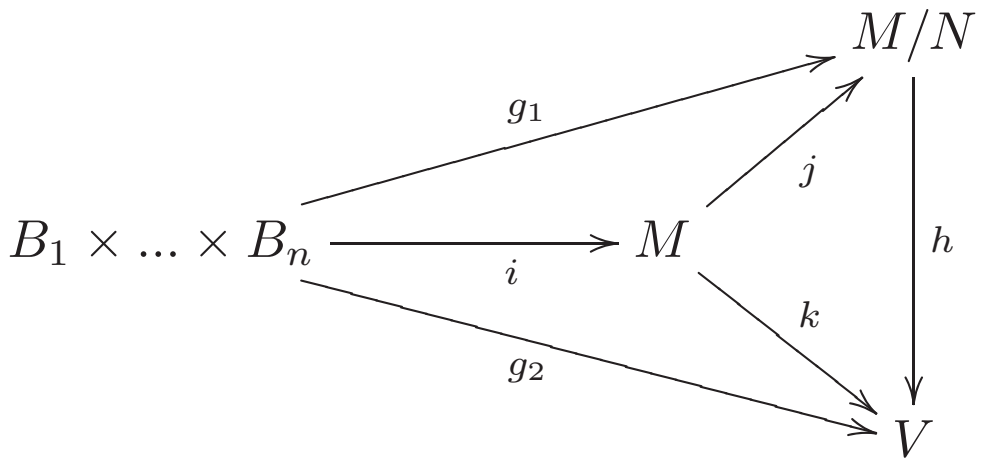

Since $\operatorname{Im} g_1$ generates $M/N$, than map $h$ is uniquely determined. □

According to proof of theorem 3.3.5

$$B_1 \otimes ... \otimes B_n = M/N$$

If $d_i \in A_i$, we write

$$j \circ (d_1, ..., d_n) = d_1 \otimes ... \otimes d_n \tag{3.3.20}$$

From equalities (3.3.15), (3.3.20), it follows that

$$g_1 \circ (d_1, ..., d_n) = d_1 \otimes ... \otimes d_n \tag{3.3.21}$$

THEOREM 3.3.9. *The map*

$$(x_1, ..., x_n) \in B_1 \times ... \times B_n \to x_1 \otimes ... \otimes x_n \in B_1 \otimes ... \otimes B_n$$

*is polymorphism.*

PROOF. The theorem follows from definitions 3.1.3, 3.3.1. □

THEOREM 3.3.10. *Let $B_1$, ..., $B_n$ be $\Omega_2$-algebras. Let*

$$f : B_1 \times ... \times B_n \to B_1 \otimes ... \otimes B_n$$

*be reduced polymorphism defined by equality*

$$f \circ (b_1, ..., b_n) = b_1 \otimes ... \otimes b_n \tag{3.3.22}$$

*Let*

$$g : B_1 \times ... \times B_n \to V$$

*be reduced polymorphism into $\Omega$-algebra $V$. There exists morphism of representations*

$$h : B_1 \otimes ... \otimes B_n \to V$$

*such that the diagram*

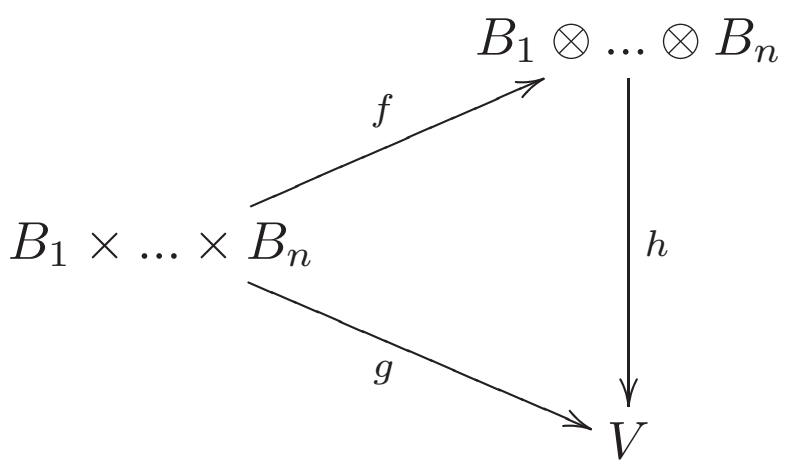

*is commutative.*

PROOF. equality (3.3.22) follows from equalities (3.3.3) and (3.3.20). An existence of the map $h$ follows from the definition 3.3.1 and constructions made in the proof of the theorem 3.3.5. □

THEOREM 3.3.11. *Let*

$$b_k \in B_k \quad k = 1, ..., n \quad b_{i \cdot 1}, ..., b_{i \cdot p} \in B_i \quad \omega \in \Omega_2(p) \quad a \in A$$

*Tensor product is distributive over operation $\omega$*

$$\begin{aligned} & b_1 \otimes ... \otimes (b_{i \cdot 1} ... b_{i \cdot p} \omega) \otimes ... \otimes b_n \\ = & (b_1 \otimes ... \otimes b_{i \cdot 1} \otimes ... \otimes b_n) ... (b_1 \otimes ... \otimes b_{i \cdot p} \otimes ... \otimes b_n)\omega \end{aligned} \tag{3.3.23}$$

*The representation of multiplicative $\Omega_1$-group $A$ in tensor product is defined by equality*

$$b_1 \otimes ... \otimes (f_i(a) \circ b_i) \otimes ... \otimes b_n = f(a) \circ (b_1 \otimes ... \otimes b_i \otimes ... \otimes b_n) \tag{3.3.24}$$

PROOF. The equality (3.3.23) follows from the equality (3.3.16) and from the definition (3.3.21). The equality (3.3.24) follows from the equality (3.3.17) and from the definition (3.3.21). □

## 3.4. Associativity of Tensor Product

Let $A$ be multiplicative $\Omega_1$-group. Let $B_1$, $B_2$, $B_3$ be $\Omega_2$-algebras. Let, for $k = 1, 2, 3$,

$$f_k : A \xrightarrow{*} B_k$$

be effective representation of multiplicative $\Omega_1$-group $A$ in $\Omega_2$-algebra $B_k$.

LEMMA 3.4.1. *For given value of $x_3 \in B_3$, the map*

$$h_{12} : (B_1 \otimes B_2) \times B_3 \to B_1 \otimes B_2 \otimes B_3 \tag{3.4.1}$$

*defined by equality*

$$h_{12}(x_1 \otimes x_2, x_3) = x_1 \otimes x_2 \otimes x_3 \tag{3.4.2}$$

*is reduced morphism of the representation $B_1 \otimes B_2$ into the representation $B_1 \otimes B_2 \otimes B_3$ .*

PROOF. According to the theorem 3.3.9, for given value of $x_3 \in B_3$, the map

$$(x_1, x_2, x_3) \in B_1 \times B_2 \times B_3 \to x_1 \otimes x_2 \otimes x_3 \in B_1 \otimes B_2 \otimes B_3 \tag{3.4.3}$$

is polymorphism with respect to $x_1 \in B_1$, $x_2 \in B_2$. Therefore, for given value of $x_3 \in B_3$, the lemma follows from the theorem 3.3.10. □

LEMMA 3.4.2. *For given value of $x_{12} \in B_1 \otimes B_2$ the map $h_{12}$ is reduced morphism of the representation $B_3$ into the representation $B_1 \otimes B_2 \otimes B_3$ .*

PROOF. According to the theorem 3.3.9 and the equality (3.3.21), for given value of $x_1 \in B_1$, $x_2 \in B_2$, the map

$$(x_1 \otimes x_2, x_3) \in B_1 \times B_2 \times B_3 \to x_1 \otimes x_2 \otimes x_3 \in B_1 \otimes B_2 \otimes B_3 \tag{3.4.4}$$

is morphism with respect to $x_3 \in B_3$. Therefore, the theorem follows from theorems 3.1.10, 3.1.11. □

LEMMA 3.4.3. *There exists reduced morphism of representations*

$$h : (B_1 \otimes B_2) \otimes B_3 \to B_1 \otimes B_2 \otimes B_3$$

PROOF. According to lemmas 3.4.1, 3.4.2 and to the definition 3.1.3, the map $h_{12}$ is reduced polymorphism of representations. The lemma follows from the theorem 3.3.10. □

LEMMA 3.4.4. *There exists reduced morphism of representations*

$$g : B_1 \otimes B_2 \otimes B_3 \to (B_1 \otimes B_2) \otimes B_3$$

Proof. The map

$$(x_1, x_2, x_3) \in B_1 \times B_2 \times B_3 \to (x_1 \otimes x_2) \otimes x_3 \in (B_1 \otimes B_2) \otimes B_3$$

is polymorphism with respect to $x_1 \in B_1$, $x_2 \in B_2$, $x_3 \in B_3$. Therefore, the lemma follows from the theorem 3.3.10. □

Theorem 3.4.5.

$$(B_1 \otimes B_2) \otimes B_3 = B_1 \otimes (B_2 \otimes B_3) = B_1 \otimes B_2 \otimes B_3 \tag{3.4.5}$$

Proof. According to lemma 3.4.3, there exists reduced morphism of representations

$$h : (B_1 \otimes B_2) \otimes B_3 \to B_1 \otimes B_2 \otimes B_3$$

According to lemma 3.4.4, there exists reduced morphism of representations

$$g : B_1 \otimes B_2 \otimes B_3 \to (B_1 \otimes B_2) \otimes B_3$$

Therefore, reduced morphisms of representations $h$, $g$ are isomorphisms. Therefore, the following equality is true

$$(B_1 \otimes B_2) \otimes B_3 = B_1 \otimes B_2 \otimes B_3 \tag{3.4.6}$$

We prove similarly the equality

$$B_1 \otimes (B_2 \otimes B_3) = B_1 \otimes B_2 \otimes B_3$$

□

Remark 3.4.6. *It is evident that structures of $\Omega_2$-algebras $(B_1 \otimes B_2) \otimes B_3$, $B_1 \otimes B_2 \otimes B_3$ are little different. We write down the equality* (3.4.6) *based on the convention 3.3.3 and this allows us to speak about associativity of tensor product of representations.* □

CHAPTER 4

# $D$-Module

## 4.1. Module over Commutative Ring

THEOREM 4.1.1. *Let ring $D$ has unit $e$. Representation*

$$f : D \xrightarrow{*} V$$

*of the ring $D$ in an Abelian group $A$ is* **effective** *iff $a = 0$ follows from equality $f(a) = 0$.*

PROOF. We define the sum of transformations $f$ and $g$ of an Abelian group according to rule

$$(f + g)(a) = f(a) + g(a)$$

Therefore, considering the representation of the ring $D$ in the Abelian group $A$, we assume

$$f(a + b)(x) = f(a)(x) + f(b)(x)$$

Suppose $a, b \in R$ cause the same transformation. Then

$$f(a) \circ m = f(b) \circ m \tag{4.1.1}$$

for any $m \in A$. From the equality (4.1.1) it follows that $a - b$ generates zero transformation

$$f(a - b) \circ m = 0$$

Element $e+a-b$ generates an identity transformation. Therefore, the representation $f$ is effective iff $a = b$. □

DEFINITION 4.1.2. *Effective representation of commutative ring $D$ in an Abelian group $V$*

$$f : D \xrightarrow{*} V \qquad f(d) : v \to d\,v \tag{4.1.2}$$

*is called* **module over ring $D$** *or* **$D$-module**. □

THEOREM 4.1.3. *Following conditions hold for $D$-module:*

- **associative law**

$$(ab)m = a(bm) \tag{4.1.3}$$

- **distributive law**

$$a(m + n) = am + an \tag{4.1.4}$$

$$(a + b)m = am + bm \tag{4.1.5}$$

- **unitarity law**

$$1m = m \tag{4.1.6}$$

*for any $a, b \in D$, $m, n \in V$.*

PROOF. Since transformation $a$ is endomorphism of the Abelian group, we obtain the equality (4.1.4). Since representation (4.1.2) is homomorphism of the additive group of ring $D$, we obtain the equality (4.1.5). Since representation (4.1.2) is left-side representation of the multiplicative group of ring $D$, we obtain equalities (4.1.3) and (4.1.6). □

THEOREM 4.1.4. *The set of vectors generated by the set of vectors* $v=(v_{\boldsymbol{i}}\in V,\boldsymbol{i}\in\boldsymbol{I})$ *has form*

$$J(v)=\left\{w:w=\sum_{\boldsymbol{i}\in\boldsymbol{I}}c^{\boldsymbol{i}}v_{\boldsymbol{i}},c^{\boldsymbol{i}}\in D\right\}\tag{4.1.7}$$

PROOF. We prove the theorem by induction based on the theorem [4]-2.6.4.

Let $k=0$. According to the theorem [4]-2.6.4, $X_0=v$. For any $v_{\boldsymbol{k}}\in v$, let $c^{\boldsymbol{i}}=\delta^{\boldsymbol{i}}_{\boldsymbol{k}}$. Then

$$v_{\boldsymbol{k}}=\sum_{\boldsymbol{i}\in\boldsymbol{I}}c^{\boldsymbol{i}}v_{\boldsymbol{i}}\tag{4.1.8}$$

$v_{\boldsymbol{k}}\in J(v)$ follows from (4.1.7), (4.1.8).

Let $X_{k-1}\subseteq J(v)$.

- Let $w_1, w_2\in X_{k-1}$. Since $V$ is Abelian group, then according to the statement [4]-2.6.4.3, $w_1+w_2\in X_k$. According to statements [4]-(2.6.1), (4.1.7), there exist $D$-numbers $c^{\boldsymbol{i}}$, $d^{\boldsymbol{i}}$, $\boldsymbol{i}\in\boldsymbol{I}$, such that

$$\begin{aligned}w_1&=\sum_{\boldsymbol{i}\in\boldsymbol{I}}c^{\boldsymbol{i}}v_{\boldsymbol{i}}\\ w_2&=\sum_{\boldsymbol{i}\in\boldsymbol{I}}d^{\boldsymbol{i}}v_{\boldsymbol{i}}\end{aligned}\tag{4.1.9}$$

  Since $V$ is Abelian group, then from the equality (4.1.9) it follows that

$$w_1+w_2=\sum_{\boldsymbol{i}\in\boldsymbol{I}}c^{\boldsymbol{i}}v_{\boldsymbol{i}}+\sum_{\boldsymbol{i}\in\boldsymbol{I}}d^{\boldsymbol{i}}v_{\boldsymbol{i}}=\sum_{\boldsymbol{i}\in\boldsymbol{I}}(c^{\boldsymbol{i}}v_{\boldsymbol{i}}+d^{\boldsymbol{i}}v_{\boldsymbol{i}})\tag{4.1.10}$$

  The equality

$$w_1+w_2=\sum_{\boldsymbol{i}\in\boldsymbol{I}}(c^{\boldsymbol{i}}+d^{\boldsymbol{i}})v_{\boldsymbol{i}}\tag{4.1.11}$$

  follows from equalities (4.1.5), (4.1.10). From the equality (4.1.11), it follows that $w_1+w_2\in J(v)$.
- Let $w\in X_{k-1}$. According to the statement [4]-2.6.4.4, for any $D$-number $a$, $aw\in X_k$. According to statements [4]-(2.6.1), (4.1.7), there exist $D$-numbers $c^{\boldsymbol{i}}$, $\boldsymbol{i}\in\boldsymbol{I}$, such that

$$w=\sum_{\boldsymbol{i}\in\boldsymbol{I}}c^{\boldsymbol{i}}v_{\boldsymbol{i}}\tag{4.1.12}$$

  From the equality (4.1.12) it follows that

$$aw=a\sum_{\boldsymbol{i}\in\boldsymbol{I}}c^{\boldsymbol{i}}v_{\boldsymbol{i}}=\sum_{\boldsymbol{i}\in\boldsymbol{I}}a(c^{\boldsymbol{i}}v_{\boldsymbol{i}})=\sum_{\boldsymbol{i}\in\boldsymbol{I}}(ac^{\boldsymbol{i}})v_{\boldsymbol{i}}\tag{4.1.13}$$

  From the equality (4.1.13), it follows that $aw\in J(v)$.

□

Convention 4.1.5. *We will use Einstein summation convention in which repeated index (one above and one below) implies summation with respect to repeated index. In this case we assume that we know the set of summation index and do not use summation symbol*

$$c^i v_i = \sum_{i\in I} c^i v_i$$

□

Definition 4.1.6. *Let* $v = (v_i \in V, i \in I)$ *be set of vectors. The expression* $c^i v_i$ *is called* **linear composition of vectors** $v_i$*. A vector*

$$w = c^i v_i$$

*is called* **vector linearly dependent on vectors** $v_i$*.* □

Theorem 4.1.7. *Let* $D$ *be field. Since equality*

$$c^i v_i = 0 \tag{4.1.14}$$

*implies existence of index* $i = j$ *such that* $c^j \neq 0$*, then the vector* $v_j$ *linearly depends on rest of vectors* $v$*.*

Proof. The theorem follows from the equality

$$v_j = \sum_{i\in I\setminus\{j\}} \frac{c^i}{c^j} v_i$$

and from the definition 4.1.6. □

It is evident that for any set of vectors $v_i$

$$0 = c^i v_i \quad c^i = 0$$

Definition 4.1.8. *Vectors*[4.1] $e_i$, $i \in I$, *of* $D$*-module* $A$ *are* **linearly independent** *if* $c = 0$ *follows from the equality*

$$c^i e_i = 0$$

*Otherwise vectors* $e_i$, $i \in I$, *are* **linearly dependent**. □

Theorem 4.1.9. *Let* $D$ *be field. A set of vectors* $\overline{\overline{e}} = (e_i, i \in I)$ *is a* **basis of** $D$**-vector space** $V$*, if vectors* $e_i$ *are linearly independent and any vector* $v \in V$ *linearly depends on vectors* $e_i$*.*

Proof. Let $\overline{\overline{e}} = (e_i, i \in I)$ be basis of $D$-vector space $V$. According to the definition [4]-2.7.1 and to theorems [4]-2.6.4, 4.1.4, any vector $v \in V$ is linear composition of vectors $e_i$

$$v = v^i e_i \tag{4.1.15}$$

From the equality (4.1.15), it follows that the set of vectors $v$, $e_i$, $i \in I$, is not linearly independent.

Consider equality

$$c^i e_i = 0 \tag{4.1.16}$$

According to the theorem 4.1.7, since

$$c^j \neq 0 \tag{4.1.17}$$

[4.1] I follow to the definition in [1], p. 130.

then the vector $e_j$ linearly depends on rest of vectors $e$. Therefore the set of vectors $e_i$, $i \in I \setminus \{j\}$, generates $D$-vector space $V$. According to the definition [4]-2.7.1, the statement (4.1.17) is not true. According to the definition 4.1.8, vectors $e_i$ are linearly independent. □

Definition 4.1.10. *A set of vectors* $\overline{\overline{e}} = (e_i, i \in I)$ *is called*[4.2] *a* **basis of $D$-module** $V$*, if arbitrary vector* $v \in V$ *is linear combination of vectors of the basis and arbitrary vector of the basis cannot be represented as a linear combination of the remaining vectors of basis.* $A$ *is* **free module over ring** $D$*, if* $A$ *has basis over ring* $D$*.*[4.3] □

Definition 4.1.11. *Let* $\overline{\overline{e}}$ *be the basis of* $D$*-module* $A$ *and* $A$*-number* $a$ *has expansion*

$$a = a^i e_i$$

*with respect to the basis* $\overline{\overline{e}}$*.* $D$*-numbers* $a^i$ *are called* **coordinates** *of* $A$*-number* $a$ *with respect to the basis* $\overline{\overline{e}}$*.* □

## 4.2. Linear Map of $D$-Module

Definition 4.2.1. *Let* $A_1$*,* $A_2$ *be* $D$*-modules. Morphism of representations*

$$f : A_1 \to A_2$$

*of* $D$*-module* $A_1$ *into* $D$*-module* $A_2$ *is called* **linear map** *of* $D$*-module* $A_1$ *into* $D$*-module* $A_2$*. Let us denote* $\mathcal{L}(D; A_1 \to A_2)$ *set of linear maps of* $D$*-module* $A_1$ *into* $D$*-module* $A_2$*.* □

Theorem 4.2.2. *Linear map*

$$f : A_1 \to A_2$$

*of* $D$*-module* $A_1$ *into* $D$*-module* $A_2$ *satisfies to equalities*[4.4]

$$f \circ (a + b) = f \circ a + f \circ b \tag{4.2.1}$$

$$f \circ (pa) = p(f \circ a) \tag{4.2.2}$$

$$a, b \in A_1 \quad p \in D$$

Proof. From definition 4.2.1 and theorem [7]-3.3.1 it follows that the map $f$ is a homomorphism of the Abelian group $A_1$ into the Abelian group $A_2$ (the equality (4.2.1)). The equality (4.2.2) follows from the equality [7]-(3.3.3). □

Theorem 4.2.3. *Let* $A_1$*,* $A_2$ *be* $D$*-modules. The map*

$$f + g : A_1 \to A_2 \quad f, g \in \mathcal{L}(D; A_1 \to A_2) \tag{4.2.3}$$

*defined by equation*

$$(f + g) \circ x = f \circ x + g \circ x \tag{4.2.4}$$

*is called* **sum of maps** $f$ *and* $g$ *and is linear map. The set* $\mathcal{L}(D; A_1; A_2)$ *is an Abelian group relative sum of maps.*

---

[4.2] The definition 4.1.10 follows from the theorem [4]-2.7.2 and the remark [4]-2.7.3.

[4.3] I follow to the definition in [1], p. 135.

[4.4] In some books (for instance, [1], p. 119) the theorem 4.2.2 is considered as a definition.

Proof. According to the theorem 4.2.2

$$f\circ(v+w)=f\circ v+f\circ w \tag{4.2.5}$$

$$f\circ(d\,v)=d\,(f\circ v) \tag{4.2.6}$$

$$g\circ(v+w)=g\circ v+g\circ w \tag{4.2.7}$$

$$g\circ(d\,v)=d\,(g\circ v) \tag{4.2.8}$$

The equality

$$\begin{aligned}(f+g)\circ(v+w)&=f\circ(v+w)+g\circ(v+w)\\&=f\circ v+f\circ w+g\circ v+g\circ w\\&=(f+g)\circ v+(f+g)\circ w\end{aligned} \tag{4.2.9}$$

follows from the equalities (4.2.4), (4.2.5), (4.2.7). The equality

$$\begin{aligned}(f+g)\circ(d\,v)&=f\circ(d\,v)+g\circ(d\,v)\\&=d\,f\circ v+d\,g\circ v\\&=d\,(f\circ v+g\circ v)\\&=d\,((f+g)\circ v)\end{aligned} \tag{4.2.10}$$

follows from the equalities (4.2.4), (4.2.6), (4.2.8). From equalities (4.2.9), (4.2.10) and from the theorem 4.2.2, it follows that the map (4.2.3) is linear map of $D$-modules.

From the equality (4.2.4), it follows that the map

$$0:v\in A_1\to 0\in A_2$$

is zero of addition

$$(0+f)\circ v=0\circ v+f\circ v=f\circ v$$

From the equality (4.2.4), it follows that the map

$$-f:v\in A_1\to-(f\circ v)\in A_2$$

is map inversed to map $f$

$$(f+(-f))\circ v=f\circ v+(-f)\circ v=f\circ v-f\circ v=0=0\circ v$$

From the equality

$$\begin{aligned}(f+g)\circ x&=f\circ x+g\circ x\\&=g\circ x+f\circ x\\&=(g+f)\circ x\end{aligned}$$

it follows that sum of maps is commutative. □

Theorem 4.2.4. *Let $A_1$, $A_2$ be $D$-modules. The map*

$$d\,f:A_1\to A_2\quad d\in D\quad f\in\mathcal{L}(D;A_1\to A_2) \tag{4.2.11}$$

*defined by equality*

$$(d\,f)\circ x=d(f\circ x) \tag{4.2.12}$$

*is linear map and is called* **product of map $f$ over scalar** $d$*. The representation*

$$a:f\in\mathcal{L}(D;A_1\to A_2)\to af\in\mathcal{L}(D;A_1\to A_2) \tag{4.2.13}$$

*of ring $D$ in Abelian group* $\mathcal{L}(D;A_1\to A_2)$ *generates structure of $D$-module.*

Proof. According to the theorem 4.2.2

$$f\circ(v+w)=f\circ v+f\circ w \tag{4.2.14}$$

$$f\circ(d\,v)=d\,(f\circ v) \tag{4.2.15}$$

The equality

$$\begin{aligned}(d\,f)\circ(v+w)&=d(f\circ(v+w))\\&=d(f\circ v+f\circ w)=d(f\circ v)+d(f\circ w)\\&=(d\,f)\circ v+(d\,f)\circ w\end{aligned} \tag{4.2.16}$$

follows from the equalities (4.2.12), (4.2.14). The equality

$$\begin{aligned}(c\,f)\circ(d\,v)&=c(f\circ(d\,v))=cd\,(f\circ v)=d(c(f\circ v))\\&=d\,((c\,f)\circ v)\end{aligned} \tag{4.2.17}$$

follows from the equalities (4.2.12), (4.2.15). From equalities (4.2.16), (4.2.17) and from the theorem 4.2.2, it follows that the map (4.2.11) is linear map of $D$-modules.

The equality

$$(p+q)f=pf+qf \tag{4.2.18}$$

follows from the equality

$$\begin{aligned}((a+b)f)\circ v&=(a+b)(f\circ v)=a(f\circ v)+b(f\circ v)=(af)\circ v+(bf)\circ v\\&=(af+bf)\circ v\end{aligned}$$

The equality

$$p(qf)=(pq)f \tag{4.2.19}$$

follows from the equality

$$((ab)f)\circ v=(ab)(f\circ v)=a(b(f\circ v))=a((bf)\circ v)=(a(bf))\circ v$$

From equalities (4.2.18), (4.2.19), it follows that the map (4.2.13) is representation of ring $D$ in Abelian group $\mathcal{L}(D;A_1\to A_2)$. According to the definition 4.1.2 and the theorem 4.2.3, Abelian group $\mathcal{L}(D;A_1\to A_2)$ is $D$-module. □

Theorem 4.2.5. *Let map*

$$f:A_1\to A_2$$

*be linear map of $D$-module $A_1$ into $D$-module $A_2$. Then*

$$f\circ 0=0$$

Proof. The theorem follows from the equality

$$f\circ(a+0)=f\circ a+f\circ 0$$

□

## 4.3. Polylinear Map of $D$-Module

Definition 4.3.1. *Let $D$ be the commutative ring. Reduced polymorphism of $D$-modules $A_1$, ..., $A_n$ into $D$-module $S$*

$$f : A_1 \times ... \times A_n \to S$$

*is called* **polylinear map** *of $D$-modules $A_1$, ..., $A_n$ into $D$-module $S$. Let us denote $\mathcal{L}(D; A_1 \times ... \times A_n \to S)$ set of polylinear maps of $D$-modules $A_1$, ..., $A_n$ into $D$-module $S$. Let us denote $\mathcal{L}(D; A^n \to S)$ set of $n$-linear maps of $D$-module $A$ $(A_1 = ... = A_n = A)$ into $D$-module $S$.* □

Theorem 4.3.2. *Let $D$ be the commutative ring. The polylinear map of $D$-modules $A_1$, ..., $A_n$ into $D$-module $S$*

$$f : A_1 \times ... \times A_n \to S$$

*satisfies to equalities*

$$f \circ (a_1, ..., a_i + b_i, ..., a_n) = f \circ (a_1, ..., a_i, ..., a_n) + f \circ (a_1, ..., b_i, ..., a_n)$$
$$f \circ (a_1, ..., pa_i, ..., a_n) = pf \circ (a_1, ..., a_i, ..., a_n)$$

$$f \circ (a_1, ..., a_i + b_i, ..., a_n) = f \circ (a_1, ..., a_i, ..., a_n) + f \circ (a_1, ..., b_i, ..., a_n)$$
$$f \circ (a_1, ..., pa_i, ..., a_n) = pf \circ (a_1, ..., a_i, ..., a_n)$$

$$1 \le i \le n \quad a_i, b_i \in A_i \quad p \in D$$

Proof. The theorem follows from definitions 3.1.1, 4.2.1 and from the theorem 4.2.2. □

Theorem 4.3.3. *Let $D$ be the commutative ring. Let $A_1$, ..., $A_n$, $S$ be $D$-modules. The map*

$$f + g : A_1 \times ... \times A_n \to S \quad f, g \in \mathcal{L}(D; A_1 \times ... \times A_n \to S) \tag{4.3.1}$$

*defined by equation*

$$(f + g) \circ (a_1, ..., a_n) = f \circ (a_1, ..., a_n) + g \circ (a_1, ..., a_n) \tag{4.3.2}$$

*is called* **sum of polylinear maps** *$f$ and $g$ and is polylinear map. The set $\mathcal{L}(D; A_1 \times ... \times A_n \to S)$ is an Abelian group relative sum of maps.*

Proof. According to the theorem 4.3.2

$$f \circ (a_1, ..., a_i + b_i, ..., a_n) = f \circ (a_1, ..., a_i, ..., a_n) + f \circ (a_1, ..., b_i, ..., a_n) \tag{4.3.3}$$

$$f \circ (a_1, ..., pa_i, ..., a_n) = pf \circ (a_1, ..., a_i, ..., a_n) \tag{4.3.4}$$

$$g \circ (a_1, ..., a_i + b_i, ..., a_n) = g \circ (a_1, ..., a_i, ..., a_n) + g \circ (a_1, ..., b_i, ..., a_n) \tag{4.3.5}$$

$$g \circ (a_1, ..., pa_i, ..., a_n) = pg \circ (a_1, ..., a_i, ..., a_n) \tag{4.3.6}$$

The equality

$$\begin{aligned}&(f+g)\circ(x_1,...,x_i+y_i,...,x_n)\\&=f\circ(x_1,...,x_i+y_i,...,x_n)+g\circ(x_1,...,x_i+y_i,...,x_n)\\&=f\circ(x_1,...,x_i,...,x_n)+f\circ(x_1,...,y_i,...,x_n)\\&+g\circ(x_1,...,x_i,...,x_n)+g\circ(x_1,...,y_i,...,x_n)\\&=(f+g)\circ(x_1,...,x_i,...,x_n)+(f+g)\circ(x_1,...,y_i,...,x_n)\end{aligned}\tag{4.3.7}$$

follows from the equalities (4.3.2), (4.3.3), (4.3.5). The equality

$$\begin{aligned}&(f+g)\circ(x_1,...,px_i,...,x_n)\\&=f\circ(x_1,...,px_i,...,x_n)+g\circ(x_1,...,px_i,...,x_n)\\&=pf\circ(x_1,...,x_i,...,x_n)+pg\circ(x_1,...,x_i,...,x_n)\\&=p(f\circ(x_1,...,x_i,...,x_n)+g\circ(x_1,...,x_i,...,x_n))\\&=p(f+g)\circ(x_1,...,x_i,...,x_n)\end{aligned}\tag{4.3.8}$$

follows from the equalities (4.3.2), (4.3.4), (4.3.6). From equalities (4.3.7), (4.3.8) and from the theorem 4.3.2, it follows that the map (4.3.1) is linear map of $D$-modules.

Let $f, g, h\in\mathcal{L}(D;A_1\times...\times A_2\to S)$. For any $a=(a_1,...,a_n)$, $a_1\in A_1$, ..., $a_n\in A_n$,

$$\begin{aligned}(f+g)\circ a&=f\circ a+g\circ a=g\circ a+f\circ a\\&=(g+f)\circ a\\((f+g)+h)\circ a&=(f+g)\circ a+h\circ a=(f\circ a+g\circ a)+h\circ a\\&=f\circ a+(g\circ a+h\circ a)=f\circ a+(g+h)\circ a\\&=(f+(g+h))\circ a\end{aligned}$$

Therefore, sum of polylinear maps is commutative and associative.

From the equality (4.3.2), it follows that the map

$$0:v\in A_1\times...\times A_n\to 0\in S$$

is zero of addition

$$(0+f)\circ(a_1,...,a_n)=0\circ(a_1,...,a_n)+f\circ(a_1,...,a_n)=f\circ(a_1,...,a_n)$$

From the equality (4.3.2), it follows that the map

$$-f:(a_1,...,a_n)\in A_1\times...\times A_n\to-(f\circ(a_1,...,a_n))\in S$$

is map inversed to map $f$

$$f+(-f)=0$$

because

$$\begin{aligned}(f+(-f))\circ(a_1,...,a_n)&=f\circ(a_1,...,a_n)+(-f)\circ(a_1,...,a_n)\\&=f\circ(a_1,...,a_n)-f\circ(a_1,...,a_n)\\&=0=0\circ(a_1,...,a_n)\end{aligned}$$

From the equality

$$\begin{aligned}(f+g)\circ(a_1,...,a_n)&=f\circ(a_1,...,a_n)+g\circ(a_1,...,a_n)\\&=g\circ(a_1,...,a_n)+f\circ(a_1,...,a_n)\\&=(g+f)\circ(a_1,...,a_n)\end{aligned}$$

it follows that sum of maps is commutative. Therefore, the set $\mathcal{L}(D;A_1\times...\times A_n\rightarrow S)$ is an Abelian group. □

THEOREM 4.3.4. *Let $D$ be the commutative ring. Let $A_1$, ..., $A_n$, $S$ be $D$-modules. The map*

$$df:A_1\times...\times A_n\rightarrow S\quad d\in D\quad f\in\mathcal{L}(D;A_1\times...\times A_n\rightarrow S)\tag{4.3.9}$$

*defined by equality*

$$(df)\circ(a_1,...,a_n)=d(f\circ(a_1,...,a_n))\tag{4.3.10}$$

*is polylinear map and is called* **product of map $f$ over scalar** $d$. *The representation*

$$a:f\in\mathcal{L}(D;A_1\times...\times A_n\rightarrow S)\rightarrow af\in\mathcal{L}(D;A_1\times...\times A_n\rightarrow S)\tag{4.3.11}$$

*of ring $D$ in Abelian group $\mathcal{L}(D;A_1\times...\times A_n\rightarrow S)$ generates structure of $D$-module.*

PROOF. According to the theorem 4.3.2

$$f\circ(a_1,...,a_i+b_i,...,a_n)=f\circ(a_1,...,a_i,...,a_n)+f\circ(a_1,...,b_i,...,a_n)\tag{4.3.12}$$

$$f\circ(a_1,...,pa_i,...,a_n)=pf\circ(a_1,...,a_i,...,a_n)\tag{4.3.13}$$

The equality

$$\begin{aligned}&(pf)\circ(x_1,...,x_i+y_i,...,x_n)\\&=p\ f\circ(x_1,...,x_i+y_i,...,x_n)\\&=p\ (f\circ(x_1,...,x_i,...,x_n)+f\circ(x_1,...,y_i,...,x_n))\\&=p(f\circ(x_1,...,x_i,...,x_n))+p(f\circ(x_1,...,y_i,...,x_n))\\&=(pf)\circ(x_1,...,x_i,...,x_n)+(pf)\circ(x_1,...,y_i,...,x_n)\end{aligned}\tag{4.3.14}$$

follows from equalities (4.3.10), (4.3.12). The equality

$$\begin{aligned}&(pf)\circ(x_1,...,qx_i,...,x_n)\\&=p(f\circ(x_1,...,qx_i,...,x_n))=pq(f\circ(x_1,...,x_i,...,x_n))\\&=qp(f\circ(x_1,...,x_n))=q(pf)\circ(x_1,...,x_n)\end{aligned}\tag{4.3.15}$$

follows from equalities (4.3.10), (4.3.13). From equalities (4.3.14), (4.3.15) and from the theorem 4.3.2, it follows that the map (4.3.9) is polylinear map of $D$-modules.

The equality

$$(p+q)f=pf+qf\tag{4.3.16}$$

follows from the equality

$$\begin{aligned}((p+q)f)\circ(x_1,...,x_n)&=(p+q)(f\circ(x_1,...,x_n))\\&=p(f\circ(x_1,...,x_n))+q(f\circ(x_1,...,x_n))\\&=(pf)\circ(x_1,...,x_n)+(qf)\circ(x_1,...,x_n)\end{aligned}$$

The equality

$$p(qf)=(pq)f \tag{4.3.17}$$

follows from the equality

$$\begin{aligned}(p(qf))\circ(x_1,...,x_n)=&p\ (qf)\circ(x_1,...,x_n)=p\ (q\ f\circ(x_1,...,x_n))\\ =&(pq)\ f\circ(x_1,...,x_n)=((pq)f)\circ(x_1,...,x_n)\end{aligned}$$

From equalities (4.3.16) (4.3.17) it follows that the map (4.3.11) is representation of ring $D$ in Abelian group $\mathcal{L}(D;A_1\times...\times A_n\rightarrow S)$ . Since specified representation is effective, then, according to the definition 4.1.2 and the theorem 4.3.3, Abelian group $\mathcal{L}(D;A_1\rightarrow A_2)$ is $D$-module. □

## 4.4. $D$-module $\mathcal{L}(D;A\rightarrow B)$

THEOREM 4.4.1.

$$\mathcal{L}(D;A^p\rightarrow\mathcal{L}(D;A^q\rightarrow B))=\mathcal{L}(D;A^{p+q}\rightarrow B) \tag{4.4.1}$$

PROOF. □

THEOREM 4.4.2. *Let*

$$\overline{\overline{e}}=\{e_i:i\in I\}$$

*be basis of $D$-module $A$. The set*

$$\overline{\overline{h}}=\{h^i\in\mathcal{L}(D;A\rightarrow D):i\in I,h^i\circ e_j=\delta^i_j\} \tag{4.4.2}$$

*is the basis of $D$-module* $\mathcal{L}(D;A\rightarrow D)$.

PROOF.

LEMMA 4.4.3. *Maps $h^i$ are linear independent.*
PROOF. Let there exist $D$-numbers $c_i$ such that

$$c_ih^i=0$$

Then for any $A$-number $e_j$

$$0=c_ih^i\circ e_j=c_i\delta^i_j=c_j$$

The lemma follows from the definition 4.1.8. ⊙

LEMMA 4.4.4. *The map* $f\in\mathcal{L}(D;A\rightarrow D)$, *is linear composition of maps $h^i$.*

PROOF. For any $A$-number $a$,

$$a=a^ie_i \tag{4.4.3}$$

the equality

$$h^i\circ a=h^i\circ(a^je_j)=a^j(h^i\circ e_j)=a^j\delta^i_j=a^i \tag{4.4.4}$$

follows from equalities (4.4.2), (4.4.3) and from the theorem 4.2.2. The equality

$$f\circ a=f\circ(a^ie_i)=a^i(f\circ e_i)=(f\circ e_i)(h^i\circ a) \tag{4.4.5}$$

follows from the equality (4.4.4). The equality

$$f=(f\circ e_i)h^i$$

follows from equalities (4.2.4), (4.2.12), (4.4.5). ⊙

The theorem follows from lemmas 4.4.3, 4.4.4 and from the definition 4.1.10. □

THEOREM 4.4.5. *Let $D$ be commutative ring. Let*

$$\overline{\overline{e}}_i=\{e_{i\cdot\boldsymbol{i}}:\boldsymbol{i}\in\boldsymbol{I}_i\}$$

*be basis of $D$-module $A_i$, $i=1$, ..., $n$. Let*

$$\overline{\overline{e}}_B=\{e_{B\cdot\boldsymbol{i}}:\boldsymbol{i}\in\boldsymbol{I}\}$$

*be basis of $D$-module $B$. The set*

$$\begin{aligned}\overline{\overline{h}}=\{\,&h^{\boldsymbol{i}_1...\boldsymbol{i}_n}_{\boldsymbol{i}}\in\mathcal{L}(D;A_1\times...\times A_n\to B):\boldsymbol{i}\in\boldsymbol{I},\boldsymbol{i}_i\in\boldsymbol{I}_i,i=1,...,n,\\&h^{\boldsymbol{i}_1...\boldsymbol{i}_n}_{\boldsymbol{i}}\circ(e_{1\cdot\boldsymbol{j}_1},...,e_{n\cdot\boldsymbol{j}_n})=\delta^{\boldsymbol{i}_1}_{\boldsymbol{j}_1}...\delta^{\boldsymbol{i}_n}_{\boldsymbol{j}_n}e_{B\cdot\boldsymbol{i}}\}\end{aligned}\tag{4.4.6}$$

*is the basis of $D$-module $\mathcal{L}(D;A_1\times...\times A_n\to B)$.*

PROOF.

LEMMA 4.4.6. *Maps $h^{\boldsymbol{i}_1...\boldsymbol{i}_n}_{\boldsymbol{i}}$ are linear independent.*

PROOF. Let there exist $D$-numbers $c^{\boldsymbol{i}}_{\boldsymbol{i}_1...\boldsymbol{i}_n}$ such that

$$c^{\boldsymbol{i}}_{\boldsymbol{i}_1...\boldsymbol{i}_n}h^{\boldsymbol{i}_1...\boldsymbol{i}_n}_{\boldsymbol{i}}=0$$

Then for any set of indices $\boldsymbol{j}_1$, ..., $\boldsymbol{j}_n$

$$0=c^{\boldsymbol{i}}_{\boldsymbol{i}_1...\boldsymbol{i}_n}h^{\boldsymbol{i}_1...\boldsymbol{i}_n}_{\boldsymbol{i}}\circ(e_{1\cdot\boldsymbol{j}_1},...,e_{n\cdot\boldsymbol{j}_n})=c^{\boldsymbol{i}}_{\boldsymbol{i}_1...\boldsymbol{i}_n}\delta^{\boldsymbol{i}_1}_{\boldsymbol{j}_1}...\delta^{\boldsymbol{i}_n}_{\boldsymbol{j}_n}e_{B\cdot\boldsymbol{i}}=c^{\boldsymbol{i}}_{\boldsymbol{j}_1...\boldsymbol{j}_n}e_{B\cdot\boldsymbol{i}}$$

Therefore, $c^{\boldsymbol{i}}_{\boldsymbol{j}_1...\boldsymbol{j}_n}=0$. The lemma follows from the definition 4.1.8. ⊙

LEMMA 4.4.7. *The map $f\in\mathcal{L}(D;A_1\times...\times A_n\to B)$ is linear composition of maps $h^{\boldsymbol{i}_1...\boldsymbol{i}_n}_{\boldsymbol{i}}$.*

PROOF. For any $A_1$-number $a_1$

$$a_1=a_1^{\boldsymbol{i}_1}e_{1\cdot\boldsymbol{i}_1}\tag{4.4.7}$$

..., for any $A_n$-number $a_n$

$$a_n=a_n^{\boldsymbol{i}_n}e_{n\cdot\boldsymbol{i}_n}\tag{4.4.8}$$

the equality

$$\begin{aligned}h^{\boldsymbol{i}_1...\boldsymbol{i}_n}_{\boldsymbol{i}}\circ(a_1,...,a_n)&=h^{\boldsymbol{i}_1...\boldsymbol{i}_n}_{\boldsymbol{i}}\circ(a_1^{\boldsymbol{j}_1}e_{1\cdot\boldsymbol{j}_1},...,a_n^{\boldsymbol{j}_n}e_{n\cdot\boldsymbol{j}_n})\\&=a_1^{\boldsymbol{j}_1}...a_n^{\boldsymbol{j}_n}(h^{\boldsymbol{i}_1...\boldsymbol{i}_n}_{\boldsymbol{i}}\circ(e_{1\cdot\boldsymbol{j}_1},...,e_{n\cdot\boldsymbol{j}_n}))\\&=a_1^{\boldsymbol{j}_1}...a_n^{\boldsymbol{j}_n}\delta^{\boldsymbol{i}_1}_{\boldsymbol{j}_1}...\delta^{\boldsymbol{i}_n}_{\boldsymbol{j}_n}e_{B\cdot\boldsymbol{i}}\\&=a_1^{\boldsymbol{i}_1}...a_n^{\boldsymbol{i}_n}e_{B\cdot\boldsymbol{i}}\end{aligned}\tag{4.4.9}$$

follows from equalities (4.4.6), (4.4.7), (4.4.8) and from the theorem 4.2.2. The equality

$$\begin{aligned}f\circ(a_1,...,a_n)&=f\circ(a_1^{\boldsymbol{j}_1}e_{1\cdot\boldsymbol{j}_1},...,a_n^{\boldsymbol{j}_n}e_{n\cdot\boldsymbol{j}_n})\\&=a_1^{\boldsymbol{j}_1}...a_n^{\boldsymbol{j}_n}f\circ(e_{1\cdot\boldsymbol{j}_1}...e_{n\cdot\boldsymbol{j}_n})\end{aligned}\tag{4.4.10}$$

follows from equalities (4.4.7), (4.4.8). Since

$$f\circ(e_{1\cdot\boldsymbol{j}_1}...e_{n\cdot\boldsymbol{j}_n})\in B$$

then

$$f\circ(e_{1\cdot\boldsymbol{j}_1}...e_{n\cdot\boldsymbol{j}_n})=f^{\boldsymbol{i}}_{\boldsymbol{j}_1...\boldsymbol{j}_n}e_{\boldsymbol{i}}\tag{4.4.11}$$

The equality

$$f\circ(a_1,...,a_n)=a_1^{i_1}...a_n^{i_n}f^{i}_{i_1...i_n}e_i=f^{i}_{i_1...i_n}h_i^{i_1...i_n}\circ(a_1,...,a_n)\tag{4.4.12}$$

follows from equalities (4.4.9), (4.4.10), (4.4.11). The equality

$$f=f^{i}_{i_1...i_n}h_i^{i_1...i_n}$$

follows from equalities (4.3.2), (4.3.10), (4.4.12). ⊙

The theorem follows from lemmas 4.4.6, 4.4.7 and from the definition 4.1.10. □

THEOREM 4.4.8. *Let $A_1$, ..., $A_n$, $B$ be free modules over commutative ring $D$. $D$-module $\mathcal{L}(D;A_1\times...\times A_n\rightarrow B)$ is free $D$-module.*

PROOF. The theorem follows from the theorem 4.4.5. □

## 4.5. Tensor Product of $D$-Modules

THEOREM 4.5.1. *The commutative ring $D$ is Abelian multiplicative $\Omega$-group.*

PROOF. Let the product $\circ$ in the ring $D$ be defined according to rule

$$a\circ b=ab$$

Since the product in the ring is distributive over addition, the theorem follows from definitions 3.1.5, 3.1.8. □

THEOREM 4.5.2. *There exists* **tensor product** $A_1\otimes...\otimes A_n$ *of $D$-modules $A_1$, ..., $A_n$.*

PROOF. The theorem follows from the definition 4.1.2 and from theorems 3.3.5, 4.5.1. □

THEOREM 4.5.3. *Let $D$ be the commutative ring. Let $A_1$, ..., $A_n$ be $D$-modules. Tensor product is distributive over sum*

$$\begin{aligned}&a_1\otimes...\otimes(a_i+b_i)\otimes...\otimes a_n\\=&a_1\otimes...\otimes a_i\otimes...\otimes a_n+a_1\otimes...\otimes b_i\otimes...\otimes a_n\\&a_i,b_i\in A_i\end{aligned}\tag{4.5.1}$$

*The representation of the ring $D$ in tensor product is defined by equality*

$$\begin{aligned}&a_1\otimes...\otimes(ca_i)\otimes...\otimes a_n=c(a_1\otimes...\otimes a_i\otimes...\otimes a_n)\\&a_i\in A_i\quad c\in D\end{aligned}\tag{4.5.2}$$

PROOF. The equality (4.5.1) follows from the equality (3.3.23). The equality (4.5.2) follows from the equality (3.3.24). □

THEOREM 4.5.4. *Let $A_1$, ..., $A_n$ be modules over commutative ring $D$. Let*

$$f:A_1\times...\times A_n\rightarrow A_1\otimes...\otimes A_n$$

*be polylinear map defined by the equality*

$$f\circ(d_1,...,d_n)=d_1\otimes...\otimes d_n\tag{4.5.3}$$

*Let*

$$g:A_1\times...\times A_n\rightarrow V$$

*be polylinear map into $D$-module $V$. There exists a linear map*

$$h : A_1 \otimes ... \otimes A_n \to V$$

*such that the diagram*

(4.5.4)

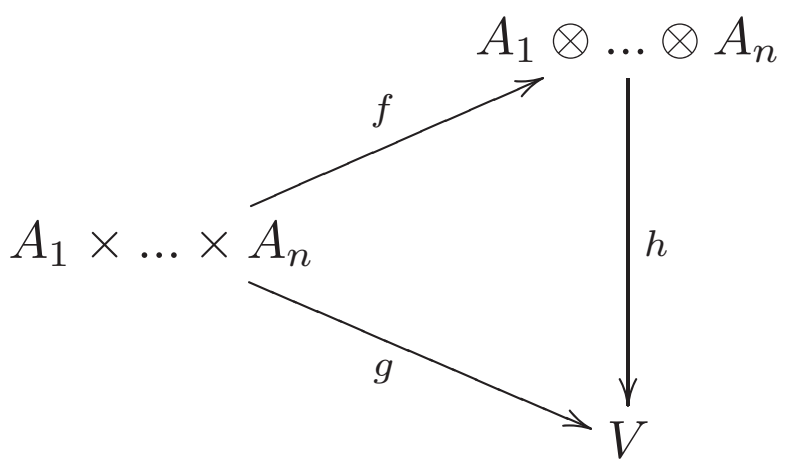

*is commutative. The map $h$ is defined by the equality*

$$h(a_1 \otimes ... \otimes a_n) = g(a_1, ..., a_n) \tag{4.5.5}$$

Proof. The theorem follows from the theorem 3.3.10 and from definitions 4.2.1, 4.3.1. □

Theorem 4.5.5. *The map*

$$(v_1, ..., v_n) \in V_1 \times ... \times V_n \to v_1 \otimes ... \otimes v_n \in V_1 \otimes ... \otimes V_n$$

*is polylinear map.*

Proof. The theorem follows from the theorem 3.3.9 and from the definition 4.3.1. □

Theorem 4.5.6. *Tensor product $A_1 \otimes ... \otimes A_n$ of free finite dimensional modules $A_1$, ..., $A_n$ over the commutative ring $D$ is free finite dimensional module.*

*Let $\overline{\overline{e}}_i$ be the basis of module $A_i$ over ring $D$. We can represent any tensor $a \in A_1 \otimes ... \otimes A_n$ in the following form*

$$a = a^{i_1...i_n} e_{1\cdot i_1} \otimes ... \otimes e_{n\cdot i_n} \tag{4.5.6}$$

*Expression $a^{i_1...i_n}$ is called* **standard component of tensor**.

Proof. Vector $a_i \in A_i$ has expansion

$$a_i = a_i^k \overline{e}_{i\cdot k}$$

relative to basis $\overline{\overline{e}}_i$. From equalities (4.5.1), (4.5.2), it follows

$$a_1 \otimes ... \otimes a_n = a_1^{i_1} ... a_n^{i_n} e_{1\cdot i_1} \otimes ... \otimes e_{n\cdot i_n}$$

Since set of tensors $a_1 \otimes ... \otimes a_n$ is the generating set of module $A_1 \otimes ... \otimes A_n$, then we can write tensor $a \in A_1 \otimes ... \otimes A_n$ in form

$$a = a^s a_{s\cdot 1}{}^{i_1} ... a_{s\cdot n}{}^{i_n} e_{1\cdot i_1} \otimes ... \otimes e_{n\cdot i_n} \tag{4.5.7}$$

where $a^s$, $a_{s\cdot 1}{}^{i_1}$, ..., $a_{s\cdot n}{}^{i_n} \in F$. Let

$$a^s a_{s\cdot 1}{}^{i_1} ... a_{s\cdot n}{}^{i_n} = a^{i_1...i_n}$$

Then equality (4.5.7) has form (4.5.6).

Therefore, set of tensors $e_{1\cdot i_1} \otimes ... \otimes e_{n\cdot i_n}$ is the generating set of module $A_1 \otimes ... \otimes A_n$. Since the dimension of module $A_i$, $i = 1$, ..., $n$, is finite, then the set of tensors $e_{1\cdot i_1} \otimes ... \otimes e_{n\cdot i_n}$ is finite. Therefore, the set of tensors $e_{1\cdot i_1} \otimes ... \otimes e_{n\cdot i_n}$ contains a basis of module $A_1 \otimes ... \otimes A_n$, and the module $A_1 \otimes ... \otimes A_n$ is free module over the ring $D$. □

CHAPTER 5

# $D$-Algebra

## 5.1. Algebra over Commutative Ring

DEFINITION 5.1.1. *Let $D$ be commutative ring. $D$-module $A$ is called* **algebra over ring** *$D$ or* **$D$-algebra***, if we defined product*[5.1] *in $A$*

$$vw = C \circ (v, w) \tag{5.1.1}$$

*where $C$ is bilinear map*

$$C : A \times A \to A$$

*If $A$ is free $D$-module, then $A$ is called* **free algebra over ring** *$D$.* □

THEOREM 5.1.2. *The multiplication in the algebra $A$ is distributive over addition.*

PROOF. The statement of the theorem follows from the chain of equations

$$(a+b)c = f \circ (a+b, c) = f \circ (a, c) + f \circ (b, c) = ac + bc$$
$$a(b+c) = f \circ (a, b+c) = f \circ (a, b) + f \circ (a, c) = ab + ac$$

□

The multiplication in algebra can be neither commutative nor associative. Following definitions are based on definitions given in [10], p. 13.

DEFINITION 5.1.3. *The* **commutator**

$$[a, b] = ab - ba$$

*measures commutativity in $D$-algebra $A$. $D$-algebra $A$ is called* **commutative***, if*

$$[a, b] = 0$$

□

DEFINITION 5.1.4. *The* **associator**

$$(a, b, c) = (ab)c - a(bc) \tag{5.1.2}$$

*measures associativity in $D$-algebra $A$. $D$-algebra $A$ is called* **associative***, if*

$$(a, b, c) = 0$$

□

[5.1] I follow the definition given in [10], p. 1, [8], p. 4. The statement which is true for any $D$-module, is true also for $D$-algebra.

THEOREM 5.1.5. *Let $A$ be algebra over commutative ring $D$.*[5.2]

$$a(b,c,d)+(a,b,c)d=(ab,c,d)-(a,bc,d)+(a,b,cd) \tag{5.1.3}$$

*for any $a, b, c, d \in A$.*

PROOF. The equation (5.1.3) follows from the chain of equations

$$\begin{aligned}a(b,c,d)+(a,b,c)d &= a((bc)d-b(cd))+((ab)c-a(bc))d\\ &= a((bc)d)-a(b(cd))+((ab)c)d-(a(bc))d\\ &= ((ab)c)d-(ab)(cd)+(ab)(cd)\\ &+a((bc)d)-a(b(cd))-(a(bc))d\\ &= (ab,c,d)-(a(bc))d+a((bc)d)+(ab)(cd)-a(b(cd))\\ &= (ab,c,d)-(a,(bc),d)+(a,b,cd)\end{aligned}$$

□

DEFINITION 5.1.6. *The set*[5.3]

$$N(A)=\{a\in A:\forall b,c\in A,(a,b,c)=(b,a,c)=(b,c,a)=0\}$$

*is called the* **nucleus of an $D$-algebra** $A$. □

DEFINITION 5.1.7. *The set*[5.4]

$$Z(A)=\{a\in A:a\in N(A),\forall b\in A,ab=ba\}$$

*is called the* **center of an $D$-algebra** $A$. □

THEOREM 5.1.8. *Let $D$ be commutative ring. If $D$-algebra $A$ has unit, then there exits an isomorphism $f$ of the ring $D$ into the center of the algebra $A$.*

PROOF. Let $e\in A$ be the unit of the algebra $A$. Then $f\circ a=ae$. □

Let $\overline{\overline{e}}$ be the basis of free algebra $A$ over ring $D$. If algebra $A$ has unit, then we assume that $e_0$ is the unit of algebra $A$.

THEOREM 5.1.9. *Let $\overline{\overline{e}}$ be the basis of free algebra $A$ over ring $D$. Let*

$$a=a^ie_i \quad b=b^ie_i \quad a,b\in A$$

*We can get the product of $a$, $b$ according to rule*

$$(ab)^k=C^k_{ij}a^ib^j \tag{5.1.4}$$

*where $C^k_{ij}$ are* **structural constants** *of algebra $A$ over ring $D$. The product of basis vectors in the algebra $A$ is defined according to rule*

$$e_ie_j=C^k_{ij}e_k \tag{5.1.5}$$

PROOF. The equation (5.1.5) is corollary of the statement that $\overline{\overline{e}}$ is the basis of the algebra $A$. Since the product in the algebra is a bilinear map, then we can write the product of $a$ and $b$ as

$$ab=a^ib^je_ie_j \tag{5.1.6}$$

[5.2] The statement of the theorem is based on the equation [10]-(2.4).
[5.3] The definition is based on the similar definition in [10], p. 13
[5.4] The definition is based on the similar definition in [10], p. 14

From equations (5.1.5), (5.1.6), it follows that

$$ab = a^i b^j C^k_{ij} e_k \tag{5.1.7}$$

Since $\overline{\overline{e}}$ is a basis of the algebra $A$, then the equation (5.1.4) follows from the equation (5.1.7). □

Theorem 5.1.10. *Since the algebra $A$ is commutative, then*

$$C^p_{ij} = C^p_{ji} \tag{5.1.8}$$

*Since the algebra $A$ is associative, then*

$$C^p_{ij} C^q_{pk} = C^q_{ip} C^p_{jk} \tag{5.1.9}$$

Proof. For commutative algebra, the equation (5.1.8) follows from equation

$$e_i e_j = e_j e_i$$

For associative algebra, the equation (5.1.9) follows from equation

$$(e_i e_j) e_k = e_i (e_j e_k)$$

□

Theorem 5.1.11. *The representation*

$$f_{2,3} : A \overset{*}{\dashrightarrow} A \tag{5.1.10}$$

*of $D$-module $A$ in $D$-module $A$ is equivalent to structure of $D$-algebra $A$.*

Proof.

- Let the structure of $D$-algebra $A$ defined in $D$-module $A$, be generated by product

$$v\,w = C \circ (v, w)$$

According to definitions 5.1.1 and 4.3.1, **left shift of $D$-module** $A$ defined by equation

$$l \circ v : w \in A \to v\,w \in A \tag{5.1.11}$$

is linear map. According to the definition 4.2.1, the map $l(v)$ is endomorphism of $D$-module $A$.

The equation

$$(l \circ (v_1 + v_2)) \circ w = (v_1 + v_2) w = v_1 w + v_2 w = (l \circ v_1) \circ w + (l \circ v_2) \circ w \tag{5.1.12}$$

follows from the definition 4.3.1 and from the equation (5.1.11). According to the theorem 4.2.3, the equation

$$l \circ (v_1 + v_2) = l \circ v_1 + l \circ v_2 \tag{5.1.13}$$

follows from equation (5.1.12). The equation

$$(l \circ (dv)) \circ w = (dv) w = d(vw) = d((l \circ v) \circ w) \tag{5.1.14}$$

follows from the definition 4.3.1 and from the equation (5.1.11). According to the theorem 4.2.3, the equation

$$l \circ (dv) = d(l \circ v) \tag{5.1.15}$$

follows from equation (5.1.14). From equations (5.1.13), (5.1.15), it follows that the map

$$f_{2,3} : v \to l \circ v$$

is the representation of $D$-module $A$ in $D$-module $A$

(5.1.16)
$$f_{2,3} : A \overset{*}{\longrightarrow} A \qquad f_{2,3} \circ v : w \to (l \circ v) \circ w$$

- Consider the representation (5.1.10) of $D$-module $A$ in $D$-module $A$. Since map $f_{2,3} \circ v$ is endomorphism of $D$-module $A$, then

(5.1.17)
$$\begin{aligned}(f_{2,3} \circ v)(w_1 + w_2) &= (f_{2,3} \circ v) \circ w_1 + (f_{2,3} \circ v) \circ w_2\\ (f_{2,3} \circ v) \circ (dw) &= d((f_{2,3} \circ v) \circ w)\end{aligned}$$

Since the map (5.1.10) is linear map
$$f_{2,3} : A \to \mathcal{L}(D; A; A)$$
then, according to theorems 4.2.3, 4.2.4,

(5.1.18)
$$(f_{2,3} \circ (v_1 + v_2)) \circ w = (f_{2,3} \circ v_1 + f_{2,3} \circ v_2)(w) = (f_{2,3} \circ v_1) \circ w + (f_{2,3} \circ v_2) \circ w$$

(5.1.19)
$$(f_{2,3} \circ (d\,v)) \circ w = (d\,(f_{2,3} \circ v)) \circ w = d\,((f_{2,3} \circ v) \circ w)$$

From equations (5.1.17), (5.1.18), (5.1.19) and the definition 4.3.1, it follows that the map $f_{2,3}$ is bilinear map. Therefore, the map $f_{2,3}$ determines the product in $D$-module $A$ according to rule
$$vw = (f_{2,3} \circ v) \circ w$$

□

Corollary 5.1.12. *$D$ is commutative ring, $A$ is Abelian group. The diagram of representations*

(5.1.20)
$$\begin{array}{ccccc} D & \xrightarrow[*]{g_{1,2}} & A & \xrightarrow[*]{g_{2,3}} & A \\ & & & & \uparrow{\scriptstyle *\, g_{1,2}} \\ & & & & D \end{array} \qquad \begin{array}{l} g_{1,2}(d) : v \to d\,v \\ g_{2,3}(v) : w \to C \circ (v, w) \\ C \in \mathcal{L}(D; A^2 \to A) \end{array}$$

*generates the structure of $D$-algebra $A$.* □

## 5.2. Linear Homomorphism

Theorem 5.2.1. *Let diagram of representations*

(5.2.1)
$$\begin{array}{ccccc} D_1 & \xrightarrow[*]{g_{1\cdot 1,2}} & A_1 & \xrightarrow[*]{g_{1\cdot 2,3}} & A_1 \\ & & & & \uparrow{\scriptstyle *\, g_{1\cdot 1,2}} \\ & & & & D_1 \end{array} \qquad \begin{array}{l} g_{1\cdot 1,2}(d) : v \to d\,v \\ g_{1\cdot 2,3}(v) : w \to C_1 \circ (v, w) \\ C_1 \in \mathcal{L}(D_1; A_1^2 \to A_1) \end{array}$$

*describe $D_1$-algebra $A_1$. Let diagram of representations*

(5.2.2)
$$\begin{array}{ccccc} D_2 & \xrightarrow[*]{g_{2\cdot 1,2}} & A_2 & \xrightarrow[*]{g_{2\cdot 2,3}} & A_2 \\ & & & & \uparrow{\scriptstyle *\, g_{2\cdot 1,2}} \\ & & & & D_2 \end{array} \qquad \begin{array}{l} g_{2\cdot 1,2}(d) : v \to d\,v \\ g_{2\cdot 2,3}(v) : w \to C_2 \circ (v, w) \\ C_2 \in \mathcal{L}(D_2; A_2^2 \to A_2) \end{array}$$

*describe $D_2$-algebra $A_2$. Morphism of $D_1$-algebra $A_1$ into $D_2$-algebra $A_2$ is tuple of maps*
$$r_1 : D_1 \to D_2 \quad r_2 : A_1 \to A_2$$

*where the map $r_1$ is homomorphim of ring $D_1$ into ring $D_2$ and the map $r_2$ is linear map of $D_1$-algebra $A_1$ into $D_2$-algebra $A_2$ such that*

$$r_2(ab) = r_2(a)r_2(b) \tag{5.2.3}$$

Proof. According to the equation [4]-(4.2.3), morphism $(r_1, r_2)$ of representation $f_{1,2}$ satisfies to the equation

$$\begin{aligned} r_2(f_{1\cdot 1,2}(d)(a)) &= f_{2\cdot 1,2}(r_1(d))(r_2(a)) \\ r_2(d\,a) &= r_1(d)r_2(a) \end{aligned} \tag{5.2.4}$$

Therefore, the map $(r_1, r_2)$ is linear map.

According to equations [4]-(4.2.3), the morphism $(r_2, r_2)$ of representation $f_{2,3}$ satisfies to the equation[5.5]

$$r_2(f_{1\cdot 2,3}(a_2)(a_3)) = f_{2\cdot 2,3}(r_2(a_2))(r_2(a_3)) \tag{5.2.5}$$

From equations (5.2.5), (5.2.1), (5.2.2), it follows that

$$r_2(C_1(v, w)) = C_2(r_2(v), r_2(w)) \tag{5.2.6}$$

Equation (5.2.3) follows from equations (5.2.6), (5.1.1). □

Definition 5.2.2. *The morphism of representations of $D_1$-algebra $A_1$ into $D_2$-algebra $A_2$ is called* **linear homomorphism** *of $D_1$-algebra $A_1$ into $D_2$-algebra $A_2$.* □

Theorem 5.2.3. *Let $\overline{\overline{e}}_1$ be the basis of $D_1$-algebra $A_1$. Let $\overline{\overline{e}}_2$ be the basis of $D_2$-algebra $A_2$. Then linear homomorphism*[5.6] *$(r_1, r_2)$ of $D_1$-algebra $A_1$ into $D_2$-algebra $A_2$ has presentation*[5.7]

$$b = e_{2*}{}^*r_{2*}{}^*r_1(a) = e_{2\cdot i}r_{2\cdot}{}^i_j r_1(a^j) \tag{5.2.7}$$

$$b = r_{2*}{}^*r_1(a) \tag{5.2.8}$$

*relative to selected bases. Here*

- *$a$ is coordinate matrix of vector $\overline{a}$ relative the basis $\overline{\overline{e}}_1$.*
- *$b$ is coordinate matrix of vector*

  $$\overline{b} = r_2(\overline{a})$$

  *relative the basis $\overline{\overline{e}}_2$.*
- *$r_2$ is coordinate matrix of set of vectors $(r_2(e_{1\cdot i}))$ relative the basis $\overline{\overline{e}}_2$. The matrix $r_2$ is called* **matrix of linear homomorphism** *relative bases $\overline{\overline{e}}_1$ and $\overline{\overline{e}}_2$.*

Proof. Vector $\overline{a} \in A_1$ has expansion

$$\overline{a} = e_{1*}{}^*a$$

relative to the basis $\overline{\overline{e}}_1$. Vector $\overline{b} \in A_2$ has expansion

$$\overline{b} = e_2{}^*{}_*b \tag{5.2.9}$$

[5.5] Since in diagrams of representations (5.2.1), (5.2.2), supports of $\Omega_2$-algebra and $\Omega_3$-algebra coincide, then morphisms of representations on levels 2 and 3 coincide also.

[5.6] This theorem is similar to the theorem [3]-5.4.3.

[5.7] We define product of matrices over commutative ring only as ${}_*{}^*$-product. However, I prefer to explicitly specify the operation, because in such case we see that this is expression with matrices. In addition, I expect to consider similar statement in non commutative case.

relative to the basis $\overline{\overline{e}}_2$.

Since $(r_1, r_2)$ is a linear homomorphism, then from (5.2.4), it follows that

$$b = r_2(a) = r_2(e_1{}^*{}_*a) = r_2(e_1)_*{}^*r_1(a) \tag{5.2.10}$$

where

$$r_1(a) = \begin{pmatrix} r_1(a^1) \\ ... \\ r_1(a^n) \end{pmatrix}$$

$r_2(e_{1\cdot i})$ is also a vector of $D$-module $A_2$ and has expansion

$$r_2(e_{1\cdot i}) = e_{2*}{}^*r_{2\cdot i} = e_{2\cdot j}\, r_{2\cdot}{}^j_i \tag{5.2.11}$$

relative to basis $\overline{\overline{e}}_2$. Combining (5.2.10) and (5.2.11) we get (5.2.7). (5.2.8) follows from comparison of (5.2.9) and (5.2.7) and theorem [3]-5.3.3. □

THEOREM 5.2.4. *Let $\overline{\overline{e}}_1$ be the basis of $D_1\star$-algebra $A_1$. Let $\overline{\overline{e}}_2$ be the basis of $D_2\star$-algebra $A_2$. If the map $r_1$ is injection, then there is relation between the matrix of linear homomorphism and structural constants*

$$r_{2\cdot}{}^l_k r_1(C_{1\cdot}{}^k_{ij}) = C_{2\cdot}{}^l_{pq} r_{2\cdot}{}^p_i r_{2\cdot}{}^q_j \tag{5.2.12}$$

PROOF. Let

$$\overline{a}, \overline{b} \in A_1 \quad \overline{a} = e_{1*}{}^*a \quad \overline{b} = e_{1*}{}^*b$$

From equations (5.1.4), (5.1.1), (5.2.1), it follows that

$$ab = e_{1\cdot k} C_{1\cdot}{}^k_{ij} a^i b^j \tag{5.2.13}$$

From equations (5.2.4), (5.2.13), it follows that

$$r_2(ab) = r_2(e_{1\cdot k}) r_1(C_{1\cdot}{}^k_{ij} a^i b^j) \tag{5.2.14}$$

Since the map $r_1$ is homomorphism of rings, then from equation (5.2.14), it follows that

$$r_2(ab) = r_2(e_{1\cdot k}) r_1(C_{1\cdot}{}^k_{ij}) r_1(a^i) r_1(b^j) \tag{5.2.15}$$

From the theorem 5.2.3 and the equation (5.2.15), it follows that

$$r_2(ab) = e_{2\cdot l}\, r_{2\cdot}{}^l_k r_1(C_{1\cdot}{}^k_{ij}) r_1(a^i) r_1(b^j) \tag{5.2.16}$$

From the equation (5.2.3) and the theorem 5.2.3, it follows that

$$r_2(ab) = r_2(a) r_2(b) = e_{2\cdot p}\, r_1(a^i) r_{2\cdot}{}^p_i e_{2\cdot q}\, r_1(b^j) r_{2\cdot}{}^q_j \tag{5.2.17}$$

From equations (5.1.4), (5.1.1), (5.2.2), (5.2.17), it follows that

$$r_2(ab) = e_{2\cdot l}\, C_{2\cdot}{}^l_{pq} r_1(a^i) r_{2\cdot}{}^p_i r_1(b^j) r_{2\cdot}{}^q_j \tag{5.2.18}$$

From equations (5.2.16), (5.2.18), it follows that

$$e_{2\cdot l}\, r_{2\cdot}{}^l_k r_1(C_{1\cdot}{}^k_{ij}) r_1(a^i) r_1(b^j) = e_{2\cdot l}\, C_{2\cdot}{}^l_{pq} r_1(a^i) r_{2\cdot}{}^p_i r_1(b^j) r_{2\cdot}{}^q_j \tag{5.2.19}$$

The equation (5.2.12) follows from the equation (5.2.19), because vectors of basis $\overline{\overline{e}}_2$ are linear independent, and $a^i$, $b^i$ (and therefore, $r_1(a^i)$, $r_1(b^i)$) are arbitrary values. □

## 5.3. Linear Automorphism of Quaternion Algebra

Determining of coordinates of linear automorphism is not a simple task. In this section we consider example of nontrivial linear automorphism of quaternion algebra.

THEOREM 5.3.1. *Coordinates of linear automorphism of quaternion algebra satisfy to the system of equations*

$$\left\{\begin{array}{lll} r_1^1=r_2^2r_3^3-r_2^3r_3^2 & r_2^1=r_3^2r_1^3-r_3^3r_1^2 & r_3^1=r_1^2r_2^3-r_1^3r_2^2\\ r_1^2=r_2^3r_3^1-r_2^1r_3^3 & r_2^2=r_3^3r_1^1-r_3^1r_1^3 & r_3^2=r_1^3r_2^1-r_1^1r_2^3\\ r_1^3=r_2^1r_3^2-r_2^2r_3^1 & r_2^3=r_3^1r_1^2-r_3^2r_1^1 & r_3^3=r_1^1r_2^2-r_1^2r_2^1 \end{array}\right. \tag{5.3.1}$$

PROOF. According to the theorems [5]-4.3.1, 5.2.4, linear automorphism of quaternion algebra satisfies to equations

$$\begin{array}{llll} r_0^l=r_0^pr_0^qC_{pq}^l & r_1^l=r_0^pr_1^qC_{pq}^l & r_2^l=r_0^pr_2^qC_{pq}^l & r_3^l=r_0^pr_3^qC_{pq}^l\\ r_1^l=r_1^pr_0^qC_{pq}^l & -r_0^l=r_1^pr_1^qC_{pq}^l & r_3^l=r_1^pr_2^qC_{pq}^l & -r_2^l=r_1^pr_3^qC_{pq}^l\\ r_2^l=r_2^pr_0^qC_{pq}^l & -r_3^l=r_2^pr_1^qC_{pq}^l & -r_0^l=r_2^pr_2^qC_{pq}^l & r_0^l=r_2^pr_3^qC_{pq}^l\\ r_3^l=r_3^pr_0^qC_{pq}^l & r_2^l=r_3^pr_1^qC_{pq}^l & -r_1^l=r_3^pr_2^qC_{pq}^l & -r_0^l=r_3^pr_3^qC_{pq}^l \end{array} \tag{5.3.2}$$

From the equation (5.3.2), it follows that

$$\begin{array}{lll} r_1^l=r_0^pr_1^qC_{pq}^l=r_2^pr_3^qC_{pq}^l & r_0^pr_1^qC_{pq}^l=r_0^pr_1^qC_{qp}^l & r_2^pr_3^qC_{pq}^l=-r_2^pr_3^qC_{qp}^l\\ r_2^l=r_0^pr_2^qC_{pq}^l=r_3^pr_1^qC_{pq}^l & r_0^pr_2^qC_{pq}^l=r_0^pr_2^qC_{qp}^l & r_3^pr_1^qC_{pq}^l=-r_1^pr_3^qC_{pq}^l\\ r_3^l=r_0^pr_3^qC_{pq}^l=r_1^pr_2^qC_{pq}^l & r_0^pr_3^qC_{pq}^l=r_0^pr_3^qC_{qp}^l & r_1^pr_2^qC_{pq}^l=-r_1^pr_2^qC_{qp}^l \end{array} \tag{5.3.3}$$

$$r_0^l=r_0^pr_0^qC_{pq}^l=-r_1^pr_1^qC_{pq}^l=-r_2^pr_2^qC_{pq}^l=-r_3^pr_3^qC_{pq}^l \tag{5.3.4}$$

If $l=0$, then from the equation

$$C_{pq}^0=C_{qp}^0$$

it follows that

$$r_i^pr_j^qC_{pq}^0=r_i^pr_j^qC_{qp}^0 \tag{5.3.5}$$

From the equation (5.3.3) for $l=0$ and the equation (5.3.5), it follows that

$$r_1^0=r_2^0=r_3^0=0 \tag{5.3.6}$$

If $l=1,2,3,$ then we can write the equation (5.3.3) in the following form

$$(5.3.7)\quad \left\{\begin{array}{l} r^l_i = r^l_0 r^0_i C^l_{l0} + r^0_0 r^l_i C^l_{0l} + r^a_0 r^b_i C^l_{ab} + r^b_0 r^a_i C^l_{ba} \\ \quad r^l_0 r^0_i C^l_{l0} + r^0_0 r^l_i C^l_{0l} + r^a_0 r^b_i C^l_{ab} + r^b_0 r^a_i C^l_{ba} \\ \quad = r^l_0 r^0_i C^l_{0l} + r^0_0 r^l_i C^l_{l0} + r^a_0 r^b_i C^l_{ba} + r^b_0 r^a_i C^l_{ab} \\ \quad i=1,2,3 \\ r^l_i = r^0_k r^l_j C^l_{0l} + r^l_k r^0_j C^l_{l0} + r^a_k r^b_j C^l_{ab} + r^b_k r^a_j C^l_{ba} \\ \quad r^0_k r^l_j C^l_{0l} + r^l_k r^0_j C^l_{l0} + r^a_k r^b_j C^l_{ab} + r^b_k r^a_j C^l_{ba} \\ \quad = -r^0_k r^l_j C^l_{l0} - r^l_k r^0_j C^l_{0l} - r^a_k r^b_j C^l_{ba} - r^b_k r^a_j C^l_{ab} \\ \quad\quad i=1 \quad k=2 \quad j=3 \\ \quad\quad i=2 \quad k=3 \quad j=1 \\ \quad\quad i=3 \quad k=1 \quad j=2 \\ \quad 0<a<b \quad a\ne l \quad b\ne l \end{array}\right.$$

From equations (5.3.7), (5.3.6) and equations

$$(5.3.8)\quad \begin{array}{c} C^l_{0l} = C^l_{l0} = 1 \\ C^l_{ab} = -C^l_{ba} \end{array}$$

it follows that

$$(5.3.9)\quad \left\{\begin{array}{l} r^l_i = r^0_0 r^l_i + r^a_0 r^b_i C^l_{ab} - r^b_0 r^a_i C^l_{ab} \\ \quad r^0_0 r^l_i + r^a_0 r^b_i C^l_{ab} - r^b_0 r^a_i C^l_{ab} \\ \quad = r^0_0 r^l_i - r^a_0 r^b_i C^l_{ab} + r^b_0 r^a_i C^l_{ab} \\ \quad i=1,2,3 \\ r^l_i = r^a_k r^b_j C^l_{ab} - r^b_k r^a_j C^l_{ab} \\ \quad r^a_k r^b_j C^l_{ab} - r^b_k r^a_j C^l_{ab} \\ \quad = r^a_k r^b_j C^l_{ab} - r^b_k r^a_j C^l_{ab} \\ \quad\quad i=1 \quad k=2 \quad j=3 \\ \quad\quad i=2 \quad k=3 \quad j=1 \\ \quad\quad i=3 \quad k=1 \quad j=2 \\ \quad 0<a<b \quad a\ne l \quad b\ne l \end{array}\right.$$

From equations (5.3.9) it follows that

(5.3.10)
$$\left\{\begin{array}{l} r^l_i = r^0_0 r^l_i \\ \qquad r^a_0 r^b_i - r^b_0 r^a_i = 0 \\ \qquad i = 1, 2, 3 \\ r^l_i = r^a_k r^b_j C^l_{ab} - r^b_k r^a_j C^l_{ab} \\ \qquad\quad i = 1 \qquad k = 2 \quad j = 3 \\ \qquad\quad i = 2 \qquad k = 3 \quad j = 1 \\ \qquad\quad i = 3 \qquad k = 1 \quad j = 2 \\ \qquad 0 < a < b \quad a \ne l \quad b \ne l \end{array}\right.$$

From the equation (5.3.10) it follows that

(5.3.11)
$$r^0_0 = 1$$

From the equation (5.3.4) for $l = 0$, it follows that

(5.3.12)
$$\begin{aligned} r^0_0 &= r^0_0 r^0_0 - r^1_0 r^1_0 - r^2_0 r^2_0 - r^3_0 r^3_0 \\ &= -r^0_i r^0_i + r^1_i r^1_i + r^2_i r^2_i + r^3_i r^3_i \\ i &= 1, 2, 3 \end{aligned}$$

From equations (5.3.6), (5.3.10), (5.3.12), it follows that

(5.3.13)
$$\begin{aligned} 0 &= r^1_0 r^1_0 + r^2_0 r^2_0 + r^3_0 r^3_0 \\ 1 &= r^1_i r^1_i + r^2_i r^2_i + r^3_i r^3_i \\ i &= 1, 2, 3 \end{aligned}$$

From equations (5.3.13) it follows that[5.8]

(5.3.14)
$$r^1_0 = r^2_0 = r^3_0 = 0$$

From the equation (5.3.4) for $l > 0$, it follows that

(5.3.15)
$$\begin{aligned} r^l_0 &= r^l_0 r^0_0 C^l_{l0} + r^0_0 r^l_0 C^l_{0l} + r^a_0 r^b_0 C^l_{ab} + r^b_0 r^a_0 C^l_{ba} \\ &= -r^l_i r^0_i C^l_{l0} - r^0_i r^l_i C^l_{0l} - r^a_i r^b_i C^l_{ab} - r^b_i r^a_i C^l_{ba} \\ i &> 0 \\ l &> 0 \quad 0 < a < b \quad a \ne l \quad b \ne l \end{aligned}$$

[5.8]Here, we rely on the fact that quaternion algebra is defined over real field. If we consider quaternion algebra over complex field, then the equation (5.3.13) defines a cone in the complex space. Correspondingly, we have wider choice of coordinates of linear automorphism.

Equations (5.3.15) are identically true by equations (5.3.6), (5.3.14), (5.3.8). From equations (5.3.14), (5.3.10), it follows that

$$(5.3.16)\quad \begin{cases} r^{l}_{i} = r^{a}_{k} r^{b}_{j} C^{l}_{ab} - r^{b}_{k} r^{a}_{j} C^{l}_{ab} \\ \quad i=1 \quad k=2 \quad j=3 \\ \quad i=2 \quad k=3 \quad j=1 \\ \quad i=3 \quad k=1 \quad j=2 \\ \quad l>0 \quad 0<a<b \quad a\neq l \quad b\neq l \end{cases}$$

Equations (5.3.1) follow from equations (5.3.16). □

Example 5.3.2. *It is evident that coordinates*

$$r^{i}_{j} = \delta^{i}_{j}$$

*satisfy the equation* (5.3.1). *It can be verified directly that coordinates of map*

$$r^{0}_{0} = 1 \quad r^{1}_{2} = 1 \quad r^{2}_{3} = 1 \quad r^{3}_{1} = 1$$

*also satisfy the equation* (5.3.1). *The matrix of coordinates of this map has form*

$$r = \begin{pmatrix} 1 & 0 & 0 & \\ 0 & 0 & 1 & 0 \\ 0 & 0 & 0 & 1 \\ 0 & 1 & 0 & 0 \end{pmatrix}$$

*According to the theorem [5]-4.3.4, standard components of the map r have form*

$$\begin{matrix} r^{00} = \frac{1}{4} & r^{11} = -\frac{1}{4} & r^{22} = -\frac{1}{4} & r^{33} = -\frac{1}{4} \\ r^{10} = -\frac{1}{4} & r^{01} = \frac{1}{4} & r^{32} = -\frac{1}{4} & r^{23} = -\frac{1}{4} \\ r^{20} = -\frac{1}{4} & r^{31} = -\frac{1}{4} & r^{02} = \frac{1}{4} & r^{13} = -\frac{1}{4} \\ r^{30} = -\frac{1}{4} & r^{21} = -\frac{1}{4} & r^{12} = -\frac{1}{4} & r^{03} = \frac{1}{4} \end{matrix}$$

*Therefore, the map r has form*

$$r(a) = a^{0} + a^{2} i + a^{3} j + a^{1} k$$

$$\begin{aligned} r(a) = \frac{1}{4}(&a - iai - jaj - kak - ia + ai - kaj - jak \\ &-ja - kai + aj - iak - ka - jai - iaj + ak) \end{aligned}$$

□

CHAPTER 6

# Linear Map of Algebra

## 6.1. Linear Map of Algebra

Definition 6.1.1. *Let $A_1$ and $A_2$ be algebras over ring $D$. The linear map of the $D$-module $A_1$ into the $D$-module $A_2$ is called* **linear map** *of $D$-algebra $A_1$ into $D$-algebra $A_2$. Let us denote $\mathcal{L}(D; A_1 \to A_2)$ set of linear maps of $D$-algebra $A_1$ into $D$-algebra $A_2$.* □

Definition 6.1.2. *Let $A_1$, ..., $A_n$, $S$ be $D$-algebras. Polylinear map*

$$f : A_1 \times ... \times A_n \to S$$

*of $D$-modules $A_1$, ..., $A_n$ into $D$-module $S$ is called* **polylinear map** *of $D$-algebras $A_1$, ..., $A_n$ into $D$-algebra $S$. Let us denote $\mathcal{L}(D; A_1 \times ... \times A_n \to S)$ set of polylinear maps of $D$-algebras $A_1$, ..., $A_n$ into $D$-algebra $S$. Let us denote $\mathcal{L}(D; A^n \to S)$ set of $n$-linear maps of $D$-algebra $A$ ($A_1 = ... = A_n = A$) into $D$-algebra $S$.* □

Theorem 6.1.3. *Tensor product $A_1 \otimes ... \otimes A_n$ of $D$-algebras $A_1$, ..., $A_n$ is $D$-algebra.*

Proof. According to the definition 5.1.1 and to the theorem 4.5.2, tensor product $A_1 \otimes ... \otimes A_n$ of $D$-algebras $A_1$, ..., $A_n$ is $D$-module.

Consider the map

$$* : (A_1 \times ... \times A_n) \times (A_1 \times ... \times A_n) \to A_1 \otimes ... \otimes A_n \tag{6.1.1}$$

defined by the equation

$$(a_1, ..., a_n) * (b_1, ..., b_n) = (a_1 b_1) \otimes ... \otimes (a_n b_n) \tag{6.1.2}$$

For given values of variables $b_1$, ..., $b_n$, the map (6.1.1) is polylinear map with respect to variables $a_1$, ..., $a_n$. According to the theorem 4.5.4, there exists a linear map

$$* (b_1, ..., b_n) : A_1 \otimes ... \otimes A_n \to A_1 \otimes ... \otimes A_n \tag{6.1.3}$$

defined by the equation

$$(a_1 \otimes ... \otimes a_n) * (b_1, ..., b_n) = (a_1 b_1) \otimes ... \otimes (a_n b_n) \tag{6.1.4}$$

Since we can present any tensor $a \in A_1 \otimes ... \otimes A_n$ as sum of tensors $a_1 \otimes ... \otimes a_n$, then, for given tensor $a \in A_1 \otimes ... \otimes A_n$, the map (6.1.3) is polylinear map of variables $b_1$, ..., $b_n$. According to the theorem 4.5.4, there exists a linear map

$$* (a) : A_1 \otimes ... \otimes A_n \to A_1 \otimes ... \otimes A_n \tag{6.1.5}$$

defined by the equation

$$(a_1 \otimes ... \otimes a_n) * (b_1 \otimes ... \otimes b_n) = (a_1 b_1) \otimes ... \otimes (a_n b_n) \tag{6.1.6}$$

Therefore, the equation (6.1.6) defines bilinear map

$$*:(A_1\otimes...\otimes A_n)\times(A_1\otimes...\otimes A_n)\rightarrow A_1\otimes...\otimes A_n \tag{6.1.7}$$

Bilinear map (6.1.7) generates the product in $D$-module $A_1\otimes...\otimes A_n$. □

In case of tensor product of $D$-algebras $A_1$, $A_2$ we consider product defined by the equation

$$(a_1\otimes a_2)\circ(b_1\otimes b_2)=(a_1b_1)\otimes(b_2a_2) \tag{6.1.8}$$

THEOREM 6.1.4. *Let $\overline{\overline{e}}_i$ be the basis of the algebra $A_i$ over the ring $D$. Let $B_{i\cdot}{}^{j}_{kl}$ be structural constants of the algebra $A_i$ relative the basis $\overline{\overline{e}}_i$. Structural constants of the tensor product $A_1\otimes...\otimes A_n$ relative to the basis $e_{1\cdot i_1}\otimes...\otimes e_{n\cdot i_n}$ have form*

$$C^{\cdot j_1...j_n}_{\cdot k_1...k_n\cdot l_1...l_n}=C_{1\cdot}{}^{j_1}_{k_1l_1}...C_{n\cdot}{}^{j_n}_{k_nl_n} \tag{6.1.9}$$

PROOF. Direct multiplication of tensors $e_{1\cdot i_1}\otimes...\otimes e_{n\cdot i_n}$ has form

$$\begin{aligned}&(e_{1\cdot k_1}\otimes...\otimes e_{n\cdot k_n})(e_{1\cdot l_1}\otimes...\otimes e_{n\cdot l_n})\\ =&(\overline{e}_{1\cdot k_1}\overline{e}_{1\cdot l_1})\otimes...\otimes(\overline{e}_{n\cdot k_n}\overline{e}_{n\cdot l_n})\\ =&(\overline{e}_{1\cdot k_1}\overline{e}_{1\cdot l_1})\otimes...\otimes(\overline{e}_{n\cdot k_n}\overline{e}_{n\cdot l_n})\\ =&(C_{1\cdot}{}^{j_1}_{k_1l_1}\overline{e}_{1\cdot j_1})\otimes...\otimes(C_{n\cdot}{}^{j_n}_{k_nl_n}\overline{e}_{n\cdot j_n})\\ =&C_{1\cdot}{}^{j_1}_{k_1l_1}...C_{n\cdot}{}^{j_n}_{k_nl_n}\overline{e}_{1\cdot j_1}\otimes...\otimes\overline{e}_{n\cdot j_n}\end{aligned} \tag{6.1.10}$$

According to the definition of structural constants

$$(e_{1\cdot k_1}\otimes...\otimes e_{n\cdot k_n})(e_{1\cdot l_1}\otimes...\otimes e_{n\cdot l_n})=C^{\cdot j_1...j_n}_{\cdot k_1...k_n\cdot l_1...l_n}(e_{1\cdot j_1}\otimes...\otimes e_{n\cdot j_n}) \tag{6.1.11}$$

The equation (6.1.9) follows from comparison (6.1.10), (6.1.11).

From the chain of equations

$$\begin{aligned}&(a_1\otimes...\otimes a_n)(b_1\otimes...\otimes b_n)\\ =&(a_1^{k_1}\overline{e}_{1\cdot k_1}\otimes...\otimes a_n^{k_n}\overline{e}_{n\cdot k_n})(b_1^{l_1}\overline{e}_{1\cdot l_1}\otimes...\otimes b_n^{l_n}\overline{e}_{n\cdot l_n})\\ =&a_1^{k_1}...a_n^{k_n}b_1^{l_1}...b_n^{l_n}(\overline{e}_{1\cdot k_1}\otimes...\otimes\overline{e}_{n\cdot k_n})(\overline{e}_{1\cdot l_1}\otimes...\otimes\overline{e}_{n\cdot l_n})\\ =&a_1^{k_1}...a_n^{k_n}b_1^{l_1}...b_n^{l_n}C^{\cdot j_1...j_n}_{\cdot k_1...k_n\cdot l_1...l_n}(e_{1\cdot j_1}\otimes...\otimes e_{n\cdot j_n})\\ =&a_1^{k_1}...a_n^{k_n}b_1^{l_1}...b_n^{l_n}C_{1\cdot}{}^{j_1}_{k_1l_1}...C_{n\cdot}{}^{j_n}_{k_nl_n}(e_{1\cdot j_1}\otimes...\otimes e_{n\cdot j_n})\\ =&(a_1^{k_1}b_1^{l_1}\ C_{1\cdot}{}^{j_1}_{k_1l_1}\overline{e}_{1\cdot j_1})\otimes...\otimes(a_n^{k_n}b_n^{l_n}C_{n\cdot}{}^{j_n}_{k_nl_n}\overline{e}_{n\cdot j_n})\\ =&(a_1b_1)\otimes...\otimes(a_nb_n)\end{aligned}$$

it follows that definition of product (6.1.11) with structural constants (6.1.9) agreed with the definition of product (6.1.6). □

THEOREM 6.1.5. *For tensors $a$, $b\in A_1\otimes...\otimes A_n$, standard components of product satisfy to equation*

$$(ab)^{j_1...j_n}=C^{\cdot j_1...j_n}_{\cdot k_1...k_n\cdot l_1...l_n}a^{k_1...k_n}b^{l_1...l_n} \tag{6.1.12}$$

PROOF. According to the definition

$$ab=(ab)^{j_1...j_n}e_{1\cdot j_1}\otimes...\otimes e_{n\cdot j_n} \tag{6.1.13}$$

At the same time

$$
\begin{aligned}
ab= & a^{k_1...k_n}e_{1\cdot k_1}\otimes ...\otimes e_{n\cdot k_n}b^{k_1...k_n}e_{1\cdot l_1}\otimes ...\otimes e_{n\cdot l_n}\\
= & a^{k_1...k_n}b^{k_1...k_n}C^{\cdot j_1...j_n}_{\cdot k_1...k_n\cdot l_1...l_n}e_{1\cdot j_1}\otimes ...\otimes e_{n\cdot j_n}
\end{aligned}
\tag{6.1.14}
$$

The equation (6.1.12) follows from equations (6.1.13), (6.1.14). □

Theorem 6.1.6. *If the algebra $A_i$, $i=1$, ..., $n$, is associative, then the tensor product $A_1\otimes ...\otimes A_n$ is associative algebra.*

Proof. Since

$$
\begin{aligned}
&((e_{1\cdot i_1}\otimes ...\otimes e_{n\cdot i_n})(e_{1\cdot j_1}\otimes ...\otimes e_{n\cdot j_n}))(e_{1\cdot k_1}\otimes ...\otimes e_{n\cdot k_n})\\
=&((\overline{e}_{1\cdot i_1}\overline{e}_{1\cdot j_1})\otimes ...\otimes (\overline{e}_{n\cdot i_n}\overline{e}_{1\cdot j_n}))(e_{1\cdot k_1}\otimes ...\otimes e_{n\cdot k_n})\\
=&((\overline{e}_{1\cdot i_1}\overline{e}_{1\cdot j_1})\overline{e}_{1\cdot k_1})\otimes ...\otimes ((\overline{e}_{n\cdot i_n}\overline{e}_{1\cdot j_n})\overline{e}_{1\cdot k_n})\\
=&(\overline{e}_{1\cdot i_1}(\overline{e}_{1\cdot j_1}\overline{e}_{1\cdot k_1}))\otimes ...\otimes (\overline{e}_{n\cdot i_n}(\overline{e}_{1\cdot j_n}\overline{e}_{1\cdot k_n}))\\
=&(e_{1\cdot i_1}\otimes ...\otimes e_{n\cdot i_n})((\overline{e}_{1\cdot j_1}\overline{e}_{1\cdot k_1})\otimes ...\otimes (\overline{e}_{1\cdot j_n}\overline{e}_{1\cdot k_n}))\\
=&(e_{1\cdot i_1}\otimes ...\otimes e_{n\cdot i_n})((e_{1\cdot j_1}\otimes ...\otimes e_{n\cdot j_n})(e_{1\cdot k_1}\otimes ...\otimes e_{n\cdot k_n}))
\end{aligned}
$$

then

$$
\begin{aligned}
(ab)c=&a^{i_1...i_n}b^{j_1...j_n}c^{k_1...k_n}\\
&((e_{1\cdot i_1}\otimes ...\otimes e_{n\cdot i_n})(e_{1\cdot j_1}\otimes ...\otimes e_{n\cdot j_n}))(e_{1\cdot k_1}\otimes ...\otimes e_{n\cdot k_n})\\
=&a^{i_1...i_n}b^{j_1...j_n}c^{k_1...k_n}\\
&(e_{1\cdot i_1}\otimes ...\otimes e_{n\cdot i_n})((e_{1\cdot j_1}\otimes ...\otimes e_{n\cdot j_n})(e_{1\cdot k_1}\otimes ...\otimes e_{n\cdot k_n}))\\
=&a(bc)
\end{aligned}
$$

□

Theorem 6.1.7. *Let $A$ be algebra over commutative ring $D$. There exists a linear map*

$$h:a\otimes b\in A\otimes A\rightarrow ab\in A$$

Proof. The theorem is corollary of the definition 5.1.1 and the theorem 4.5.4. □

Theorem 6.1.8. *Let map*

$$f:A_1\rightarrow A_2$$

*be linear map of $D$-algebra $A_1$ into $D$-algebra $A_2$. Then maps $af$, $fb$, $a$, $b\in A_2$, defined by equations*

$$(af)\circ x=a(f\circ x)$$
$$(fb)\circ x=(f\circ x)b$$

*are linear.*

Proof. Statement of theorem follows from chains of equations

$$\begin{aligned}(af)\circ(x+y)&=a(f\circ(x+y))=a(f\circ x+f\circ y)=a(f\circ x)+a(f\circ y)\\&=(af)\circ x+(af)\circ y\\(af)\circ(px)&=a(f\circ(px))=ap(f\circ x)=pa(f\circ x)\\&=p((af)\circ x)\\(fb)\circ(x+y)&=(f\circ(x+y))b=(f\circ x+f\circ y)\ b=(f\circ x)b+(f\circ y)b\\&=(fb)\circ x+(fb)\circ y\\(fb)\circ(px)&=(f\circ(px))b=p(f\circ x)b\\&=p((fb)\circ x)\end{aligned}$$

□

## 6.2. Algebra $\mathcal{L}(D;A\to A)$

Theorem 6.2.1. *Let $A$, $B$, $C$ be algebras over commutative ring $D$. Let $f$ be linear map from $D$-algebra $A$ into $D$-algebra $B$. Let $g$ be linear map from $D$-algebra $B$ into $D$-algebra $C$. The map $g\circ f$ defined by diagram*

(6.2.1)

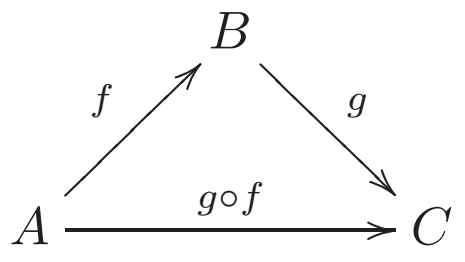

*is linear map from $D$-algebra $A$ into $D$-algebra $C$.*

Proof. The proof of the theorem follows from chains of equations

$$\begin{aligned}(g\circ f)\circ(a+b)&=g\circ(f\circ(a+b))=g\circ(f\circ a+f\circ b)\\&=g\circ(f\circ a)+g\circ(f\circ b)=(g\circ f)\circ a+(g\circ f)\circ b\\(g\circ f)\circ(pa)&=g\circ(f\circ(pa))=g\circ(p\ f\circ a)=p\ g\circ(f\circ a)\\&=p\ (g\circ f)\circ a\end{aligned}$$

□

Theorem 6.2.2. *Let $A$, $B$, $C$ be algebras over the commutative ring $D$. Let $f$ be a linear map from $D$-algebra $A$ into $D$-algebra $B$. The map $f$ generates a linear map*

$$f^*:g\in\mathcal{L}(D;B\to C)\to g\circ f\in\mathcal{L}(D;A\to C)\tag{6.2.2}$$

(6.2.3)

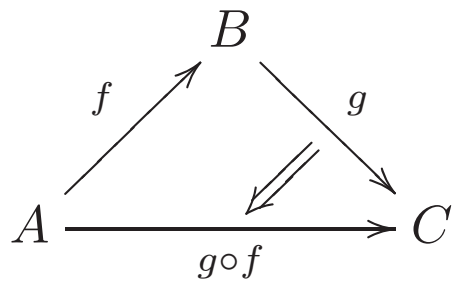

PROOF. The proof of the theorem follows from chains of equations[6.1]

$$\begin{aligned}((g_1+g_2)\circ f)\circ a &= (g_1+g_2)\circ(f\circ a) = g_1\circ(f\circ a)+g_2\circ(f\circ a)\\ &= (g_1\circ f)\circ a+(g_2\circ f)\circ a\\ &= (g_1\circ f+g_2\circ f)\circ a\\ ((pg)\circ f)\circ a &= (pg)\circ(f\circ a) = p\ g\circ(f\circ a) = p\ (g\circ f)\circ a\\ &= (p(g\circ f))\circ a\end{aligned}$$

□

THEOREM 6.2.3. *Let $A$, $B$, $C$ be algebras over the commutative ring $D$. Let $g$ be a linear map from $D$-algebra $B$ into $D$-algebra $C$. The map $g$ generates a linear map*

$$g_* : f \in \mathcal{L}(D; A \to B) \to g\circ f \in \mathcal{L}(D; A \to C) \tag{6.2.6}$$

(6.2.7)

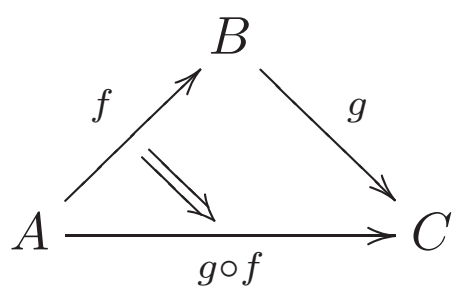

PROOF. The proof of the theorem follows from chains of equations[6.2]

$$\begin{aligned}(g\circ(f_1+f_2))\circ a &= g\circ((f_1+f_2)\circ a) = g\circ(f_1\circ a+f_2\circ a)\\ &= g\circ(f_1\circ a)+g\circ(f_2\circ a) = (g\circ f_1)\circ a+(g\circ f_2)\circ a\\ &= (g\circ f_1+g\circ f_2)\circ a\\ (g\circ(pf))\circ a &= g\circ((pf)\circ a) = g\circ(p\ (f\circ a)) = p\ g\circ(f\circ a)\\ &= p\ (g\circ f)\circ a = (p(g\circ f))\circ a\end{aligned}$$

□

THEOREM 6.2.4. *Let $A$, $B$, $C$ be algebras over the commutative ring $D$. The map*

$$\circ : (g,f) \in \mathcal{L}(D; B \to C)\times\mathcal{L}(D; A \to B) \to g\circ f \in \mathcal{L}(D; A \to C) \tag{6.2.10}$$

*is bilinear map.*

PROOF. The theorem follows from theorems 6.2.2, 6.2.3. □

---

[6.1] We use following definitions of operations over maps

$$(f+g)\circ a = f\circ a+g\circ a \tag{6.2.4}$$

$$(pf)\circ a = p\,f\circ a \tag{6.2.5}$$

[6.2] We use following definitions of operations over maps

$$(f+g)\circ a = f\circ a+g\circ a \tag{6.2.8}$$

$$(pf)\circ a = p\,f\circ a \tag{6.2.9}$$

THEOREM 6.2.5. *Let $A$ be algebra over commutative ring $D$. $D$-module $\mathcal{L}(D;A\to A)$ equiped by product*

$$\circ : (g,f)\in\mathcal{L}(D;A\to A)\times\mathcal{L}(D;A\to A)\to g\circ f\in\mathcal{L}(D;A\to A) \tag{6.2.11}$$

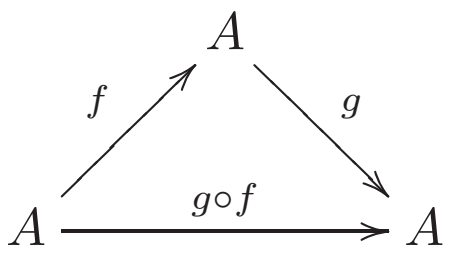

*is algebra over $D$.*

PROOF. The theorem follows from definition 5.1.1 and theorem 6.2.4. □

## 6.3. Linear Map into Associative Algebra

THEOREM 6.3.1. *Consider $D$-algebras $A_1$ and $A_2$. For given map $f\in\mathcal{L}(D;A_1\to A_2)$, the map*

$$g:A_2\times A_2\to\mathcal{L}(D;A_1\to A_2)$$
$$g(a,b)\circ f=afb$$

*is bilinear map.*

PROOF. The statement of theorem follows from chains of equations

$$\begin{aligned}((a_1+a_2)fb)\circ x&=(a_1+a_2)\ f\circ x\ b=a_1\ f\circ x\ b+a_2\ f\circ x\ b\\&=(a_1fb)\circ x+(a_2fb)\circ x=(a_1fb+a_2fb)\circ x\\((pa)fb)\circ x&=(pa)\ f\circ x\ b=p(a\ f\circ x\ b)=p((afb)\circ x)=(p(afb))\circ x\\(af(b_1+b_2))\circ x&=a\ f\circ x\ (b_1+b_2)=a\ f\circ x\ b_1+a\ f\circ x\ b_2\\&=(afb_1)\circ x+(afb_2)\circ x=(afb_1+afb_2)\circ x\\(af(pb))\circ x&=a\ f\circ x\ (pb)=p(a\ f\circ x\ b)=p((afb)\circ x)=(p(afb))\circ x\end{aligned}$$

□

THEOREM 6.3.2. *Consider $D$-algebras $A_1$ and $A_2$. For given map $f\in\mathcal{L}(D;A_1\to A_2)$, there exists linear map*

$$h:A_2\otimes A_2\to\mathcal{L}(D;A_1\to A_2)$$

*defined by the equation*

$$(a\otimes b)\circ f=afb \tag{6.3.1}$$

PROOF. The statement of the theorem is corollary of theorems 4.5.4, 6.3.1. □

THEOREM 6.3.3. *Consider $D$-algebras $A_1$ and $A_2$. Let us define product in algebra $A_2\otimes A_2$ according to rule*

$$(c\otimes d)\circ(a\otimes b)=(ca)\otimes(bd) \tag{6.3.2}$$

*A linear map*

$$h:A_2\otimes A_2\to{}^*\mathcal{L}(D;A_1\to A_2) \tag{6.3.3}$$

*defined by the equation*

$$(a\otimes b)\circ f=afb\quad a,b\in A_2\quad f\in\mathcal{L}(D;A_1\to A_2) \tag{6.3.4}$$

*is representation*[6.3] *of algebra* $A_2\otimes A_2$ *in module* $\mathcal{L}(D;A_1\rightarrow A_2)$.

PROOF. According to theorem 6.1.8, map (6.3.4) is transformation of module $\mathcal{L}(D;A_1\rightarrow A_2)$. For a given tensor $c\in A_2\otimes A_2$, a transformation $h(c)$ is a linear transformation of module $\mathcal{L}(D;A_1\rightarrow A_2)$, because

$$\begin{aligned}((a\otimes b)\circ(f_1+f_2))\circ x&=(a(f_1+f_2)b)\circ x=a((f_1+f_2)\circ x)b\\&=a(f_1\circ x+f_2\circ x)b=a(f_1\circ x)b+a(f_2\circ x)b\\&=(af_1b)\circ x+(af_2b)\circ x\\&=(a\otimes b)\circ f_1\circ x+(a\otimes b)\circ f_2\circ x\\&=((a\otimes b)\circ f_1+(a\otimes b)\circ f_2)\circ x\\((a\otimes b)\circ(pf))\circ x&=(a(pf)b)\circ x=a((pf)\circ x)b\\&=a(p\ f\circ x)b=pa(f\circ x)b\\&=p\ (afb)\circ x=p\ ((a\otimes b)\circ f)\circ x\\&=(p((a\otimes b)\circ f))\circ x\end{aligned}$$

According to theorem 6.3.2, map (6.3.4) is linear map.

Let $f\in\mathcal{L}(D;A_1\rightarrow A_2)$, $a\otimes b$, $c\otimes d\in A_2\otimes A_2$. According to the theorem 6.3.2

$$(a\otimes b)\circ f=afb\in\mathcal{L}(D;A_1\rightarrow A_2)$$

Therefore, according to the theorem 6.3.2

$$(c\otimes d)\circ((a\otimes b)\circ f)=c(afb)d$$

Since the product in algebra $A_2$ is associative, then

$$(c\otimes d)\circ((a\otimes b)\circ f)=c(afb)d=(ca)f(bd)=(ca\otimes bd)\circ f$$

Therefore, since we define the product in algebra $A_2\otimes A_2$ according to equation (6.3.2), then the map (6.3.3) is morphism of algebras. According to the definition [4]-2.1.2, map (6.3.4) is a representation of the algebra $A_2\otimes A_2$ in the module $\mathcal{L}(D;A_1\rightarrow A_2)$. □

THEOREM 6.3.4. *Consider $D$-algebra $A$. Let us define product in algebra $A\otimes A$ according to rule (6.3.2). A representation of algebra $A\otimes A$*

$$h:A\otimes A\rightarrow {}^*\mathcal{L}(D;A\rightarrow A)\tag{6.3.5}$$

*in module $\mathcal{L}(D;A\rightarrow A)$ defined by the equation*

$$(a\otimes b)\circ f=afb\quad a,b\in A\quad f\in\mathcal{L}(D;A\rightarrow A)\tag{6.3.6}$$

*allows us to identify tensor $d\in A\otimes A$ and map $d\circ\delta\in\mathcal{L}(D;A\rightarrow A)$ where $\delta\in\mathcal{L}(D;A\rightarrow A)$ is identity map.*

PROOF. According to the theorem 6.3.2, the map $f\in\mathcal{L}(D;A\rightarrow A)$ and the tensor $d\in A\otimes A$ generate the map

$$x\rightarrow(d\circ f)\circ x\tag{6.3.7}$$

If we assume $f=\delta$, $d=a\otimes b$, then the equation (6.3.7) gets form

$$((a\otimes b)\circ\delta)\circ x=(a\delta b)\circ x=a\ (\delta\circ x)\ b=axb\tag{6.3.8}$$

---

[6.3]See the definition of representation of $\Omega$-algebra in the definition [4]-2.1.2.

If we assume

$$((a \otimes b) \circ \delta) \circ x = (a \otimes b) \circ (\delta \circ x) = (a \otimes b) \circ x \tag{6.3.9}$$

then comparison of equations (6.3.8) and (6.3.9) gives a basis to identify the action of the tensor $a \otimes b$ and transformation $(a \otimes b) \circ \delta$. □

From the theorem 6.3.4, it follows that we can consider the map (6.3.4) as the product of maps $a \otimes b$ and $f$. The tensor $a \in A_2 \otimes A_2$ is **nonsingular**, if there exists the tensor $b \in A_2 \otimes A_2$ such that $a \circ b = 1 \otimes 1$.

DEFINITION 6.3.5. *Consider*[6.4] *the representation of algebra $A_2 \otimes A_2$ in the module $\mathcal{L}(D; A_1 \to A_2)$. The set*

$$(A_2 \otimes A_2) \circ f = \{g = d \circ f : d \in A_2 \otimes A_2\}$$

*is called* **orbit of linear map** $f \in \mathcal{L}(D; A_1 \to A_2)$. □

THEOREM 6.3.6. *Consider $D$-algebra $A_1$ and associative $D$-algebra $A_2$. Consider the representation of algebra $A_2 \otimes A_2$ in the module $\mathcal{L}(D; A_1; A_2)$. The map*

$$h : A_1 \to A_2$$

*generated by the map*

$$f : A_1 \to A_2$$

*has form*

$$h = (a_{s\cdot 0} \otimes a_{s\cdot 1}) \circ f = a_{s\cdot 0} f a_{s\cdot 1} \tag{6.3.10}$$

PROOF. We can represent any tensor $a \in A_2 \otimes A_2$ in the form

$$a = a_{s\cdot 0} \otimes a_{s\cdot 1}$$

According to the theorem 6.3.3, the map (6.3.4) is linear. This proofs the statement of the theorem. □

THEOREM 6.3.7. *Let $A_2$ be algebra with unit $e$. Let $a \in A_2 \otimes A_2$ be a nonsingular tensor. Orbits of linear maps $f \in \mathcal{L}(D; A_1 \to A_2)$ and $g = a \circ f$ coincide*

$$(A_2 \otimes A_2) \circ f = (A_2 \otimes A_2) \circ g \tag{6.3.11}$$

PROOF. If $h \in (A_2 \otimes A_2) \circ g$, then there exists $b \in A_2 \otimes A_2$ such that $h = b \circ g$. In that case

$$h = b \circ (a \circ f) = (b \circ a) \circ f \tag{6.3.12}$$

Therefore, $h \in (A_2 \otimes A_2) \circ f$,

$$(A_2 \otimes A_2) \circ g \subset (A_2 \otimes A_2) \circ f \tag{6.3.13}$$

Since $a$ is nonsingular tensor, then

$$f = a^{-1} \circ g \tag{6.3.14}$$

If $h \in (A_2 \otimes A_2) \circ f$, then there exists $b \in A_2 \otimes A_2$ such that

$$h = b \circ f \tag{6.3.15}$$

From equations (6.3.14), (6.3.15), it follows that

$$h = b \circ (a^{-1} \circ g) = (b \circ a^{-1}) \circ g$$

[6.4] The definition is made by analogy with the definition [4]-3.1.8.

Therefore, $h \in (A_2 \otimes A_2) \circ g$,

$$(A_2 \otimes A_2) \circ f \subset (A_2 \otimes A_2) \circ g \tag{6.3.16}$$

(6.3.11) follows from equations (6.3.13), (6.3.16). □

From the theorem 6.3.7, it also follows that if $g = a \circ f$ and $a \in A_2 \otimes A_2$ is a singular tensor, then relationship (6.3.13) is true. However, the main result of the theorem 6.3.7 is that the representations of the algebra $A_2 \otimes A_2$ in module $\mathcal{L}(D; A_1 \to A_2)$ generates an equivalence in the module $\mathcal{L}(D; A_1 \to A_2)$. If we successfully choose the representatives of each equivalence class, then the resulting set will be generating set of considered representation.[6.5]

## 6.4. Linear Map into Free Finite Dimensional Associative Algebra

THEOREM 6.4.1. *Let $A_1$ be free $D$-module. Let $A_2$ be free finite dimensional associative $D$-algebra. Let $\overline{\overline{e}}$ be basis of $D$ module $A_2$. Let $\overline{\overline{I}}$ be the basis of left $A_2 \otimes A_2$-module $\mathcal{L}(D; A_1 \to A_2)$.*[6.6]

6.4.1.1: *The map*

$$f : A_1 \to A_2$$

*has the following expansion*

$$f = f^k \circ I_k \tag{6.4.1}$$

*where*

$$f^k = f^k_{s_k \cdot 0} \otimes f^k_{s_k \cdot 1} \qquad f^k \in A_2 \otimes A_2$$

6.4.1.2: *The map $f$ has the standard representation*

$$f = f^{k \cdot ij}(e_i \otimes e_j) \circ I_k = f^{k \cdot ij} e_i I_k e_j \tag{6.4.2}$$

PROOF. Since $\overline{\overline{I}}$ is the basis of left $A_2 \otimes A_2$-module $\mathcal{L}(D; A_1 \to A_2)$, then according to the definition [4]-2.7.1 and the theorem 4.1.4, there exists expansion

$$f = f^k \circ I_k \qquad f^k \in A_2 \otimes A_2 \tag{6.4.3}$$

of the linear map $f$ with respect to the basis $\overline{\overline{I}}$. According to the definition (3.3.20),

$$f^k = f^k_{s_k \cdot 0} \otimes f^k_{s_k \cdot 1} \tag{6.4.4}$$

The equality (6.4.1) follows from equalities (6.4.3), (6.4.4). According to theorem 4.5.6, the standard representation of the tensor $f^k$ has form

$$f^k = f^{k \cdot ij} e_i \otimes e_j \tag{6.4.5}$$

The equation (6.4.2) follows from equations (6.4.1), (6.4.5). □

---

[6.5] Generating set of representation is defined in definition [4]-2.6.5.

[6.6] If $D$-module $A_1$ or $D$-module $A_2$ is not free $D$-nodule, then we may consider the set

$$\overline{\overline{I}} = \{I_k \in \mathcal{L}(D; A_1 \to A_2) : k = 1, ..., n\}$$

of linear independent linear maps. The theorem is true for any linear map

$$f : A_1 \to A_2$$

generated by the set of linear maps $\overline{\overline{I}}$.

THEOREM 6.4.2. *Let $A_1$ be free $D$-module. Let $A_2$ be free associative $D$-algebra. Let $\overline{\overline{I}}$ be the basis of left $A_2\otimes A_2$-module $\mathcal{L}(D;A_1\to A_2)$. For any map $I_k\in\overline{\overline{I}}$, there exists set of linear maps*

$$I_k^l : A_1\otimes A_1\to A_2\otimes A_2$$

*of $D$-module $A_1\otimes A_1$ into $D$-module $A_2\otimes A_2$ such that*

$$I_k\circ a\circ x=(I_k^l\circ a)\circ I_l\circ x \tag{6.4.6}$$

*The map $I_k^l$ is called* **conjugation transformation**.

PROOF. According to the theorem 6.2.1, for any tensor $a\in A_1\otimes A_1$, the map

$$x\to I_k\circ a\circ x \tag{6.4.7}$$

is linear. According to the statement 6.4.1.1, there exists expansion

$$I_k\circ a\circ x=b^l\circ I_l\circ x\qquad b^l\in A_2\times A_2 \tag{6.4.8}$$

Let

$$b^l=I_k^l\circ a \tag{6.4.9}$$

The equality (6.4.6) follows from equalities (6.4.8), (6.4.9). From equalities

$$\begin{aligned}(I_k^l\circ(a_1+a_2))\circ I_l\circ x&=I_k\circ(a_1+a_2)\circ x\\&=I_k\circ a_1\circ x+I_k\circ a_2\circ x\\&=(I_k^l\circ a_1)\circ I_l\circ x+(I_k^l\circ a_1)\circ I_l\circ x\end{aligned}$$

$$\begin{aligned}(I_k^l\circ(da))\circ I_l\circ x&=I_k\circ(da)\circ x=I_k\circ(d(a\circ x))\\&=d(I_k\circ a\circ x)=d((I_k^l\circ a)\circ I_l\circ x)\\&=(d(I_k^l\circ a))\circ I_l\circ x\end{aligned}$$

it follows that the map $I_k^l$ is linear map. □

THEOREM 6.4.3. *Let $A_1$ be free $D$-module. Let $A_2$, $A_2$ be free associative $D$-algebras. Let $\overline{\overline{I}}$ be the basis of left $A_2\otimes A_2$-module $\mathcal{L}(D;A_1\to A_2)$. Let $\overline{\overline{J}}$ be the basis of left $A_3\otimes A_3$-module $\mathcal{L}(D;A_2\to A_3)$.*

6.4.3.1: *The set of maps*

$$\overline{\overline{K}}=\{K_{lk}:K_{lk}=J_l\circ I_k,J_l\in\overline{\overline{J}},I_k\in\overline{\overline{I}}\} \tag{6.4.10}$$

*is the basis of left $A_3\otimes A_3$-module $\mathcal{L}(D;A_1\to A_2\to A_3)$.*

6.4.3.2: *Let*

$$f=f^k\circ I_k \tag{6.4.11}$$

*be expansion of linear map*

$$f:A_1\to A_2$$

*with respect to the basis $\overline{\overline{I}}$. Let*

$$g=g^l\circ J_l \tag{6.4.12}$$

*be expansion of linear map*

$$g:A_2\to A_3$$

*with respect to the basis* $\overline{\overline{J}}$. *Then linear map*

$$h = g \circ f \tag{6.4.13}$$

*has expansion*

$$h = h^{lk} \circ K_{lk} \tag{6.4.14}$$

*with respect to the basis* $\overline{\overline{K}}$ *where*

$$h^{lk} = g^l \circ (J^k_m \circ f^m) \tag{6.4.15}$$

Proof. The equality

$$h \circ a = g \circ f \circ a = g^l \circ J_l \circ f^k \circ I_k \circ a \tag{6.4.16}$$

follows from equalities (6.4.11), (6.4.12), (6.4.13). The equality

$$\begin{aligned} h \circ a = g \circ f \circ a &= g^l \circ (J^m_l \circ f^k) \circ J_m \circ I_k \circ a \\ &= g^l \circ (J^m_l \circ f^k) \circ K_{mk} \circ a \end{aligned} \tag{6.4.17}$$

follows from equalities (6.4.10), (6.4.16) and from the theorem 6.4.2. From the equality (6.4.17) it follows that set of maps $\overline{\overline{K}}$ generates left $A_3 \otimes A_3$-module $\mathcal{L}(D; A_1 \to A_2 \to A_3)$. From the equality

$$a^{lk} K_{lk} = (a^{lk} \circ J_l) \circ I_k = 0$$

it follows that

$$a^{lk} \circ J_l = 0$$

and, therefore, $a^{lk} = 0$. Therefore, the set $\overline{\overline{K}}$ is the basis of left $A_3 \otimes A_3$-module $\mathcal{L}(D; A_1 \to A_2 \to A_3)$. □

Theorem 6.4.4. *Let $A$ be free associative $D$-algebra. Let left $A \otimes A$-module $\mathcal{L}(D; A \to A)$ is generated by the identity map $I_0 = \delta$. Let*

$$f = f_{s\cdot 0} \otimes f_{s\cdot 1} \tag{6.4.18}$$

*be expansion of linear map*

$$f : A \to A$$

*Let*

$$g = g_{t\cdot 0} \otimes g_{t\cdot 1} \tag{6.4.19}$$

*be expansion of linear map*

$$g : A \to A$$

*Then linear map*

$$h = g \circ f \tag{6.4.20}$$

*has expansion*

$$h = h_{ts\cdot 0} \otimes h_{ts\cdot 1} \tag{6.4.21}$$

*where*

$$\begin{aligned} h_{ts\cdot 0} &= g_{t\cdot 0} f_{s\cdot 0} \\ h_{ts\cdot 1} &= f_{s\cdot 1} g_{t\cdot 1} \end{aligned} \tag{6.4.22}$$

Proof. The equality

$$\begin{aligned} h\circ a &= g\circ f\circ a \\ &= (g_{t\cdot 0}\otimes g_{t\cdot 1})\circ(f_{s\cdot 0}\otimes f_{s\cdot 1})\circ a \\ &= (g_{t\cdot 0}\otimes g_{t\cdot 1})\circ(f_{s\cdot 0}af_{s\cdot 1}) \\ &= g_{t\cdot 0}f_{s\cdot 0}af_{s\cdot 1}g_{t\cdot 1} \end{aligned} \tag{6.4.23}$$

follows from equalities (6.4.18), (6.4.19), (6.4.20). The equality (6.4.22) follows from the equality (6.4.23). □

Theorem 6.4.5. *Let $\overline{\overline{e}}_1$ be basis of the free finite dimensional $D$-module $A_1$. Let $\overline{\overline{e}}_2$ be basis of the free finite dimensional associative $D$-algebra $A_2$. Let $C^{p}_{kl}$ be structural constants of algebra $A_2$. Let $\overline{\overline{I}}$ be the basis of left $A_2\otimes A_2$-module $\mathcal{L}(D;A_1\to A_2)$ and $I_{k\cdot}{}^{j}_{i}$ be coordinates of map $I_k$ with respect to bases $\overline{\overline{e}}_1$ and $\overline{\overline{e}}_2$. Coordinates $f^{k}_{l}$ of the map $f\in\mathcal{L}(D;A_1\to A_2)$ and its standard components $f^{k\cdot ij}$ are connected by the equation*

$$f^{k}_{l}=f^{k\cdot ij}I_{k\cdot}{}^{m}_{l}C^{p}_{im}C^{k}_{pj} \tag{6.4.24}$$

Proof. Relative to bases $\overline{\overline{e}}_1$ and $\overline{\overline{e}}_2$, linear maps $f$ and $I_k$ have form

$$f\circ x=f^{i}_{j}x^{j}e_{2\cdot i} \tag{6.4.25}$$

$$I_k\circ x=I_{k\cdot}{}^{i}_{j}x^{j}e_{2\cdot i} \tag{6.4.26}$$

The equality

$$\begin{aligned} f^{k}_{l}x^{l}e_{2\cdot k} &= f^{k\cdot ij}e_{2\cdot i}I_{k\cdot}{}^{m}_{l}x^{l}e_{2\cdot m}\overline{e}_{2\cdot j} \\ &= f^{k\cdot ij}I_{k\cdot}{}^{m}_{l}x^{l}C^{p}_{im}C^{k}_{pj}e_{2\cdot k} \end{aligned} \tag{6.4.27}$$

follows from equalities (6.4.2), (6.4.25), (6.4.26). Since vectors $\overline{e}_{2\cdot k}$ are linear independent and $x^{i}$ are arbitrary, then the equation (6.4.24) follows from the equation (6.4.27). □

Theorem 6.4.6. *Let $D$ be field. Let $\overline{\overline{e}}_1$ be basis of the free finite dimensional $D$-algebra $A_1$. Let $\overline{\overline{e}}_2$ be basis of the free finite dimensional associative $D$-algebra $A_2$. Let $C_{2\cdot}{}^{p}_{kl}$ be structural constants of algebra $A_2$. Consider matrix*

$$\mathcal{B}=\left(\mathcal{C}^{\cdot k}{}_{m\cdot ij}\right)=\left(C_{2\cdot}{}^{p}_{im}C_{2\cdot}{}^{k}_{pj}\right) \tag{6.4.28}$$

*whose rows and columns are indexed by ${}^{\cdot k}_{m}$ and ${}_{\cdot ij}$, respectively. If matrix $\mathcal{B}$ is nonsingular, then, for given coordinates of linear transformation $g^{l}_{k}$ and for map $f=\delta$, the system of linear equations (6.4.24) with standard components of this transformation $g^{kr}$ has the unique solution.*

*If matrix $\mathcal{B}$ is singular, then the equation*

$$\operatorname{rank}\begin{pmatrix}\mathcal{C}^{\cdot k}{}_{m\cdot ij} & g^{k}_{m}\end{pmatrix}=\operatorname{rank}\mathcal{C} \tag{6.4.29}$$

*is the condition for the existence of solutions of the system of linear equations (6.4.24). In such case the system of linear equations (6.4.24) has infinitely many solutions and there exists linear dependence between values $g^{k}_{m}$.*

Proof. The statement of the theorem is corollary of the theory of linear equations over field. □

Theorem 6.4.7. *Let $A$ be free finite dimensional associative algebra over the field $D$. Let $\overline{\overline{e}}$ be the basis of the algebra $A$ over the field $D$. Let $C^{p}_{kl}$ be structural constants of algebra $A$. Let matrix (6.4.28) be singular. Let the linear map $f \in \mathcal{L}(D; A \to A)$ be nonsingular. If coordinates of linear transformations $f$ and $g$ satisfy to the equation*

$$\operatorname{rank}\begin{pmatrix}\mathcal{C}\cdot^{k}_{m\cdot ij} & g^{k}_{m} & f^{k}_{m}\end{pmatrix} = \operatorname{rank}\mathcal{C} \tag{6.4.30}$$

*then the system of linear equations*

$$g^{k}_{l} = f^{m}_{l} g^{ij} C^{p}_{im} C^{k}_{pj} \tag{6.4.31}$$

*has infinitely many solutions.*

Proof. According to the equation (6.4.30) and the theorem 6.4.6, the system of linear equations

$$f^{k}_{l} = f^{ij} C^{p}_{il} C^{k}_{pj} \tag{6.4.32}$$

has infinitely many solutions corresponding to linear map

$$f = f^{ij}\overline{e}_i \otimes \overline{e}_j \tag{6.4.33}$$

According to the equation (6.4.30) and the theorem 6.4.6, the system of linear equations

$$g^{k}_{l} = g^{ij} C^{p}_{il} C^{k}_{pj} \tag{6.4.34}$$

has infinitely many solutions corresponding to linear map

$$g = g^{ij}\overline{e}_i \otimes \overline{e}_j \tag{6.4.35}$$

Maps $f$ and $g$ are generated by the map $\delta$. According to the theorem 6.3.7, the map $f$ generates the map $g$. This proves the statement of the theorem. □

Theorem 6.4.8. *Let $A$ be free finite dimensional associative algebra over the field $D$. The representation of algebra $A \otimes A$ in algebra $\mathcal{L}(D; A \to A)$ has finite* **basis** $\overline{\overline{I}}$.

6.4.8.1: *The linear map $f \in \mathcal{L}(D; A \to A)$ has form*

$$f = (a_{k\cdot s_k\cdot 0} \otimes a_{k\cdot s_k\cdot 1}) \circ I_k = \sum_{k} a_{k\cdot s_k\cdot 0} I_k a_{k\cdot s_k\cdot 1} \tag{6.4.36}$$

6.4.8.2: *Its standard representation has form*

$$f = a^{k\cdot ij}(e_i \otimes e_j) \circ I_k = a^{k\cdot ij} e_i I_k e_j \tag{6.4.37}$$

Proof. From the theorem 6.4.7, it follows that if matrix $\mathcal{B}$ is singular and the map $f$ satisfies to the equation

$$\operatorname{rank}\begin{pmatrix}\mathcal{C}\cdot^{k}_{m\cdot ij} & f^{k}_{m}\end{pmatrix} = \operatorname{rank}\mathcal{C} \tag{6.4.38}$$

then the map $f$ generates the same set of maps that is generated by the map $\delta$. Therefore, to build the basis of representation of the algebra $A \otimes A$ in the module $\mathcal{L}(D; A \to A)$, we must perform the following construction.

The set of solutions of system of equations (6.4.31) generates a free submodule $\mathcal{L}$ of the module $\mathcal{L}(D; A \to A)$. We build the basis $(\overline{h}_1, ..., \overline{h}_k)$ of the submodule $\mathcal{L}$. Then we supplement this basis by linearly independent vectors $\overline{h}_{k+1}$, ..., $\overline{h}_m$, that do not belong to the submodule $\mathcal{L}$ so that the set of vectors $\overline{h}_1$, ..., $\overline{h}_m$ forms a basis of the module $\mathcal{L}(D;A_1 \to A_2)$. The set of orbits $(A\otimes A)\circ\delta$, $(A\otimes A)\circ\overline{h}_{k+1}$, ..., $(A\otimes A)\circ\overline{h}_m$ generates the module $\mathcal{L}(D;A \to A)$. Since the set of orbits is finite, we can choose the orbits so that they do not intersect. For each orbit we can choose a representative which generates the orbit. □

EXAMPLE 6.4.9. *For complex field, the algebra $\mathcal{L}(R;C \to C)$ has basis*

$$I_0\circ z = z$$
$$I_1\circ z = \overline{z}$$

*For quaternion algebra, the algebra $\mathcal{L}(R;H \to H)$ has basis*

$$I_0\circ z = z$$

□

## 6.5. Linear Map into Nonassociative Algebra

Since the product is nonassociative, we may assume that action of $a$, $b\in A$ over the map $f$ may have form either $a(fb)$, or $(af)b$. However this assumption leads us to a rather complex structure of the linear map. To better understand how complex the structure of the linear map, we begin by considering the left and right shifts in nonassociative algebra.

THEOREM 6.5.1. *Let*

$$l(a)\circ x = ax \tag{6.5.1}$$

*be map of left shift. Then*

$$l(a)\circ l(b) = l(ab)-(a,b)_1 \tag{6.5.2}$$

*where we introduced linear map*

$$(a,b)_1\circ x = (a,b,x)$$

PROOF. From the equations (5.1.2), (6.5.1), it follows that

$$\begin{aligned}(l(a)\circ l(b))\circ x &= l(a)\circ(l(b)\circ x)\\ &= a(bx) = (ab)x-(a,b,x)\\ &= l(ab)\circ x-(a,b)_1\circ x\end{aligned} \tag{6.5.3}$$

The equation (6.5.2) follows from equation (6.5.3). □

THEOREM 6.5.2. *Let*

$$r(a)\circ x = xa \tag{6.5.4}$$

*be map of right shift. Then*

$$r(a)\circ r(b) = r(ba)+(b,a)_2 \tag{6.5.5}$$

*where we introduced linear map*

$$(b,a)_2\circ x = (x,b,a)$$

Proof. From the equations (5.1.2), (6.5.4) it follows that

$$(r(a)\circ r(b))\circ x = r(a)\circ(r(b)\circ x) = (xb)a = x(ba)+(x,b,a) = r(ba)\circ x+(x,b,a) \tag{6.5.6}$$

The equation (6.5.5) follows from equation (6.5.6). □

Let

$$f:A\rightarrow A \quad f=(ax)b$$

be linear map of the algebra $A$. According to the theorem 6.1.8, the map

$$g:A\rightarrow A \quad g=(cf)d$$

is also a linear map. However, it is not obvious whether we can write the map $g$ as a sum of terms of type $(ax)b$ and $a(xb)$.

If $A$ is free finite dimensional algebra, then we can assume that the linear map has the standard representation like[6.7]

$$f\circ x = f^{ij}\ (\overline{e}_i x)\overline{e}_j \tag{6.5.8}$$

In this case we can use the theorem 6.4.8 for maps into nonassociative algebra.

Theorem 6.5.3. *Let $\overline{\overline{e}}_1$ be basis of the free finite dimensional $D$-algebra $A_1$. Let $\overline{\overline{e}}_2$ be basis of the free finite dimensional nonassociative $D$-algebra $A_2$. Let $C_{2\cdot}{}^{p}_{kl}$ be structural constants of algebra $A_2$. Let the map*

$$g=a\circ f \tag{6.5.9}$$

*generated by the map $f\in\mathcal{L}(D;A_1\rightarrow A_2)$ through the tensor $a\in A_2\otimes A_2$, has the standard representation*

$$g=a^{ij}(\overline{e}_i\otimes\overline{e}_j)\circ f=a^{ij}(\overline{e}_i f)\overline{e}_j \tag{6.5.10}$$

*Coordinates of the map (6.5.9) and its standard components are connected by the equation*

$$g^k_l = f^m_l g^{ij} C_{2\cdot}{}^{p}_{im} C_{2\cdot}{}^{k}_{pj} \tag{6.5.11}$$

Proof. Relative to bases $\overline{\overline{e}}_1$ and $\overline{\overline{e}}_2$, linear maps $f$ and $g$ have form

$$f\circ x = f^i_j x^j \overline{e}_{2\cdot i} \tag{6.5.12}$$

$$g\circ x = g^i_j x^j \overline{e}_{2\cdot i} \tag{6.5.13}$$

From equations (6.5.12), (6.5.13), (6.5.10) it follows that

$$g^k_l x^l \overline{e}_{2\cdot k} = a^{ij}(\overline{e}_{2\cdot i}(f^m_l x^l \overline{e}_{2\cdot m}))\overline{e}_{2\cdot j} = a^{ij} f^m_l x^l C_{2\cdot}{}^{p}_{im} C_{2\cdot}{}^{k}_{pj} \overline{e}_{2\cdot k} \tag{6.5.14}$$

[6.7] The choice is arbitrary. We may consider the standard representation like

$$f\circ x = f^{ij}\overline{e}_i(x\overline{e}_j)$$

Then the equation (6.5.11) has form

$$g^k_l = f^m_l g^{ij} C_{2\cdot}{}^{k}_{ip} C_{2\cdot}{}^{p}_{mj} \tag{6.5.7}$$

I chose the expression (6.5.8) because order of the factors corresponds to the order chosen in the theorem 6.4.8.

Since vectors $\overline{e}_{2\cdot k}$ are linear independent and $x^i$ are arbitrary, then the equation (6.5.11) follows from the equation (6.5.14). □

THEOREM 6.5.4. *Let $A$ be free finite dimensional nonassociative algebra over the ring $D$. The representation of algebra $A\otimes A$ in algebra $\mathcal{L}(D;A\rightarrow A)$ has finite basis $\overline{\overline{I}}$.*

(1) *The linear map $f\in\mathcal{L}(D;A\rightarrow A)$ has form*

$$f=(a_{k\cdot s_k\cdot 0}\otimes a_{k\cdot s_k\cdot 1})\circ I_k=(a_{k\cdot s_k\cdot 0}I_k)a_{k\cdot s_k\cdot 1} \tag{6.5.15}$$

(2) *Its standard representation has form*

$$f=a^{k\cdot ij}(\overline{e}_i\otimes\overline{e}_j)\circ I_k=a^{k\cdot ij}(\overline{e}_iI_k)\overline{e}_j \tag{6.5.16}$$

PROOF. Consider matrix (6.4.28). If matrix $\mathcal{B}$ is nonsingular, then, for given coordinates of linear transformation $g^l_k$ and for map $f=\delta$, the system of linear equations (6.5.11) with standard components of this transformation $g^{kr}$ has the unique solution. If matrix $\mathcal{B}$ is singular, then according to the theorem 6.4.8 there exists finite basis $\overline{\overline{I}}$ generating the set of linear maps. □

Unlike the case of an associative algebra, the set of generators $I$ in the theorem 6.5.4 is not minimal. From the equation (6.5.2) it follows that the equation (6.3.12) does not hold. Therefore, orbits of maps $I_k$ do not generate an equivalence relation in the algebra $L(A;A)$. Since we consider only maps like $(aI_k)b$, then it is possible that for $k\neq l$ the map $I_k$ generates the map $I_l$, if we consider all possible operations in the algebra $A$. Therefore, the set of generators $I_k$ of nonassociative algebra $A$ does not play such a critical role as conjugation in complex field. The answer to the question of how important it is the map $I_k$ in nonassociative algebra requires additional research.

## 6.6. Polylinear Map into Associative Algebra

THEOREM 6.6.1. *Let $A_1$, ..., $A_n$, $A$ be associative $D$-algebras. Let*

$$f_i\in\mathcal{L}(D;A_i\rightarrow A)\quad i=1,...,n$$

$$a_j\in A\quad j=0,...,n$$

*For given transposition $\sigma$ of $n$ variables, the map*

$$\begin{aligned}&((a_0,...,a_n,\sigma)\circ(f_1,...,f_n))\circ(x_1,...,x_n)\\ &=(a_0\sigma(f_1)a_1...a_{n-1}\sigma(f_n)a_n)\circ(x_1,...,x_n)\\ &=a_0\sigma(f_1\circ x_1)a_1...a_{n-1}\sigma(f_n\circ x_n)a_n\end{aligned} \tag{6.6.1}$$

*is $n$-linear map into algebra $A$.*

Proof. The statement of theorem follows from chains of equations

$$\begin{aligned}&((a_0,...,a_n,\sigma)\circ(f_1,...,f_n))\circ(x_1,...,x_i+y_i,...,x_n)\\ =&a_0\sigma(f_1\circ x_1)a_1...\sigma(f_i\circ(x_i+y_i))...a_{n-1}\sigma(f_n\circ x_n)a_n\\ =&a_0\sigma(f_1\circ x_1)a_1...\sigma(f_i\circ x_i+f_i\circ y_i)...a_{n-1}\sigma(f_n\circ x_n)a_n\\ =&a_0\sigma(f_1\circ x_1)a_1...\sigma(f_i\circ x_i)...a_{n-1}\sigma(f_n\circ x_n)a_n\\ +&a_0\sigma(f_1\circ x_1)a_1...\sigma(f_i\circ y_i)...a_{n-1}\sigma(f_n\circ x_n)a_n\\ =&((a_0,...,a_n,\sigma)\circ(f_1,...,f_n))\circ(x_1,...,x_i,...,x_n)\\ +&((a_0,...,a_n,\sigma)\circ(f_1,...,f_n))\circ(x_1,...,y_i,...,x_n)\end{aligned}$$

$$\begin{aligned}&((a_0,...,a_n,\sigma)\circ(f_1,...,f_n))\circ(x_1,...,px_i,...,x_n)\\ =&a_0\sigma(f_1\circ x_1)a_1...\sigma(f_i\circ(px_i))...a_{n-1}\sigma(f_n\circ x_n)a_n\\ =&a_0\sigma(f_1\circ x_1)a_1...\sigma(p(f_i\circ x_i))...a_{n-1}\sigma(f_n\circ x_n)a_n\\ =&p(a_0\sigma(f_1\circ x_1)a_1...\sigma(f_i\circ x_i)...a_{n-1}\sigma(f_n\circ x_n)a_n)\\ =&p(((a_0,...,a_n,\sigma)\circ(f_1,...,f_n))\circ(x_1,...,x_i,...,x_n))\end{aligned}$$

□

In the equation (6.6.1), as well as in other expressions of polylinear map, we have convention that map $f_i$ has variable $x_i$ as argument.

Theorem 6.6.2. *Let $A_1$, ..., $A_n$, $A$ be associative $D$-algebras. For given set of maps*

$$f_i\in\mathcal{L}(D;A_i\to A)\quad i=1,...,n$$

*the map*

$$h:A^{n+1}\to\mathcal{L}(D;A_1\times...\times A_n\to A)$$

*defined by equation*

$$(a_0,...,a_n,\sigma)\circ(f_1,...,f_n)=a_0\sigma(f_1)a_1...a_{n-1}\sigma(f_n)a_n$$

*is $n+1$-linear map into $D$-module $\mathcal{L}(D;A_1\times...\times A_n\to A)$.*

Proof. The statement of theorem follows from chains of equations

$$\begin{aligned}&((a_0,...,a_i+b_i,...a_n,\sigma)\circ(f_1,...,f_n))\circ(x_1,...,x_n)\\ =&a_0\sigma(f_1\circ x_1)a_1...(a_i+b_i)...a_{n-1}\sigma(f_n\circ x_n)a_n\\ =&a_0\sigma(f_1\circ x_1)a_1...a_i...a_{n-1}\sigma(f_n\circ x_n)a_n+a_0\sigma(f_1\circ x_1)a_1...b_i...a_{n-1}\sigma(f_n\circ x_n)a_n\\ =&((a_0,...,a_i,...,a_n,\sigma)\circ(f_1,...,f_n))\circ(x_1,...,x_n)\\ +&((a_0,...,b_i,...,a_n,\sigma)\circ(f_1,...,f_n))\circ(x_1,...,x_n)\\ =&((a_0,...,a_i,...,a_n,\sigma)\circ(f_1,...,f_n)+(a_0,...,b_i,...,a_n,\sigma)\circ(f_1,...,f_n))\circ(x_1,...,x_n)\end{aligned}$$

$$\begin{aligned}&((a_0,...,pa_i,...a_n,\sigma)\circ(f_1,...,f_n))\circ(x_1,...,x_n)\\ =&a_0\sigma(f_1\circ x_1)a_1...pa_i...a_{n-1}\sigma(f_n\circ x_n)a_n\\ =&p(a_0\sigma(f_1\circ x_1)a_1...a_i...a_{n-1}\sigma(f_n\circ x_n)a_n)\\ =&p(((a_0,...,a_i,...,a_n,\sigma)\circ(f_1,...,f_n))\circ(x_1,...,x_n))\\ =&(p((a_0,...,a_i,...,a_n,\sigma)\circ(f_1,...,f_n)))\circ(x_1,...,x_n)\end{aligned}$$

□

Theorem 6.6.3. *Let $A_1$, ..., $A_n$, $A$ be associative $D$-algebras. For given set of maps*

$$f_i\in\mathcal{L}(D;A_i\to A)\quad i=1,...,n$$

*there exists linear map*

$$h:A^{\otimes n+1}\times S_n\to\mathcal{L}(D;A_1\times...\times A_n\to A)$$

*defined by the equation*

$$\begin{aligned}(a_0\otimes...\otimes a_n,\sigma)\circ(f_1,...,f_n)&=(a_0,...,a_n,\sigma)\circ(f_1,...,f_n)\\&=a_0\sigma(f_1)a_1...a_{n-1}\sigma(f_n)a_n\end{aligned}\tag{6.6.2}$$

Proof. The statement of the theorem is corollary of theorems 4.5.4, 6.6.2. □

Theorem 6.6.4. *Let $A_1$, ..., $A_n$, $A$ be associative $D$-algebras. For given tensor $a\in A^{\otimes n+1}$ and given transposition $\sigma\in S_n$ the map*

$$h:\prod_{i=1}^n\mathcal{L}(D;A_i\to A)\to\mathcal{L}(D;A_1\times...\times A_n\to A)$$

*defined by equation*

$$(a_0\otimes...\otimes a_n,\sigma)\circ(f_1,...,f_n)=a_0\sigma(f_1)a_1...a_{n-1}\sigma(f_n)a_n$$

*is $n$-linear map into $D$-module $\mathcal{L}(D;A_1\times...\times A_n\to A)$.*

Proof. The statement of theorem follows from chains of equations

$$\begin{aligned}&((a_0\otimes...\otimes a_n,\sigma)\circ(f_1,...,f_i+g_i,...,f_n))\circ(x_1,...,x_n)\\=&(a_0\sigma(f_1)a_1...\sigma(f_i+g_i)...a_{n-1}\sigma(f_n)a_n)\circ(x_1,...,x_n)\\=&a_0\sigma(f_1\circ x_1)a_1...\sigma((f_i+g_i)\circ x_i)...a_{n-1}\sigma(f_n\circ x_n)a_n\\=&a_0\sigma(f_1\circ x_1)a_1...\sigma(f_i\circ x_i+g_i\circ x_i)...a_{n-1}\sigma(f_n\circ x_n)a_n\\=&a_0\sigma(f_1\circ x_1)a_1...\sigma(f_i\circ x_i)...a_{n-1}\sigma(f_n\circ x_n)a_n\\+&a_0\sigma(f_1\circ x_1)a_1...\sigma(g_i\circ x_i)...a_{n-1}\sigma(f_n\circ x_n)a_n\\=&(a_0\sigma(f_1)a_1...\sigma(f_i)...a_{n-1}\sigma(f_n)a_n)\circ(x_1,...,x_n)\\+&(a_0\sigma(f_1)a_1...\sigma(g_i)...a_{n-1}\sigma(f_n)a_n)\circ(x_1,...,x_n)\\=&((a_0\otimes...\otimes a_n,\sigma)\circ(f_1,...,f_i,...,f_n))\circ(x_1,...,x_n)\\+&((a_0\otimes...\otimes a_n,\sigma)\circ(f_1,...,g_i,...,f_n))\circ(x_1,...,x_n)\\=&((a_0\otimes...\otimes a_n,\sigma)\circ(f_1,...,f_i,...,f_n)\\+&(a_0\otimes...\otimes a_n,\sigma)\circ(f_1,...,g_i,...,f_n))\circ(x_1,...,x_n)\end{aligned}$$

$$\begin{aligned}&((a_0\otimes...\otimes a_n,\sigma)\circ(f_1,...,pf_i,...,f_n))\circ(x_1,...,x_n)\\=&(a_0\sigma(f_1)a_1,...\sigma(pf_i)...a_{n-1}\sigma(f_n)a_n)\circ(x_1,...,x_n)\\=&a_0\sigma(f_1\circ x_1)a_1,...\sigma((pf_i)\circ x_i)...a_{n-1}\sigma(f_n\circ x_n)a_n\\=&a_0\sigma(f_1\circ x_1)a_1,...\sigma(p(f_i\circ x_i))...a_{n-1}\sigma(f_n\circ x_n)a_n\\=&p(a_0\sigma(f_1\circ x_1)a_1,...\sigma(f_i\circ x_i)...a_{n-1}\sigma(f_n\circ x_n)a_n)\\=&p(((a_0\otimes...\otimes a_n,\sigma)\circ(f_1,...,f_i,...,f_n))\circ(x_1,...,x_n))\\=&(p((a_0\otimes...\otimes a_n,\sigma)\circ(f_1,...,f_i,...,f_n)))\circ(x_1,...,x_n)\end{aligned}$$

□

Theorem 6.6.5. *Let $A_1$, ..., $A_n$, $A$ be associative $D$-algebras. For given tensor $a \in A^{\otimes n+1}$ and given transposition $\sigma \in S_n$ there exists linear map*

$$h : \mathcal{L}(D; A_1 \to A) \otimes ... \otimes \mathcal{L}(D; A_n \to A) \to \mathcal{L}(D; A_1 \times ... \times A_n \to A)$$

*defined by the equation*

$$(a_0 \otimes ... \otimes a_n, \sigma) \circ (f_1 \otimes ... \otimes f_n) = (a_0 \otimes ... \otimes a_n, \sigma) \circ (f_1, ..., f_n) \tag{6.6.3}$$

Proof. The statement of the theorem is corollary of theorems 4.5.4, 6.6.4. □

Theorem 6.6.6. *Let $A$ be associative $D$-algebra. Polylinear map*

$$f : A^n \to A, a = f \circ (a_1, ..., a_n) \tag{6.6.4}$$

*generated by maps $I_{s\cdot 1}$, ..., $I_{s\cdot n} \in \mathcal{L}(D; A \to A)$ has form*

$$a = f^n_{s\cdot 0}\ \sigma_s(I_{s\cdot 1} \circ a_1)\ f^n_{s\cdot 1}\ ...\ \sigma_s(I_{s\cdot n} \circ a_n)\ f^n_{s\cdot n} \tag{6.6.5}$$

*where $\sigma_s$ is a transposition of set of variables $\{a_1, ..., a_n\}$*

$$\sigma_s = \begin{pmatrix} a_1 & ... & a_n \\ \sigma_s(a_1) & ... & \sigma_s(a_n) \end{pmatrix}$$

Proof. We prove statement by induction on $n$.

When $n = 1$ the statement of theorem follows from the statement 6.4.8.1. In such case we may identify[6.8]

$$f^1_{s\cdot p} = f_{s\cdot p} \quad p = 0, 1$$

Let statement of theorem be true for $n = k-1$. Then it is possible to represent map (6.6.4) as

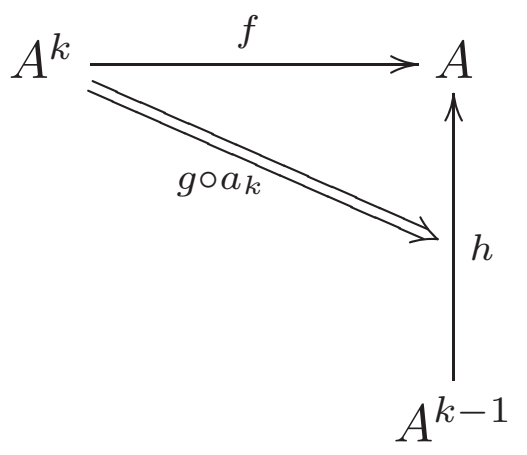

$$a = f \circ (a_1, ..., a_k) = (g \circ a_k) \circ (a_1, ..., a_{k-1})$$

According to statement of induction polylinear map $h$ has form

$$a = h^{k-1}_{t\cdot 0}\ \sigma_t(I_{1\cdot t} \circ a_1)\ h^{k-1}_{t\cdot 1}\ ...\ \sigma_t(I_{k-1\cdot t} \circ a_{k-1})\ h^{k-1}_{t\cdot k-1}$$

According to construction $h = g \circ a_k$. Therefore, expressions $h_{t\cdot p}$ are functions of $a_k$. Since $g \circ a_k$ is linear map of $a_k$, then only one expression $h_{t\cdot p}$ is linear map of $a_k$, and rest expressions $h_{t\cdot q}$ do not depend on $a_k$.

[6.8] In representation (6.6.5) we will use following rules.

- If range of any index is set consisting of one element, then we will omit corresponding index.
- If $n = 1$, then $\sigma_s$ is identical transformation. We will not show such transformation in the expression.

Without loss of generality, assume $p=0$. According to the equation (6.3.10) for given $t$

$$h^{k-1}_{t\cdot 0}=g_{tr\cdot 0}\ I_{k\cdot r}\circ a_k\ g_{tr\cdot 1}$$

Assume $s=tr$. Let us define transposition $\sigma_s$ according to rule

$$\sigma_s=\sigma_{tr}=\begin{pmatrix}a_k & a_1 & ... & a_{k-1}\\ a_k & \sigma_t(a_1) & ... & \sigma_t(a_{k-1})\end{pmatrix}$$

Suppose

$$f^k_{tr\cdot q+1}=h^{k-1}_{t\cdot q}\qquad q=1,...,k-1$$
$$f^k_{tr\cdot q}=g_{tr\cdot q}\qquad q=0,1$$

We proved step of induction. □

Definition 6.6.7. *Expression* $f s\cdot p^n$ *in equation* (6.6.5) *is called* **component of polylinear map** $f$. □

Theorem 6.6.8. *Consider* $D$*-algebras* $A$. *A linear map*

$$h:A^{\otimes n+1}\times S_n\rightarrow {}^*\mathcal{L}(D;A^n\rightarrow A)$$

*defined by the equation*

(6.6.6)
$$(a_0\otimes...\otimes a_n,\sigma)\circ(f_1\otimes...\otimes f_n)=a_0\sigma(f_1)a_1...a_{n-1}\sigma(f_n)a_n$$
$$a_0,...,a_n\in A\quad \sigma\in S_n\quad f_1,...,f_n\in\mathcal{L}(D;A_n\rightarrow A)$$

*is representation*[6.9] *of algebra* $A^{\otimes n+1}\times S^n$ *in* $D$*-module* $\mathcal{L}(D;A^n\rightarrow A)$.

Proof. According to the theorems 6.4.8, 6.6.6, we can represent $n$-linear map as sum of terms (6.6.1), where $f_i$, $i=1$, ..., $n$, are generators of representation (6.3.3). Let us write the term $s$ of the expression (6.6.5) as

(6.6.7)
$$b_1\sigma(I_{1\cdot s}\circ x_1)c_1b_2...c_{n-1}b_n\sigma(I_{n\cdot s}\circ x_n)c_n$$

where

$$b_1=f^n_{s\cdot 0}\quad b_2=...=b_n=e\quad c_1=f^n_{s\cdot 1}\quad ...\quad c_n=f^n_{s\cdot n}$$

Let us assume

$$f_i=\sigma^{-1}(b_i)I_{i\cdot s}\sigma^{-1}(c_i)\quad i=1,...,n$$

in equation (6.6.7). Therefore, according to theorem 6.6.3, map (6.6.6) is transformation of module $\mathcal{L}(D;A^n\rightarrow A)$. For a given tensor $c\in A^{\otimes n+1}$ and given transposition $\sigma\in S_n$, a transformation $h(c,\sigma)$ is a linear transformation of module $\mathcal{L}(D;A^n\rightarrow A)$ according to the theorem 6.6.5. According to the theorem 6.6.3, map (6.6.6) is linear map. According to the definition [4]-2.1.2, map (6.6.6) is a representation of the algebra $A^{\otimes n+1}\times S^n$ in the module $\mathcal{L}(D;A^n\rightarrow A)$. □

Theorem 6.6.9. *Consider* $D$*-algebra* $A$. *A representation*

$$h:A^{\otimes n+1}\times S_n\rightarrow {}^*\mathcal{L}(D;A^n\rightarrow A)$$

*of algebra* $A^{\otimes n+1}$ *in module* $\mathcal{L}(D;A^n\rightarrow A)$ *defined by the equation* (6.6.6) *allows us to identify tensor* $d\in A^{\otimes n+1}$ *and transposition* $\sigma\in S^n$ *with map*

(6.6.8)
$$(d,\sigma)\circ(f_1,...,f_n)\quad f_i=\delta\in\mathcal{L}(D;A\rightarrow A)$$

[6.9] See the definition of representation of $\Omega$-algebra in the definition [4]-2.1.2.

*where* $\delta \in \mathcal{L}(D; A \to A)$ *is identity map.*

PROOF. If we assume $f_i = \delta, \ \ d = a_0 \otimes ... \otimes a_n$ in the equation (6.6.2), then the equation (6.6.2) gets form

$$\begin{aligned}((a_0 \otimes ... \otimes a_n, \sigma) \circ (\delta, ..., \delta)) \circ (x_1, ..., x_n) &= a_0 \ (\delta \circ x_1) \ ... \ (\delta \circ x_n) \ a_n \\ &= a_0 \ x_1 \ ... \ x_n \ a_n\end{aligned} \tag{6.6.9}$$

If we assume

$$\begin{aligned}&((a_0 \otimes ... \otimes a_n, \sigma) \circ (\delta, ..., \delta)) \circ (x_1, ..., x_n) \\ =& (a_0 \otimes ... \otimes a_n, \sigma) \circ (\delta \circ x_1, ..., \delta \circ x_n) \\ =& (a_0 \otimes ... \otimes a_n, \sigma) \circ (x_1, ..., x_n)\end{aligned} \tag{6.6.10}$$

then comparison of equations (6.6.9) and (6.6.10) gives a basis to identify the action of the tensor $d = a_0 \otimes ... \otimes a_n$ and transposition $\sigma \in S^n$ with map (6.6.8). □

Instead of notation $(a_0 \otimes ...a_n, \sigma)$, we also use notation

$$a_0 \otimes_{\sigma(1)} ... \otimes_{\sigma(n)} a_n$$

when we want to show order of arguments in expression. For instance, the following expressions are equivalent

$$(a_0 \otimes a_1 \otimes a_2 \otimes a_3, (2, 1, 3)) \circ (x_1, x_2, x_3) = a_0 x_2 a_1 x_1 a_2 x_3 a_3$$
$$(a_0 \otimes_2 a_1 \otimes_1 a_2 \otimes_3 a_3) \circ (x_1, x_2, x_3) = a_0 x_2 a_1 x_1 a_2 x_3 a_3$$

## 6.7. Polylinear Map into Free Finite Dimensional Associative Algebra

THEOREM 6.7.1. *Let* $A$ *be free finite dimensional associative algebra over the ring* $D$. *Let* $\overline{\overline{I}}$ *be basis of algebra* $\mathcal{L}(D; A \to A)$. *Let* $\overline{\overline{e}}$ *be the basis of the algebra* $A$ *over the ring* $D$. **Standard representation of polylinear map** *into associative algebra has form*

$$f \circ (a_1, ..., a_n) = f^{i_0 ... i_n}_{t \cdot k_1 ... k_n} \ \overline{e}_{i_0} \ \sigma_t(I_{k_1} \circ a_1) \ \overline{e}_{i_1} \ ... \ \sigma_t(I_{k_n} \circ a_n) \ \overline{e}_{i_n} \tag{6.7.1}$$

*Index* $t$ *enumerates every possible transpositions* $\sigma_t$ *of the set of variables* $\{a_1, ..., a_n\}$. *Expression* $f^{i_0 ... i_n}_{t \cdot k_1 ... k_n}$ *in equation* (6.7.1) *is called* **standard component of polylinear map** $f$.

PROOF. We change index $s$ in the equation (6.6.5) so as to group the terms with the same set of generators $I_k$. Expression (6.6.5) gets form

$$a = f^n_{k_1 ... k_n \cdot s \cdot 0} \ \sigma_s(I_{k_1 \cdot s} \circ a_1) \ f^n_{k_1 ... k_n \cdot s \cdot 1} \ ... \ \sigma_s(I_{k_n \cdot s} \circ a_n) \ f^n_{k_1 ... k_n \cdot s \cdot n} \tag{6.7.2}$$

We assume that the index $s$ takes values depending on $k_1$, ..., $k_n$. Components of polylinear map $f$ have expansion

$$f^n_{k_1 ... k_n \cdot s \cdot p} = \overline{e}_i f^{ni}_{k_1 ... k_n \cdot s \cdot p} \tag{6.7.3}$$

relative to basis $\overline{\overline{e}}$. If we substitute (6.7.3) into (6.6.5), we get

$$a = f^{nj_1}_{k_1 ... k_n \cdot s \cdot 0} \ \overline{e}_{j_1} \ \sigma_s(I_{k_1} \circ a_1) \ f^{nj_2}_{k_1 ... k_n \cdot s \cdot 1} \ \overline{e}_{j_2} \ ... \ \sigma_s(I_{k_n} \circ a_n) \ f^{nj_n}_{k_1 ... k_n \cdot s \cdot n} \ \overline{e}_{j_n} \tag{6.7.4}$$

Let us consider expression

$$f^{j_0 ... j_n}_{t \cdot k_1 ... k_n} = f^{nj_1}_{k_1 ... k_n \cdot s \cdot 0} \ ... f^{nj_n}_{k_1 ... k_n \cdot s \cdot n} \tag{6.7.5}$$

The right-hand side is supposed to be the sum of the terms with the index $s$, for which the transposition $\sigma_s$ is the same. Each such sum has a unique index $t$. If we substitute expression (6.7.5) into equation (6.7.4) we get equation (6.7.1). □

THEOREM 6.7.2. *Let $A$ be free finite dimensional associative algebra over the ring $D$. Let $\overline{\overline{e}}$ be the basis of the algebra $A$ over the ring $D$. Polylinear map (6.6.4) can be represented as $D$-valued form of degree $n$ over ring $D$*[6.10]

$$f(a_1,...,a_n)=a_1^{i_1}...a_n^{i_n}f_{i_1...i_n} \tag{6.7.6}$$

*where*

$$\begin{aligned} a_j&=\overline{e}_i a_j^i\\ f_{i_1...i_n}&=f\circ(\overline{e}_{i_1},...,\overline{e}_{i_n})\end{aligned} \tag{6.7.7}$$

PROOF. According to the definition 6.1.2, the equation (6.7.6) follows from the chain of equations

$$f\circ(a_1,...,a_n)=f\circ(\overline{e}_{i_1}a_1^{i_1},...,\overline{e}_{i_n}a_n^{i_n})=a_1^{i_1}...a_n^{i_n}f\circ(\overline{e}_{i_1},...,\overline{e}_{i_n})$$

Let $\overline{\overline{e}}'$ be another basis. Let

$$\overline{e}'_i=\overline{e}_j h_i^j \tag{6.7.8}$$

be transformation, map basis $\overline{\overline{e}}$ into basis $\overline{\overline{e}}'$. From equations (6.7.8) and (6.7.7) it follows

$$\begin{aligned} f'_{i_1...i_n}&=f\circ(\overline{e}'_{i_1},...,\overline{e}'_{i_n})\\ &=f\circ(\overline{e}_{j_1}h_{i_1}^{j_1},...,\overline{e}'_{j_n}h_{i_n}^{j_n})\\ &=h_{i_1}^{j_1}...h_{i_n}^{j_n}f\circ(\overline{e}_{j_1},...,\overline{e}_{j_n})\\ &=h_{i_1}^{j_1}...h_{i_n}^{j_n}f_{j_1...j_n}\end{aligned} \tag{6.7.9}$$

□

Polylinear map (6.6.4) is **symmetric**, if

$$f\circ(a_1,...,a_n)=f\circ(\sigma(a_1),...,\sigma(a_n))$$

for any transposition $\sigma$ of set $\{a_1,...,a_n\}$.

THEOREM 6.7.3. *If polyadditive map $f$ is symmetric, then*

$$f_{i_1,...,i_n}=f_{\sigma(i_1),...,\sigma(i_n)} \tag{6.7.10}$$

PROOF. Equation (6.7.10) follows from equation

$$\begin{aligned} a_1^{i_1}...a_n^{i_n}f_{i_1...i_n}&=f\circ(a_1,...,a_n)\\ &=f\circ(\sigma(a_1),...,\sigma(a_n))\\ &=a_1^{i_1}...a_n^{i_n}f_{\sigma(i_1)...\sigma(i_n)}\end{aligned}$$

□

[6.10] We proved the theorem by analogy with theorem in [2], p. 107, 108

Polylinear map (6.6.4) is **skew symmetric**, if

$$f\circ(a_1,...,a_n)=|\sigma|f\circ(\sigma(a_1),...,\sigma(a_n))$$

for any transposition $\sigma$ of set $\{a_1,...,a_n\}$. Here

$$|\sigma|=\begin{cases}1 & \text{transposition }\sigma\text{ even}\\ -1 & \text{transposition }\sigma\text{ odd}\end{cases}$$

THEOREM 6.7.4. *If polylinear map $f$ is skew symmetric, then*

$$f_{i_1,...,i_n}=|\sigma|f_{\sigma(i_1),...,\sigma(i_n)}\tag{6.7.11}$$

PROOF. Equation (6.7.11) follows from equation

$$\begin{aligned}a_1^{i_1}...a_n^{i_n}f_{i_1...i_n}=&f\circ(a_1,...,a_n)\\=&|\sigma|f\circ(\sigma(a_1),...,\sigma(a_n))\\=&a_1^{i_1}...a_n^{i_n}|\sigma|f_{\sigma(i_1)...\sigma(i_n)}\end{aligned}$$

□

THEOREM 6.7.5. *Let $A$ be free finite dimensional associative algebra over the ring $D$. Let $\overline{\overline{I}}$ be basis of algebra $\mathcal{L}(D;A\rightarrow A)$. Let $\overline{\overline{e}}$ be the basis of the algebra $A$ over the ring $D$. Let polylinear over ring $D$ map $f$ be generated by set of maps $(I_{k_1},...,I_{k_n})$. Coordinates of the map $f$ and its components relative basis $\overline{\overline{e}}$ satisfy to the equation*

$$f_{l_1...l_n}=f^{i_0...i_n}_{t\cdot k_1...k_n}I_{k_1}{}^{j_1}_{\cdot l_1}...I_{k_n}{}^{j_n}_{\cdot l_n}C^{k_1}_{i_0\sigma_t(j_1)}C^{l_1}_{k_1i_1}...B^{k_n}_{l_{n-1}\sigma_t(j_n)}C^{l_n}_{k_ni_n}\overline{e}_{l_n}\tag{6.7.12}$$

$$f^p_{l_1...l_n}=f^{i_0...i_n}_{t\cdot k_1...k_n}I_{k_1}{}^{j_1}_{\cdot l_1}...I_{k_n}{}^{j_n}_{\cdot l_n}C^{k_1}_{i_0\sigma_t(j_1)}C^{l_1}_{k_1i_1}...C^{k_n}_{l_{n-1}\sigma_t(j_n)}C^{p}_{k_ni_n}\tag{6.7.13}$$

PROOF. In equation (6.7.1), we assume

$$I_{k_i}\circ a_i=\overline{e}_{j_i}I_{k_i}{}^{j_i}_{\cdot l_i}a_i^{l_i}$$

Then equation (6.7.1) gets form

(6.7.14)

$$\begin{aligned}f\circ(a_1,...,a_n)=&f^{i_0...i_n}_{t\cdot k_1...k_n}\ \overline{e}_{i_0}\sigma_t(a_1^{l_1}I_{k_1}{}^{j_1}_{\cdot l_1}\overline{e}_{j_1})\overline{e}_{i_1}...\sigma_t(a_n^{l_n}I_{k_n}{}^{j_n}_{\cdot l_n}\overline{e}_{j_n})\overline{e}_{i_n}\\=&a_1^{l_1}...a_n^{l_n}f^{i_0...i_n}_{t\cdot k_1...k_n}I_{k_1}{}^{j_1}_{\cdot l_1}...I_{k_n}{}^{j_n}_{\cdot l_n}\overline{e}_{i_0}\sigma_t(\overline{e}_{j_1})\overline{e}_{i_1}...\sigma_t(\overline{e}_{j_n})\overline{e}_{i_n}\\=&a_1^{l_1}...a_n^{l_n}f^{i_0...i_n}_{t\cdot k_1...k_n}I_{k_1}{}^{j_1}_{\cdot l_1}...I_{k_n}{}^{j_n}_{\cdot l_n}C^{k_1}_{i_0\sigma_t(j_1)}C^{l_1}_{k_1i_1}\\&...C^{k_n}_{l_{n-1}\sigma_t(j_n)}C^{l_n}_{k_ni_n}\overline{e}_{l_n}\end{aligned}$$

From equation (6.7.6) it follows that

$$f\circ(a_1,...,a_n)=\overline{e}_pf^p_{i_1...i_n}a_1^{i_1}...a_n^{i_n}\tag{6.7.15}$$

Equation (6.7.12) follows from comparison of equations (6.7.14) and (6.7.6). Equation (6.7.13) follows from comparison of equations (6.7.14) and (6.7.15). □

# Index

# Special Symbols and Notations

# Линейное отображение $D$-алгебры

Александр Клейн

Aleks_Kleyn@MailAPS.org
http://AleksKleyn.dyndns-home.com:4080/
http://sites.google.com/site/AleksKleyn/
http://arxiv.org/a/kleyn_a_1
http://AleksKleyn.blogspot.com/

Аннотация. Модуль - это эффективное представление кольца в абелевой группе. Линейное отображение модуля над коммутативным кольцом - это морфизм соответствующего представления. Это определение является центральной темой предлагаемой книги.

Чтобы рассмотреть это определение с более общей точки зрения, в первой половине книги я рассмотрел декартово произведение представлений. Полиморфизм представлений - это отображение декартова произведения представлений, которое является морфизмом представлений по каждому отдельному аргументу. Приведенный морфизм представлений позволяет упростить изучение морфизмов представлений. Однако представление должно удостоверять определённым требованиям для того, чтобы существовал приведенный полиморфизм представлений. Возможно, что абелевая группа является единственной Ω-алгеброй, представление в которой допускает полиморфизм представлений. Однако сегодня это утверждение не доказано.

Мультипликативная Ω-группа - это Ω-алгебра, в которой определено произведение. Определение тензорного произведения представлений абелевой мультипликативной Ω-группы опирается на свойства приведенного полиморфизма представлений абелевой мультипликативной Ω-группы.

Так как алгебра - это модуль, в котором определено произведение, то мы можем применить эту теории к изучению линейных отображений алгебры. Например, множество линейных преобразований $D$-алгебры $A$ можно рассматривать как представление алгебры $A\otimes A$ в алгебре $A$.

# Оглавление

Глава 1

# Предисловие

## 1.1. Предисловие к изданию 1

Существует несколько эквивалентных определений модуля. Для меня наиболее интересно определение модуля как эффективное представление коммутативного кольца в абелевой группе. Это позволяет рассмотреть конструкции линейной алгебры с более общей точки зрения и понять какие конструкции верны в других алгебраических теориях. Например, мы можем рассматривать линейное отображение как морфизм представления, т. е. отображение, которое сохраняет структуру представления.

Модуль, в котором определено произведение, называется алгеброй. В зависимости от решаемой задачи, мы можем рассматривать различные алгебраические структуры на алгебре. Соответственно меняются отображения, сохраняющие структуру алгебры.

Если нас не интересует представление кольца $D$ в $D$-алгебре $A$, то алгебра $A$ является кольцом. Отображение, сохраняющее структуру алгебры как кольца, называется гомоморфизмом алгебры. Если нас не интересует произведение в $D$-алгебре $A$, то мы рассматриваем $D$-алгебру как модуль. Отображение, сохраняющее структуру алгебры как модуля, называется линейным отображением $D$-алгебры. Линейные отображения являются важным инструментом при изучении математического анализа. Отображение, сохраняющее структуру представления, называется линейным гомоморфизмом. В разделе 5.3 я показал существование нетривиального линейного гомоморфизма.

Февраль, 2015

## 1.2. Предисловие к изданию 2

Основной замысел этой книги - попытка перенести в представление $\Omega_1$-алгебры (которую мы будем называть $\Omega_1$-алгеброй преобразований $\Omega_2$-алгебры или просто $\Omega_1$-алгеброй преобразований) некоторые понятия, с которыми мы хорошо знакомы в линейной алгебре. Поскольку линейное отображение - это приведенный морфизм модуля, если мы рассматриваем модуль как представление кольца в абелевой группе, то определение полиморфизма представлений является естественным обобщением полилинейного отображения.

Изучая линейное отображение над некоммутативным кольцом с делением $D$, мы приходим к выводу что в равенстве

$$f(ax) = af(x)$$

$a$ не может быть произвольным элементом кольца $D$, но должно принадлежать центру кольца $D$. Это позволяет сделать теорию линейных отображений интересной и богатой.

Поэтому, также как в случае полилинейного отображения, мы требуем, чтобы преобразование, которое $\Omega_1$-алгебра порождает в одном сомножителе, можно было бы перенести на любой другой сомножитель. Вообще говоря, это требование не представляется жёстким, но выполнение этого требования позволяет рассматривать $\Omega_1$-алгебру преобразований как абелевую мультипликативную $\Omega_1$-группу.

Тензорное произведение векторных пространств ассоциативно. Рассмотрение ассоциативности тензорного произведения представлений порождает новое ограничение на полиморфизм представлений. Любые две операции $\Omega_2$-алгебры должны удовлетворять равенству (3.1.19). Возможно, что абелевая группа является единственной $\Omega_2$-алгеброй, представление в которой допускает полиморфизм представлений. Однако сегодня это утверждение не доказано.

Если мы верим в некоторое утверждение, которое мы не можем доказать или опровергнуть, то открытие противоположного утверждения может привести к интересным открытиям. На протяжении тысячелетий математики пытались доказать пятый постулат Евклида. Лобачевский и Больяй доказали, что существует геометрия, где этот постулат неверен. Математики пытались найти решение алгебраического уравнения, используя радикалы. Однако Галуа доказал, что это невозможно сделать, если степень уравнения выше 4.

Я долгое время полагал, что из существования приведенного морфизма представления следует существование приведенного полиморфизма представления. Теорема 3.1.13 оказалась для меня полной неожиданностью.

Я полагаю, что текст о полиморфизме и тензорном произведении представлений важен.

- Написав серию статей, посвящённых математическому анализу над банаховыми алгебрами, я полагал, что однажды я смогу написать подобную статью о математическом анализе в произвольном представлении. Меня смущала необходимость отказаться от сложения как оценки как мало расстояние между отображениями. Конечно, я могу утверждать, что морфизм представления находится в окрестности рассматриваемого отображения. Однако сложение является существенной компонентой построения. Отсутствие приведенного полиморфизма лишает меня возможности определить производные второго порядка.
- Аналогично тому, как мы рассматриваем линейное отображение $D$-алгебры, мы можем рассматривать линейное отображение эффективного представления в абелевой $\Omega_2$-группе. Структура линейного отображения эффективного представления в абелевой $\Omega_2$-группе сложнее чем структура линейного отображения $D$-алгебры. Однако это задача интересна для меня, и я надеюсь вернуться к ней в будущем.

Май, 2015

## 1.3. Соглашения

Соглашение 1.3.1. *Мы будем пользоваться соглашением Эйнштейна о сумме, в котором повторяющийся индекс (один вверху и один внизу) подразумевает сумму по повторяющемуся индексу. В этом случае предполагается* *известным множество индекса суммирования и знак суммы опускается*

$$c^i v_i = \sum_{i \in I} c^i v_i$$

□

Соглашение 1.3.2. *В выражении вида*

$$a_{i \cdot 0} x a_{i \cdot 1}$$

*предполагается сумма по индексу $i$.* □

Соглашение 1.3.3. *Пусть $A$ - свободная алгебра с конечным или счётным базисом. При разложении элемента алгебры $A$ относительно базиса $\overline{\overline{e}}$ мы пользуемся одной и той же корневой буквой для обозначения этого элемента и его координат. В выражении $a^2$ не ясно - это компонента разложения элемента $a$ относительно базиса или это операция возведения в степень. Для облегчения чтения текста мы будем индекс элемента алгебры выделять цветом. Например,*

$$a = a^{\boldsymbol{i}} e_{\boldsymbol{i}}$$

□

Соглашение 1.3.4. *Если свободная конечномерная алгебра имеет единицу, то мы будем отождествлять вектор базиса $\overline{e}_{\mathbf{0}}$ с единицей алгебры.* □

Глава 2

# Произведение представлений

## 2.1. Декартово произведение универсальных алгебр

Определение 2.1.1. *Пусть $\mathcal{A}$ - категория. Пусть $\{B_i, i \in I\}$ - множество объектов из $\mathcal{A}$. Объект*

$$P = \prod_{i \in I} B_i$$

*и множество морфизмов*

$$\{\ f_i : P \longrightarrow B_i\ , i \in I\}$$

*называется* **произведением объектов** $\{B_i, i \in I\}$ **в категории** $\mathcal{A}$[2.1], *если для любого объекта $R$ и множество морфизмов*

$$\{\ g_i : R \longrightarrow B_i\ , i \in I\}$$

*существует единственный морфизм*

$$h : R \longrightarrow P$$

*такой, что диаграмма*

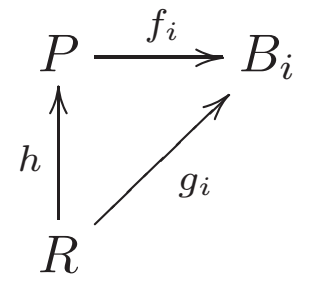

$$f_i \circ h = g_i$$

*коммутативна для всех $i \in I$.*

*Если $|I| = n$, то для произведения объектов $\{B_i, i \in I\}$ в $\mathcal{A}$ мы так же будем пользоваться записью*

$$P = \prod_{i=1}^{n} B_i = B_1 \times ... \times B_n$$

□

Пример 2.1.2. *Пусть $\mathcal{S}$ - категория множеств.*[2.2] *Согласно определению 2.1.1, декартово произведение*

$$A = \prod_{i \in I} A_i$$

*семейства множеств $(A_i, i \in I)$ и семейство проекций на $i$-й множитель*

$$p_i : A \to A_i$$

[2.1] Определение дано согласно [1], страница 45.

[2.2] Смотри также пример в [1], страница 45.

*являются произведением в категории* $\mathcal{S}$. □

Теорема 2.1.3. *Произведение существует в категории* $\mathcal{A}$ $\Omega$*-алгебр. Пусть* $\Omega$*-алгебра* $A$ *и семейство морфизмов*

$$p_i : A \to A_i \quad i \in I$$

*является произведением в категории* $\mathcal{A}$*. Тогда*

2.1.3.1: *Множество* $A$ *является декартовым произведением семейства множеств* $(A_i, i \in I)$

2.1.3.2: *Гомоморфизм* $\Omega$*-алгебры*

$$p_i : A \to A_i$$

*является проекцией на* $i$*-й множитель.*

2.1.3.3: *Любое* $A$*-число* $a$ *может быть однозначно представлено в виде кортежа* $(p_i(a), i \in I)$ $A_i$*-чисел.*

2.1.3.4: *Пусть* $\omega \in \Omega$ *-* $n$*-арная операция. Тогда операция* $\omega$ *определена покомпонентно*

$$a_1 ... a_n \omega = (a_{1i} ... a_{ni} \omega, i \in I) \tag{2.1.1}$$

*где* $a_1 = (a_{1i}, i \in I)$, ..., $a_n = (a_{ni}, i \in I)$ .

Доказательство. Пусть

$$A = \prod_{i \in I} A_i$$

декартово произведение семейства множеств $(A_i, i \in I)$ и, для каждого $i \in I$, отображение

$$p_i : A \to A_i$$

является проекцией на $i$-й множитель. Рассмотрим диаграмму морфизмов в категории множеств $\mathcal{S}$

(2.1.2)

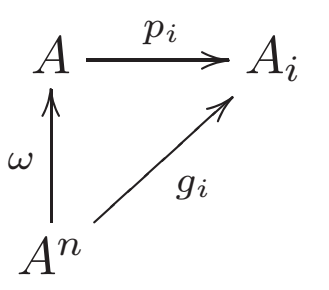

$$p_i \circ \omega = g_i$$

где отображение $g_i$ определено равенством

$$g_i(a_1, ..., a_n) = p_i(a_1) ... p_i(a_n) \omega$$

Согласно определению 2.1.1, отображение $\omega$ определено однозначно из множества диаграмм (2.1.2)

$$a_1 ... a_n \omega = (p_i(a_1) ... p_i(a_n) \omega, i \in I) \tag{2.1.3}$$

Равенство (2.1.1) является следствием равенства (2.1.3). □

Определение 2.1.4. *Если* $\Omega$*-алгебра* $A$ *и семейство морфизмов*

$$p_i : A \to A_i \quad i \in I$$

*является произведением в категории* $\mathcal{A}$*, то* $\Omega$*-алгебра* $A$ *называется* **прямым** *или* **декартовым произведением** $\Omega$**-алгебр** $(A_i, i \in I)$ . □

Теорема 2.1.5. *Пусть множество $A$ является декартовым произведением множеств $(A_i, i \in I)$ и множество $B$ является декартовым произведением множеств $(B_i, i \in I)$. Для каждого $i \in I$, пусть*

$$f_i : A_i \to B_i$$

*является отображением множества $A_i$ в множество $B_i$. Для каждого $i \in I$, рассмотрим коммутативную диаграмму*

(2.1.4)

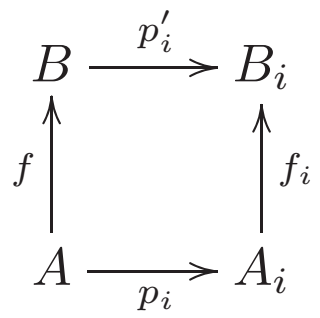

*где отображения $p_i$, $p'_i$ являются проекцией на $i$-й множитель. Множество коммутативных диаграмм (2.1.4) однозначно определяет отображение*

$$f : A_1 \to A_2$$

$$f(a_i, i \in I) = (f_i(a_i), i \in I)$$

Доказательство. Для каждого $i \in I$, рассмотрим коммутативную диаграмму

(2.1.5)

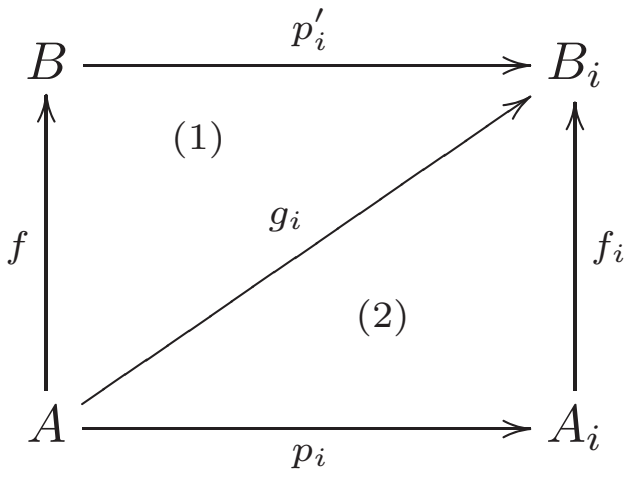

Пусть $a \in A$. Согласно утверждению 2.1.3.3, $A$-число $a$ может быть представлено в виде кортежа $A_i$-чисел

$$a = (a_i, i \in I) \quad a_i = p_i(a) \in A_i \tag{2.1.6}$$

Пусть

$$b = f(a) \in B \tag{2.1.7}$$

Согласно утверждению 2.1.3.3, $B$-число $b$ может быть представлено в виде кортежа $B_i$-чисел

$$b = (b_i, i \in I) \quad b_i = p'_i(b) \in B_i \tag{2.1.8}$$

Из коммутативности диаграммы (1) и из равенств (2.1.7), (2.1.8) следует, что

$$b_i = g_i(b) \tag{2.1.9}$$

Из коммутативности диаграммы (2) и из равенства (2.1.6) следует, что

$$b_i = f_i(a_i)$$

□

Теорема 2.1.6. *Пусть* $\Omega$*-алгебра* $A$ *является декартовым произведением* $\Omega$*-алгебр* $(A_i, i \in I)$ *и* $\Omega$*-алгебра* $B$ *является декартовым произведением* $\Omega$*-алгебр* $(B_i, i \in I)$ *. Для каждого* $i \in I$*, пусть отображение*

$$f_i : A_i \to B_i$$

*является гомоморфизмом* $\Omega$*-алгебры. Тогда отображение*

$$f : A_1 \to A_2$$

*определённое равенством*

$$f(a_i, i \in I) = (f_i(a_i), i \in I) \tag{2.1.10}$$

*является гомоморфизмом* $\Omega$*-алгебры.*

Доказательство. Пусть $\omega \in \Omega$ - n-арная операция. Пусть $a_1 = (a_{1i}, i \in I)$, ..., $a_n = (a_{ni}, i \in I)$ и $b_1 = (b_{1i}, i \in I)$, ..., $b_n = (b_{ni}, i \in I)$ . Из равенств (2.1.1), (2.1.10) следует, что

$$\begin{aligned} f(a_1...a_n\omega) &= f(a_{1i}...a_{ni}\omega, i \in I) \\ &= (f_i(a_{1i}...a_{ni}\omega), i \in I) \\ &= ((f_i(a_{1i}))...(f_i(a_{ni})), i \in I) \\ &= (b_{1i}...b_{ni}\omega, i \in I) \end{aligned}$$

$$f(a_1)...f(a_n)\omega = b_1...b_n\omega = (b_{1i}...b_{ni}\omega, i \in I)$$

□

## 2.2. Декартово произведение представлений

Лемма 2.2.1. *Пусть*

$$A = \prod_{i \in I} A_i$$

*декартово произведение семейства* $\Omega_1$*-алгебр* $(A_i, i \in I)$ *. Для каждого* $i \in I$*, пусть множество* ${}^*A_i$ *является* $\Omega_2$*-алгеброй. Тогда множество*

$$^\circ A = \{f \in {}^*A : f(a_i, i \in I) = (f_i(a_i), i \in I)\} \tag{2.2.1}$$

*является декартовым произведением* $\Omega_2$*-алгебр* ${}^*A_i$*.*

Доказательство. Согласно определению (2.2.1), мы можем представить отображение $f \in {}^\circ A$ в виде кортежа

$$f = (f_i, i \in I)$$

отображений $f_i \in {}^*A_i$. Согласно определению (2.2.1),

$$(f_i, i \in I)(a_i, i \in I) = (f_i(a_i), i \in I)$$

Пусть $\omega \in \Omega_2$ - n-арная операция. Мы определим операцию $\omega$ на множестве ${}^\circ A$ равенством

$$((f_{1i}, i \in I)...(f_{ni}, i \in I)\omega)(a_i, i \in I) = ((f_{1i}(a_i))...(f_{ni}(a_i))\omega, i \in I)$$

□

Определение 2.2.2. *Пусть* $\mathcal{A}_1$ *- категория* $\Omega_1$*-алгебр. Пусть* $\mathcal{A}_2$ *- категория* $\Omega_2$*-алгебр. Мы определим* **категорию** $(\mathcal{A}_1*)\mathcal{A}_2$ **левосторонних представлений**. *Объектами этой категории являются левосторонние представления* $\Omega_1$*-алгебры в* $\Omega_2$*-алгебре. Морфизмами этой категории являются морфизмы соответствующих представлений.* □

Теорема 2.2.3. *В категории* $(\mathcal{A}_1*)\mathcal{A}_2$ *существует произведение однотранзитивных левосторонних представлений* $\Omega_1$*-алгебры в* $\Omega_2$*-алгебре.*

Доказательство. Для $j = 1, 2,$ пусть

$$P_j = \prod_{i\in I} B_{ji}$$

произведение семейства $\Omega_j$-алгебр $\{B_{ji}, i \in I\}$ и для любого $i \in I$ отображение

$$t_{ji} : P_j \longrightarrow B_{ji}$$

является проекцией на множитель $i$. Для каждого $i \in I$, пусть

$$h_i : B_{1i} \xrightarrow{*} B_{2i}$$

однотранзитивное $B_{1i}*$-представление в $\Omega_2$-алгебре $B_{2i}$.

Пусть $b_1 \in P_1$. Согласно утверждению 2.1.3.3, $P_1$-число $b_1$ может быть представлено в виде кортежа $B_{1i}$-чисел

(2.2.2) $$b_1 = (b_{1i}, i \in I) \quad b_{1i} = t_{1i}(b_1) \in B_{1i}$$

Пусть $b_2 \in P_2$. Согласно утверждению 2.1.3.3, $P_2$-число $b_2$ может быть представлено в виде кортежа $B_{2i}$-чисел

(2.2.3) $$b_2 = (b_{2i}, i \in I) \quad b_{2i} = t_{2i}(b_2) \in B_{2i}$$

Лемма 2.2.4. *Для каждого* $i \in I$*, рассмотрим диаграмму отображений*

(2.2.4)

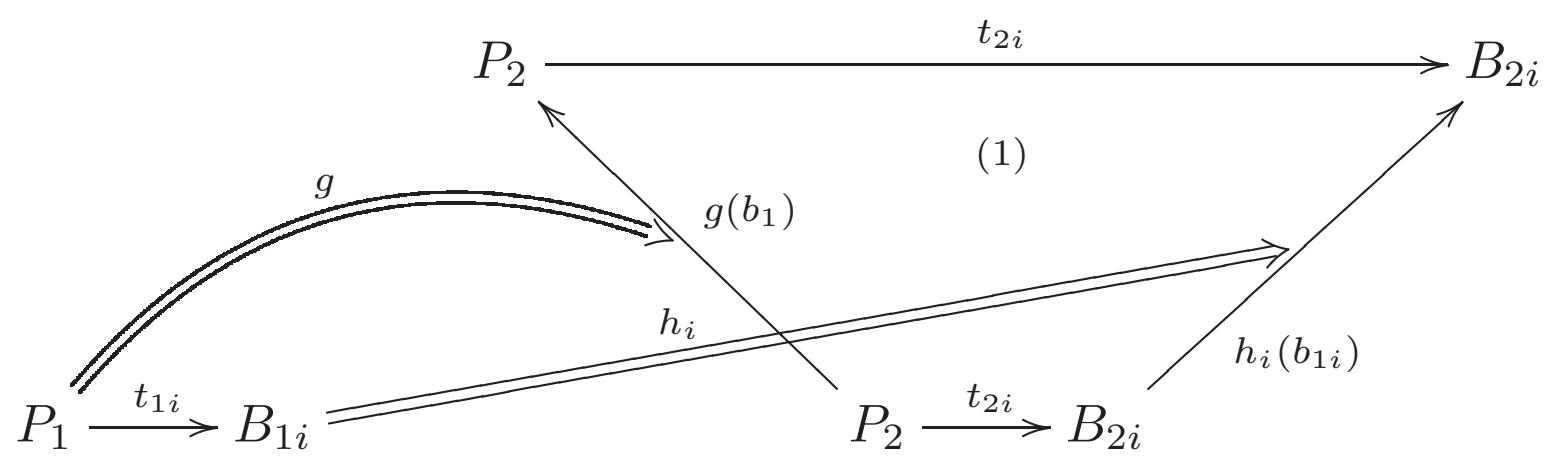

*Пусть отображение*

$$g : P_1 \to {}^*P_2$$

*определено равенством*

(2.2.5) $$g(b_1) \circ b_2 = (h_i(b_{1i}) \circ b_{2i}, i \in I)$$

*Тогда отображение* $g$ *является однотранзитивным* $P_1*$*-представлением в* $\Omega_2$*-алгебре* $P_2$

$$g : P_1 \xrightarrow{*} P_2$$

*Отображение* $(t_{1i}, t_{2i})$ *является морфизмом представления* $g$ *в представление* $h_i$.

Доказательство.

2.2.4.1: Согласно определениям [4]-2.1.1, [4]-2.1.2, отображение $h_i(b_{1i})$ является гомоморфизмом $\Omega_2$-алгебры $B_{2i}$. Согласно теореме 2.1.6, из коммутативности диаграммы (1) для каждого $i \in I$, следует, что отображение

$$g(b_1) : P_2 \to P_2$$

определённое равенством (2.2.5) является гомоморфизмом $\Omega_2$-алгебры $P_2$.

2.2.4.2: Согласно определению [4]-2.1.2, множество ${}^*B_{2i}$ является $\Omega_1$-алгеброй. Согласно лемме 2.2.1, множество ${}^\circ P_2 \subseteq {}^*P_2$ является $\Omega_1$-алгеброй.

2.2.4.3: Согласно определению [4]-2.1.2, отображение

$$h_i : B_{1i} \to {}^*B_{2i}$$

является гомоморфизмом $\Omega_1$-алгебры. Согласно теореме 2.1.6, отображение

$$g : P_1 \to {}^*P_2$$

определённое равенством

$$g(b_1) = (h_i(b_{1i}), i \in I)$$

является гомоморфизмом $\Omega_1$-алгебры.

Согласно утверждениям 2.2.4.1, 2.2.4.3 и определению [4]-2.1.2, отображение $g$ является $P_1*$-представлением в $\Omega_2$-алгебре $P_2$.

Пусть $b_{21}, b_{22} \in P_2$. Согласно утверждению 2.1.3.3, $P_2$-числа $b_{21}, b_{22}$ могут быть представлены в виде кортежей $B_{2i}$-чисел

$$\begin{array}{ll} b_{21} = (b_{21i}, i \in I) & b_{21i} = t_{2i}(b_{21}) \in B_{2i} \\ b_{22} = (b_{22i}, i \in I) & b_{22i} = t_{2i}(b_{22}) \in B_{2i} \end{array} \tag{2.2.6}$$

Согласно теореме [4]-2.1.9, поскольку представление $h_i$ однотранзитивно, то существует единственное $B_{1i}$-число $b_{1i}$ такое, что

$$b_{22i} = h_i(b_{1i}) \circ b_{21i}$$

Согласно определениям (2.2.2), (2.2.5), (2.2.6), существует единственное $P_1$-число $b_1$ такое, что

$$b_{22} = g(b_1) \circ b_{21}$$

Согласно теореме [4]-2.1.9, представление $g$ однотранзитивно.

Из коммутативности диаграммы (1) и определения [4]-2.2.2, следует, что отображение $(t_{1i}, t_{2i})$ является морфизмом представления $g$ в представление $h_i$. ⊙

Пусть

$$d_2 = g(b_1) \circ b_2 \quad d_2 = (d_{2i}, i \in I) \tag{2.2.7}$$

Из равенств (2.2.5), (2.2.7) следует, что

$$d_{2i} = h_i(b_{1i}) \circ b_{2i} \tag{2.2.8}$$

Для $j = 1, 2$, пусть $R_j$ - другой объект категории $\mathcal{A}_j$. Для любого $i \in I$, пусть отображение

$$r_{1i} : R_1 \longrightarrow B_{1i}$$

является морфизмом из $\Omega_1$-алгебра $R_1$ в $\Omega_1$-алгебру $B_{1i}$. Согласно определению 2.1.1, существует единственный морфизм $\Omega_1$-алгебры

$$s_1 : R_1 \longrightarrow P_1$$

такой, что коммутативна диаграмма

(2.2.9)

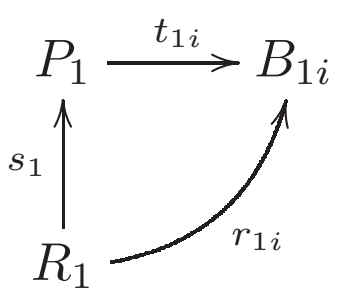

$$t_{1i} \circ s_1 = r_{1i}$$

Пусть $a_1 \in R_1$. Пусть

(2.2.10) $$b_1 = s_1(a_1) \in P_1$$

Из коммутативности диаграммы (2.2.9) и утверждений (2.2.10), (2.2.2) следует, что

(2.2.11) $$b_{1i} = r_{1i}(a_1)$$

Пусть

$$f : R_1 \overset{*}{\longrightarrow} R_2$$

однотранзитивное $R_1*$-представление в $\Omega_2$-алгебре $R_2$. Согласно теореме [4]-2.2.10, морфизм $\Omega_2$-алгебры

$$r_{2i} : R_2 \longrightarrow B_{2i}$$

такой, что отображение $(r_{1i}, r_{2i})$ является морфизмом представлений из $f$ в $h_i$, определён однозначно с точностью до выбора образа $R_2$-числа $a_2$. Согласно замечанию [4]-2.2.6, в диаграмме отображений

(2.2.12)

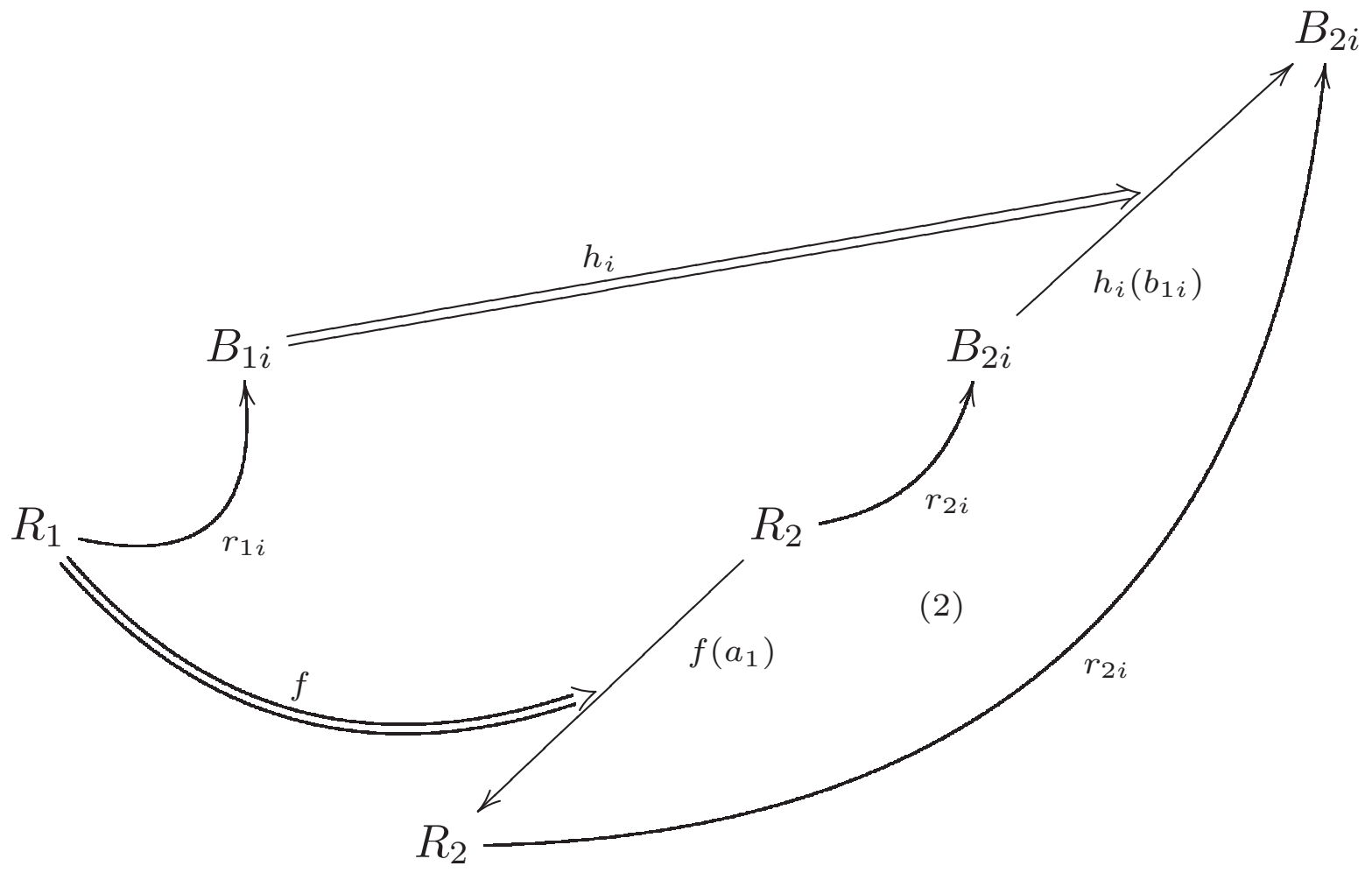

диаграмма (2) коммутативна. Согласно определению 2.1.1, существует единственный морфизм $\Omega_2$-алгебры

$$s_2 : R_2 \longrightarrow P_2$$

такой, что коммутативна диаграмма

(2.2.13)

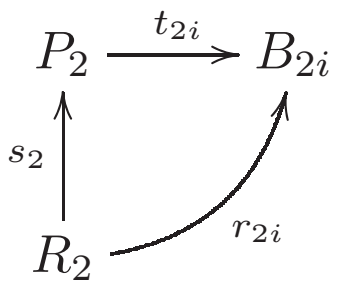

$$t_{2i} \circ s_2 = r_{2i}$$

Пусть $a_2 \in R_2$. Пусть

$$b_2 = s_2(a_2) \in P_2 \tag{2.2.14}$$

Из коммутативности диаграммы (2.2.13) и утверждений (2.2.14), (2.2.3) следует, что

$$b_{2i} = r_{2i}(a_2) \tag{2.2.15}$$

Пусть

$$c_2 = f(a_1) \circ a_2 \tag{2.2.16}$$

Из коммутативности диаграммы (2) и равенств (2.2.8), (2.2.15), (2.2.16) следует, что

$$d_{2i} = r_{2i}(c_2) \tag{2.2.17}$$

Из равенств (2.2.8), (2.2.17) следует, что

$$d_2 = s_2(c_2) \tag{2.2.18}$$

что согласуется с коммутативносью диаграммы (2.2.13).

Для каждого $i \in I$, мы объединим диаграммы отображений (2.2.4), (2.2.9), (2.2.13), (2.2.12)

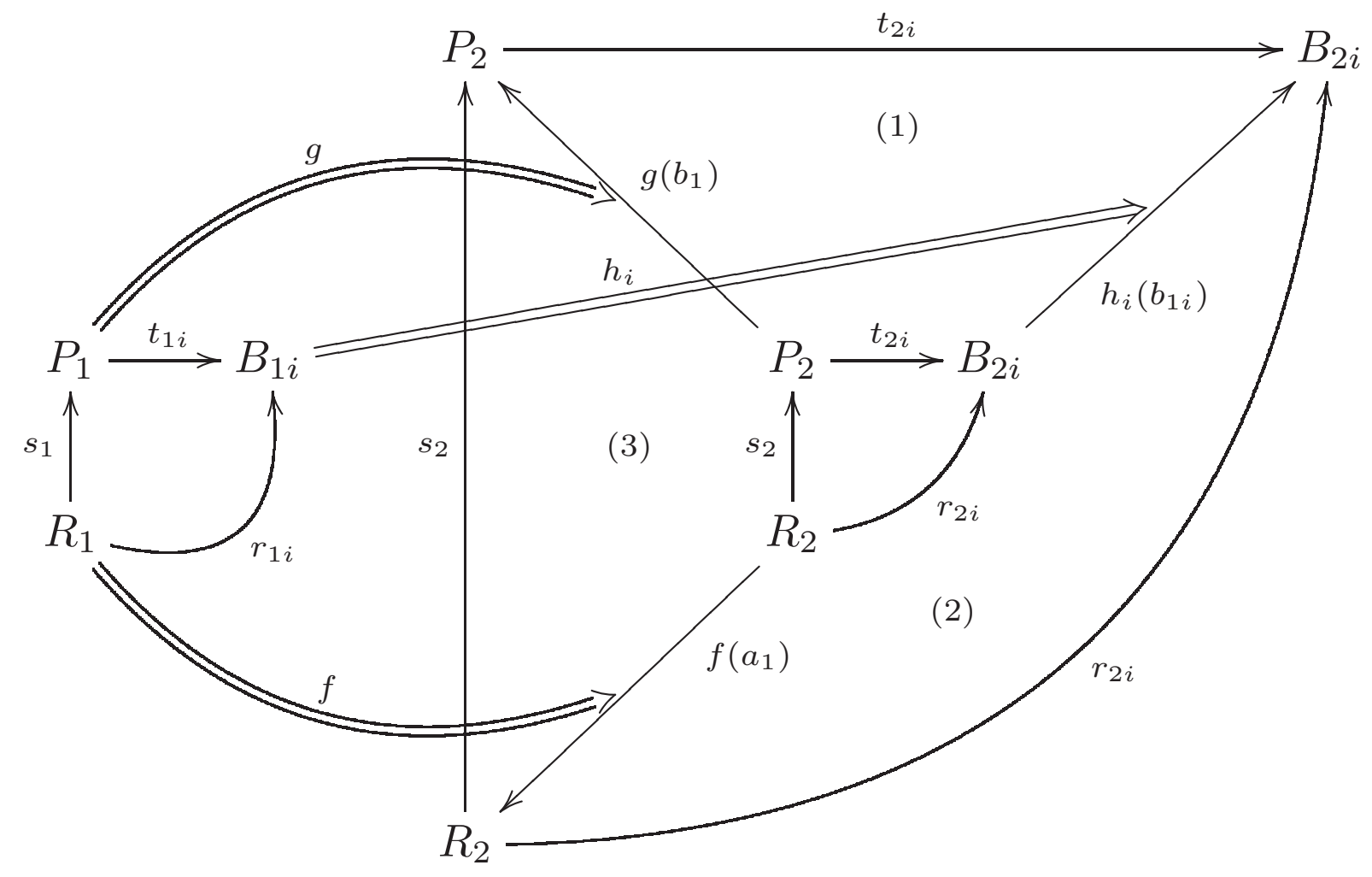

Из равенств (2.2.7) (2.2.14) и из равенств (2.2.16), (2.2.18), следует коммутативность диаграммы (3). Следовательно, отображение $(s_1, s_2)$ является морфизмом представлений из $f$ в $g$., Согласно теореме [4]-2.2.10, морфизм $(s_1, s_2)$ определён однозначно, так как мы требуем (2.2.18).

Согласно определению 2.1.1, представление $g$ и семейство морфизмов представления $((t_{1i}, t_{2i}), i \in I)$ является произведением в категории $(\mathcal{A}_1*)\mathcal{A}_2$ . □

## 2.3. Приведенное декартово произведение представлений

Определение 2.3.1. *Пусть $A_1$ - $\Omega_1$-алгебра. Пусть $\mathcal{A}_2$ - категория $\Omega_2$-алгебр. Мы определим* **категорию** $(A_1*)\mathcal{A}_2$ **левосторонних представлений**. *Объектами этой категории являются левосторонние представления $\Omega_1$-алгебры $A_1$ в $\Omega_2$-алгебре. Морфизмами этой категории являются приведенные морфизмы соответствующих представлений.* □

Теорема 2.3.2. *В категории $(A_1*)\mathcal{A}_2$ существует произведение эффективных левосторонних представлений $\Omega_1$-алгебры $A_1$ в $\Omega_2$-алгебре и это произведение является эффективным левосторонним представлением $\Omega_1$-алгебры $A_1$.*

Доказательство. Пусть

$$A_2 = \prod_{i \in I} A_{2i}$$

произведение семейства $\Omega_2$-алгебр $\{A_{2i}, i \in I\}$ и для любого $i \in I$ отображение

$$t_i : A_2 \longrightarrow A_{2i}$$

является проекцией на множитель $i$. Для каждого $i \in I$, пусть

$$h_i : A_1 \overset{*}{\longrightarrow} A_{2i}$$

эффективное $A_1*$-представление в $\Omega_2$-алгебре $A_{2i}$.

Пусть $b_1 \in A_1$. Пусть $b_2 \in A_2$. Согласно утверждению 2.1.3.3, $A_2$-число $b_2$ может быть представлено в виде кортежа $A_{2i}$-чисел

$$b_2 = (b_{2i}, i \in I) \quad b_{2i} = t_i(b_2) \in A_{2i} \tag{2.3.1}$$

Лемма 2.3.3. *Для каждого $i \in I$, рассмотрим диаграмму отображений*

(2.3.2)

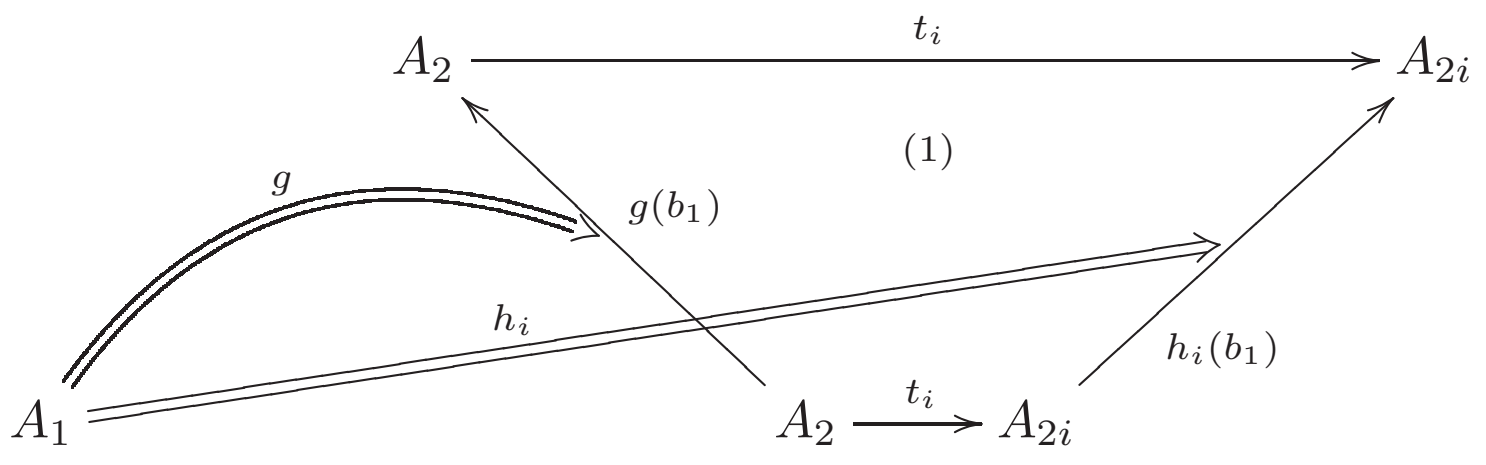

*Пусть отображение*

$$g : A_1 \to {}^* A_2$$

*определено равенством*

$$g(b_1) \circ b_2 = (h_i(b_1) \circ b_{2i}, i \in I) \tag{2.3.3}$$

*Тогда отображение $g$ является эффективным $A_1*$-представлением в $\Omega_2$-алгебре $A_2$*

$$g : A_1 \overset{*}{\longrightarrow} A_2$$

*Отображение $t_i$ является приведенным морфизмом представления $g$ в представление $h_i$.*

Доказательство.

2.3.3.1: Согласно определениям [4]-2.1.1, [4]-2.1.2, отображение $h_i(b_1)$ является гомоморфизмом $\Omega_2$-алгебры $A_{2i}$. Согласно теореме 2.1.6, из коммутативности диаграммы (1) для каждого $i \in I$, следует, что отображение

$$g(b_1) : A_2 \to A_2$$

определённое равенством (2.3.3) является гомоморфизмом $\Omega_2$-алгебры $A_2$.

2.3.3.2: Согласно определению [4]-2.1.2, множество $^*A_{2i}$ является $\Omega_1$-алгеброй. Согласно лемме 2.2.1, множество $^\circ A_2 \subseteq {}^*A_2$ является $\Omega_1$-алгеброй.

2.3.3.3: Согласно определению [4]-2.1.2, отображение

$$h_i : A_1 \to {}^*A_{2i}$$

является гомоморфизмом $\Omega_1$-алгебры. Согласно теореме 2.1.6, отображение

$$g : A_1 \to {}^*A_2$$

определённое равенством

$$g(b_1) = (h_i(b_1), i \in I)$$

является гомоморфизмом $\Omega_1$-алгебры.

Согласно утверждениям 2.3.3.1, 2.3.3.3 и определению [4]-2.1.2, отображение $g$ является $A_1*$-представлением в $\Omega_2$-алгебре $A_2$.

Для любого $i \in I$, согласно определению [4]-2.1.6, $A_1$-число $a_1$ порождает единственное преобразование

$$b_{22i} = h_i(b_1) \circ b_{21i} \tag{2.3.4}$$

Пусть $b_{21}, b_{22} \in A_2$. Согласно утверждению 2.1.3.3, $A_2$-числа $b_{21}, b_{22}$ могут быть представлены в виде кортежей $A_{2i}$-чисел

$$\begin{aligned} b_{21} &= (b_{21i}, i \in I) \quad b_{21i} = t_i(b_{21}) \in A_{2i} \\ b_{22} &= (b_{22i}, i \in I) \quad b_{22i} = t_i(b_{22}) \in A_{2i} \end{aligned} \tag{2.3.5}$$

Согласно определению (2.3.3) представления $g$, из равенств (2.3.4), (2.3.5) следует, что $A_1$-число $a_1$ порождает единственное преобразование

$$b_{22} = (h_i(b_1) \circ b_{21i}, i \in I) = g(b_1) \circ b_{21} \tag{2.3.6}$$

Согласно определению [4]-2.1.6, представление $g$ эффективно.

Из коммутативности диаграммы (1) и определения [4]-2.2.2, следует, что отображение $t_i$ является приведенным морфизмом представления $g$ в представление $h_i$. ⊙

Пусть

$$d_2 = g(b_1) \circ b_2 \quad d_2 = (d_{2i}, i \in I) \tag{2.3.7}$$

Из равенств (2.3.3), (2.3.7) следует, что

$$d_{2i} = h_i(b_1) \circ b_{2i} \tag{2.3.8}$$

Пусть $R_2$ - другой объект категории $\mathcal{A}_2$. Пусть

$$f : A_1 \xrightarrow{*} R_2$$

эффективное $A_1*$-представление в $\Omega_2$-алгебре $R_2$. Для любого $i \in I$, пусть существует морфизм

$$r_i : R_2 \longrightarrow A_{2i}$$

представлений из $f$ в $h_i$. Согласно замечанию [4]-2.2.6, в диаграмме отображений

(2.3.9)

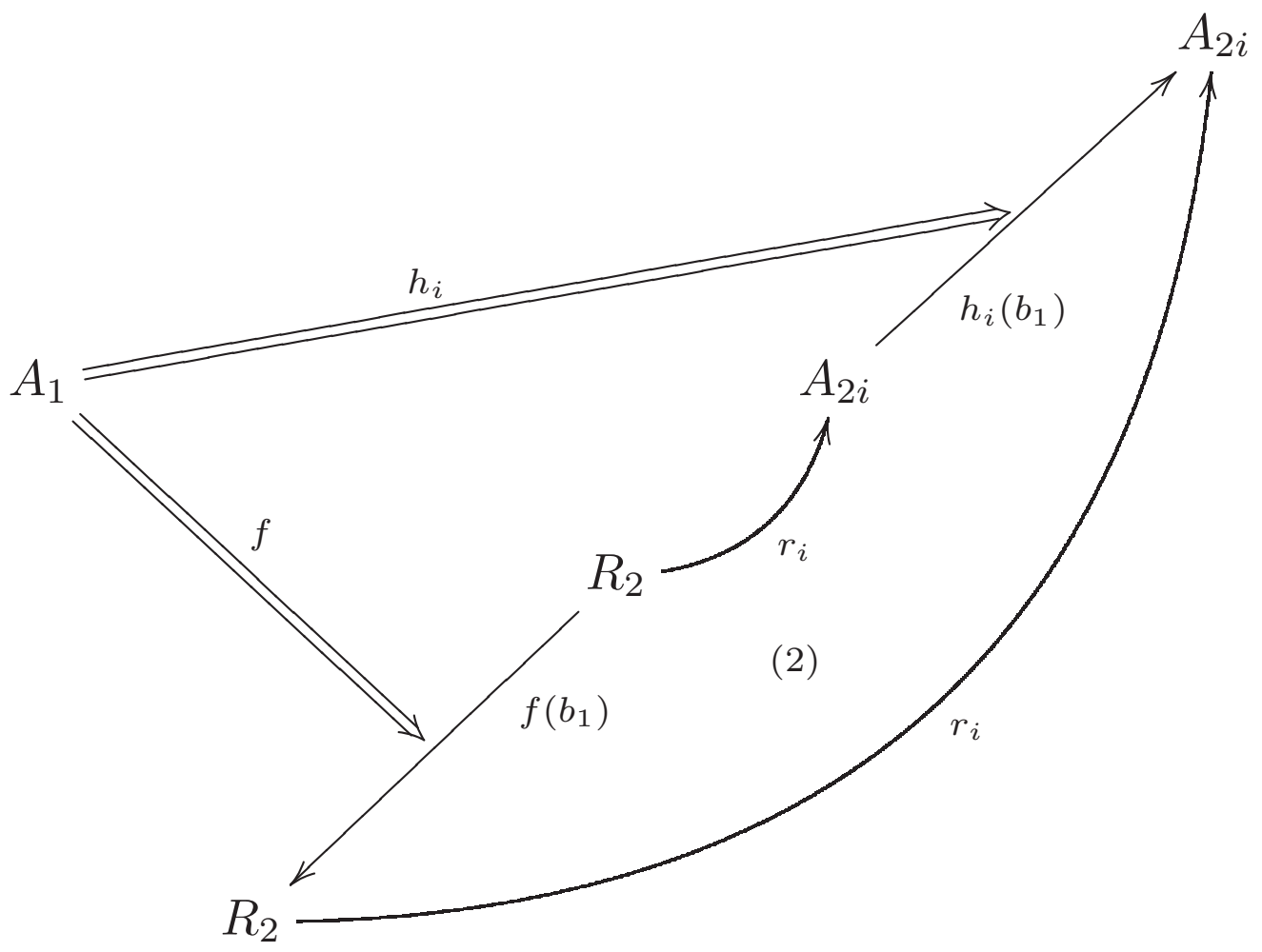

диаграмма (2) коммутативна. Согласно определению 2.1.1, существует единственный морфизм $\Omega_2$-алгебры

$$s : R_2 \longrightarrow A_2$$

такой, что коммутативна диаграмма

(2.3.10)

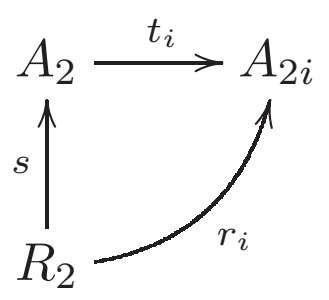

$$t_i \circ s = r_i$$

Пусть $a_2 \in R_2$. Пусть

$$b_2 = s(a_2) \in A_2 \tag{2.3.11}$$

Из коммутативности диаграммы (2.3.10) и утверждений (2.3.11), (2.3.1) следует, что

$$b_{2i} = r_i(a_2) \tag{2.3.12}$$

Пусть

$$c_2 = f(a_1) \circ a_2 \tag{2.3.13}$$

Из коммутативности диаграммы (2) и равенств (2.3.8), (2.3.12), (2.3.13) следует, что

$$d_{2i} = r_i(c_2) \tag{2.3.14}$$

Из равенств (2.3.8), (2.3.14) следует, что

$$d_2 = s(c_2) \tag{2.3.15}$$

что согласуется с коммутативносью диаграммы (2.3.10).

Для каждого $i \in I$, мы объединим диаграммы отображений (2.3.2), (2.3.10), (2.3.9)

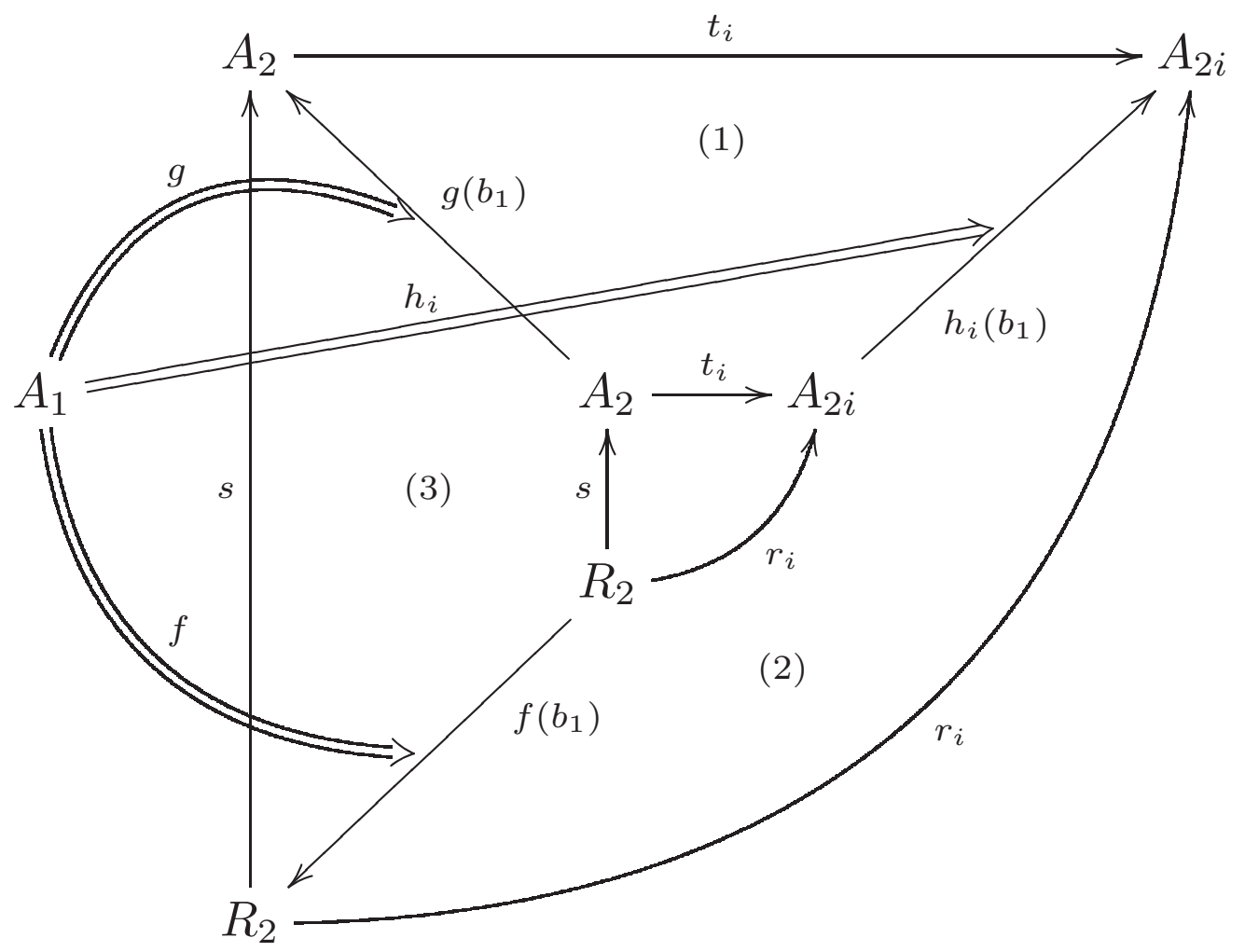

Из равенств (2.3.7), (2.3.11) и из равенств (2.3.13), (2.3.15), следует коммутативность диаграммы (3). Следовательно, отображение $s$ является приведенным морфизмом представлений из $f$ в $g$. Согласно замечанию [4]-2.3.2, отображение $s$ является гомоморфизмом $\Omega_2$ алгебры. Согласно теореме 2.1.3 и определению 2.1.1, приведенный морфизм $s$ определён однозначно.

Согласно определению 2.1.1, представление $g$ и семейство морфизмов представления $(t_i, i \in I)$ является произведением в категории $(A_1*)\mathcal{A}_2$ . $\square$

Глава 3

# Тензорное произведение представлений

## 3.1. Полиморфизм представлений

Определение 3.1.1. *Пусть* $A_1$, ..., $A_n$ - $\Omega_1$*-алгебры. Пусть* $B_1$, ..., $B_n$, $B$ - $\Omega_2$*-алгебры. Пусть для любого* $k$, $k=1$, ..., $n$,

$$f_k : A_k \longrightarrow\!\!\!\!\!\!\!*\;\;\; B_k$$

*представление* $\Omega_1$*-алгебры* $A_k$ *в* $\Omega_2$*-алгебре* $B_k$. *Пусть*

$$f : A \longrightarrow\!\!\!\!\!\!\!*\;\;\; B$$

*представление* $\Omega_1$*-алгебры* $A$ *в* $\Omega_2$*-алгебре* $B$. *Отображение*

$$r : A_1 \times ... \times A_n \to A \quad R : B_1 \times ... \times B_n \to B$$

*называется* **полиморфизмом представлений** $f_1$, ..., $f_n$ *в представление* $f$, *если для любого* $k$, $k=1$, ..., $n$,*при условии, что все переменные кроме переменных* $a_k \in A_k$, $b_k \in B_k$ *имеют заданное значение, отображение* $(r, R)$ *является морфизмом представления* $f_k$ *в представление* $f$.

*Если* $f_1 = ... = f_n$, *то мы будем говорить, что отображение* $(r, R)$ *является полиморфизмом представления* $f_1$ *в представление* $f$.

*Если* $f_1 = ... = f_n = f$, *то мы будем говорить, что отображение* $(r, R)$ *является полиморфизмом представления* $f$. □

Теорема 3.1.2. *Пусть отображение* $(r, R)$ *является полиморфизмом представлений* $f_1$, ..., $f_n$ *в представление* $f$. *Отображение* $(r, R)$ *удовлетворяет равенству*

$$R(f_1(a_1)(m_1), ..., f_n(a_n)(m_n)) = f(r(a_1, ..., a_n))(R(m_1, ..., m_n)) \tag{3.1.1}$$

*Пусть* $\omega_1 \in \Omega_1(p)$. *Для любого* $k$, $k=1$, ..., $n$,*отображение* $r$ *удовлетворяет равенству*

$$\begin{aligned} &r(a_1, ..., a_{k\cdot 1}...a_{k\cdot p}\omega_1, ..., a_n) \\ =\, &r(a_1, ..., a_{k\cdot 1}, ..., a_n)...r(a_1, ..., a_{k\cdot p}, ..., a_n)\omega_1 \end{aligned} \tag{3.1.2}$$

*Пусть* $\omega_2 \in \Omega_2(p)$. *Для любого* $k$, $k=1$, ..., $n$,*отображение* $R$ *удовлетворяет равенству*

$$\begin{aligned} &R(m_1, ..., m_{k\cdot 1}...m_{k\cdot p}\omega_2, ..., m_n) \\ =\, &R(m_1, ..., m_{k\cdot 1}, ..., m_n)...R(m_1, ..., m_{k\cdot p}, ..., m_n)\omega_2 \end{aligned} \tag{3.1.3}$$

Доказательство. Равенство (3.1.1) следует из определения 3.1.1 и равенства [4]-(2.2.4). Равенство (3.1.2) следует из утверждения, что для любого $k$, $k=1, ..., n$,при условии, что все переменные кроме переменной $x_k \in A_k$ имеют заданное значение, отображение $r$ является гомоморфизмом $\Omega_1$-алгебры $A_k$ в $\Omega_1$-алгебру $A$. Равенство (3.1.3) следует из утверждения, что для любого $k$, $k=1, ..., n$,при условии, что все переменные кроме переменной $m_k \in B_k$ имеют заданное значение, отображение $R$ является гомоморфизмом $\Omega_2$-алгебры $B_k$ в $\Omega_2$-алгебру $B$. □

Определение 3.1.3. *Пусть $A, B_1, ..., B_n, B$ - универсальные алгебры. Пусть для любого $k$, $k=1, ..., n$,*

$$f_k : A \xrightarrow{*} B_k$$

*эффективное представление $\Omega_1$-алгебры $A_k$ в $\Omega_2$-алгебре $B_k$. Пусть*

$$f : A \xrightarrow{*} B$$

*эффективное представление $\Omega_1$-алгебры $A$ в $\Omega_2$-алгебре $B$. Отображение*

$$r_2 : B_1 \times ... \times B_n \to B$$

*называется* **приведенным полиморфизмом представлений** *$f_1, ..., f_n$ в представление $f$, если для любого $k$, $k=1, ..., n$,при условии, что все переменные кроме переменной $x_k \in B_k$ имеют заданное значение, отображение $R$ является приведенным морфизмом представления $f_k$ в представление $f$.*

*Если $f_1 = ... = f_n$, то мы будем говорить, что отображение $R$ является приведенным полиморфизмом представления $f_1$ в представление $f$.*

*Если $f_1 = ... = f_n = f$, то мы будем говорить, что отображение $R$ является приведенным полиморфизмом представления $f$.* □

Теорема 3.1.4. *Пусть отображение $R$ является приведенным полиморфизмом эффективных представлений $f_1, ..., f_n$ в эффективное представление $f$. Для любого $k$, $k=1, ..., n$,отображение $R$ удовлетворяет равенству*

$$R(m_1, ..., f_k(a) \circ m_k, ..., m_n) = f(a) \circ R(m_1, ..., m_n) \tag{3.1.4}$$

*Пусть $\omega_2 \in \Omega_2(p)$. Для любого $k$, $k=1, ..., n$,отображение $R$ удовлетворяет равенству*

$$\begin{aligned} &R(m_1, ..., m_{k\cdot 1}...m_{k\cdot p}\omega_2, ..., m_n) \\ =&R(m_1, ..., m_{k\cdot 1}, ..., m_n)...R(m_1, ..., m_{k\cdot p}, ..., m_n)\omega_2 \end{aligned} \tag{3.1.5}$$

Доказательство. Равенство (3.1.4) следует из определения 3.1.3 и равенства [4]-(2.3.3). Равенство (3.1.5) следует из утверждения, что для любого $k$, $k=1, ..., n$,при условии, что все переменные кроме переменной $m_k \in B_k$ имеют заданное значение, отображение $R$ является гомоморфизмом $\Omega_2$-алгебры $B_k$ в $\Omega_2$-алгебру $B$. □

Мы также будем говорить, что отображение $(r, R)$ является полиморфизмом представлений в $\Omega_2$-алгебрах $B_1, ..., B_n$ в представление в $\Omega_2$-алгебре $B$. Аналогично, мы будем говорить, что отображение $R$ является приведенным полиморфизмом представлений в $\Omega_2$-алгебрах $B_1, ..., B_n$ в представление в $\Omega_2$-алгебре $B$.

Сравнение определений 3.1.1 и 3.1.3 показывает, что существует разница между этими двумя формами полиморфизма. Особенно хорошо это различие видно при сравнении равенств (3.1.1) и (3.1.4). Если мы хотим иметь возможность выразить приведенный полиморфизм представлений через полиморфизм представлений, то мы должны потребовать, два условия:

(1) Представление $f$ универсальной алгебры содержит тождественное преобразование $\delta$. Следовательно, существует $e \in A$ такой, что $f(e) = \delta$. Не нарушая общности, мы положим, что выбор $e \in A$ не зависит от того, какое представление $f_1$, ..., $f_n$ мы рассматриваем.
(2) Для любого $k$, $k = 1$, ..., $n$,

$$r(a_1, ..., a_n) = a_k \quad a_i = e \quad i \neq k \tag{3.1.6}$$

Тогда, при условии $a_i = e$, $i \neq k$, равенство (3.1.1) имеет вид

$$R(m_1, ..., f_k(a_k) \circ m_k, ..., m_n) = f(r(e, ..., a_k, ..., e)) \circ R(m_1, ..., m_n) \tag{3.1.7}$$

Очевидно, что равенство (3.1.7) совпадает с равенством (3.1.4).

Похожая задача появляется при анализе приведенного полиморфизма представлений. Пользуясь равенством (3.1.4), мы можем записать выражение

$$R(m_1, ..., f_k(a_k) \circ m_k, ..., f_l(a_l) \circ m_l, ..., m_n) \tag{3.1.8}$$

либо в виде

$$\begin{aligned}&R(m_1, ..., f_k(a_k) \circ m_k, ..., f_l(a_l) \circ m_l, ..., m_n)\\ =&f(a_k) \circ R(m_1, ..., m_k, ..., f_l(a_l) \circ m_l, ..., m_n)\\ =&f(a_k) \circ (f(a_l) \circ R(m_1, ..., m_k, ..., m_l, ..., m_n))\\ =&(f(a_k) \circ f(a_l)) \circ R(m_1, ..., m_k, ..., m_l, ..., m_n)\end{aligned} \tag{3.1.9}$$

либо в виде

$$\begin{aligned}&R(m_1, ..., f_k(a_k) \circ m_k, ..., f_l(a_l) \circ m_l, ..., m_n)\\ =&f(a_l) \circ R(m_1, ..., f_k(a_k) \circ m_k, ..., m_l, ..., m_n)\\ =&f(a_l) \circ (f(a_k) \circ R(m_1, ..., m_k, ..., m_l, ..., m_n))\\ =&(f(a_l) \circ f(a_k)) \circ R(m_1, ..., m_k, ..., m_l, ..., m_n)\end{aligned} \tag{3.1.10}$$

Отображения $f(a_k)$, $f(a_l)$ являются гомоморфизмами $\Omega_2$-алгебры $B$. Следовательно, отображение $f(a_k) \circ f(a_l)$ является гомоморфизмом $\Omega_2$-алгебры $B$. Однако, не всякая $\Omega_1$-алгебра $A$ имеет такой $a$ (зависящий от $a_k$ и $a_l$), что

$$f(a) = f(a_k) \circ f(a_l)$$

Если представление $f$ однотранзитивно и для любых $A$-чисел $a$, $b$ существует $A$-число $c$ такое, что

$$f(c) = f(a) \circ f(b) \tag{3.1.11}$$

то равенство (3.1.11) определяет $A$-число $c$ единственным образом. Следовательно, мы можем определить умножение

$$c_1 = a_1 * b_1$$

таким образом, что

$$f(a * b) = f(a) \circ f(b) \tag{3.1.12}$$

Определение 3.1.5. *Пусть произведение*

$$c_1 = a_1 * b_1$$

*является операцией* $\Omega_1$*-алгебры* $A$. *Положим* $\Omega = \Omega_1 \setminus \{*\}$. *Для любой операции* $\omega \in \Omega(p)$, *умножение дистрибутивно относительно операция* $\omega$

$$a * (b_1 ... b_n \omega) = (a * b_1) ... (a * b_n) \omega \tag{3.1.13}$$

$$(b_1 ... b_n \omega) * a = (b_1 * a) ... (b_n * a) \omega \tag{3.1.14}$$

$\Omega_1$*-алгебра* $A$ *называется* **мультипликативной $\Omega$-группой**.[3.1] □

Определение 3.1.6. *Пусть* $A$, $B$ *- мультипликативные* $\Omega$*-группы. Отображение*

$$f : A \to B$$

*называется* **мультипликативным**, *если*

$$f(a * b) = f(a) \circ f(b)$$

□

Теорема 3.1.7. *Однотранзитивное представление мультипликативной* $\Omega$*-группы является мультипликативным отображением.*

Доказательство. Утверждение является следствием равенства (3.1.12) и определения 3.1.6. □

Однако утверждение теоремы 3.1.7 недостаточно, чтобы доказать равенство выражений (3.1.9) и (3.1.10).

Определение 3.1.8. *Если*

$$a * b = b * a \tag{3.1.15}$$

*то мультипликативная* $\Omega$*-группа называется* **абелевой**. □

Теорема 3.1.9. *Пусть*

$$f : A \xrightarrow{*} M$$

*эффективное представление абелевой мультипликативной* $\Omega$*-группы* $A$. *Тогда*

$$\begin{aligned} & f(a_k) \circ (f(a_l) \circ R(m_1, ..., m_k, ..., m_l, ..., m_n)) \\ = & f(a_l) \circ (f(a_k) \circ R(m_1, ..., m_k, ..., m_l, ..., m_n)) \end{aligned} \tag{3.1.16}$$

Доказательство. Из равенств (3.1.9), (3.1.10), (3.1.12), (3.1.15) следует, что

$$\begin{aligned} & f(a_k) \circ (f(a_l) \circ R(m_1, ..., m_k, ..., m_l, ..., m_n)) \\ = & (f(a_k) \circ f(a_l)) \circ R(m_1, ..., m_k, ..., m_l, ..., m_n) \\ = & f(a_k * a_l) \circ R(m_1, ..., m_k, ..., m_l, ..., m_n) \\ = & f(a_l * a_k) \circ R(m_1, ..., m_k, ..., m_l, ..., m_n) \\ = & (f(a_l) \circ f(a_k)) \circ R(m_1, ..., m_k, ..., m_l, ..., m_n) \\ = & f(a_l) \circ (f(a_k) \circ R(m_1, ..., m_k, ..., m_l, ..., m_n)) \end{aligned} \tag{3.1.17}$$

[3.1] Определение мультипликативной $\Omega$-группы похоже на определение [6]-2.1.3 $\Omega$-группы. Однако, $\Omega$-группа предполагает сложение как групповую операцию. Для нас существенно, что групповая операция мультипликативной $\Omega$-группы является произведением. Кроме того, операция $\omega$ $\Omega$-группы дистрибутивна относительно сложения. В мультипликативной $\Omega$-группе, произведение дистрибутивно относительно операции $\omega$.

Равенство (3.1.16) является следствием равенства (3.1.17). □

Теорема 3.1.10. *Пусть $A$ - абелевая мультипликативная $\Omega$-группа. Пусть $R$ - приведенный полиморфизм эффективных представлений $f_1$, ..., $f_n$ в эффективное представление $f$. Тогда для любых $k$, $l$, $k=1$, ..., $n$, $l=1$, ..., $n$,*

$$\begin{aligned}&R(m_1,...,f_k(a)\circ m_k,...,m_l,...,m_n)\\=&R(m_1,...,m_k,...,f_l(a)\circ m_l,...,m_n)\end{aligned} \tag{3.1.18}$$

Доказательство. Равенство (3.1.18) непосредственно следует из равенства (3.1.4). □

Теорема 3.1.11. *Пусть*

$$A \xrightarrow{*} B_1 \qquad A \xrightarrow{*} B_2 \qquad A \xrightarrow{*} B$$

*эффективные представления абелевой мультиплиткативной $\Omega_1$-группы $A$ в $\Omega_2$-алгебрах $B_1$, $B_2$, $B$. Допустим $\Omega_2$-алгебра имеет 2 операции, а именно $\omega_1\in\Omega_2(p)$, $\omega_2\in\Omega_2(q)$. Необходимым условием существования приведенного полиморфизма*

$$R:B_1\times B_2\rightarrow B$$

*является равенство*

$$(a_{1\cdot 1}...a_{1\cdot q}\omega_2)...(a_{p\cdot 1}...a_{p\cdot q}\omega_2)\omega_1=(a_{1\cdot 1}...a_{p\cdot 1}\omega_1)...(a_{1\cdot q}...a_{p\cdot q}\omega_1)\omega_2 \tag{3.1.19}$$

Доказательство. Пусть $a_1$, ..., $a_p\in B_1$, $b_1$, ..., $b_q\in B_2$. Согласно равенству (3.1.5), выражение

$$R(a_1...a_p\omega_1,b_1...b_q\omega_2) \tag{3.1.20}$$

может иметь 2 значения

$$\begin{aligned}&R(a_1...a_p\omega_1,b_1...b_q\omega_2)\\=&R(a_1,b_1...b_q\omega_2)...R(a_p,b_1...b_q\omega_2)\omega_1\\=&(R(a_1,b_1)...R(a_1,b_q)\omega_2)...(R(a_p,b_1)...R(a_p,b_q)\omega_2)\omega_1\end{aligned} \tag{3.1.21}$$

$$\begin{aligned}&R(a_1...a_p\omega_1,b_1...b_q\omega_2)\\=&R(a_1...a_p\omega_1,b_1)...R(a_1...a_p\omega_1,b_q)\omega_2\\=&(R(a_1,b_1)...R(a_p,b_1)\omega_1)...(R(a_1,b_q)...R(a_p,b_q)\omega_1)\omega_2\end{aligned} \tag{3.1.22}$$

Из равенств (3.1.21), (3.1.21) следует, что

$$\begin{aligned}&(R(a_1,b_1)...R(a_1,b_q)\omega_2)...(R(a_p,b_1)...R(a_p,b_q)\omega_2)\omega_1\\=&(R(a_1,b_1)...R(a_p,b_1)\omega_1)...(R(a_1,b_q)...R(a_p,b_q)\omega_1)\omega_2\end{aligned} \tag{3.1.23}$$

Следовательно, выражение (3.1.20) определенно корректно тогда и только тогда, когда равенство (3.1.23) верно. Положим

$$a_{i\cdot j}=R(a_i,b_j)\in A \tag{3.1.24}$$

Равенство (3.1.19) является следствием равенств (3.1.23), (3.1.24). □

Теорема 3.1.12. *Существует приведенный полиморфизм эффективного представления абелевой мультиплиткативной $\Omega$-группы в абелевой группе.*

Доказательство. Поскольку операция сложения в абелевой группе коммутативна и ассоциативна, то теорема является следствием теоремы 3.1.11. □

Теорема 3.1.13. *Не существует приведенный полиморфизм эффективного представления абелевой мультиплиткативной $\Omega$-группы в кольце.*

Доказательство. В кольце определены две операции: сложение, которое коммутативно и ассоциативно, и произведение, которое дистрибутивно относительно сложения. Согласно теореме 3.1.11, если существует полиморфизм эффективного представления в кольцо, то сложение и произведение должны удовлетворять равенству

$$a_{1\cdot 1}a_{2\cdot 1}+a_{1\cdot 2}a_{2\cdot 2}=(a_{1\cdot 1}+a_{1\cdot 2})(a_{2\cdot 1}+a_{2\cdot 2}) \tag{3.1.25}$$

Однако правая часть равенства (3.1.25) имеет вид

$$\begin{aligned}(a_{1\cdot 1}+a_{1\cdot 2})(a_{2\cdot 1}+a_{2\cdot 2})&=(a_{1\cdot 1}+a_{1\cdot 2})a_{2\cdot 1}+(a_{1\cdot 1}+a_{1\cdot 2})a_{2\cdot 2}\\&=a_{1\cdot 1}a_{2\cdot 1}+a_{1\cdot 2}a_{2\cdot 1}+a_{1\cdot 1}a_{2\cdot 2}+a_{1\cdot 2}a_{2\cdot 2}\end{aligned}$$

Следовательно, равенство (3.1.25) не верно. □

Вопрос 3.1.14. *Возможно, что полиморфизм представлений существует только для эффективного представления в Абелевая группе. Однако это утверждение пока не доказано.* □

## 3.2. Конгруэнция

Теорема 3.2.1. *Пусть $N$ - отношение эквивалентности на множестве $A$. Рассмотрим категорию $\mathcal{A}$ объектами которой являются отображения*[3.2]

$$f_1:A\rightarrow S_1\quad \ker f_1\supseteq N$$

$$f_2:A\rightarrow S_2\quad \ker f_2\supseteq N$$

*Мы определим морфизм $f_1\rightarrow f_2$ как отображение $h:S_1\rightarrow S_2$, для которого коммутативна диаграмма*

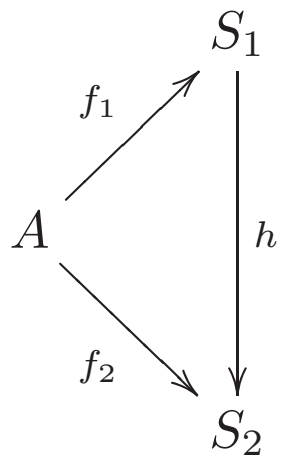

*Отображение*

$$\mathrm{nat}\,N:A\rightarrow A/N$$

*является универсально отталкивающим в категории $\mathcal{A}$.*[3.3]

[3.2] Утверждение леммы аналогично утверждению на с. [1]-94.

[3.3] Определение универсального объекта смотри в определении на с. [1]-47.

Доказательство. Рассмотрим диаграмму

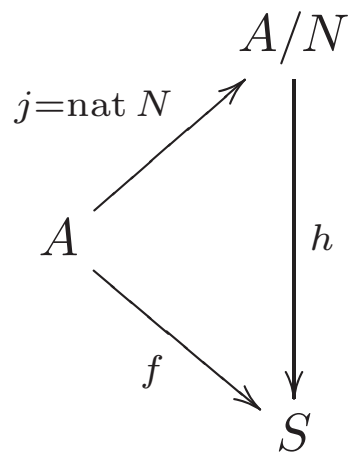

(3.2.1) $$\ker f \supseteq N$$

Из утверждения (3.2.1) и равенства

$$j(a_1) = j(a_2)$$

следует

$$f(a_1) = f(a_2)$$

Следовательно, мы можем однозначно определить отображение $h$ с помощью равенства

$$h(j(b)) = f(b)$$

□

Теорема 3.2.2. *Пусть*

$$f : A \xrightarrow{*} B$$

*представление* $\Omega_1$*-алгебры* $A$ *в* $\Omega_2$*-алгебре* $B$*. Пусть* $N$ *- такая конгруэнция*[3.4] *на* $\Omega_2$*-алгебре* $B$*, что любое преобразование* $h \in {}^*B$ *согласованно с конгруэнцией* $N$*. Существует представление*

$$f_1 : A \xrightarrow{*} B/N$$

$\Omega_1$*-алгебры* $A$ *в* $\Omega_2$*-алгебре* $B/N$ *и отображение*

$$\operatorname{nat} N : B \to B/N$$

*является приведенным морфизмом представления* $f$ *в представление* $f_1$

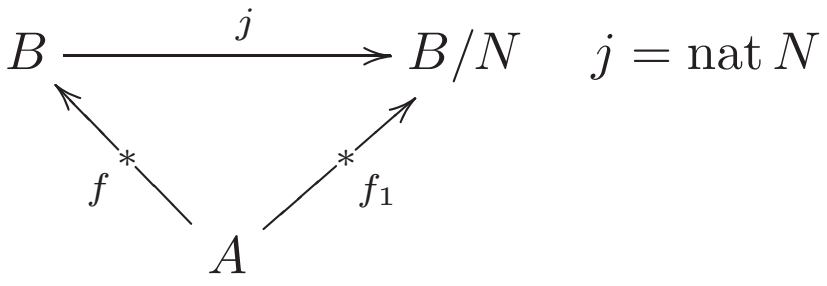

Доказательство. Любой элемент множества $B/N$ мы можем представить в виде $j(a)$, $a \in B$.

Согласно теореме [9]-II.3.5, мы можем определить единственную структуру $\Omega_2$-алгебры на множестве $B/N$. Если $\omega \in \Omega_2(p)$, то мы определим операцию $\omega$ на множестве $B/N$ согласно равенству (3) на странице [9]-73

(3.2.2) $$j(b_1)...j(b_p)\omega = j(b_1...b_p\omega)$$

Также как в доказательстве теоремы [4]-2.2.16, мы можем определить представление

$$f_1 : A \xrightarrow{*} B/N$$

[3.4] Смотри определение конгруэнции на с. [9]-71.

с помощью равенства

(3.2.3) $$f_1(a)\circ j(b)=j(f(a)\circ b)$$

Равенство (3.2.3) можно представить с помощью диаграммы

(3.2.4)

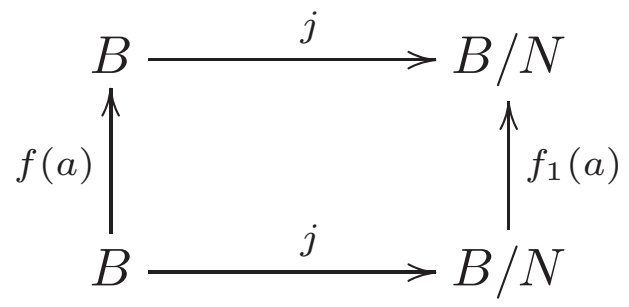

Пусть $\omega\in\Omega_2(p)$. Так как отображения $f(a)$ и $j$ являются гомоморфизмами $\Omega_2$-алгебры, то

(3.2.5)
$$\begin{aligned}f_1(a)\circ(j(b_1)...j(b_p)\omega)&=f_1(a)\circ j(b_1...b_p\omega)\\&=j(f(a)\circ(b_1...b_p\omega))\\&=j((f(a)\circ b_1)...(f(a)\circ b_p)\omega)\\&=j(f(a)\circ b_1)...j(f(a)\circ b_p)\omega\\&=(f_1(a)\circ j(b_1))...(f_1(a)\circ j(b_p))\omega\end{aligned}$$

Из равенства (3.2.5) следует, что отображение $f_1(a)$ является гомоморфизмом $\Omega_2$-алгебры. Из равенства (3.2.3) и [4]-2.3.2, следует, что отображение $j$ является приведенным морфизмом представления $f$ в представление $f_1$. □

Теорема 3.2.3. *Пусть*
$$f:A\xrightarrow{*}B$$
*представление $\Omega_1$-алгебры $A$ в $\Omega_2$-алгебре $B$. Пусть $N$ - такая конгруэнция на $\Omega_2$-алгебре $B$, что любое преобразование $h\in{}^*B$ согласованно с конгруэнцией $N$. Рассмотрим категорию $\mathcal{A}$ объектами которой являются приведенные морфизмы представлений*[3.5]

$$R_1:B\to S_1\quad \ker R_1\supseteq N$$
$$R_2:B\to S_2\quad \ker R_2\supseteq N$$

*где $S_1$, $S_2$ - $\Omega_2$-алгебры и*

$$g_1:A\xrightarrow{*}S_1\qquad g_2:A\xrightarrow{*}S_2$$

*представления $\Omega_1$-алгебры $A$. Мы определим морфизм $R_1\to R_2$ как приведенный морфизм представлений $h:S_1\to S_2$, для которого коммутативна диаграмма*

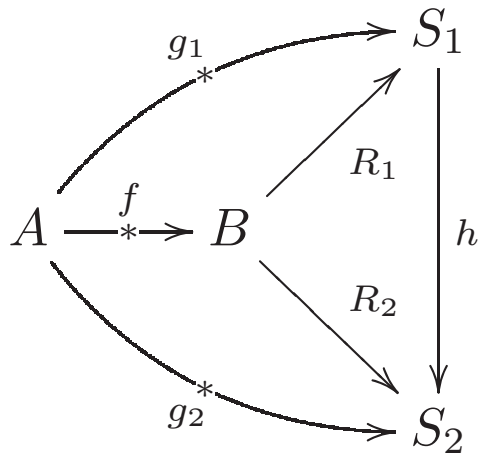

[3.5]Утверждение леммы аналогично утверждению на с. [1]-94.

*Приведенный морфизм* nat $N$ *представления* $f$ *в представление* $f_1$ *(теорема 3.2.2) является универсально отталкивающим в категории* $\mathcal{A}$.[3.6]

Доказательство. Существование и единственность отображения $h$, для которого коммутативна диаграмма

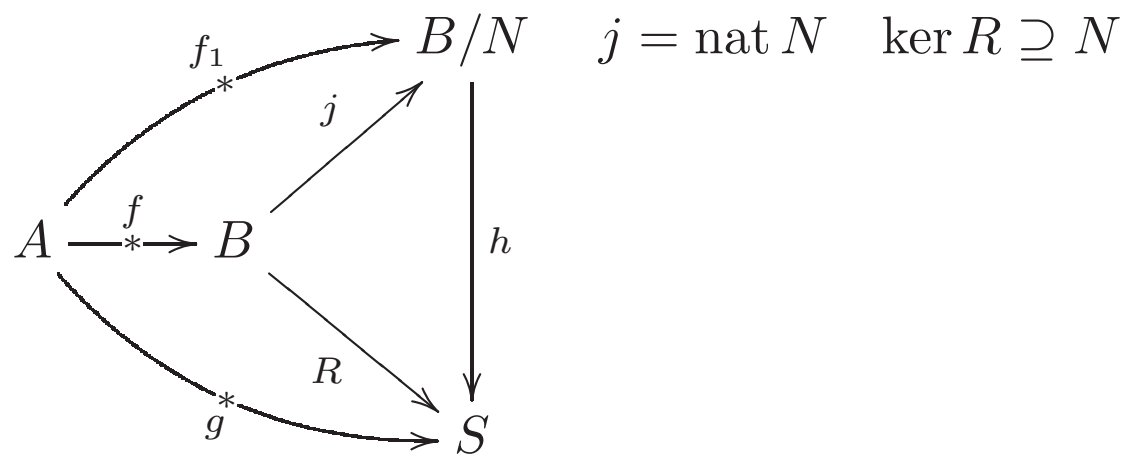

следует из теоремы 3.2.1. Следовательно, мы можем однозначно определить отображение $h$ с помощью равенства

$$h(j(b)) = R(b) \tag{3.2.6}$$

Пусть $\omega \in \Omega_2(p)$. Так как отображения $R$ и $j$ являются гомоморфизмами $\Omega_2$-алгебры, то

$$\begin{aligned} h(j(b_1)...j(b_p)\omega) &= h(j(b_1...b_p\omega)) \\ &= R(b_1...b_p\omega) \\ &= R(b_1)...R(b_p)\omega \\ &= h(j(b_1))...h(j(b_p))\omega \end{aligned} \tag{3.2.7}$$

Из равенства (3.2.7) следует, что отображение $h$ является гомоморфизмом $\Omega_2$-алгебры.

Так как отображение $R$ является приведенным морфизмом представления $f$ в представление $g$, то верно равенство

$$g(a)(R(b)) = R(f(a)(b)) \tag{3.2.8}$$

Из равенства (3.2.6) следует

$$g(a)(h(j(b))) = g(a)(R(b)) \tag{3.2.9}$$

Из равенств (3.2.8), (3.2.9) следует

$$g(a)(h(j(b))) = R(f(a)(b)) \tag{3.2.10}$$

Из равенств (3.2.6), (3.2.10) следует

$$g(a)(h(j(b))) = h(j(f(a)(b))) \tag{3.2.11}$$

Из равенств (3.2.3), (3.2.11) следует

$$g(a)(h(j(b))) = h(f_1(a)(j(b))) \tag{3.2.12}$$

Из равенства (3.2.12) следует, что отображение $h$ является приведенным морфизмом представления $f_1$ в представление $g$. □

[3.6] Определение универсального объекта смотри в определении на с. [1]-47.

## 3.3. Тензорное произведение представлений

Определение 3.3.1. *Пусть $A$ является абелевой мультипликативной $\Omega_1$-группой. Пусть $B_1$, ..., $B_n$ - $\Omega_2$-алгебры.*[3.7] *Пусть для любого $k$, $k = 1$, ..., $n$,*

$$f_k : A \longrightarrow^{*} B_k$$

*эффективное представление мультипликативной $\Omega_1$-группы $A$ в $\Omega_2$-алгебре $B_k$. Рассмотрим категорию $\mathcal{A}$ объектами которой являются приведенные полиморфизмы представлений $f_1$, ..., $f_n$*

$$r_1 : B_1 \times ... \times B_n \longrightarrow S_1 \qquad r_2 : B_1 \times ... \times B_n \longrightarrow S_2$$

*где $S_1$, $S_2$ - $\Omega_2$-алгебры и*

$$g_1 : A \longrightarrow^{*} S_1 \qquad g_2 : A \longrightarrow^{*} S_2$$

*эффективные представления мультипликативной $\Omega_1$-группы $A$. Мы определим морфизм $R_1 \to R_2$ как приведенный морфизм представлений $h : S_1 \to S_2$, для которого коммутативна диаграмма*

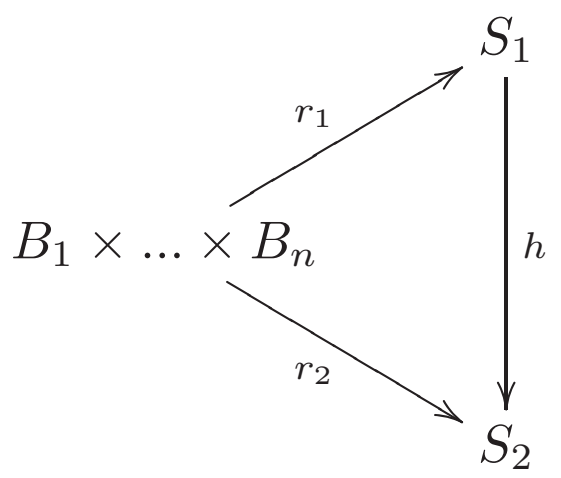

*Универсальный объект $B_1 \otimes ... \otimes B_n$ категории $\mathcal{A}$ называется* **тензорным произведением** *представлений $B_1$, ..., $B_n$.* □

Теорема 3.3.2. *Если тензорное произведение эффективных представлений существует, то тензорное произведение определено однозначно с точностью до изоморфизма представлений.*

Доказательство. Пусть $A$ является абелевой мультипликативной $\Omega_1$-группой. Пусть $B_1$, ..., $B_n$ - $\Omega_2$-алгебры. Пусть для любого $k$, $k = 1$, ..., $n$,

$$f_k : A \longrightarrow^{*} B_k$$

эффективное представление мультипликативной $\Omega_1$-группы $A$ в $\Omega_2$-алгебре $B_k$. Пусть эффективные представления

$$g_1 : A \longrightarrow^{*} S_1 \qquad g_2 : A \longrightarrow^{*} S_2$$

[3.7] Я определяю тензорное произведение представлений универсальной алгебры по аналогии с определением в [1], с. 456 - 458.

являются тензорным произведением представлений $B_1$, ..., $B_n$. Из коммутативности диаграммы

(3.3.1)

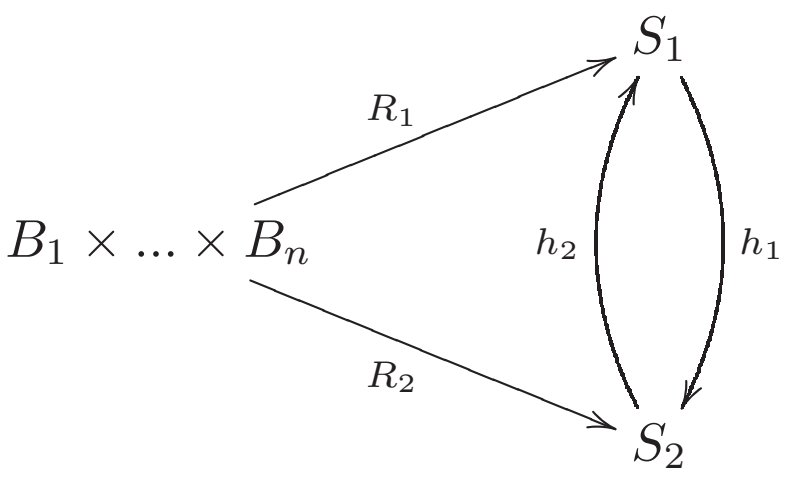

следует, что

$$
\begin{aligned}
R_1 &= h_2 \circ h_1 \circ R_1 \\
R_2 &= h_1 \circ h_2 \circ R_2
\end{aligned}
\tag{3.3.2}
$$

Из равенств (3.3.2) следует, что морфизмы представления $h_1 \circ h_2$, $h_2 \circ h_1$ являются тождественными отображениями. Следовательно, морфизмы представления $h_1$, $h_2$ являются изоморфизмами. □

Соглашение 3.3.3. *Алгебры $S_1$, $S_2$ могут быть различными множествами. Однако они неразличимы для нас, если мы рассматриваем их как изоморфные представления. В этом случае мы будем писать $S_1 = S_2$.* □

Определение 3.3.4. *Тензорное произведение*

$$B^{\otimes n} = B_1 \otimes ... \otimes B_n \quad B_1 = ... = B_n = B$$

*называется* **тензорной степенью** *представления $B$.* □

Теорема 3.3.5. *Если существует полиморфизм представлений, то тензорное произведение представлений существует.*

Доказательство. Пусть

$$f : A \longrightarrow\!\!\!\!\!*\;\; M$$

представление $\Omega_1$-алгебры $A$, порождённое декартовым произведением $B_1 \times ... \times B_n$ множеств $B_1$, ..., $B_n$.[3.8] Инъекция

$$i : B_1 \times ... \times B_n \longrightarrow M$$

определена по правилу[3.9]

$$i \circ (b_1, ..., b_n) = (b_1, ..., b_n) \tag{3.3.3}$$

Пусть $N$ - отношение эквивалентности, порождённое равенствами[3.10]

$$(b_1, ..., b_{i\cdot 1}...b_{i\cdot p}\omega, ..., b_n) = (b_1, ..., b_{i\cdot 1}, ..., b_n)...(b_1, ..., b_{i\cdot p}, ..., b_n)\omega \tag{3.3.4}$$

$$(b_1, ..., f_i(a) \circ b_i, ..., b_n) = f(a) \circ (b_1, ..., b_i, ..., b_n) \tag{3.3.5}$$

---

[3.8] Согласно теоремам 2.1.3, 2.3.2, множество, порождённое приведенным декартовым произведением представлений $B_1$, ..., $B_n$ совпадает с декартовым произведением $B_1 \times ... \times B_n$ множеств $B_1$, ..., $B_n$. В этом месте доказательства нас не интересует алгебраическая структра на множестве $B_1 \times ... \times B_n$.

[3.9] Равенство (3.3.3) утверждает, что мы отождествляем базис представления $M$ с множеством $B_1 \times ... \times B_n$.

[3.10] Я рассматриваю формирование элементов представления из элементов базиса согласно теореме [4]-2.6.4. Теорема 3.3.11 требует выполнения условий (3.3.4), (3.3.5).

$$b_k \in B_k \quad k=1,...,n \quad b_{i\cdot 1},...,b_{i\cdot p} \in B_i \quad \omega \in \Omega_2(p) \quad a \in A$$

Лемма 3.3.6. *Пусть $\omega \in \Omega_2(p)$. Тогда*

$$\begin{aligned}&f(c)\circ(b_1,...,b_{i\cdot 1}...b_{i\cdot p}\omega,...,b_n)\\=&f(c)\circ((b_1,...,b_{i\cdot 1},...,b_n)...(b_1,...,b_{i\cdot p},...,b_n)\omega)\end{aligned} \tag{3.3.6}$$

Доказательство. Из равенства (3.3.5) следует

$$f(c)\circ(b_1,...,b_{i\cdot 1}...b_{i\cdot p}\omega,...,b_n)=(b_1,...,f_i(c)\circ(b_{i\cdot 1}...b_{i\cdot p}\omega),...,b_n) \tag{3.3.7}$$

Так как $f_i(c)$ - эндоморфизм $\Omega_2$-алгебры $B_i$, то из равенства (3.3.7) следует

$$f(c)\circ(b_1,...,b_{i\cdot 1}...b_{i\cdot p}\omega,...,b_n)=(b_1,...,(f_i(c)\circ b_{i\cdot 1})...(f_i(c)\circ b_{i\cdot p})\omega,...,b_n) \tag{3.3.8}$$

Из равенств (3.3.8), (3.3.4) следует

$$\begin{aligned}&f(c)\circ(b_1,...,b_{i\cdot 1}...b_{i\cdot p}\omega,...,b_n)\\=&(b_1,...,f_i(c)\circ b_{i\cdot 1},...,b_n)...(b_1,...,f_i(c)\circ b_{i\cdot p},...,b_n)\omega\end{aligned} \tag{3.3.9}$$

Из равенств (3.3.9), (3.3.5) следует

$$\begin{aligned}&f(c)\circ(b_1,...,b_{i\cdot 1}...b_{i\cdot p}\omega,...,b_n)\\=&(f(c)\circ(b_1,...,b_{i\cdot 1},...,b_n))...(f(c)\circ(b_1,...,b_{i\cdot p},...,b_n))\omega\end{aligned} \tag{3.3.10}$$

Так как $f(c)$ - эндоморфизм $\Omega_2$-алгебры $B$, то равенство (3.3.6) следует из равенства (3.3.10). ⊙

Лемма 3.3.7.

$$f(c)\circ(b_1,...,f_i(a)\circ b_i,...,b_n)=f(c)\circ(f(a)\circ(b_1,...,b_i,...,b_n)) \tag{3.3.11}$$

Доказательство. Из равенства (3.3.5) следует, что

$$\begin{aligned}f(c)\circ(b_1,...,f_i(a)\circ b_i,...,b_n)&=(b_1,...,f_i(c)\circ(f_i(a)\circ b_i),...,b_n)\\&=(b_1,...,(f_i(c)\circ f_i(a))\circ b_i,...,b_n)\\&=(f(c)\circ f(a))\circ(b_1,...,b_i,...,b_n)\\&=f(c)\circ(f(a)\circ(b_1,...,b_i,...,b_n))\end{aligned} \tag{3.3.12}$$

Равенство (3.3.11) следует из равенства (3.3.12). ⊙

Лемма 3.3.8. *Для любого $c \in A$ эндоморфизм $f(c)$ $\Omega_2$-алгебры $M$ согласовано с эквивалентностью $N$.*

Доказательство. Утверждение леммы следует из лемм 3.3.6, 3.3.7 и определения [4]-2.2.14. ⊙

Из леммы 3.3.8 и теоремы [7]-3.2.15, следует, что на множестве ${}^*M/N$ определена $\Omega_1$-алгебра. Рассмотрим диаграмму

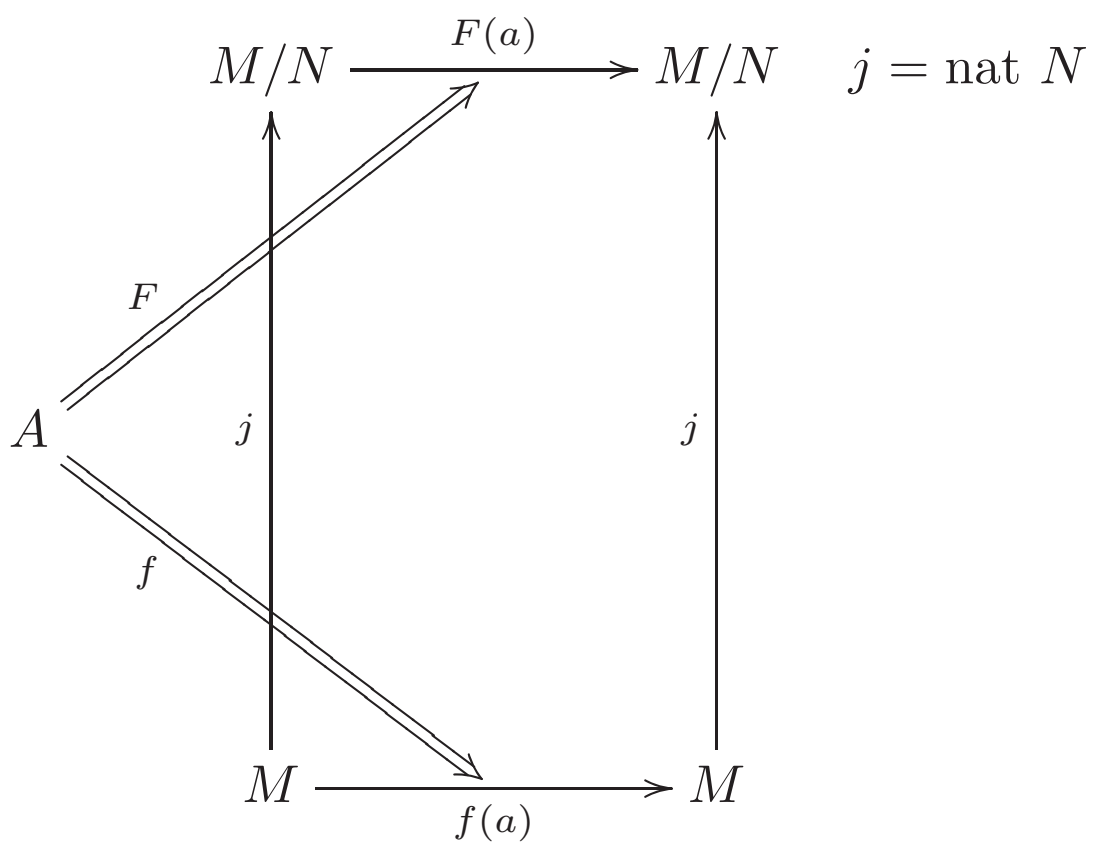

Согласно лемме 3.3.8, из условия

$$j\circ b_1=j\circ b_2$$

следует

$$j\circ(f(a)\circ b_1)=j\circ(f(a)\circ b_2)$$

Следовательно, преобразование $F(a)$ определено корректно и

$$F(a)\circ j=j\circ f(a) \tag{3.3.13}$$

Если $\omega\in\Omega_1(p)$, то мы положим

$$(F(a_1)...F(a_p)\omega)\circ(J\circ b)=J\circ((f(a_1)...f(a_p)\omega)\circ b)$$

Следовательно, отображение $F$ является представлением $\Omega_1$-алгебры $A$. Из (3.3.13) следует, что $j$ является приведенным морфизмом представлений $f$ и $F$.

Рассмотрим коммутативную диаграмму

(3.3.14)

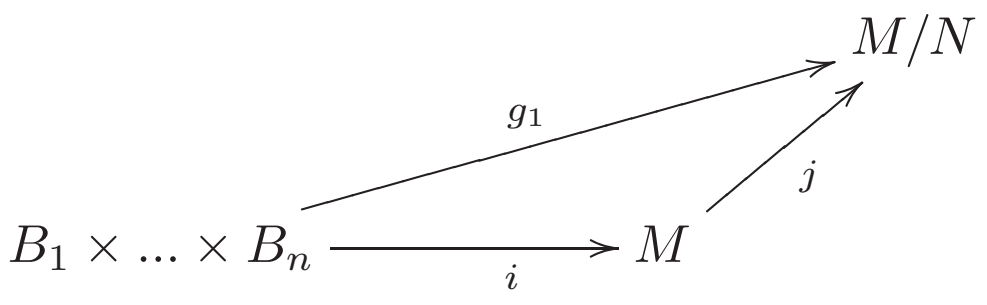

Из коммутативности диаграммы (3.3.14) и равенства (3.3.3) следует, что

$$g_1\circ(b_1,...,b_n)=j\circ(b_1,...,b_n) \tag{3.3.15}$$

Из равенств (3.3.3), (3.3.4), (3.3.5) следует

$$\begin{aligned}&g_1\circ(b_1,...,b_{i\cdot 1}...b_{i\cdot p}\omega,...,b_n)\\=&(g_1\circ(b_1,...,b_{i\cdot 1},...,b_n))...(g_1\circ(b_1,...,b_{i\cdot p},...,b_n))\omega\end{aligned} \tag{3.3.16}$$

$$g_1\circ(b_1,...,f_i(a)\circ b_i,...,b_n)=f(a)\circ(g_1\circ(b_1,...,b_i,...,b_n)) \tag{3.3.17}$$

Из равенств (3.3.16) и (3.3.17) следует, что отображение $g_1$ является приведенным полиморфизмом представлений $f_1$, ..., $f_n$.

Поскольку $B_1\times...\times B_n$ - базис представления $M$ $\Omega_1$-алгебры $A$, то, согласно теореме [4]-2.7.7, для любого представления

$$A \xrightarrow{*} V$$

и любого приведенного полиморфизма

$$g_2 : B_1 \times ... \times B_n \longrightarrow V$$

существует единственный морфизм представлений $k : M \to V$, для которого коммутативна следующая диаграмма

(3.3.18)

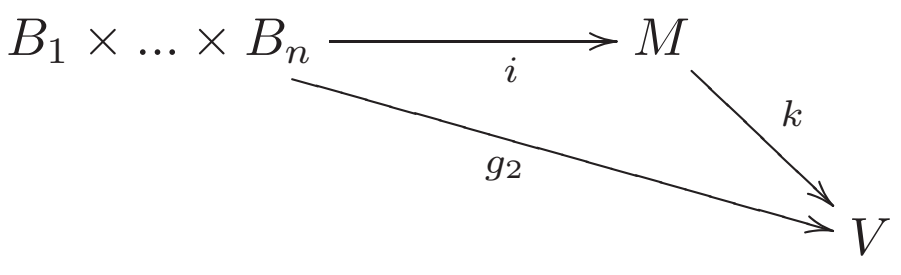

Так как $g_2$ - приведенный полиморфизм, то $\ker k \supseteq N$.

Согласно теореме 3.2.3 отображение $j$ универсально в категории морфизмов представления $f$, ядро которых содержит $N$. Следовательно, определён морфизм представлений

$$h : M/N \to V$$

для которого коммутативна диаграмма

(3.3.19)

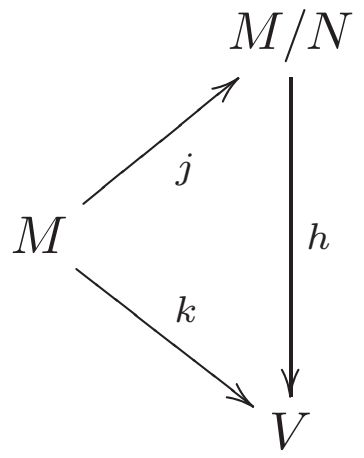

Объединяя диаграммы (3.3.14), (3.3.18), (3.3.19), получим коммутативную диаграмму

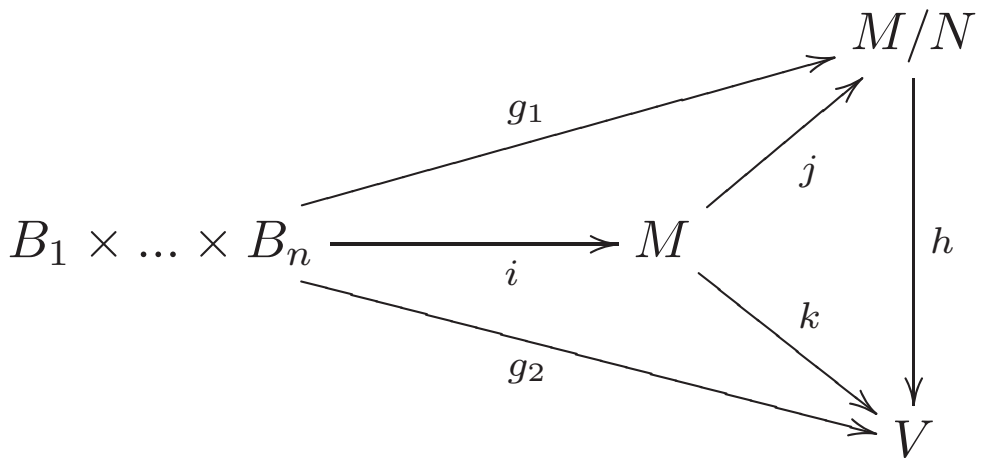

Так как $\operatorname{Im} g_1$ порождает $M/N$, то отображение $h$ однозначно определено. □

Согласно доказательству теоремы 3.3.5

$$B_1 \otimes ... \otimes B_n = M/N$$

Для $d_i \in A_i$ будем записывать

(3.3.20) $$j \circ (d_1, ..., d_n) = d_1 \otimes ... \otimes d_n$$

Из равенств (3.3.15), (3.3.20) следует, что

(3.3.21) $$g_1 \circ (d_1, ..., d_n) = d_1 \otimes ... \otimes d_n$$

Теорема 3.3.9. *Отображение*

$$(x_1, ..., x_n) \in B_1 \times ... \times B_n \to x_1 \otimes ... \otimes x_n \in B_1 \otimes ... \otimes B_n$$

*является полиморфизмом.*

Доказательство. Теорема является следствием определений 3.1.3, 3.3.1. □

Теорема 3.3.10. *Пусть* $B_1$, ..., $B_n$ - $\Omega_2$*-алгебры. Пусть*

$$f : B_1 \times ... \times B_n \to B_1 \otimes ... \otimes B_n$$

*приведенный полиморфизм, определённый равенством*

$$f \circ (b_1, ..., b_n) = b_1 \otimes ... \otimes b_n \tag{3.3.22}$$

*Пусть*

$$g : B_1 \times ... \times B_n \to V$$

*приведенный полиморфизм в* $\Omega_2$*-алгебру* $V$*. Существует морфизм представлений*

$$h : B_1 \otimes ... \otimes B_n \to V$$

*такой, что диаграмма*

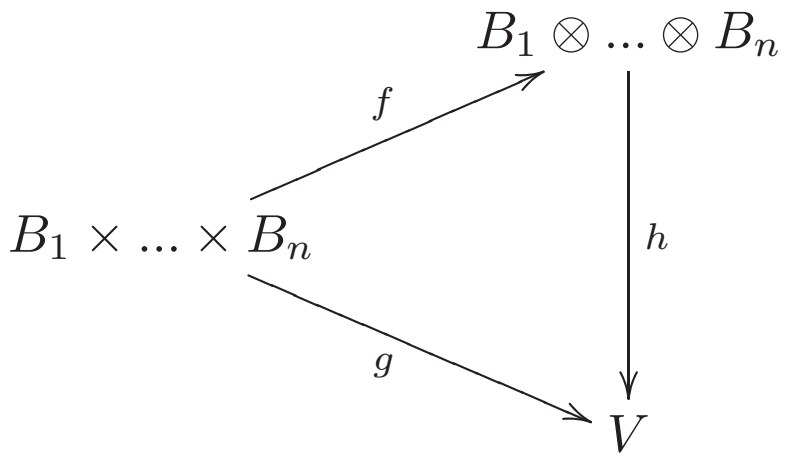

*коммутативна.*

Доказательство. Равенство (3.3.22) следует из равенств (3.3.3) и (3.3.20). Существование отображения $h$ следует из определения 3.3.1 и построений, выполненных при доказательстве теоремы 3.3.5. □

Теорема 3.3.11. *Пусть*

$$b_k \in B_k \quad k = 1, ..., n \quad b_{i\cdot 1}, ..., b_{i\cdot p} \in B_i \quad \omega \in \Omega_2(p) \quad a \in A$$

*Тензорное произведение дистрибутивно относительно операции* $\omega$

$$\begin{aligned} & b_1 \otimes ... \otimes (b_{i\cdot 1}...b_{i\cdot p}\omega) \otimes ... \otimes b_n \\ =& (b_1 \otimes ... \otimes b_{i\cdot 1} \otimes ... \otimes b_n)...(b_1 \otimes ... \otimes b_{i\cdot p} \otimes ... \otimes b_n)\omega \end{aligned} \tag{3.3.23}$$

*Представление мультипликативной* $\Omega_1$*-группы* $A$ *в тензорном произведении определено равенством*

$$b_1 \otimes ... \otimes (f_i(a) \circ b_i) \otimes ... \otimes b_n = f(a) \circ (b_1 \otimes ... \otimes b_i \otimes ... \otimes b_n) \tag{3.3.24}$$

Доказательство. Равенство (3.3.23) является следствием равенства (3.3.16) и определения (3.3.21). Равенство (3.3.24) является следствием равенства (3.3.17) и определения (3.3.21). □

## 3.4. Ассоциативность тензорного произведения

Пусть $A$ является мультипликативной $\Omega_1$-группой. Пусть $B_1, B_2, B_3$ - $\Omega_2$-алгебры. Пусть для $k = 1, 2, 3$

$$f_k : A \xrightarrow{*} B_k$$

эффективное представление мультипликативной $\Omega_1$-группы $A$ в $\Omega_2$-алгебре $B_k$.

Лемма 3.4.1. *Для заданного значения $x_3 \in B_3$, отображение*

$$h_{12} : (B_1 \otimes B_2) \times B_3 \to B_1 \otimes B_2 \otimes B_3 \tag{3.4.1}$$

*определённое равенством*

$$h_{12}(x_1 \otimes x_2, x_3) = x_1 \otimes x_2 \otimes x_3 \tag{3.4.2}$$

*является приведенным морфизмом представления $B_1 \otimes B_2$ в представление $B_1 \otimes B_2 \otimes B_3$ .*

Доказательство. Согласно теореме 3.3.9, для заданного значения $x_3 \in B_3$, отображение

$$(x_1, x_2, x_3) \in B_1 \times B_2 \times B_3 \to x_1 \otimes x_2 \otimes x_3 \in B_1 \otimes B_2 \otimes B_3 \tag{3.4.3}$$

является полиморфизмом по переменным $x_1 \in B_1$, $x_2 \in B_2$. Следовательно, для заданного значения $x_3 \in B_3$, лемма является следствием теоремы 3.3.10. □

Лемма 3.4.2. *Для заданного значения $x_{12} \in B_1 \otimes B_2$ отображение $h_{12}$ является приведенным морфизмом представления $B_3$ в представление $B_1 \otimes B_2 \otimes B_3$ .*

Доказательство. Согласно теореме 3.3.9 и равенству (3.3.21), для заданного значения $x_1 \in B_1$, $x_2 \in B_2$, отображение

$$(x_1 \otimes x_2, x_3) \in B_1 \times B_2 \times B_3 \to x_1 \otimes x_2 \otimes x_3 \in B_1 \otimes B_2 \otimes B_3 \tag{3.4.4}$$

является морфизмом по переменной $x_3 \in B_3$. Следовательно, теорема является следсвием теорем 3.1.10, 3.1.11. □

Лемма 3.4.3. *Существует приведенный морфизм представлений*

$$h : (B_1 \otimes B_2) \otimes B_3 \to B_1 \otimes B_2 \otimes B_3$$

Доказательство. Согласно леммам 3.4.1, 3.4.2 и определению 3.1.3, отображение $h_{12}$ является приведенным полиморфизмом представлений. Утверждение леммы является следствием теоремы 3.3.10. □

Лемма 3.4.4. *Существует приведенный морфизм представлений*

$$g : B_1 \otimes B_2 \otimes B_3 \to (B_1 \otimes B_2) \otimes B_3$$

Доказательство. Отображение

$$(x_1, x_2, x_3) \in B_1 \times B_2 \times B_3 \to (x_1 \otimes x_2) \otimes x_3 \in (B_1 \otimes B_2) \otimes B_3$$

является полиморфизмом по переменным $x_1 \in B_1$, $x_2 \in B_2$, $x_3 \in B_3$. Следовательно, лемма является следствием теоремы 3.3.10. □

Теорема 3.4.5.

$$(B_1 \otimes B_2) \otimes B_3 = B_1 \otimes (B_2 \otimes B_3) = B_1 \otimes B_2 \otimes B_3 \tag{3.4.5}$$

Доказательство. Согласно лемме 3.4.3, существует приведенный морфизм представлений

$$h:(B_1\otimes B_2)\otimes B_3\rightarrow B_1\otimes B_2\otimes B_3$$

Согласно лемме 3.4.4, существует приведенный морфизм представлений

$$g:B_1\otimes B_2\otimes B_3\rightarrow (B_1\otimes B_2)\otimes B_3$$

Следовательно, приведенные морфизмы представлений $h$, $g$ являются изоморфизмами, откуда следует равенство

$$(B_1\otimes B_2)\otimes B_3=B_1\otimes B_2\otimes B_3 \tag{3.4.6}$$

Аналогично мы можем доказать равенство

$$B_1\otimes (B_2\otimes B_3)=B_1\otimes B_2\otimes B_3$$

□

Замечание 3.4.6. *Очевидно, что структура $\Omega_2$-алгебр $(B_1\otimes B_2)\otimes B_3$, $B_1\otimes B_2\otimes B_3$ слегка различна. Мы записываем равенство (3.4.6), опираясь на соглашение 3.3.3 и это позволяет нам говорить об ассоциативности тензорного произведения представлений.* □

Глава 4

# $D$-модуль

## 4.1. Модуль над коммутативным кольцом

Теорема 4.1.1. *Пусть кольцо $D$ имеет единицу $e$. Представление*

$$f : D \overset{*}{\longrightarrow} V$$

*кольца $D$ в абелевой группе $A$* **эффективно** *тогда и только тогда, когда из равенства $f(a) = 0$ следует $a = 0$.*

Доказательство. Сумма преобразований $f$ и $g$ абелевой группы определяется согласно правилу

$$(f + g)(a) = f(a) + g(a)$$

Поэтому, рассматривая представление кольца $D$ в абелевой группе $A$, мы полагаем

$$f(a + b)(x) = f(a)(x) + f(b)(x)$$

Если $a, b \in R$ порождают одно и то же преобразование, то

$$f(a) \circ m = f(b) \circ m \tag{4.1.1}$$

для любого $m \in A$. Из равенства (4.1.1) следует, что $a - b$ порождает нулевое преобразование

$$f(a - b) \circ m = 0$$

Элемент $e + a - b$ порождает тождественное преобразование. Следовательно, представление $f$ эффективно тогда и только тогда, когда $a = b$. □

Определение 4.1.2. *Эффективное представление коммутативного кольца $D$ в абелевой группе $V$*

$$f : D \overset{*}{\longrightarrow} V \qquad f(d) : v \to d\,v \tag{4.1.2}$$

*называется* **модулем над кольцом** $D$ *или* $D$**-модулем**. □

Теорема 4.1.3. *Элементы $D$-модуля $V$ удовлетворяют соотношениям*

- **закону ассоциативности**

$$(ab)m = a(bm) \tag{4.1.3}$$

- **закону дистрибутивности**

$$a(m + n) = am + an \tag{4.1.4}$$

$$(a + b)m = am + bm \tag{4.1.5}$$

- **закону унитарности**

$$1m = m \tag{4.1.6}$$

*для любых $a, b \in D$, $m, n \in V$.*

Доказательство. Равенство (4.1.4) следует из утверждения, что преобразование $a$ является эндоморфизмом абелевой группы. Равенство (4.1.5) следует из утверждения, что представление (4.1.2) является гомоморфизмом аддитивной группы кольца $D$. Равенства (4.1.3) и (4.1.6) следуют из утверждения, что представление (4.1.2) является левосторонним представлением мультипликативной группы кольца $D$. □

Теорема 4.1.4. *Множество векторов, порождённое множеством векторов* $v=(v_i\in V, i\in I)$ *имеет вид*

$$J(v)=\left\{w:w=\sum_{i\in I}c^iv_i, c^i\in D\right\}\tag{4.1.7}$$

Доказательство. Мы докажем теорему по индукции, опираясь на теорему [4]-2.6.4.

Пусть $k=0$. Согласно теореме [4]-2.6.4, $X_0=v$. Для произвольного $v_k\in v$, положим $c^i=\delta^i_k$. Тогда

$$v_k=\sum_{i\in I}c^iv_i\tag{4.1.8}$$

$v_k\in J(v)$ следует из (4.1.7), (4.1.8).

Пусть $X_{k-1}\subseteq J(v)$.

- Пусть $w_1$, $w_2\in X_{k-1}$. Так как $V$ является абелевой группой, то согласно утверждению [4]-2.6.4.3, $w_1+w_2\in X_k$. Согласно утверждениям [4]-(2.6.1), (4.1.7), существуют $D$-числа $c^i$, $d^i$, $i\in I$, такие, что

$$\begin{aligned}w_1&=\sum_{i\in I}c^iv_i\\ w_2&=\sum_{i\in I}d^iv_i\end{aligned}\tag{4.1.9}$$

  Так как $V$ является абелевой группой, то из равенства (4.1.9) следует, что

$$w_1+w_2=\sum_{i\in I}c^iv_i+\sum_{i\in I}d^iv_i=\sum_{i\in I}(c^iv_i+d^iv_i)\tag{4.1.10}$$

  Равенство

$$w_1+w_2=\sum_{i\in I}(c^i+d^i)v_i\tag{4.1.11}$$

  является следствием равенств (4.1.5), (4.1.10). Из равенства (4.1.11) следует, что $w_1+w_2\in J(v)$.
- Пусть $w\in X_{k-1}$. Согласно утверждению [4]-2.6.4.4, для любого $D$-числа $a$, $aw\in X_k$. Согласно утверждениям [4]-(2.6.1), (4.1.7), существуют $D$-числа $c^i$, $i\in I$, такие, что

$$w=\sum_{i\in I}c^iv_i\tag{4.1.12}$$

  Из равенства (4.1.12) следует, что

$$aw=a\sum_{i\in I}c^iv_i=\sum_{i\in I}a(c^iv_i)=\sum_{i\in I}(ac^i)v_i\tag{4.1.13}$$

Из равенства (4.1.13) следует, что $aw \in J(v)$.

□

Соглашение 4.1.5. *Мы будем пользоваться соглашением Эйнштейна о сумме, в котором повторяющийся индекс (один вверху и один внизу) подразумевает сумму по повторяющемуся индексу. В этом случае предполагается известным множество индекса суммирования и знак суммы опускается*

$$c^i v_i = \sum_{i \in I} c^i v_i$$

□

Определение 4.1.6. *Пусть* $v = (v_i \in V, i \in I)$ *- множество векторов. Выражение* $c^i v_i$ *называется* **линейной комбинацией векторов** $v_i$. *Вектор*

$$w = c^i v_i$$

*называется* **вектором, линейно зависимым от векторов** $v_i$. □

Теорема 4.1.7. *Пусть $D$ - поле. Если уравнение*

$$c^i v_i = 0 \tag{4.1.14}$$

*предполагает существования индекса* $i = j$ *такого, что* $c^j \neq 0$, *то вектор* $v_j$ *линейно зависит от остальных векторов* $v$.

Доказательство. Теорема является следствием равенства

$$v_j = \sum_{i \in I \setminus \{j\}} \frac{c^i}{c^j} v_i$$

и определения 4.1.6. □

Очевидно, что для любого множества векторов $v_i$,

$$0 = c^i v_i \quad c^i = 0$$

Определение 4.1.8. *Векторы*[4.1] $e_i$, $i \in I$, *$D$-модуля $A$* **линейно независимы**, *если $c = 0$ следует из уравнения*

$$c^i e_i = 0$$

*В противном случае, векторы* $e_i$, $i \in I$, **линейно зависимы**. □

Теорема 4.1.9. *Пусть $D$ - поле. Множество векторов* $\overline{\overline{e}} = (e_i, i \in I)$ *является* **базисом $D$-векторного пространства** $V$, *если векторы $e_i$ линейно независимы и любой вектор $v \in V$ линейно зависит от векторов $e_i$.*

Доказательство. Пусть $\overline{\overline{e}} = (e_i, i \in I)$ - базис $D$-векторного пространства $V$. Согласно определению [4]-2.7.1 и теоремам [4]-2.6.4, 4.1.4, произвольный вектор $v \in V$ является линейной комбинацией векторов $e_i$

$$v = v^i e_i \tag{4.1.15}$$

Из равенства (4.1.15) следует, что множество векторов $v$, $e_i$, $i \in I$, не является линейно независимым.

[4.1] Я следую определению в [1], с. 100.

Рассмотрим равенство

$$c^i e_i = 0 \tag{4.1.16}$$

Согласно теореме 4.1.7, если

$$c^j \neq 0 \tag{4.1.17}$$

то вектор $e_j$ линейно зависит от остальных векторов $e$. Следовательно, множество векторов $e_i$, $i \in I \setminus \{j\}$, порождает $D$-векторное пространство $V$. Согласно определению [4]-2.7.1, утверждение (4.1.17) неверно. Согласно определению 4.1.8, векторы $e_i$ линейно независимы. □

Определение 4.1.10. *Множество векторов* $\overline{\overline{e}} = (e_i, i \in I)$ *называется*[4.2] **базисом $D$-модуля** $V$, *если произвольный вектор* $v \in V$ *является линейной комбинацией векторов базиса и произвольный вектор базиса нельзя представить в виде линейной комбинации остальных векторов базиса.* $A$ - **свободный модуль над кольцом** $D$, *если* $A$ *имеет базис над кольцом* $D$.[4.3] □

Определение 4.1.11. *Пусть* $\overline{\overline{e}}$ *- базис* $D$*-модуля* $A$, *и* $A$*-число* $a$ *имеет разложение*

$$a = a^i e_i$$

*относительно базиса* $\overline{\overline{e}}$. $D$*-числа* $a^i$ *называются* **координатами** $A$*-числа* $a$ *относительно базиса* $\overline{\overline{e}}$. □

## 4.2. Линейное отображение $D$-модуля

Определение 4.2.1. *Пусть* $A_1$, $A_2$ *-* $D$*-модули. Морфизм представлений*

$$f : A_1 \to A_2$$

$D$*-модуля* $A_1$ *в* $D$*-модуль* $A_2$ *называется* **линейным отображением** $D$*-модуля* $A_1$ *в* $D$*-модуль* $A_2$. *Обозначим* $\mathcal{L}(D; A_1 \to A_2)$ *множество линейных отображений* $D$*-модуля* $A_1$ *в* $D$*-модуль* $A_2$. □

Теорема 4.2.2. *Линейное отображение*

$$f : A_1 \to A_2$$

$D$*-модуля* $A_1$ *в* $D$*-модуль* $A_2$ *удовлетворяет равенствам*[4.4]

$$f \circ (a + b) = f \circ a + f \circ b \tag{4.2.1}$$

$$f \circ (pa) = p(f \circ a) \tag{4.2.2}$$

$$a, b \in A_1 \quad p \in D$$

Доказательство. Из определения 4.2.1 и теоремы [7]-3.2.19 следует, что отображение $f$ является гомоморфизмом абелевой группы $A_1$ в абелеву группу $A_2$ (равенство (4.2.1)). Равенство (4.2.2) следует из равенства [7]-(3.2.45). □

---

[4.2] Определение 4.1.10 является следствием теоремы [4]-2.7.2 и замечания [4]-2.7.3.

[4.3] Я следую определению в [1], с. 103.

[4.4] В некоторых книгах (например, [1], с. 94) теорема 4.2.2 рассматривается как определение.

Теорема 4.2.3. *Пусть $A_1$, $A_2$ - $D$-модули. Отображение*

$$f+g : A_1 \to A_2 \quad f, g \in \mathcal{L}(D; A_1 \to A_2) \tag{4.2.3}$$

*определённое равенством*

$$(f+g)\circ x = f\circ x + g\circ x \tag{4.2.4}$$

*называется* **суммой отображений** $f$ *и* $g$ *и является линейным отображением. Множество* $\mathcal{L}(D; A_1; A_2)$ *является абелевой группой относительно суммы отображений.*

Доказательство. Согласно теореме 4.2.2

$$f\circ(v+w) = f\circ v + f\circ w \tag{4.2.5}$$

$$f\circ(d\,v) = d\,(f\circ v) \tag{4.2.6}$$

$$g\circ(v+w) = g\circ v + g\circ w \tag{4.2.7}$$

$$g\circ(d\,v) = d\,(g\circ v) \tag{4.2.8}$$

Равенство

$$\begin{aligned}(f+g)\circ(v+w) &= f\circ(v+w) + g\circ(v+w)\\ &= f\circ v + f\circ w + g\circ v + g\circ w\\ &= (f+g)\circ v + (f+g)\circ w\end{aligned} \tag{4.2.9}$$

является следствием равенств (4.2.4), (4.2.5), (4.2.7). Равенство

$$\begin{aligned}(f+g)\circ(d\,v) &= f\circ(d\,v) + g\circ(d\,v)\\ &= d\,f\circ v + d\,g\circ v\\ &= d\,(f\circ v + g\circ v)\\ &= d\,((f+g)\circ v)\end{aligned} \tag{4.2.10}$$

является следствием равенств (4.2.4), (4.2.6), (4.2.8). Из равенств (4.2.9), (4.2.10) и теоремы 4.2.2 следует, что отображение (4.2.3) является линейным отображением $D$-модулей.

Из равенства (4.2.4) следует, что отображение

$$0 : v \in A_1 \to 0 \in A_2$$

является нулём операции сложения

$$(0+f)\circ v = 0\circ v + f\circ v = f\circ v$$

Из равенства (4.2.4) следует, что отображение

$$-f : v \in A_1 \to -(f\circ v) \in A_2$$

является отображением, обратным отображению $f$

$$(f+(-f))\circ v = f\circ v + (-f)\circ v = f\circ v - f\circ v = 0 = 0\circ v$$

Из равенства

$$\begin{aligned}(f+g)\circ x &= f\circ x + g\circ x\\ &= g\circ x + f\circ x\\ &= (g+f)\circ x\end{aligned}$$

следует, что сумма отображений коммутативно. □

Теорема 4.2.4. *Пусть $A_1$, $A_2$ - $D$-модули. Отображение*

$$d\,f : A_1 \to A_2 \quad d \in D \quad f \in \mathcal{L}(D; A_1 \to A_2) \tag{4.2.11}$$

*определённое равенством*

$$(d\,f) \circ x = d(f \circ x) \tag{4.2.12}$$

*называется* **произведением отображения** $f$ **на скаляр** $d$ *и является линейным отображением. Представление*

$$a : f \in \mathcal{L}(D; A_1 \to A_2) \to af \in \mathcal{L}(D; A_1 \to A_2) \tag{4.2.13}$$

*кольца $D$ в абелевой группе $\mathcal{L}(D; A_1 \to A_2)$ порождает структуру $D$-модуля.*

Доказательство. Согласно теореме 4.2.2

$$f \circ (v + w) = f \circ v + f \circ w \tag{4.2.14}$$

$$f \circ (d\,v) = d\,(f \circ v) \tag{4.2.15}$$

Равенство

$$\begin{aligned}(d\,f) \circ (v + w) &= d(f \circ (v + w)) \\ &= d(f \circ v + f \circ w) = d(f \circ v) + d(f \circ w) \\ &= (d\,f) \circ v + (d\,f) \circ w\end{aligned} \tag{4.2.16}$$

является следствием равенств (4.2.12), (4.2.14). Равенство

$$\begin{aligned}(c\,f) \circ (d\,v) &= c(f \circ (d\,v)) = cd\,(f \circ v) = d(c(f \circ v)) \\ &= d\,((c\,f) \circ v)\end{aligned} \tag{4.2.17}$$

является следствием равенств (4.2.12), (4.2.15). Из равенств (4.2.16), (4.2.17) и теоремы 4.2.2 следует, что отображение (4.2.11) является линейным отображением $D$-модулей.

Равенство

$$(p + q)f = pf + qf \tag{4.2.18}$$

является следствием равенства

$$\begin{aligned}((a + b)f) \circ v &= (a + b)(f \circ v) = a(f \circ v) + b(f \circ v) = (af) \circ v + (bf) \circ v \\ &= (af + bf) \circ v\end{aligned}$$

Равенство

$$p(qf) = (pq)f \tag{4.2.19}$$

является следствием равенства

$$((ab)f) \circ v = (ab)(f \circ v) = a(b(f \circ v)) = a((bf) \circ v) = (a(bf)) \circ v$$

Из равенств (4.2.18), (4.2.19), следует, что отображение (4.2.13) является представлением кольца $D$ в абелевой группе $\mathcal{L}(D; A_1 \to A_2)$. Согласно определению 4.1.2 и теореме 4.2.3, абелевая группа $\mathcal{L}(D; A_1 \to A_2)$ является $D$-модулем. □

Теорема 4.2.5. *Пусть отображение*

$$f: A_1 \to A_2$$

*является линейным отображением $D$-модуля $A_1$ в $D$-модуль $A_2$. Тогда*

$$f \circ 0 = 0$$

Доказательство. Следствие равенства

$$f \circ (a + 0) = f \circ a + f \circ 0$$

□

## 4.3. Полилинейное отображение $D$-модуля

Определение 4.3.1. *Пусть $D$ - коммутативное кольцо. Приведенный полиморфизм $D$-модулей $A_1$, ..., $A_n$ в $D$-модуль $S$*

$$f: A_1 \times ... \times A_n \to S$$

*называется* **полилинейным отображением** *$D$-модулей $A_1$, ..., $A_n$ в $D$-модуль $S$. Обозначим $\mathcal{L}(D; A_1 \times ... \times A_n \to S)$ множество полилинейных отображений $D$-модулей $A_1$, ..., $A_n$ в $D$-модуль $S$. Обозначим $\mathcal{L}(D; A^n \to S)$ множество $n$-линейных отображений $D$-модуля $A$ $(A_1 = ... = A_n = A)$ в $D$-модуль $S$.* □

Теорема 4.3.2. *Пусть $D$ - коммутативное кольцо. Полилинейное отображение $D$-модулей $A_1$, ..., $A_n$ в $D$-модуль $S$*

$$f: A_1 \times ... \times A_n \to S$$

*удовлетворяет равенствам*

$$f \circ (a_1, ..., a_i + b_i, ..., a_n) = f \circ (a_1, ..., a_i, ..., a_n) + f \circ (a_1, ..., b_i, ..., a_n)$$
$$f \circ (a_1, ..., pa_i, ..., a_n) = pf \circ (a_1, ..., a_i, ..., a_n)$$

$$f \circ (a_1, ..., a_i + b_i, ..., a_n) = f \circ (a_1, ..., a_i, ..., a_n) + f \circ (a_1, ..., b_i, ..., a_n)$$
$$f \circ (a_1, ..., pa_i, ..., a_n) = pf \circ (a_1, ..., a_i, ..., a_n)$$

$$1 \le i \le n \quad a_i, b_i \in A_i \quad p \in D$$

Доказательство. Теорема является следствием определений 3.1.1, 4.2.1 и теоремы 4.2.2. □

Теорема 4.3.3. *Пусть $D$ - коммутативное кольцо. Пусть $A_1$, ..., $A_n$, $S$ - $D$-модули. Отображение*

(4.3.1) $$f + g: A_1 \times ... \times A_n \to S \quad f, g \in \mathcal{L}(D; A_1 \times ... \times A_n \to S)$$

*определённое равенством*

(4.3.2) $$(f + g) \circ (a_1, ..., a_n) = f \circ (a_1, ..., a_n) + g \circ (a_1, ..., a_n)$$

*называется* **суммой полилинейных отображений** *$f$ и $g$ и является полилинейным отображением. Множество $\mathcal{L}(D; A_1 \times ... \times A_n \to S)$ является абелевой группой относительно суммы отображений.*

Доказательство. Согласно теореме 4.3.2

$$f\circ(a_1,...,a_i+b_i,...,a_n)=f\circ(a_1,...,a_i,...,a_n)+f\circ(a_1,...,b_i,...,a_n)\tag{4.3.3}$$

$$f\circ(a_1,...,pa_i,...,a_n)=pf\circ(a_1,...,a_i,...,a_n)\tag{4.3.4}$$

$$g\circ(a_1,...,a_i+b_i,...,a_n)=g\circ(a_1,...,a_i,...,a_n)+g\circ(a_1,...,b_i,...,a_n)\tag{4.3.5}$$

$$g\circ(a_1,...,pa_i,...,a_n)=pg\circ(a_1,...,a_i,...,a_n)\tag{4.3.6}$$

Равенство

$$\begin{aligned}&(f+g)\circ(x_1,...,x_i+y_i,...,x_n)\\&=f\circ(x_1,...,x_i+y_i,...,x_n)+g\circ(x_1,...,x_i+y_i,...,x_n)\\&=f\circ(x_1,...,x_i,...,x_n)+f\circ(x_1,...,y_i,...,x_n)\\&+g\circ(x_1,...,x_i,...,x_n)+g\circ(x_1,...,y_i,...,x_n)\\&=(f+g)\circ(x_1,...,x_i,...,x_n)+(f+g)\circ(x_1,...,y_i,...,x_n)\end{aligned}\tag{4.3.7}$$

является следствием равенств (4.3.2), (4.3.3), (4.3.5). Равенство

$$\begin{aligned}&(f+g)\circ(x_1,...,px_i,...,x_n)\\&=f\circ(x_1,...,px_i,...,x_n)+g\circ(x_1,...,px_i,...,x_n)\\&=pf\circ(x_1,...,x_i,...,x_n)+pg\circ(x_1,...,x_i,...,x_n)\\&=p(f\circ(x_1,...,x_i,...,x_n)+g\circ(x_1,...,x_i,...,x_n))\\&=p(f+g)\circ(x_1,...,x_i,...,x_n)\end{aligned}\tag{4.3.8}$$

является следствием равенств (4.3.2), (4.3.4), (4.3.6). Из равенств (4.3.7), (4.3.8) и теоремы 4.3.2 следует, что отображение (4.3.1) является полилинейным отображением $D$-модулей.

Пусть $f, g, h\in\mathcal{L}(D;A_1\times...\times A_2\to S)$. Для любого $a=(a_1,...,a_n)$, $a_1\in A_1$, ..., $a_n\in A_n$,

$$\begin{aligned}(f+g)\circ a&=f\circ a+g\circ a=g\circ a+f\circ a\\&=(g+f)\circ a\\((f+g)+h)\circ a&=(f+g)\circ a+h\circ a=(f\circ a+g\circ a)+h\circ a\\&=f\circ a+(g\circ a+h\circ a)=f\circ a+(g+h)\circ a\\&=(f+(g+h))\circ a\end{aligned}$$

Следовательно, сумма полилинейных отображений коммутативна и ассоциативна.

Из равенства (4.3.2) следует, что отображение

$$0:v\in A_1\times...\times A_n\to 0\in S$$

является нулём операции сложения

$$(0+f)\circ(a_1,...,a_n)=0\circ(a_1,...,a_n)+f\circ(a_1,...,a_n)=f\circ(a_1,...,a_n)$$

Из равенства (4.3.2) следует, что отображение

$$-f:(a_1,...,a_n)\in A_1\times...\times A_n\to-(f\circ(a_1,...,a_n))\in S$$

является отображением, обратным отображению $f$

$$f+(-f)=0$$

так как

$$\begin{aligned}(f+(-f))\circ(a_1,...,a_n)&=f\circ(a_1,...,a_n)+(-f)\circ(a_1,...,a_n)\\&=f\circ(a_1,...,a_n)-f\circ(a_1,...,a_n)\\&=0=0\circ(a_1,...,a_n)\end{aligned}$$

Из равенства

$$\begin{aligned}(f+g)\circ(a_1,...,a_n)&=f\circ(a_1,...,a_n)+g\circ(a_1,...,a_n)\\&=g\circ(a_1,...,a_n)+f\circ(a_1,...,a_n)\\&=(g+f)\circ(a_1,...,a_n)\end{aligned}$$

следует, что сумма отображений коммутативно. Следовательно, множество $\mathcal{L}(D;A_1\times...\times A_n\to S)$ является абелевой группой. □

Теорема 4.3.4. *Пусть $D$ - коммутативное кольцо. Пусть $A_1$, ..., $A_n$, $S$ - $D$-модули. Отображение*

(4.3.9) $$d\,f:A_1\times...\times A_n\to S\quad d\in D\quad f\in\mathcal{L}(D;A_1\times...\times A_n\to S)$$

*определённое равенством*

(4.3.10) $$(d\,f)\circ(a_1,...,a_n)=d(f\circ(a_1,...,a_n))$$

*называется* **произведением отображения** $f$ **на скаляр** $d$ *и является полилинейным отображением. Представление*

(4.3.11) $$a:f\in\mathcal{L}(D;A_1\times...\times A_n\to S)\to af\in\mathcal{L}(D;A_1\times...\times A_n\to S)$$

*кольца $D$ в абелевой группе $\mathcal{L}(D;A_1\times...\times A_n\to S)$ порождает структуру $D$-модуля.*

Доказательство. Согласно теореме 4.3.2

(4.3.12) $$f\circ(a_1,...,a_i+b_i,...,a_n)=f\circ(a_1,...,a_i,...,a_n)+f\circ(a_1,...,b_i,...,a_n)$$

(4.3.13) $$f\circ(a_1,...,pa_i,...,a_n)=pf\circ(a_1,...,a_i,...,a_n)$$

Равенство

(4.3.14) $$\begin{aligned}&(pf)\circ(x_1,...,x_i+y_i,...,x_n)\\=&p\ f\circ(x_1,...,x_i+y_i,...,x_n)\\=&p\ (f\circ(x_1,...,x_i,...,x_n)+f\circ(x_1,...,y_i,...,x_n))\\=&p(f\circ(x_1,...,x_i,...,x_n))+p(f\circ(x_1,...,y_i,...,x_n))\\=&(pf)\circ(x_1,...,x_i,...,x_n)+(pf)\circ(x_1,...,y_i,...,x_n)\end{aligned}$$

является следствием равенств (4.3.10), (4.3.12). Равенство

(4.3.15) $$\begin{aligned}&(pf)\circ(x_1,...,qx_i,...,x_n)\\=&p(f\circ(x_1,...,qx_i,...,x_n))=pq(f\circ(x_1,...,x_i,...,x_n))\\=&qp(f\circ(x_1,...,x_n))=q(pf)\circ(x_1,...,x_n)\end{aligned}$$

является следствием равенств (4.3.10), (4.3.13). Из равенств (4.3.14), (4.3.15) и теоремы 4.3.2 следует, что отображение (4.3.9) является полилинейным отображением $D$-модулей.

Равенство

(4.3.16) $$(p+q)f=pf+qf$$

является следствием равенства

$$\begin{aligned}((p+q)f)\circ(x_1,...,x_n)=&(p+q)(f\circ(x_1,...,x_n))\\=&p(f\circ(x_1,...,x_n))+q(f\circ(x_1,...,x_n))\\=&(pf)\circ(x_1,...,x_n)+(qf)\circ(x_1,...,x_n)\end{aligned}$$

Равенство

$$p(qf)=(pq)f \tag{4.3.17}$$

является следствием равенства

$$\begin{aligned}(p(qf))\circ(x_1,...,x_n)=&p\ (qf)\circ(x_1,...,x_n)=p\ (q\ f\circ(x_1,...,x_n))\\=&(pq)\ f\circ(x_1,...,x_n)=((pq)f)\circ(x_1,...,x_n)\end{aligned}$$

Из равенств (4.3.16), (4.3.17), следует, что отображение (4.3.11) является представлением кольца $D$ в абелевой группе $\mathcal{L}(D;A_1\times...\times A_n\to S)$ . Так как указанное представление эффективно, то, согласно определению 4.1.2 и теореме 4.3.3, абелевая группа $\mathcal{L}(D;A_1\to A_2)$ является $D$-модулем. □

## 4.4. $D$-модуль $\mathcal{L}(D;A\to B)$

Теорема 4.4.1.

$$\mathcal{L}(D;A^p\to\mathcal{L}(D;A^q\to B))=\mathcal{L}(D;A^{p+q}\to B) \tag{4.4.1}$$

Доказательство. □

Теорема 4.4.2. *Пусть*

$$\overline{\overline{e}}=\{e_i:i\in I\}$$

*базис $D$-модуля $A$. Множество*

$$\overline{\overline{h}}=\{h^i\in\mathcal{L}(D;A\to D):i\in I,h^i\circ e_j=\delta^i_j\} \tag{4.4.2}$$

*является базисом $D$-модуля* $\mathcal{L}(D;A\to D)$.

Доказательство.

Лемма 4.4.3. *Отображения $h^i$ линейно независимы.*

Доказательство. Пусть существуют $D$-числа $c_i$ такие, что

$$c_ih^i=0$$

Тогда для любого $A$-числа $e_j$

$$0=c_ih^i\circ e_j=c_i\delta^i_j=c_j$$

Лемма является следствием определения 4.1.8. ⊙

Лемма 4.4.4. *Отображение $f\in\mathcal{L}(D;A\to D)$, является линейной композицией отображений $h^i$.*

Доказательство. Для любого $A$-числа $a$,

$$a=a^ie_i \tag{4.4.3}$$

равенство

$$h^i\circ a=h^i\circ(a^je_j)=a^j(h^i\circ e_j)=a^j\delta^i_j=a^i \tag{4.4.4}$$

является следствием равенств (4.4.2), (4.4.3) и теоремы 4.2.2. Равенство

$$f\circ a=f\circ(a^ie_i)=a^i(f\circ e_i)=(f\circ e_i)(h^i\circ a) \tag{4.4.5}$$

является следствием равенства (4.4.4). Равенство

$$f=(f\circ e_{\boldsymbol{i}})h^{\boldsymbol{i}}$$

является следствием равенств (4.2.4), (4.2.12), (4.4.5). ⊙

Теорема является следствием лемм 4.4.3, 4.4.4 и определения 4.1.10. □

Теорема 4.4.5. *Пусть $D$ - коммутативное кольцо. Пусть*

$$\overline{\overline{e}}_i=\{e_{i\cdot\boldsymbol{i}}:\boldsymbol{i}\in\boldsymbol{I}_i\}$$

*базис $D$-модуля $A_i$, $i=1,\ ...,\ n$. Пусть*

$$\overline{\overline{e}}_B=\{e_{B\cdot\boldsymbol{i}}:\boldsymbol{i}\in\boldsymbol{I}\}$$

*базис $D$-модуля $B$. Множество*

$$\begin{aligned}\overline{\overline{h}}=\{\,&h^{\boldsymbol{i}_1...\boldsymbol{i}_n}_{\boldsymbol{i}}\in\mathcal{L}(D;A_1\times...\times A_n\to B):\boldsymbol{i}\in\boldsymbol{I},\boldsymbol{i}_i\in\boldsymbol{I}_i,i=1,...,n,\\&h^{\boldsymbol{i}_1...\boldsymbol{i}_n}_{\boldsymbol{i}}\circ(e_{1\cdot\boldsymbol{j}_1},...,e_{n\cdot\boldsymbol{j}_n})=\delta^{\boldsymbol{i}_1}_{\boldsymbol{j}_1}...\delta^{\boldsymbol{i}_n}_{\boldsymbol{j}_n}e_{B\cdot\boldsymbol{i}}\}\end{aligned}\tag{4.4.6}$$

*является базисом $D$-модуля $\mathcal{L}(D;A_1\times...\times A_n\to B)$.*

Доказательство.

Лемма 4.4.6. *Отображения $h^{\boldsymbol{i}_1...\boldsymbol{i}_n}_{\boldsymbol{i}}$ линейно независимы.*

Доказательство. Пусть существуют $D$-числа $c^{\boldsymbol{i}}_{\boldsymbol{i}_1...\boldsymbol{i}_n}$ такие, что

$$c^{\boldsymbol{i}}_{\boldsymbol{i}_1...\boldsymbol{i}_n}h^{\boldsymbol{i}_1...\boldsymbol{i}_n}_{\boldsymbol{i}}=0$$

Тогда для любого набора индексов $\boldsymbol{j}_1,\ ...,\ \boldsymbol{j}_n$

$$0=c^{\boldsymbol{i}}_{\boldsymbol{i}_1...\boldsymbol{i}_n}h^{\boldsymbol{i}_1...\boldsymbol{i}_n}_{\boldsymbol{i}}\circ(e_{1\cdot\boldsymbol{j}_1},...,e_{n\cdot\boldsymbol{j}_n})=c^{\boldsymbol{i}}_{\boldsymbol{i}_1...\boldsymbol{i}_n}\delta^{\boldsymbol{i}_1}_{\boldsymbol{j}_1}...\delta^{\boldsymbol{i}_n}_{\boldsymbol{j}_n}e_{B\cdot\boldsymbol{i}}=c^{\boldsymbol{i}}_{\boldsymbol{j}_1...\boldsymbol{j}_n}e_{B\cdot\boldsymbol{i}}$$

Следовательно, $c^{\boldsymbol{i}}_{\boldsymbol{j}_1...\boldsymbol{j}_n}=0$. Лемма является следствием определения 4.1.8. ⊙

Лемма 4.4.7. *Отображение $f\in\mathcal{L}(D;A_1\times...\times A_n\to B)$ является линейной композицией отображений $h^{\boldsymbol{i}_1...\boldsymbol{i}_n}_{\boldsymbol{i}}$.*

Доказательство. Для любого $A_1$-числа $a_1$

$$a_1=a_1^{\boldsymbol{i}_1}e_{1\cdot\boldsymbol{i}_1}\tag{4.4.7}$$

..., для любого $A_n$-числа $a_n$

$$a_n=a_n^{\boldsymbol{i}_n}e_{n\cdot\boldsymbol{i}_n}\tag{4.4.8}$$

равенство

$$\begin{aligned}h^{\boldsymbol{i}_1...\boldsymbol{i}_n}_{\boldsymbol{i}}\circ(a_1,...,a_n)&=h^{\boldsymbol{i}_1...\boldsymbol{i}_n}_{\boldsymbol{i}}\circ(a_1^{\boldsymbol{j}_1}e_{1\cdot\boldsymbol{j}_1},...,a_n^{\boldsymbol{j}_n}e_{n\cdot\boldsymbol{j}_n})\\&=a_1^{\boldsymbol{j}_1}...a_n^{\boldsymbol{j}_n}(h^{\boldsymbol{i}_1...\boldsymbol{i}_n}_{\boldsymbol{i}}\circ(e_{1\cdot\boldsymbol{j}_1},...,e_{n\cdot\boldsymbol{j}_n}))\\&=a_1^{\boldsymbol{j}_1}...a_n^{\boldsymbol{j}_n}\delta^{\boldsymbol{i}_1}_{\boldsymbol{j}_1}...\delta^{\boldsymbol{i}_n}_{\boldsymbol{j}_n}e_{B\cdot\boldsymbol{i}}\\&=a_1^{\boldsymbol{i}_1}...a_n^{\boldsymbol{i}_n}e_{B\cdot\boldsymbol{i}}\end{aligned}\tag{4.4.9}$$

является следствием равенств (4.4.6), (4.4.7), (4.4.8) и теоремы 4.2.2. Равенство

$$\begin{aligned}f\circ(a_1,...,a_n)&=f\circ(a_1^{\boldsymbol{j}_1}e_{1\cdot\boldsymbol{j}_1},...,a_n^{\boldsymbol{j}_n}e_{n\cdot\boldsymbol{j}_n})\\&=a_1^{\boldsymbol{j}_1}...a_n^{\boldsymbol{j}_n}f\circ(e_{1\cdot\boldsymbol{j}_1}...e_{n\cdot\boldsymbol{j}_n})\end{aligned}\tag{4.4.10}$$

является следствием равенств (4.4.7), (4.4.8). Так как

$$f\circ(e_{1\cdot\boldsymbol{j}_1}...e_{n\cdot\boldsymbol{j}_n})\in B$$

то

(4.4.11) $$f\circ(e_{1\cdot j_1}...e_{n\cdot j_n})=f^{i}_{j_1...j_n}e_i$$

Равенство

(4.4.12) $$\begin{aligned}f\circ(a_1,...,a_n)&=a_1^{i_1}...a_n^{i_n}f^{i}_{i_1...i_n}e_i\\&=f^{i}_{i_1...i_n}h_i^{i_1...i_n}\circ(a_1,...,a_n)\end{aligned}$$

является следствием равенств (4.4.9), (4.4.10), (4.4.11). Равенство

$$f=f^{i}_{i_1...i_n}h_i^{i_1...i_n}$$

является следствием равенств (4.3.2), (4.3.10), (4.4.12). ⊙

Теорема является следствием лемм 4.4.6, 4.4.7 и определения 4.1.10. □

Теорема 4.4.8. *Пусть $A_1$, ..., $A_n$, $B$ - свободные модули над коммутативным кольцом $D$. $D$-модуль $\mathcal{L}(D;A_1\times...\times A_n\to B)$ является свободным $D$-модулем.*

Доказательство. Теорема является следствием теоремы 4.4.5. □

## 4.5. Тензорное произведение $D$-модулей

Теорема 4.5.1. *Коммутативное кольцо $D$ является абелевой мультипликативной $\Omega$-группой.*

Доказательство. Мы полагаем, что произведение $\circ$ в кольце $D$ определён согласно правилу

$$a\circ b=ab$$

Так как произведение в кольце дистрибутивно относительно сложения, теорема является следствием определений 3.1.5, 3.1.8. □

Теорема 4.5.2. **Тензорное произведение** $A_1\otimes...\otimes A_n$ *$D$-модулей $A_1$, ..., $A_n$ существует.*

Доказательство. Теорема является следствием определения 4.1.2 и теорем 3.3.5, 4.5.1. □

Теорема 4.5.3. *Пусть $D$ - коммутативное кольцо. Пусть $A_1$, ..., $A_n$ - $D$-модули. Тензорное произведение дистрибутивно относительно сложения*

(4.5.1) $$\begin{aligned}&a_1\otimes...\otimes(a_i+b_i)\otimes...\otimes a_n\\=&a_1\otimes...\otimes a_i\otimes...\otimes a_n+a_1\otimes...\otimes b_i\otimes...\otimes a_n\\&a_i,b_i\in A_i\end{aligned}$$

*Представление кольца $D$ в тензорном произведении определено равенством*

(4.5.2) $$\begin{aligned}&a_1\otimes...\otimes(ca_i)\otimes...\otimes a_n=c(a_1\otimes...\otimes a_i\otimes...\otimes a_n)\\&a_i\in A_i\quad c\in D\end{aligned}$$

Доказательство. Равенство (4.5.1) является следствием равенства (3.3.23). Равенство (4.5.2) является следствием равенства (3.3.24). □

Теорема 4.5.4. *Пусть $A_1$, ..., $A_n$ - модули над коммутативным кольцом $D$. Пусть*

$$f : A_1 \times ... \times A_n \to A_1 \otimes ... \otimes A_n$$

*полилинейное отображение, определённое равенством*

$$f \circ (d_1, ..., d_n) = d_1 \otimes ... \otimes d_n \tag{4.5.3}$$

*Пусть*

$$g : A_1 \times ... \times A_n \to V$$

*полилинейное отображение в $D$-модуль $V$. Существует линейное отображение*

$$h : A_1 \otimes ... \otimes A_n \to V$$

*такое, что диаграмма*

(4.5.4)

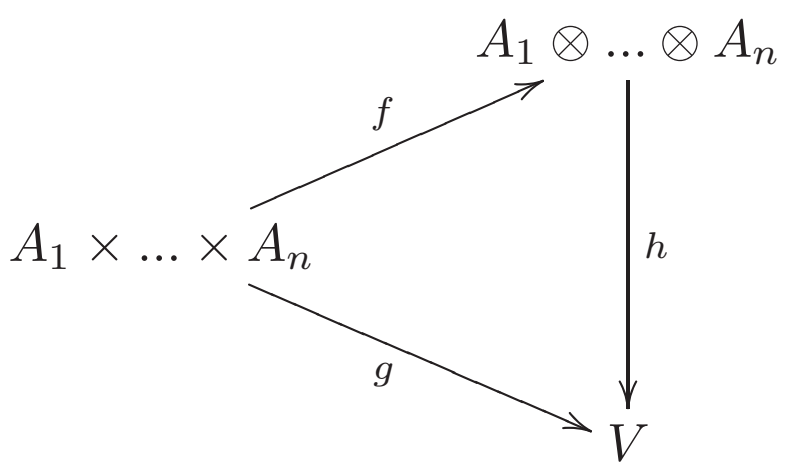

*коммутативна. Отображение $h$ определено равенством*

$$h(a_1 \otimes ... \otimes a_n) = g(a_1, ..., a_n) \tag{4.5.5}$$

Доказательство. Теорема является следствием теоремы 3.3.10 и определений 4.2.1, 4.3.1. □

Теорема 4.5.5. *Отображение*

$$(v_1, ..., v_n) \in V_1 \times ... \times V_n \to v_1 \otimes ... \otimes v_n \in V_1 \otimes ... \otimes V_n$$

*является полилинейным отображением.*

Доказательство. Теорема является следствием теоремы 3.3.9 и определения 4.3.1. □

Теорема 4.5.6. *Тензорное произведение $A_1 \otimes ... \otimes A_n$ свободных конечномерных модулей $A_1$, ..., $A_n$ над коммутативным кольцом $D$ является свободным конечномерным модулем.*

*Пусть $\overline{\overline{e}}_i$ - базис модуля $A_i$ над кольцом $D$. Произвольный тензор $a \in A_1 \otimes ... \otimes A_n$ можно представить в виде*

$$a = a^{i_1...i_n} e_{1\cdot i_1} \otimes ... \otimes e_{n\cdot i_n} \tag{4.5.6}$$

*Мы будем называть выражение $a^{i_1...i_n}$* **стандартной компонентой тензора**.

Доказательство. Вектор $a_i \in A_i$ имеет разложение

$$a_i = a_i^k \overline{e}_{i\cdot k}$$

относительно базиса $\overline{\overline{e}}_i$. Из равенств (4.5.1), (4.5.2) следует

$$a_1 \otimes ... \otimes a_n = a_1^{i_1} ... a_n^{i_n} e_{1\cdot i_1} \otimes ... \otimes e_{n\cdot i_n}$$

Так как множество тензоров $a_1\otimes...\otimes a_n$ является множеством образующих модуля $A_1\otimes...\otimes A_n$, то тензор $a\in A_1\otimes...\otimes A_n$ можно записать в виде

$$a=a^s a_{s\cdot}{}^{i_1}_{1}...a_{s\cdot}{}^{i_n}_{n} e_{1\cdot i_1}\otimes...\otimes e_{n\cdot i_n} \tag{4.5.7}$$

где $a^s$, $a_{s\cdot}{}^{i_1}_{1}$, ..., $a_{s\cdot}{}^{i_n}_{n}\in F$. Положим

$$a^s a_{s\cdot}{}^{i_1}_{1}...a_{s\cdot}{}^{i_n}_{n}=a^{i_1...i_n}$$

Тогда равенство (4.5.7) примет вид (4.5.6).

Следовательно, множество тензоров $e_{1\cdot i_1}\otimes...\otimes e_{n\cdot i_n}$ является множеством образующих модуля $A_1\otimes...\otimes A_n$. Так как размерность модуля $A_i$, $i=1$, ..., $n$, конечна, то конечно множество тензоров $e_{1\cdot i_1}\otimes...\otimes e_{n\cdot i_n}$. Следовательно, множество тензоров $e_{1\cdot i_1}\otimes...\otimes e_{n\cdot i_n}$ содержит базис модуля $A_1\otimes...\otimes A_n$, и модуль $A_1\otimes...\otimes A_n$ является свободным модулем над кольцом $D$. □

Глава 5

# $D$-алгебра

## 5.1. Алгебра над коммутативным кольцом

Определение 5.1.1. *Пусть $D$ - коммутативное кольцо. $D$-модуль $A$ называется* **алгеброй над кольцом** $D$ *или* $D$**-алгеброй**, *если определена операция произведения*[5.1] *в* $A$

$$vw = C \circ (v, w) \tag{5.1.1}$$

*где $C$ - билинейное отображение*

$$C : A \times A \to A$$

*Если $A$ является свободным $D$-модулем, то $A$ называется* **свободной алгеброй над кольцом** $D$. □

Теорема 5.1.2. *Произведение в алгебре $A$ дистрибутивно по отношению к сложению.*

Доказательство. Утверждение теоремы следует из цепочки равенств

$$(a+b)c = f \circ (a+b, c) = f \circ (a, c) + f \circ (b, c) = ac + bc$$
$$a(b+c) = f \circ (a, b+c) = f \circ (a, b) + f \circ (a, c) = ab + ac$$

□

Произведение в алгебре может быть ни коммутативным, ни ассоциативным. Следующие определения основаны на определениях, данным в [10], с. 13.

Определение 5.1.3. **Коммутатор**

$$[a, b] = ab - ba$$

*служит мерой коммутативности в $D$-алгебре $A$. $D$-алгебра $A$ называется* **коммутативной**, *если*

$$[a, b] = 0$$

□

Определение 5.1.4. **Ассоциатор**

$$(a, b, c) = (ab)c - a(bc) \tag{5.1.2}$$

*служит мерой ассоциативности в $D$-алгебре $A$. $D$-алгебра $A$ называется* **ассоциативной**, *если*

$$(a, b, c) = 0$$

□

[5.1] Я следую определению, приведенному в [10], с. 1, [8], с. 4. Утверждение, верное для произвольного $D$-модуля, верно также для $D$-алгебры.

Теорема 5.1.5. *Пусть $A$ - алгебра над коммутативным кольцом $D$.*[5.2]

$$a(b,c,d)+(a,b,c)d=(ab,c,d)-(a,bc,d)+(a,b,cd) \tag{5.1.3}$$

*для любых $a$, $b$, $c$, $d\in A$.*

Доказательство. Равенство (5.1.3) следует из цепочки равенств

$$\begin{aligned}a(b,c,d)+(a,b,c)d&=a((bc)d-b(cd))+((ab)c-a(bc))d\\&=a((bc)d)-a(b(cd))+((ab)c)d-(a(bc))d\\&=((ab)c)d-(ab)(cd)+(ab)(cd)\\&+a((bc)d)-a(b(cd))-(a(bc))d\\&=(ab,c,d)-(a(bc))d+a((bc)d)+(ab)(cd)-a(b(cd))\\&=(ab,c,d)-(a,(bc),d)+(a,b,cd)\end{aligned}$$

□

Определение 5.1.6. **Ядро $D$-алгебры** $A$ *- это множество*[5.3]

$$N(A)=\{a\in A:\forall b,c\in A,(a,b,c)=(b,a,c)=(b,c,a)=0\}$$

□

Определение 5.1.7. **Центр $D$-алгебры** $A$ *- это множество*[5.4]

$$Z(A)=\{a\in A:a\in N(A),\forall b\in A,ab=ba\}$$

□

Теорема 5.1.8. *Пусть $D$ - коммутативное кольцо. Если $D$-алгебра $A$ имеет единицу, то существует изоморфизм $f$ кольца $D$ в центр алгебры $A$.*

Доказательство. Пусть $e\in A$ - единица алгебры $A$. Положим $f\circ a=ae$. □

Пусть $\overline{\overline{e}}$ - базис свободной алгебры $A$ над кольцом $D$. Если алгебра $A$ имеет единицу, положим $e_0$ - единица алгебры $A$.

Теорема 5.1.9. *Пусть $\overline{\overline{e}}$ - базис свободной алгебры $A$ над кольцом $D$. Пусть*

$$a=a^ie_i\quad b=b^ie_i\quad a,b\in A$$

*Произведение $a$, $b$ можно получить согласно правилу*

$$(ab)^k=C^k_{ij}a^ib^j \tag{5.1.4}$$

*где $C^k_{ij}$* - **структурные константы** *алгебры $A$ над кольцом $D$. Произведение базисных векторов в алгебре $A$ определено согласно правилу*

$$e_ie_j=C^k_{ij}e_k \tag{5.1.5}$$

---

[5.2] Утверждение теоремы опирается на равенство [10]-(2.4).

[5.3] Определение дано на базе аналогичного определения в [10], с. 13

[5.4] Определение дано на базе аналогичного определения в [10], с. 14

Доказательство. Равенство (5.1.5) является следствием утверждения, что $\overline{\overline{e}}$ является базисом алгебры $A$. Так как произведение в алгебре является билинейным отображением, то произведение $a$ и $b$ можно записать в виде

$$ab = a^i b^j e_i e_j \tag{5.1.6}$$

Из равенств (5.1.5), (5.1.6), следует

$$ab = a^i b^j C_{ij}^k e_k \tag{5.1.7}$$

Так как $\overline{\overline{e}}$ является базисом алгебры $A$, то равенство (5.1.4) следует из равенства (5.1.7). □

Теорема 5.1.10. *Если алгебра $A$ коммутативна, то*

$$C_{ij}^p = C_{ji}^p \tag{5.1.8}$$

*Если алгебра $A$ ассоциативна, то*

$$C_{ij}^p C_{pk}^q = C_{ip}^q C_{jk}^p \tag{5.1.9}$$

Доказательство. Для коммутативной алгебры, равенство (5.1.8) следует из равенства

$$e_i e_j = e_j e_i$$

Для ассоциативной алгебры, равенство (5.1.9) следует из равенства

$$(e_i e_j) e_k = e_i (e_j e_k)$$

□

Теорема 5.1.11. *Представление*

$$f_{2,3} : A \twoheadrightarrow A \tag{5.1.10}$$

*$D$-модуля $A$ в $D$-модуле $A$ эквивалентно структуре $D$-алгебры $A$.*

Доказательство.

- Пусть в $D$-модуле $A$ определена структура $D$-алгебры $A$, порождённая произведением

  $$v\,w = C \circ (v, w)$$

  Согласно определениям 5.1.1 и 4.3.1, **левый сдвиг $D$-модуля** $A$, определённый равенством

  $$l \circ v : w \in A \to v\,w \in A \tag{5.1.11}$$

  является линейным отображением. Согласно определению 4.2.1, отображение $l \circ v$ является эндоморфизмом $D$-модуля $A$.

  Равенство

  $$(l \circ (v_1 + v_2)) \circ w = (v_1 + v_2)w = v_1 w + v_2 w = (l \circ v_1) \circ w + (l \circ v_2) \circ w \tag{5.1.12}$$

  является следствием определения 4.3.1 и равенства (5.1.11). Согласно теореме 4.2.3, равенство

  $$l \circ (v_1 + v_2) = l \circ v_1 + l \circ v_2 \tag{5.1.13}$$

  является следствием равенства (5.1.12). Равенство

  $$(l \circ (dv)) \circ w = (dv)w = d(vw) = d((l \circ v) \circ w) \tag{5.1.14}$$

является следствием определения 4.3.1 и равенства (5.1.11). Согласно теореме 4.2.3, равенство

$$l\circ(dv)=d(l\circ v) \tag{5.1.15}$$

является следствием равенства (5.1.14). Из равенств (5.1.13), (5.1.15) следует, что отображение

$$f_{2,3}:v\to l\circ v$$

является представлением $D$-модуля $A$ в $D$-модуле $A$

$$f_{2,3}:A \xrightarrow{*} A \qquad f_{2,3}\circ v:w\to(l\circ v)\circ w \tag{5.1.16}$$

- Рассмотрим представление (5.1.10) $D$-модуля $A$ в $D$-модуле $A$. Поскольку отображение $f_{2,3}\circ v$ является эндоморфизмом $D$-модуля $A$, то

$$\begin{aligned}(f_{2,3}\circ v)(w_1+w_2)&=(f_{2,3}\circ v)\circ w_1+(f_{2,3}\circ v)\circ w_2\\ (f_{2,3}\circ v)\circ(dw)&=d((f_{2,3}\circ v)\circ w)\end{aligned} \tag{5.1.17}$$

Поскольку отображение (5.1.10) является линейным отображением

$$f_{2,3}:A\to\mathcal{L}(D;A;A)$$

то, согласно теоремам 4.2.3, 4.2.4,

$$(f_{2,3}\circ(v_1+v_2))\circ w=(f_{2,3}\circ v_1+f_{2,3}\circ v_2)(w)=(f_{2,3}\circ v_1)\circ w+(f_{2,3}\circ v_2)\circ w \tag{5.1.18}$$

$$(f_{2,3}\circ(d\,v))\circ w=(d\,(f_{2,3}\circ v))\circ w=d\,((f_{2,3}\circ v)\circ w) \tag{5.1.19}$$

Из равенств (5.1.17), (5.1.18), (5.1.19) и определения 4.3.1, следует, что отображение $f_{2,3}$ является билинейным отображением. Следовательно, отображение $f_{2,3}$ определяет произведение в $D$-модуле $A$ согласно правилу

$$vw=(f_{2,3}\circ v)\circ w$$

□

Следствие 5.1.12. *$D$ - коммутативное кольцо, $A$ - абелевая группа. Диаграмма представлений*

$$\begin{array}{ccccc} D & \xrightarrow[*]{g_{1,2}} & A & \xrightarrow[*]{g_{2,3}} & A \\ & & & & \uparrow{\scriptstyle *\,g_{1,2}} \\ & & & & D \end{array} \qquad \begin{array}{l} g_{1,2}(d):v\to d\,v \\ g_{2,3}(v):w\to C\circ(v,w) \\ \qquad C\in\mathcal{L}(D;A^2\to A)\end{array} \tag{5.1.20}$$

*порождает структуру $D$-алгебры $A$.* □

## 5.2. Линейный гомоморфизм

Теорема 5.2.1. *Пусть*

$$\begin{array}{ccccc} D_1 & \xrightarrow[*]{g_{1\cdot 1,2}} & A_1 & \xrightarrow[*]{g_{1\cdot 2,3}} & A_1 \\ & & & & \uparrow{\scriptstyle *\, g_{1\cdot 1,2}} \\ & & & & D_1 \end{array} \qquad \begin{array}{l} g_{1\cdot 1,2}(d) : v \rightarrow d\,v \\ g_{1\cdot 2,3}(v) : w \rightarrow C_1 \circ (v,w) \\ C_1 \in \mathcal{L}(D_1; A_1^2 \rightarrow A_1) \end{array} \tag{5.2.1}$$

*диаграмма представлений, описывающая $D_1$-алгебру $A_1$. Пусть*

$$\begin{array}{ccccc} D_2 & \xrightarrow[*]{g_{2\cdot 1,2}} & A_2 & \xrightarrow[*]{g_{2\cdot 2,3}} & A_2 \\ & & & & \uparrow{\scriptstyle *\, g_{2\cdot 1,2}} \\ & & & & D_2 \end{array} \qquad \begin{array}{l} g_{2\cdot 1,2}(d) : v \rightarrow d\,v \\ g_{2\cdot 2,3}(v) : w \rightarrow C_2 \circ (v,w) \\ C_2 \in \mathcal{L}(D_2; A_2^2 \rightarrow A_2) \end{array} \tag{5.2.2}$$

*диаграмма представлений, описывающая $D_2$-алгебру $A_2$. Морфизм $D_1$-алгебры $A_1$ в $D_2$-алгебру $A_2$ - это кортеж отображений*

$$r_1 : D_1 \rightarrow D_2 \quad r_2 : A_1 \rightarrow A_2$$

*где отображение $r_1$ является гомоморфизмом кольца $D_1$ в кольцо $D_2$ и отображение $r_2$ является линейным отображением $D_1$-алгебры $A_1$ в $D_2$-алгебру $A_2$ таким, что*

$$r_2(ab) = r_2(a)r_2(b) \tag{5.2.3}$$

Доказательство. Согласно равенствам [4]-(4.2.3), морфизм $(r_1, r_2)$ представления $f_{1,2}$ удовлетворяет равенству

$$\begin{array}{c} r_2(f_{1\cdot 1,2}(d)(a)) = f_{2\cdot 1,2}(r_1(d))(r_2(a)) \\ r_2(d\,a) = r_1(d)r_2(a) \end{array} \tag{5.2.4}$$

Следовательно, отображение $(r_1, r_2)$ является линейным отображением.

Согласно равенству [4]-(4.2.3), морфизм $(r_2, r_2)$ представления $f_{2,3}$ удовлетворяет равенству[5.5]

$$r_2(f_{1\cdot 2,3}(a_2)(a_3)) = f_{2\cdot 2,3}(r_2(a_2))(r_2(a_3)) \tag{5.2.5}$$

Из равенств (5.2.5), (5.2.1), (5.2.2), следует

$$r_2(C_1(v,w)) = C_2(r_2(v), r_2(w)) \tag{5.2.6}$$

Равенство (5.2.3) следует из равенств (5.2.6), (5.1.1). □

Определение 5.2.2. *Морфизм представлений $D_1$-алгебры $A_1$ в $D_2$-алгебру $A_2$ называется* **линейным гомоморфизмом** *$D_1$-алгебры $A_1$ в $D_2$-алгебру $A_2$.*
□

[5.5] Так как в диаграммах представлений (5.2.1), (5.2.2), носители $\Omega_2$-алгебры и $\Omega_3$-алгебры совпадают, то также совпадают морфизмы представлений на уровнях 2 и 3.

Теорема 5.2.3. *Пусть $\overline{\overline{e}}_1$ - базис $D_1$-алгебры $A_1$. Пусть $\overline{\overline{e}}_2$ - базис $D_2$-алгебры $A_2$. Тогда линейный гомоморфизм*[5.6] $(r_1, r_2)$ *$D_1$-алгебры $A_1$ в $D_2$-алгебру $A_2$ имеет представление*[5.7]

$$b = e_{2*}{}^{*}r_{2*}{}^{*}r_1(a) = e_{2\cdot i}r_{2\cdot}{}^{i}_{j}r_1(a^j) \tag{5.2.7}$$

$$b = r_{2*}{}^{*}r_1(a) \tag{5.2.8}$$

*относительно заданных базисов. Здесь*

- *$a$ - координатная матрица вектора $\overline{a}$ относительно базиса $\overline{\overline{e}}_1$.*
- *$b$ - координатная матрица вектора*
$$b = r_2(a)$$
*относительно базиса $\overline{\overline{e}}_2$.*
- *$r_2$ - координатная матрица множества векторов $(r_2(e_{1\cdot i}))$ относительно базиса $\overline{\overline{e}}_2$. Мы будем называть матрицу $r_2$* **матрицей линейного гомоморфизма** *относительно базисов $\overline{\overline{e}}_1$ и $\overline{\overline{e}}_2$.*

Доказательство. Вектор $\overline{a} \in A_1$ имеет разложение

$$\overline{a} = e_{1*}{}^{*}a$$

относительно базиса $\overline{\overline{e}}_1$. Вектор $\overline{b} \in A_2$ имеет разложение

$$\overline{b} = e_2{}^{*}{}_{*}b \tag{5.2.9}$$

относительно базиса $\overline{\overline{e}}_2$.

Так как $(r_1, r_2)$ - линейный гомоморфизм, то на основании (5.2.4) следует, что

$$b = r_2(a) = r_2(e_1{}^{*}{}_{*}a) = r_2(e_1)_{*}{}^{*}r_1(a) \tag{5.2.10}$$

где

$$r_1(a) = \begin{pmatrix} r_1(a^1) \\ \dots \\ r_1(a^n) \end{pmatrix}$$

$r_2(e_{1\cdot i})$ также вектор $D$-модуля $A_2$ и имеет разложение

$$r_2(e_{1\cdot i}) = e_{2*}{}^{*}r_{2\cdot i} = e_{2\cdot j}\, r_{2\cdot}{}^{j}_{i} \tag{5.2.11}$$

относительно базиса $\overline{\overline{e}}_2$. Комбинируя (5.2.10) и (5.2.11), мы получаем (5.2.7). (5.2.8) следует из сравнения (5.2.9) и (5.2.7) и теоремы [3]-5.3.3. □

Теорема 5.2.4. *Пусть $\overline{\overline{e}}_1$ - базис $D_1\star$-алгебры $A_1$. Пусть $\overline{\overline{e}}_2$ - базис $D_2\star$-алгебры $A_2$. Если отображение $r_1$ является инъекцией, то матрица линейного гомоморфизма и структурные константы связаны соотношением*

$$r_{2\cdot}{}^{l}_{k}r_1(C_{1\cdot}{}^{k}_{ij}) = C_{2\cdot}{}^{l}_{pq}r_{2\cdot}{}^{p}_{i}r_{2\cdot}{}^{q}_{j} \tag{5.2.12}$$

---

[5.6] Эта теорема аналогична теореме [3]-5.4.3.

[5.7] Произведение матриц над коммутативным кольцом определено только как ${}_*{}^*$-произведение. Однако я предпочитаю явно указывать операцию, так как в этом случае видно, что в выражении участвуют матрицы. Кроме того, я планирую рассмотреть подобную теорему в некоммутативном случае.

Доказательство. Пусть

$$\overline{a},\overline{b}\in A_1\quad \overline{a}=e_{1*}{}^*a\quad \overline{b}=e_{1*}{}^*b$$

Из равенств (5.1.4), (5.1.1), (5.2.1), следует

$$ab=e_{1\cdot k}C_{1\cdot}{}^{k}_{ij}a^ib^j \tag{5.2.13}$$

Из равенств (5.2.4), (5.2.13), следует

$$r_2(ab)=r_2(e_{1\cdot k})r_1(C_{1\cdot}{}^{k}_{ij}a^ib^j) \tag{5.2.14}$$

Поскольку отображение $r_1$ является гомоморфизмом колец, то из равенства (5.2.14), следует

$$r_2(ab)=r_2(e_{1\cdot k})r_1(C_{1\cdot}{}^{k}_{ij})r_1(a^i)r_1(b^j) \tag{5.2.15}$$

Из теоремы 5.2.3 и равенства (5.2.15), следует

$$r_2(ab)=e_{2\cdot l}\,r_{2\cdot}{}^{l}_{k}r_1(C_{1\cdot}{}^{k}_{ij})r_1(a^i)r_1(b^j) \tag{5.2.16}$$

Из равенства (5.2.3) и теоремы 5.2.3, следует

$$r_2(ab)=r_2(a)r_2(b)=e_{2\cdot p}\,r_1(a^i)r_{2\cdot}{}^{p}_{i}e_{2\cdot q}\,r_1(b^j)r_{2\cdot}{}^{q}_{j} \tag{5.2.17}$$

Из равенств (5.1.4), (5.1.1), (5.2.2), (5.2.17), следует

$$r_2(ab)=e_{2\cdot l}\,C_{2\cdot}{}^{l}_{pq}r_1(a^i)r_{2\cdot}{}^{p}_{i}r_1(b^j)r_{2\cdot}{}^{q}_{j} \tag{5.2.18}$$

Из равенств (5.2.16), (5.2.18), следует

$$e_{2\cdot l}\,r_{2\cdot}{}^{l}_{k}r_1(C_{1\cdot}{}^{k}_{ij})r_1(a^i)r_1(b^j)=e_{2\cdot l}\,C_{2\cdot}{}^{l}_{pq}r_1(a^i)r_{2\cdot}{}^{p}_{i}r_1(b^j)r_{2\cdot}{}^{q}_{j} \tag{5.2.19}$$

Равенство (5.2.12) следует из равенства (5.2.19), так как векторы базиса $\overline{\overline{e}}_2$ линейно независимы, и $a^i$, $b^i$ (а следовательно, $r_1(a^i)$, $r_1(b^i)$) - произвольные величины. □

## 5.3. Линейный автоморфизм алгебры кватернионов

Определение координат линейного автоморфизма - задача непростая. В этом разделе мы рассмотрим пример нетривиального линейного автоморфизма алгебры кватернионов.

Теорема 5.3.1. *Координаты линейного автоморфизма алгебры кватернионов удовлетворяют системе уравнений*

$$\left\{\begin{array}{lll}
r^1_1=r^2_2r^3_3-r^3_2r^2_3 & r^1_2=r^2_3r^3_1-r^3_3r^2_1 & r^1_3=r^2_1r^3_2-r^3_1r^2_2\\
r^2_1=r^3_2r^1_3-r^1_2r^3_3 & r^2_2=r^3_3r^1_1-r^1_3r^3_1 & r^2_3=r^3_1r^1_2-r^1_1r^3_2\\
r^3_1=r^1_2r^2_3-r^2_2r^1_3 & r^3_2=r^1_3r^2_1-r^2_3r^1_1 & r^3_3=r^1_1r^2_2-r^2_1r^1_2
\end{array}\right. \tag{5.3.1}$$

Доказательство. Согласно теоремам [5]-4.3.1, 5.2.4, линейный автоморфизм алгебры кватернионов удовлетворяет уравнениям

$$\begin{array}{llll}
r^l_0=r^p_0r^q_0C^l_{pq} & r^l_1=r^p_0r^q_1C^l_{pq} & r^l_2=r^p_0r^q_2C^l_{pq} & r^l_3=r^p_0r^q_3C^l_{pq}\\
r^l_1=r^p_1r^q_0C^l_{pq} & -r^l_0=r^p_1r^q_1C^l_{pq} & r^l_3=r^p_1r^q_2C^l_{pq} & -r^l_2=r^p_1r^q_3C^l_{pq}\\
r^l_2=r^p_2r^q_0C^l_{pq} & -r^l_3=r^p_2r^q_1C^l_{pq} & -r^l_0=r^p_2r^q_2C^l_{pq} & r^l_0=r^p_2r^q_3C^l_{pq}\\
r^l_3=r^p_3r^q_0C^l_{pq} & r^l_2=r^p_3r^q_1C^l_{pq} & -r^l_1=r^p_3r^q_2C^l_{pq} & -r^l_0=r^p_3r^q_3C^l_{pq}
\end{array} \tag{5.3.2}$$

Из равенства (5.3.2) следует

$$r_1^l=r_0^pr_1^qC_{pq}^l=r_2^pr_3^qC_{pq}^l \quad r_0^pr_1^qC_{pq}^l=r_0^pr_1^qC_{qp}^l \quad r_2^pr_3^qC_{pq}^l=-r_2^pr_3^qC_{qp}^l$$

(5.3.3)
$$r_2^l=r_0^pr_2^qC_{pq}^l=r_3^pr_1^qC_{pq}^l \quad r_0^pr_2^qC_{pq}^l=r_0^pr_2^qC_{qp}^l \quad r_3^pr_1^qC_{pq}^l=-r_1^pr_3^qC_{pq}^l$$

$$r_3^l=r_0^pr_3^qC_{pq}^l=r_1^pr_2^qC_{pq}^l \quad r_0^pr_3^qC_{pq}^l=r_0^pr_3^qC_{qp}^l \quad r_1^pr_2^qC_{pq}^l=-r_1^pr_2^qC_{qp}^l$$

(5.3.4)
$$r_0^l=r_0^pr_0^qC_{pq}^l=-r_1^pr_1^qC_{pq}^l=-r_2^pr_2^qC_{pq}^l=-r_3^pr_3^qC_{pq}^l$$

Если $l=0$, то из равенства

$$C_{pq}^0=C_{qp}^0$$

следует

(5.3.5)
$$r_i^pr_j^qC_{pq}^0=r_i^pr_j^qC_{qp}^0$$

Из равенства (5.3.3) для $l=0$ и равенства (5.3.5), следует

(5.3.6)
$$r_1^0=r_2^0=r_3^0=0$$

Если $l=1,2,3$, то равенство (5.3.3) можно записать в виде

(5.3.7)
$$\left\{\begin{array}{l}
r_i^l=r_0^lr_i^0C_{l0}^l+r_0^0r_i^lC_{0l}^l+r_0^ar_i^bC_{ab}^l+r_0^br_i^aC_{ba}^l\\
\quad r_0^lr_i^0C_{l0}^l+r_0^0r_i^lC_{0l}^l+r_0^ar_i^bC_{ab}^l+r_0^br_i^aC_{ba}^l\\
\quad =r_0^lr_i^0C_{0l}^l+r_0^0r_i^lC_{l0}^l+r_0^ar_i^bC_{ba}^l+r_0^br_i^aC_{ab}^l\\
\quad i=1,2,3\\
r_i^l=r_k^0r_j^lC_{0l}^l+r_k^lr_j^0C_{l0}^l+r_k^ar_j^bC_{ab}^l+r_k^br_j^aC_{ba}^l\\
\quad r_k^0r_j^lC_{0l}^l+r_k^lr_j^0C_{l0}^l+r_k^ar_j^bC_{ab}^l+r_k^br_j^aC_{ba}^l\\
\quad =-r_k^0r_j^lC_{l0}^l-r_k^lr_j^0C_{0l}^l-r_k^ar_j^bC_{ba}^l-r_k^br_j^aC_{ab}^l\\
\quad i=1 \quad k=2 \quad j=3\\
\quad i=2 \quad k=3 \quad j=1\\
\quad i=3 \quad k=1 \quad j=2\\
\quad 0<a<b \quad a\ne l \quad b\ne l
\end{array}\right.$$

Из равенств (5.3.7), (5.3.6) и равенств

(5.3.8)
$$\begin{array}{l}C_{0l}^l=C_{l0}^l=1\\ C_{ab}^l=-C_{ba}^l\end{array}$$

следует

$$
(5.3.9)\qquad \left\{\begin{array}{l}
r_i^l = r_0^0 r_i^l + r_0^a r_i^b C_{ab}^l - r_0^b r_i^a C_{ab}^l \\
\qquad r_0^0 r_i^l + r_0^a r_i^b C_{ab}^l - r_0^b r_i^a C_{ab}^l \\
\quad = r_0^0 r_i^l - r_0^a r_i^b C_{ab}^l + r_0^b r_i^a C_{ab}^l \\
\qquad i = 1, 2, 3 \\
r_i^l = r_k^a r_j^b C_{ab}^l - r_k^b r_j^a C_{ab}^l \\
\qquad r_k^a r_j^b C_{ab}^l - r_k^b r_j^a C_{ab}^l \\
\quad = r_k^a r_j^b C_{ab}^l - r_k^b r_j^a C_{ab}^l \\
\qquad\quad i = 1 \qquad k = 2 \quad j = 3 \\
\qquad\quad i = 2 \qquad k = 3 \quad j = 1 \\
\qquad\quad i = 3 \qquad k = 1 \quad j = 2 \\
\qquad 0 < a < b \quad a \ne l \quad b \ne l
\end{array}\right.
$$

Из равенств (5.3.9) следует

$$
(5.3.10)\qquad \left\{\begin{array}{l}
r_i^l = r_0^0 r_i^l \\
\qquad r_0^a r_i^b - r_0^b r_i^a = 0 \\
\qquad i = 1, 2, 3 \\
r_i^l = r_k^a r_j^b C_{ab}^l - r_k^b r_j^a C_{ab}^l \\
\qquad\quad i = 1 \qquad k = 2 \quad j = 3 \\
\qquad\quad i = 2 \qquad k = 3 \quad j = 1 \\
\qquad\quad i = 3 \qquad k = 1 \quad j = 2 \\
\qquad 0 < a < b \quad a \ne l \quad b \ne l
\end{array}\right.
$$

Из равенства (5.3.10) следует

$$
(5.3.11)\qquad r_0^0 = 1
$$

Из равенства (5.3.4) для $l = 0$ следует

$$
(5.3.12)\qquad \begin{array}{l}
r_0^0 = r_0^0 r_0^0 - r_0^1 r_0^1 - r_0^2 r_0^2 - r_0^3 r_0^3 \\
\quad = -r_i^0 r_i^0 + r_i^1 r_i^1 + r_i^2 r_i^2 + r_i^3 r_i^3 \\
i = 1, 2, 3
\end{array}
$$

Из равенств (5.3.6), (5.3.10), (5.3.12), следует

$$
(5.3.13)\qquad \begin{array}{l}
0 = r_0^1 r_0^1 + r_0^2 r_0^2 + r_0^3 r_0^3 \\
1 = r_i^1 r_i^1 + r_i^2 r_i^2 + r_i^3 r_i^3 \\
i = 1, 2, 3
\end{array}
$$

Из равенств (5.3.13) следует [5.8]

$$r^1_0 = r^2_0 = r^3_0 = 0 \tag{5.3.14}$$

Из равенства (5.3.4) для $l > 0$ следует

$$\begin{aligned} r^l_0 &= r^l_0 r^0_0 C^l_{l0} + r^0_0 r^l_0 C^l_{0l} + r^a_0 r^b_0 C^l_{ab} + r^b_0 r^a_0 C^l_{ba} \\ &= -r^l_i r^0_i C^l_{l0} - r^0_i r^l_i C^l_{0l} - r^a_i r^b_i C^l_{ab} - r^b_i r^a_i C^l_{ba} \\ & \quad i > 0 \\ & \quad l > 0 \quad 0 < a < b \quad a \ne l \quad b \ne l \end{aligned} \tag{5.3.15}$$

Равенства (5.3.15) тождественно верны в силу равенств (5.3.6), (5.3.14), (5.3.8). Из равенств (5.3.14), (5.3.10), следует

$$\left\{ \begin{array}{l} r^l_i = r^a_k r^b_j C^l_{ab} - r^b_k r^a_j C^l_{ab} \\ \qquad i = 1 \qquad k = 2 \qquad j = 3 \\ \qquad i = 2 \qquad k = 3 \qquad j = 1 \\ \qquad i = 3 \qquad k = 1 \qquad j = 2 \\ \qquad l > 0 \quad 0 < a < b \quad a \ne l \quad b \ne l \end{array} \right. \tag{5.3.16}$$

Равенства (5.3.1) следуют из равенств (5.3.16). □

Пример 5.3.2. *Очевидно, координаты*

$$r^i_j = \delta^i_j$$

*удовлетворяют уравнению* (5.3.1). *Мы можем убедиться непосредственной проверкой, что координаты отображения*

$$r^0_0 = 1 \quad r^1_2 = 1 \quad r^2_3 = 1 \quad r^3_1 = 1$$

*также удовлетворяют уравнению* (5.3.1). *Матрица координат этого отображения имеет вид*

$$r = \begin{pmatrix} 1 & 0 & 0 & \\ 0 & 0 & 1 & 0 \\ 0 & 0 & 0 & 1 \\ 0 & 1 & 0 & 0 \end{pmatrix}$$

[5.8] Мы здесь опираемся на то, что алгебра кватернионов определена над полем действительных чисел. Если рассматривать алгебру кватернионов над полем комплексных чисел, то уравнение (5.3.13) определяет конус в комплексном пространстве. Соответственно, у нас шире выбор координат линейного автоморфизма.

*Согласно теореме [5]-4.3.4, стандартные компоненты отображения* $r$ *имеют вид*

$$
\begin{array}{cccc}
r^{00} = \frac{1}{4} & r^{11} = -\frac{1}{4} & r^{22} = -\frac{1}{4} & r^{33} = -\frac{1}{4} \\
r^{10} = -\frac{1}{4} & r^{01} = \frac{1}{4} & r^{32} = -\frac{1}{4} & r^{23} = -\frac{1}{4} \\
r^{20} = -\frac{1}{4} & r^{31} = -\frac{1}{4} & r^{02} = \frac{1}{4} & r^{13} = -\frac{1}{4} \\
r^{30} = -\frac{1}{4} & r^{21} = -\frac{1}{4} & r^{12} = -\frac{1}{4} & r^{03} = \frac{1}{4}
\end{array}
$$

*Следовательно, отображение* $r$ *имеет вид*

$$r(a) = a^0 + a^2 i + a^3 j + a^1 k$$

$$
\begin{aligned}
r(a) = \frac{1}{4}(&a - iai - jaj - kak - ia + ai - kaj - jak \\
&-ja - kai + aj - iak - ka - jai - iaj + ak)
\end{aligned}
$$

□

Глава 6

# Линейное отображение алгебры

## 6.1. Линейное отображение алгебры

Определение 6.1.1. *Пусть $A_1$ и $A_2$ - алгебры над кольцом $D$. Линейное отображение $D$-модуля $A_1$ в $D$-модуль $A_2$ называется* **линейным отображением** *$D$-алгебры $A_1$ в $D$-алгебру $A_2$. Обозначим $\mathcal{L}(D; A_1 \to A_2)$ множество линейных отображений $D$-алгебры $A_1$ в $D$-алгебру $A_2$.* □

Определение 6.1.2. *Пусть $A_1$, ..., $A_n$, $S$ - $D$-алгебры. Полилинейное отображение*

$$f : A_1 \times ... \times A_n \to S$$

*$D$-модулей $A_1$, ..., $A_n$ в $D$-модуль $S$, Мы будем называть отображение* **полилинейным отображением** *$D$-алгебр $A_1$, ..., $A_n$ в $D$-модуль $S$. Обозначим $\mathcal{L}(D; A_1 \times ... \times A_n \to S)$ множество полилинейных отображений $D$-алгебр $A_1$, ..., $A_n$ в $D$-алгебру $S$. Обозначим $\mathcal{L}(D; A^n \to S)$ множество $n$-линейных отображений $D$-алгебры $A$ ($A_1 = ... = A_n = A$) в $D$-алгебру $S$.* □

Теорема 6.1.3. *Тензорное произведение $A_1 \otimes ... \otimes A_n$ $D$-алгебр $A_1$, ..., $A_n$ является $D$-алгеброй.*

Доказательство. Согласно определению 5.1.1 и теореме 4.5.2, тензорное произведение $A_1 \otimes ... \otimes A_n$ $D$-алгебр $A_1$, ..., $A_n$ является $D$-модулем.

Рассмотрим отображение

$$* : (A_1 \times ... \times A_n) \times (A_1 \times ... \times A_n) \to A_1 \otimes ... \otimes A_n \tag{6.1.1}$$

определённое равенством

$$(a_1, ..., a_n) * (b_1, ..., b_n) = (a_1 b_1) \otimes ... \otimes (a_n b_n) \tag{6.1.2}$$

Для заданных значений переменных $b_1$, ..., $b_n$, отображение (6.1.1) полилинейно по переменным $a_1$, ..., $a_n$. Согласно теореме 4.5.4, существует линейное отображение

$$* (b_1, ..., b_n) : A_1 \otimes ... \otimes A_n \to A_1 \otimes ... \otimes A_n \tag{6.1.3}$$

определённое равенством

$$(a_1 \otimes ... \otimes a_n) * (b_1, ..., b_n) = (a_1 b_1) \otimes ... \otimes (a_n b_n) \tag{6.1.4}$$

Так как произвольный тензор $a \in A_1 \otimes ... \otimes A_n$ можно представить в виде суммы тензоров вида $a_1 \otimes ... \otimes a_n$, то для заданного тензора $a \in A_1 \otimes ... \otimes A_n$ отображение (6.1.3) является полилинейным отображением переменных $b_1$, ..., $b_n$. Согласно теореме 4.5.4, существует линейное отображение

$$* (a) : A_1 \otimes ... \otimes A_n \to A_1 \otimes ... \otimes A_n \tag{6.1.5}$$

определёнwhichное равенством

$$(a_1 \otimes ... \otimes a_n) * (b_1 \otimes ... \otimes b_n) = (a_1 b_1) \otimes ... \otimes (a_n b_n) \tag{6.1.6}$$

Следовательно, равенство (6.1.6) определяет билинейное отображение

$$* : (A_1 \otimes ... \otimes A_n) \times (A_1 \otimes ... \otimes A_n) \to A_1 \otimes ... \otimes A_n \tag{6.1.7}$$

Билинейное отображение (6.1.7) определяет произведение в $D$-модуле $A_1 \otimes ... \otimes A_n$. □

В случае тензорного произведения $D$-алгебр $A_1$, $A_2$ мы будем рассматривать произведение, определённое равенством

$$(a_1 \otimes a_2) \circ (b_1 \otimes b_2) = (a_1 b_1) \otimes (b_2 a_2) \tag{6.1.8}$$

Теорема 6.1.4. *Пусть $\overline{\overline{e}}_i$ - базис алгебры $A_i$ над кольцом $D$. Пусть $B_{i\cdot}{}^{j}_{kl}$ - структурные константы алгебры $A_i$ относительно базиса $\overline{\overline{e}}_i$. Структурные константы тензорного произведение $A_1 \otimes ... \otimes A_n$ относительно базиса $e_{1\cdot i_1} \otimes ... \otimes e_{n\cdot i_n}$ имеют вид*

$$C^{\cdot j_1...j_n}_{\cdot k_1...k_n \cdot l_1...l_n} = C_{1\cdot}{}^{j_1}_{k_1 l_1} ... C_{n\cdot}{}^{j_n}_{k_n l_n} \tag{6.1.9}$$

Доказательство. Непосредственное перемножение тензоров $e_{1\cdot i_1} \otimes ... \otimes e_{n\cdot i_n}$ имеет вид

$$\begin{aligned}
&(e_{1\cdot k_1} \otimes ... \otimes e_{n\cdot k_n})(e_{1\cdot l_1} \otimes ... \otimes e_{n\cdot l_n}) \\
=&(\overline{e}_{1\cdot k_1} \overline{e}_{1\cdot l_1}) \otimes ... \otimes (\overline{e}_{n\cdot k_n} \overline{e}_{n\cdot l_n}) \\
=&(\overline{e}_{1\cdot k_1} \overline{e}_{1\cdot l_1}) \otimes ... \otimes (\overline{e}_{n\cdot k_n} \overline{e}_{n\cdot l_n}) \\
=&(C_{1\cdot}{}^{j_1}_{k_1 l_1} \overline{e}_{1\cdot j_1}) \otimes ... \otimes (C_{n\cdot}{}^{j_n}_{k_n l_n} \overline{e}_{n\cdot j_n}) \\
=&C_{1\cdot}{}^{j_1}_{k_1 l_1} ... C_{n\cdot}{}^{j_n}_{k_n l_n} \overline{e}_{1\cdot j_1} \otimes ... \otimes \overline{e}_{n\cdot j_n}
\end{aligned} \tag{6.1.10}$$

Согласно определению структурных констант

$$(e_{1\cdot k_1} \otimes ... \otimes e_{n\cdot k_n})(e_{1\cdot l_1} \otimes ... \otimes e_{n\cdot l_n}) = C^{\cdot j_1...j_n}_{\cdot k_1...k_n \cdot l_1...l_n} (e_{1\cdot j_1} \otimes ... \otimes e_{n\cdot j_n}) \tag{6.1.11}$$

Равенство (6.1.9) следует из сравнения (6.1.10), (6.1.11).

Из цепочки равенств

$$\begin{aligned}
&(a_1 \otimes ... \otimes a_n)(b_1 \otimes ... \otimes b_n) \\
=&(a_1^{k_1} \overline{e}_{1\cdot k_1} \otimes ... \otimes a_n^{k_n} \overline{e}_{n\cdot k_n})(b_1^{l_1} \overline{e}_{1\cdot l_1} \otimes ... \otimes b_n^{l_n} \overline{e}_{n\cdot l_n}) \\
=&a_1^{k_1} ... a_n^{k_n} b_1^{l_1} ... b_n^{l_n} (\overline{e}_{1\cdot k_1} \otimes ... \otimes \overline{e}_{n\cdot k_n})(\overline{e}_{1\cdot l_1} \otimes ... \otimes \overline{e}_{n\cdot l_n}) \\
=&a_1^{k_1} ... a_n^{k_n} b_1^{l_1} ... b_n^{l_n} C^{\cdot j_1...j_n}_{\cdot k_1...k_n \cdot l_1...l_n} (e_{1\cdot j_1} \otimes ... \otimes e_{n\cdot j_n}) \\
=&a_1^{k_1} ... a_n^{k_n} b_1^{l_1} ... b_n^{l_n} C_{1\cdot}{}^{j_1}_{k_1 l_1} ... C_{n\cdot}{}^{j_n}_{k_n l_n} (e_{1\cdot j_1} \otimes ... \otimes e_{n\cdot j_n}) \\
=&(a_1^{k_1} b_1^{l_1} C_{1\cdot}{}^{j_1}_{k_1 l_1} \overline{e}_{1\cdot j_1}) \otimes ... \otimes (a_n^{k_n} b_n^{l_n} C_{n\cdot}{}^{j_n}_{k_n l_n} \overline{e}_{n\cdot j_n}) \\
=&(a_1 b_1) \otimes ... \otimes (a_n b_n)
\end{aligned}$$

следует, что определение произведения (6.1.11) со структурными константами (6.1.9) согласовано с определением произведения (6.1.6). □

Теорема 6.1.5. *Для тензоров $a$, $b \in A_1 \otimes ... \otimes A_n$ стандартные компоненты произведения удовлетворяют равенству*

$$(ab)^{j_1...j_n} = C^{\cdot j_1...j_n}_{\cdot k_1...k_n \cdot l_1...l_n} a^{k_1...k_n} b^{l_1...l_n} \tag{6.1.12}$$

Доказательство. Согласно определению

$$ab = (ab)^{j_1...j_n} e_{1\cdot j_1} \otimes ... \otimes e_{n\cdot j_n} \tag{6.1.13}$$

В тоже время

$$\begin{aligned} ab &= a^{k_1...k_n} e_{1\cdot k_1} \otimes ... \otimes e_{n\cdot k_n} b^{k_1...k_n} e_{1\cdot l_1} \otimes ... \otimes e_{n\cdot l_n} \\ &= a^{k_1...k_n} b^{k_1...k_n} C^{\cdot j_1...j_n}_{\cdot k_1...k_n\cdot l_1...l_n} e_{1\cdot j_1} \otimes ... \otimes e_{n\cdot j_n} \end{aligned} \tag{6.1.14}$$

Равенство (6.1.12) следует из равенств (6.1.13), (6.1.14). □

Теорема 6.1.6. *Если алгебра $A_i$, $i = 1, ..., n$, ассоциативна, то тензорное произведение $A_1 \otimes ... \otimes A_n$ - ассоциативная алгебра.*

Доказательство. Поскольку

$$\begin{aligned} &((e_{1\cdot i_1} \otimes ... \otimes e_{n\cdot i_n})(e_{1\cdot j_1} \otimes ... \otimes e_{n\cdot j_n}))(e_{1\cdot k_1} \otimes ... \otimes e_{n\cdot k_n}) \\ =&((\overline{e}_{1\cdot i_1}\overline{e}_{1\cdot j_1}) \otimes ... \otimes (\overline{e}_{n\cdot i_n}\overline{e}_{1\cdot j_n}))(e_{1\cdot k_1} \otimes ... \otimes e_{n\cdot k_n}) \\ =&((\overline{e}_{1\cdot i_1}\overline{e}_{1\cdot j_1})\overline{e}_{1\cdot k_1}) \otimes ... \otimes ((\overline{e}_{n\cdot i_n}\overline{e}_{1\cdot j_n})\overline{e}_{1\cdot k_n}) \\ =&(\overline{e}_{1\cdot i_1}(\overline{e}_{1\cdot j_1}\overline{e}_{1\cdot k_1})) \otimes ... \otimes (\overline{e}_{n\cdot i_n}(\overline{e}_{1\cdot j_n}\overline{e}_{1\cdot k_n})) \\ =&(e_{1\cdot i_1} \otimes ... \otimes e_{n\cdot i_n})((\overline{e}_{1\cdot j_1}\overline{e}_{1\cdot k_1}) \otimes ... \otimes (\overline{e}_{1\cdot j_n}\overline{e}_{1\cdot k_n})) \\ =&(e_{1\cdot i_1} \otimes ... \otimes e_{n\cdot i_n})((e_{1\cdot j_1} \otimes ... \otimes e_{n\cdot j_n})(e_{1\cdot k_1} \otimes ... \otimes e_{n\cdot k_n})) \end{aligned}$$

то

$$\begin{aligned} (ab)c =& a^{i_1...i_n} b^{j_1...j_n} c^{k_1...k_n} \\ &((e_{1\cdot i_1} \otimes ... \otimes e_{n\cdot i_n})(e_{1\cdot j_1} \otimes ... \otimes e_{n\cdot j_n}))(e_{1\cdot k_1} \otimes ... \otimes e_{n\cdot k_n}) \\ =& a^{i_1...i_n} b^{j_1...j_n} c^{k_1...k_n} \\ &(e_{1\cdot i_1} \otimes ... \otimes e_{n\cdot i_n})((e_{1\cdot j_1} \otimes ... \otimes e_{n\cdot j_n})(e_{1\cdot k_1} \otimes ... \otimes e_{n\cdot k_n})) \\ =& a(bc) \end{aligned}$$

□

Теорема 6.1.7. *Пусть $A$ - алгебра над коммутативным кольцом $D$. Существует линейное отображение*

$$h : a \otimes b \in A \otimes A \to ab \in A$$

Доказательство. Теорема является следствием определения 5.1.1 и теоремы 4.5.4. □

Теорема 6.1.8. *Пусть отображение*

$$f : A_1 \to A_2$$

*является линейным отображением $D$-алгебры $A_1$ в $D$-алгебру $A_2$. Тогда отображения $af$, $fb$, $a$, $b \in A_2$, определённые равенствами*

$$(af) \circ x = a(f \circ x)$$

$$(fb) \circ x = (f \circ x)b$$

*также являются линейными.*

Доказательство. Утверждение теоремы следует из цепочек равенств

$$\begin{aligned}
(af)\circ(x+y) &= a(f\circ(x+y)) = a(f\circ x+f\circ y) = a(f\circ x)+a(f\circ y)\\
&= (af)\circ x+(af)\circ y\\
(af)\circ(px) &= a(f\circ(px)) = ap(f\circ x) = pa(f\circ x)\\
&= p((af)\circ x)\\
(fb)\circ(x+y) &= (f\circ(x+y))b = (f\circ x+f\circ y)\ b = (f\circ x)b+(f\circ y)b\\
&= (fb)\circ x+(fb)\circ y\\
(fb)\circ(px) &= (f\circ(px))b = p(f\circ x)b\\
&= p((fb)\circ x)
\end{aligned}$$

□

## 6.2. Алгебра $\mathcal{L}(D;A\to A)$

Теорема 6.2.1. *Пусть $A$, $B$, $C$ - алгебры над коммутативным кольцом $D$. Пусть $f$ - линейное отображение из $D$-алгебры $A$ в $D$-алгебру $B$. Пусть $g$ - линейное отображение из $D$-алгебры $B$ в $D$-алгебру $C$. Отображение $g\circ f$, определённое диаграммой*

(6.2.1)

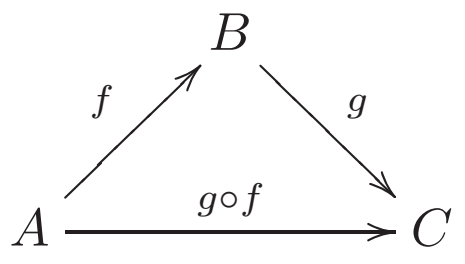

*является линейным отображением из $D$-алгебры $A$ в $D$-алгебру $C$.*

Доказательство. Доказательство теоремы следует из цепочек равенств

$$\begin{aligned}
(g\circ f)\circ(a+b) &= g\circ(f\circ(a+b)) = g\circ(f\circ a+f\circ b)\\
&= g\circ(f\circ a)+g\circ(f\circ b) = (g\circ f)\circ a+(g\circ f)\circ b\\
(g\circ f)\circ(pa) &= g\circ(f\circ(pa)) = g\circ(p\ f\circ a) = p\ g\circ(f\circ a)\\
&= p\ (g\circ f)\circ a
\end{aligned}$$

□

Теорема 6.2.2. *Пусть $A$, $B$, $C$ - алгебры над коммутативным кольцом $D$. Пусть $f$ - линейное отображение из $D$-алгебры $A$ в $D$-алгебру $B$. Отображение $f$ порождает линейное отображение*

$$f^*:g\in\mathcal{L}(D;B\to C)\to g\circ f\in\mathcal{L}(D;A\to C) \tag{6.2.2}$$

(6.2.3)

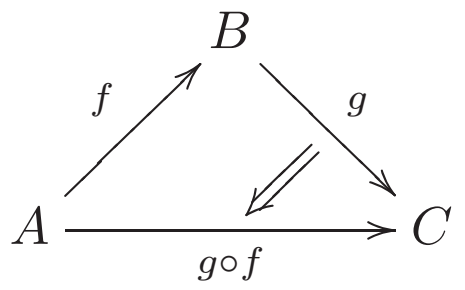

Доказательство. Доказательство теоремы следует из цепочек равенств [6.1]

$$\begin{aligned}((g_1+g_2)\circ f)\circ a&=(g_1+g_2)\circ(f\circ a)=g_1\circ(f\circ a)+g_2\circ(f\circ a)\\&=(g_1\circ f)\circ a+(g_2\circ f)\circ a\\&=(g_1\circ f+g_2\circ f)\circ a\\((pg)\circ f)\circ a&=(pg)\circ(f\circ a)=p\ g\circ(f\circ a)=p\ (g\circ f)\circ a\\&=(p(g\circ f))\circ a\end{aligned}$$

□

Теорема 6.2.3. *Пусть $A$, $B$, $C$ - алгебры над коммутативным кольцом $D$. Пусть $g$ - линейное отображение из $D$-алгебры $B$ в $D$-алгебру $C$. Отображение $g$ порождает линейное отображение*

$$g_*:f\in\mathcal{L}(D;A\to B)\to g\circ f\in\mathcal{L}(D;A\to C)\tag{6.2.6}$$

(6.2.7)

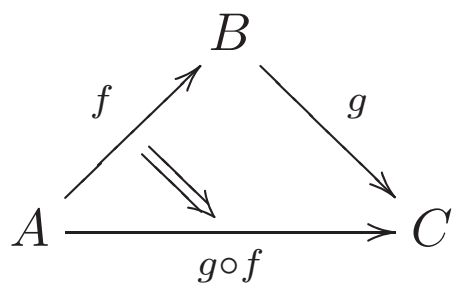

Доказательство. Доказательство теоремы следует из цепочек равенств [6.2]

$$\begin{aligned}(g\circ(f_1+f_2))\circ a&=g\circ((f_1+f_2)\circ a)=g\circ(f_1\circ a+f_2\circ a)\\&=g\circ(f_1\circ a)+g\circ(f_2\circ a)=(g\circ f_1)\circ a+(g\circ f_2)\circ a\\&=(g\circ f_1+g\circ f_2)\circ a\\(g\circ(pf))\circ a&=g\circ((pf)\circ a)=g\circ(p\ (f\circ a))=p\ g\circ(f\circ a)\\&=p\ (g\circ f)\circ a=(p(g\circ f))\circ a\end{aligned}$$

□

Теорема 6.2.4. *Пусть $A$, $B$, $C$ - алгебры над коммутативным кольцом $D$. Отображение*

$$\circ:(g,f)\in\mathcal{L}(D;B\to C)\times\mathcal{L}(D;A\to B)\to g\circ f\in\mathcal{L}(D;A\to C)\tag{6.2.10}$$

*является билинейным отображением.*

Доказательство. Теорема является следствием теорем 6.2.2, 6.2.3. □

---

[6.1] Мы пользуемся следующими определениями операций над отображениями

$$(f+g)\circ a=f\circ a+g\circ a\tag{6.2.4}$$

$$(pf)\circ a=p\,f\circ a\tag{6.2.5}$$

[6.2] Мы пользуемся следующими определениями операций над отображениями

$$(f+g)\circ a=f\circ a+g\circ a\tag{6.2.8}$$

$$(pf)\circ a=p\,f\circ a\tag{6.2.9}$$

Теорема 6.2.5. *Пусть $A$ - алгебра над коммутативным кольцом $D$. $D$-модуль $\mathcal{L}(D;A \to A)$, оснащённый произведением*

$$\circ : (g,f) \in \mathcal{L}(D;A \to A) \times \mathcal{L}(D;A \to A) \to g \circ f \in \mathcal{L}(D;A \to A) \tag{6.2.11}$$

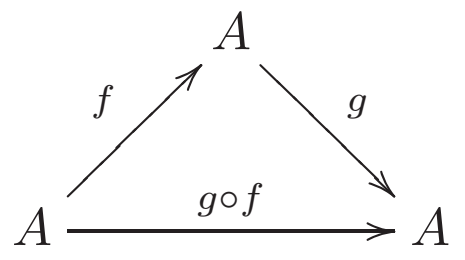

*является алгеброй над $D$.*

Доказательство. Теорема является следствием определения 5.1.1 и теоремы 6.2.4. □

## 6.3. Линейное отображение в ассоциативную алгебру

Теорема 6.3.1. *Рассмотрим $D$-алгебры $A_1$ и $A_2$. Для заданного отображения $f \in \mathcal{L}(D;A_1 \to A_2)$ отображение*

$$g : A_2 \times A_2 \to \mathcal{L}(D;A_1 \to A_2)$$

$$g(a,b) \circ f = afb$$

*является билинейным отображением.*

Доказательство. Утверждение теоремы следует из цепочек равенств

$$\begin{aligned}
((a_1+a_2)fb)\circ x &= (a_1+a_2)\ f\circ x\ b = a_1\ f\circ x\ b + a_2\ f\circ x\ b\\
&= (a_1fb)\circ x + (a_2fb)\circ x = (a_1fb+a_2fb)\circ x\\
((pa)fb)\circ x &= (pa)\ f\circ x\ b = p(a\ f\circ x\ b) = p((afb)\circ x) = (p(afb))\circ x\\
(af(b_1+b_2))\circ x &= a\ f\circ x\ (b_1+b_2) = a\ f\circ x\ b_1 + a\ f\circ x\ b_2\\
&= (afb_1)\circ x + (afb_2)\circ x = (afb_1+afb_2)\circ x\\
(af(pb))\circ x &= a\ f\circ x\ (pb) = p(a\ f\circ x\ b) = p((afb)\circ x) = (p(afb))\circ x
\end{aligned}$$

□

Теорема 6.3.2. *Рассмотрим $D$-алгебры $A_1$ и $A_2$. Для заданного отображения $f \in \mathcal{L}(D;A_1 \to A_2)$ существует линейное отображение*

$$h : A_2 \otimes A_2 \to \mathcal{L}(D;A_1 \to A_2)$$

*определённое равенством*

$$(a\otimes b)\circ f = afb \tag{6.3.1}$$

Доказательство. Утверждение теоремы является следствием теорем 4.5.4, 6.3.1. □

Теорема 6.3.3. *Рассмотрим $D$-алгебры $A_1$ и $A_2$. Определим произведение в алгебре $A_2 \otimes A_2$ согласно правилу*

$$(c\otimes d)\circ(a\otimes b) = (ca)\otimes(bd) \tag{6.3.2}$$

*Линейное отображение*

$$h : A_2 \otimes A_2 \to {}^*\mathcal{L}(D;A_1 \to A_2) \tag{6.3.3}$$

*определённое равенством*

$$(a\otimes b)\circ f = afb \quad a,b \in A_2 \quad f \in \mathcal{L}(D;A_1 \to A_2) \tag{6.3.4}$$

*является представлением*[6.3] *алгебры* $A_2\otimes A_2$ *в модуле* $\mathcal{L}(D;A_1\to A_2)$.

Доказательство. Согласно теореме 6.1.8, отображение (6.3.4) является преобразованием модуля $\mathcal{L}(D;A_1\to A_2)$. Для данного тензора $c\in A_2\otimes A_2$ преобразование $h(c)$ является линейным преобразованием модуля $\mathcal{L}(D;A_1\to A_2)$, так как

$$\begin{aligned}((a\otimes b)\circ(f_1+f_2))\circ x&=(a(f_1+f_2)b)\circ x=a((f_1+f_2)\circ x)b\\&=a(f_1\circ x+f_2\circ x)b=a(f_1\circ x)b+a(f_2\circ x)b\\&=(af_1b)\circ x+(af_2b)\circ x\\&=(a\otimes b)\circ f_1\circ x+(a\otimes b)\circ f_2\circ x\\&=((a\otimes b)\circ f_1+(a\otimes b)\circ f_2)\circ x\\((a\otimes b)\circ(pf))\circ x&=(a(pf)b)\circ x=a((pf)\circ x)b\\&=a(p\ f\circ x)b=pa(f\circ x)b\\&=p\ (afb)\circ x=p\ ((a\otimes b)\circ f)\circ x\\&=(p((a\otimes b)\circ f))\circ x\end{aligned}$$

Согласно теореме 6.3.2, отображение (6.3.4) является линейным отображением.

Пусть $f\in\mathcal{L}(D;A_1\to A_2)$, $a\otimes b$, $c\otimes d\in A_2\otimes A_2$. Согласно теореме 6.3.2

$$(a\otimes b)\circ f=afb\in\mathcal{L}(D;A_1\to A_2)$$

Следовательно, согласно теореме 6.3.2

$$(c\otimes d)\circ((a\otimes b)\circ f)=c(afb)d$$

Поскольку произведение в алгебре $A_2$ ассоциативно, то

$$(c\otimes d)\circ((a\otimes b)\circ f)=c(afb)d=(ca)f(bd)=(ca\otimes bd)\circ f$$

Следовательно, если мы определим произведение в алгебре $A_2\otimes A_2$ согласно равенству (6.3.2), то отображение (6.3.3) является морфизмом алгебр. Согласно определению [4]-2.1.2, отображение (6.3.4) является представлением алгебры $A_2\otimes A_2$ в модуле $\mathcal{L}(D;A_1\to A_2)$. □

Теорема 6.3.4. *Рассмотрим $D$-алгебру $A$. Определим произведение в алгебре $A\otimes A$ согласно правилу* (6.3.2). *Представление*

(6.3.5) $$h:A\otimes A\to{}^*\mathcal{L}(D;A\to A)$$

*алгебры $A\otimes A$ в модуле $\mathcal{L}(D;A\to A)$, определённое равенством*

(6.3.6) $$(a\otimes b)\circ f=afb\quad a,b\in A\quad f\in\mathcal{L}(D;A\to A)$$

*позволяет отождествить тензор $d\in A\otimes A$ с отображением $d\circ\delta\in\mathcal{L}(D;A\to A)$, где $\delta\in\mathcal{L}(D;A\to A)$ - тождественное отображение.*

Доказательство. Согласно теореме 6.3.2, отображение $f\in\mathcal{L}(D;A\to A)$ и тензор $d\in A\otimes A$ порождают отображение

(6.3.7) $$x\to(d\circ f)\circ x$$

Если мы положим $f=\delta$, $d=a\otimes b$, то равенство (6.3.7) приобретает вид

(6.3.8) $$((a\otimes b)\circ\delta)\circ x=(a\delta b)\circ x=a\ (\delta\circ x)\ b=axb$$

[6.3] Определение представления $\Omega$-алгебры дано в определении [4]-2.1.2.

Если мы положим

$$((a \otimes b) \circ \delta) \circ x = (a \otimes b) \circ (\delta \circ x) = (a \otimes b) \circ x \tag{6.3.9}$$

то сравнение равенств (6.3.8) и (6.3.9) даёт основание отождествить действие тензора $a \otimes b$ с преобразованием $(a \otimes b) \circ \delta$. □

Из теоремы 6.3.4 следует, что отображение (6.3.4) можно рассматривать как произведение отображений $a \otimes b$ и $f$. Тензор $a \in A_2 \otimes A_2$ **невырожден**, если существует тензор $b \in A_2 \otimes A_2$ такой, что $a \circ b = 1 \otimes 1$.

Определение 6.3.5. *Рассмотрим*[6.4] *представление алгебры $A_2 \otimes A_2$ в модуле $\mathcal{L}(D; A_1 \to A_2)$.* **Орбитой линейного отображения** $f \in \mathcal{L}(D; A_1 \to A_2)$ *называется множество*

$$(A_2 \otimes A_2) \circ f = \{g = d \circ f : d \in A_2 \otimes A_2\}$$

□

Теорема 6.3.6. *Рассмотрим $D$-алгебру $A_1$ и ассоциативную $D$-алгебру $A_2$. Рассмотрим представление алгебры $A_2 \otimes A_2$ в модуле $\mathcal{L}(D; A_1; A_2)$. Отображение*

$$h : A_1 \to A_2$$

*порождённое отображением*

$$f : A_1 \to A_2$$

*имеет вид*

$$h = (a_{s \cdot 0} \otimes a_{s \cdot 1}) \circ f = a_{s \cdot 0} f a_{s \cdot 1} \tag{6.3.10}$$

Доказательство. Произвольный тензор $a \in A_2 \otimes A_2$ можно представить в виде

$$a = a_{s \cdot 0} \otimes a_{s \cdot 1}$$

Согласно теореме 6.3.3, отображение (6.3.4) линейно. Это доказывает утверждение теоремы. □

Теорема 6.3.7. *Пусть $A_2$ - алгебра с единицей $e$. Пусть $a \in A_2 \otimes A_2$ - невырожденный тензор. Орбиты линейных отображений $f \in \mathcal{L}(D; A_1 \to A_2)$ и $g = a \circ f$ совпадают*

$$(A_2 \otimes A_2) \circ f = (A_2 \otimes A_2) \circ g \tag{6.3.11}$$

Доказательство. Если $h \in (A_2 \otimes A_2) \circ g$, то существует $b \in A_2 \otimes A_2$ такое, что $h = b \circ g$. Тогда

$$h = b \circ (a \circ f) = (b \circ a) \circ f \tag{6.3.12}$$

Следовательно, $h \in (A_2 \otimes A_2) \circ f$,

$$(A_2 \otimes A_2) \circ g \subset (A_2 \otimes A_2) \circ f \tag{6.3.13}$$

Так как $a$ - невырожденный тензор, то

$$f = a^{-1} \circ g \tag{6.3.14}$$

Если $h \in (A_2 \otimes A_2) \circ f$, то существует $b \in A_2 \otimes A_2$ такое, что

$$h = b \circ f \tag{6.3.15}$$

[6.4] Определение дано по аналогии с определением [4]-3.1.8.

Из равенств (6.3.14), (6.3.15), следует, что

$$h = b \circ (a^{-1} \circ g) = (b \circ a^{-1}) \circ g$$

Следовательно, $h \in (A_2 \otimes A_2) \circ g$,

$$(A_2 \otimes A_2) \circ f \subset (A_2 \otimes A_2) \circ g \tag{6.3.16}$$

(6.3.11) следует из равенств (6.3.13), (6.3.16). □

Из теоремы 6.3.7 также следует, что, если $g = a \circ f$ и $a \in A_2 \otimes A_2$ - вырожденный тензор, то отношение (6.3.13) верно. Однако основной результат теоремы 6.3.7 состоит в том, что представления алгебры $A_2 \otimes A_2$ в модуле $\mathcal{L}(D; A_1 \to A_2)$ порождает отношение эквивалентности в модуле $\mathcal{L}(D; A_1 \to A_2)$. Если удачно выбрать представители каждого класса эквивалентности, то полученное множество будет множеством образующих рассматриваемого представления.[6.5]

## 6.4. Линейное отображение в свободную конечно мерную ассоциативную алгебру

Теорема 6.4.1. *Пусть $A_1$ - свободный $D$-модуль. Пусть $A_2$ - свободная конечно мерная ассоциативная $D$-алгебра. Пусть $\overline{\overline{e}}$ - базис $D$-модуля $A_2$. Пусть $\overline{\overline{I}}$ - базис левого $A_2 \otimes A_2$-модуля $\mathcal{L}(D; A_1 \to A_2)$.*[6.6]

6.4.1.1: *Отображение*

$$f : A_1 \to A_2$$

*имеет следующее разложение*

$$f = f^k \circ I_k \tag{6.4.1}$$

*where*

$$f^k = f^k_{s_k\cdot 0} \otimes f^k_{s_k\cdot 1} \qquad f^k \in A_2 \otimes A_2$$

6.4.1.2: *Отображение $f$ имеет стандартное представление*

$$f = f^{k\cdot ij}(e_i \otimes e_j) \circ I_k = f^{k\cdot ij} e_i I_k e_j \tag{6.4.2}$$

Доказательство. Поскольку $\overline{\overline{I}}$ является базисом левого $A_2 \otimes A_2$-модуля $\mathcal{L}(D; A_1 \to A_2)$, то согласно определению [4]-2.7.1 и теореме 4.1.4, существует разложение

$$f = f^k \circ I_k \qquad f^k \in A_2 \otimes A_2 \tag{6.4.3}$$

линейного отображения $f$ относительно базиса $\overline{\overline{I}}$. Согласно определению (3.3.20),

$$f^k = f^k_{s_k\cdot 0} \otimes f^k_{s_k\cdot 1} \tag{6.4.4}$$

---

[6.5] Множество образующих представления определено в определении [4]-2.6.5.

[6.6] Если $D$-модуль $A_1$ или $D$-модуль $A_2$ не является свободным $D$-модулем, то мы будем рассматривать множество

$$\overline{\overline{I}} = \{I_k \in \mathcal{L}(D; A_1 \to A_2) : k = 1, ..., n\}$$

линейно независимых линейных отображений. Теорема верна для любого линейного отображения

$$f : A_1 \to A_2$$

порождённого множеством линейных отображений $\overline{\overline{I}}$.

Равенство (6.4.1) является следствием равенств (6.4.3), (6.4.4). Согласно теореме 4.5.6, стандартное представление тензора $f^k$ имеет вид

$$f^k = f^{k\cdot ij} e_i \otimes e_j \tag{6.4.5}$$

Равенство (6.4.2) следует из равенств (6.4.1), (6.4.5). □

Теорема 6.4.2. *Пусть $A_1$ - свободный $D$-модуль. Пусть $A_2$ - свободная ассоциативная $D$-алгебра. Пусть $\overline{\overline{I}}$ - базис левого $A_2\otimes A_2$-модуля $\mathcal{L}(D;A_1\to A_2)$. Для любого отображения $I_k\in\overline{\overline{I}}$, существует множество линейных отображений*

$$I^l_k : A_1\otimes A_1\to A_2\otimes A_2$$

*$D$-модуля $A_1\otimes A_1$ в $D$-модуль $A_2\otimes A_2$ таких, что*

$$I_k\circ a\circ x = (I^l_k\circ a)\circ I_l\circ x \tag{6.4.6}$$

*Отображение $I^l_k$ называется* **преобразованием сопряжения**.

Доказательство. Согласно теореме 6.2.1, для произвольного тензора $a\in A_1\otimes A_1$, отображение

$$x\to I_k\circ a\circ x \tag{6.4.7}$$

является линейным. Согласно утверждению 6.4.1.1, существует разложение

$$I_k\circ a\circ x = b^l\circ I_l\circ x \qquad b^l\in A_2\times A_2 \tag{6.4.8}$$

Положим

$$b^l = I^l_k\circ a \tag{6.4.9}$$

Равенство (6.4.6) является следствием равенств (6.4.8), (6.4.9). Из равенств

$$\begin{aligned}(I^l_k\circ(a_1+a_2))\circ I_l\circ x &= I_k\circ(a_1+a_2)\circ x\\ &= I_k\circ a_1\circ x + I_k\circ a_2\circ x\\ &= (I^l_k\circ a_1)\circ I_l\circ x + (I^l_k\circ a_1)\circ I_l\circ x\end{aligned}$$

$$\begin{aligned}(I^l_k\circ(da))\circ I_l\circ x &= I_k\circ(da)\circ x = I_k\circ(d(a\circ x))\\ &= d(I_k\circ a\circ x) = d((I^l_k\circ a)\circ I_l\circ x)\\ &= (d(I^l_k\circ a))\circ I_l\circ x\end{aligned}$$

следует, что отображение $I^l_k$ является линейным. □

Теорема 6.4.3. *Пусть $A_1$ - свободный $D$-модуль. Пусть $A_2$, $A_3$ - свободные ассоциативные $D$-алгебры. Пусть $\overline{\overline{I}}$ - базис левого $A_2\otimes A_2$-модуля $\mathcal{L}(D;A_1\to A_2)$. Пусть $\overline{\overline{J}}$ - базис левого $A_3\otimes A_3$-модуля $\mathcal{L}(D;A_2\to A_3)$.*

6.4.3.1: *Множество отображений*

$$\overline{\overline{K}} = \{K_{lk} : K_{lk} = J_l\circ I_k, J_l\in\overline{\overline{J}}, I_k\in\overline{\overline{I}}\} \tag{6.4.10}$$

*является базисом левого $A_3\otimes A_3$-модуля $\mathcal{L}(D;A_1\to A_2\to A_3)$.*

6.4.3.2: *Пусть линейное отображение*

$$f: A_1 \to A_2$$

*имеет разложение*

(6.4.11) $$f = f^k \circ I_k$$

*относительно базиса* $\overline{\overline{I}}$. *Пусть линейное отображение*

$$g: A_2 \to A_3$$

*имеет разложение*

(6.4.12) $$g = g^l \circ J_l$$

*относительно базиса* $\overline{\overline{J}}$. *Тогда линейное отображение*

(6.4.13) $$h = g \circ f$$

*имеет разложение*

(6.4.14) $$h = h^{lk} \circ K_{lk}$$

*относительно базиса* $\overline{\overline{K}}$, *где*

(6.4.15) $$h^{lk} = g^l \circ (J^k_m \circ f^m)$$

Доказательство. Равенство

(6.4.16) $$h \circ a = g \circ f \circ a = g^l \circ J_l \circ f^k \circ I_k \circ a$$

является следствием равенств (6.4.11), (6.4.12), (6.4.13). Равенство

(6.4.17) $$\begin{aligned} h \circ a = g \circ f \circ a &= g^l \circ (J^m_l \circ f^k) \circ J_m \circ I_k \circ a \\ &= g^l \circ (J^m_l \circ f^k) \circ K_{mk} \circ a \end{aligned}$$

является следствием равенств (6.4.10), (6.4.16) и теоремы 6.4.2. Из равенства (6.4.17) следует, что множество отображений $\overline{\overline{K}}$ порождает левый $A_3 \otimes A_3$-модуля
$\mathcal{L}(D; A_1 \to A_2 \to A_3)$. Из равенства

$$a^{lk} K_{lk} = (a^{lk} \circ J_l) \circ I_k = 0$$

следует, что

$$a^{lk} \circ J_l = 0$$

и, следовательно, $a^{lk} = 0$. Следовательно, множество $\overline{\overline{K}}$ является базисом левого $A_3 \otimes A_3$-модуля $\mathcal{L}(D; A_1 \to A_2 \to A_3)$. □

Теорема 6.4.4. *Пусть $A$ - свободная ассоциативная $D$-алгеба. Пусть левый $A \otimes A$-модуль $\mathcal{L}(D; A \to A)$ порождён тождественным отображением $I_0 = \delta$. Пусть линейное отображение*

$$f: A \to A$$

*имеет разложение*

(6.4.18) $$f = f_{s \cdot 0} \otimes f_{s \cdot 1}$$

*Пусть линейное отображение*

$$g: A \to A$$

*имеет разложение*

$$g = g_{t\cdot 0} \otimes g_{t\cdot 1} \tag{6.4.19}$$

*Тогда линейное отображение*

$$h = g \circ f \tag{6.4.20}$$

*имеет разложение*

$$h = h_{ts\cdot 0} \otimes h_{ts\cdot 1} \tag{6.4.21}$$

*где*

$$\begin{aligned} h_{ts\cdot 0} &= g_{t\cdot 0} f_{s\cdot 0}\\ h_{ts\cdot 1} &= f_{s\cdot 1} g_{t\cdot 1} \end{aligned} \tag{6.4.22}$$

Доказательство. Равенство

$$\begin{aligned} h\circ a &= g\circ f\circ a\\ &= (g_{t\cdot 0}\otimes g_{t\cdot 1})\circ (f_{s\cdot 0}\otimes f_{s\cdot 1})\circ a\\ &= (g_{t\cdot 0}\otimes g_{t\cdot 1})\circ (f_{s\cdot 0} a f_{s\cdot 1})\\ &= g_{t\cdot 0} f_{s\cdot 0} a f_{s\cdot 1} g_{t\cdot 1} \end{aligned} \tag{6.4.23}$$

является следствием равенств (6.4.18), (6.4.19), (6.4.20). Равенство (6.4.22) является следствием равенства (6.4.23). □

Теорема 6.4.5. *Пусть $\overline{\overline{e}}_1$ - базис свободного конечно мерного $D$-модуля $A_1$. Пусть $\overline{\overline{e}}_2$ - базис свободной конечно мерной ассоциативной $D$-алгебры $A_2$. Пусть $C^{p}_{kl}$ - структурные константы алгебры $A_2$. Пусть $\overline{\overline{I}}$ - базис левого $A_2\otimes A_2$-модуля $\mathcal{L}(D;A_1\to A_2)$ и $I_{k\cdot}{}^{j}_{i}$ - координаты отображения $I_k$ относительно базисов $\overline{\overline{e}}_1$ и $\overline{\overline{e}}_2$. Координаты $f^{k}_{l}$ отображения $f\in\mathcal{L}(D;A_1\to A_2)$ и его стандартные компоненты $f^{k\cdot ij}$ связаны равенством*

$$f^{k}_{l} = f^{k\cdot ij} I_{k\cdot}{}^{m}_{l} C^{p}_{im} C^{k}_{pj} \tag{6.4.24}$$

Доказательство. Относительно базисов $\overline{\overline{e}}_1$ и $\overline{\overline{e}}_2$, линейные отображения $f$ и $I_k$ имеют вид

$$f\circ x = f^{i}_{j} x^{j} e_{2\cdot i} \tag{6.4.25}$$

$$I_k\circ x = I_{k\cdot}{}^{i}_{j} x^{j} e_{2\cdot i} \tag{6.4.26}$$

Равенство

$$\begin{aligned} f^{k}_{l} x^{l} e_{2\cdot k} &= f^{k\cdot ij} e_{2\cdot i} I_{k\cdot}{}^{m}_{l} x^{l} e_{2\cdot m} \overline{e}_{2\cdot j}\\ &= f^{k\cdot ij} I_{k\cdot}{}^{m}_{l} x^{l} C^{p}_{im} C^{k}_{pj} e_{2\cdot k} \end{aligned} \tag{6.4.27}$$

является следствием равенств (6.4.2), (6.4.25), (6.4.26). Так как векторы $\overline{e}_{2\cdot k}$ линейно независимы и $x^{i}$ произвольны, то равенство (6.4.24) следует из равенства (6.4.27). □

Теорема 6.4.6. *Пусть $D$ является полем. Пусть $\overline{\overline{e}}_1$ - базис свободной конечно мерной $D$-алгебры $A_1$. Пусть $\overline{\overline{e}}_2$ - базис свободной конечно мерной ассоциативной $D$-алгебры $A_2$. Пусть $C_{2\cdot}{}^{p}_{kl}$ - структурные константы алгебры $A_2$. Рассмотрим матрицу*

$$\mathcal{B} = \left(\mathcal{C}^{\cdot k}_{m\cdot ij}\right) = \left(C_{2\cdot}{}^{p}_{im} C_{2\cdot}{}^{k}_{pj}\right) \tag{6.4.28}$$

*строки которой пронумерованы индексом* ${}^{\cdot k}_{m}$ *и столбцы пронумерованы индексом* ${}_{\cdot ij}$ . *Если матрица* $\mathcal{B}$ *невырождена, то для заданных координат линейного преобразования* $g^l_k$ *и для отображения* $f = \delta$ , *система линейных уравнений* (6.4.24) *относительно стандартных компонент этого преобразования* $g^{kr}$ *имеет единственное решение.*

*Если матрица* $\mathcal{B}$ *вырождена, то условием существования решения системы линейных уравнений* (6.4.24) *является равенство*

$$\mathrm{rank}\begin{pmatrix}\mathcal{C}^{\cdot k}_{m\cdot ij} & g^k_m\end{pmatrix} = \mathrm{rank}\,\mathcal{C} \tag{6.4.29}$$

*В этом случае система линейных уравнений* (6.4.24) *имеет бесконечно много решений и существует линейная зависимость между величинами* $g^k_m$.

Доказательство. Утверждение теоремы является следствием теории линейных уравнений над полем. □

Теорема 6.4.7. *Пусть* $A$ *- свободная конечно мерная ассоциативная алгебра над полем* $D$. *Пусть* $\overline{\overline{e}}$ *- базис алгебры* $A$ *над полем* $D$. *Пусть* $C^p_{kl}$ *- структурные константы алгебры* $A$. *Пусть матрица* (6.4.28) *вырождена. Пусть линейное отображение* $f \in \mathcal{L}(D; A \to A)$ *невырождено. Если координаты линейных преобразований* $f$ *и* $g$ *удовлетворяют равенству*

$$\mathrm{rank}\begin{pmatrix}\mathcal{C}^{\cdot k}_{m\cdot ij} & g^k_m & f^k_m\end{pmatrix} = \mathrm{rank}\,\mathcal{C} \tag{6.4.30}$$

*то система линейных уравнений*

$$g^k_l = f^m_l g^{ij} C^p_{im} C^k_{pj} \tag{6.4.31}$$

*имеет бесконечно много решений.*

Доказательство. Согласно равенству (6.4.30) и теореме 6.4.6, система линейных уравнений

$$f^k_l = f^{ij} C^p_{il} C^k_{pj} \tag{6.4.32}$$

имеет бесконечно много решений, соответствующих линейному отображению

$$f = f^{ij}\overline{e}_i \otimes \overline{e}_j \tag{6.4.33}$$

Согласно равенству (6.4.30) и теореме 6.4.6, система линейных уравнений

$$g^k_l = g^{ij} C^p_{il} C^k_{pj} \tag{6.4.34}$$

имеет бесконечно много решений, соответствующих линейному отображению

$$g = g^{ij}\overline{e}_i \otimes \overline{e}_j \tag{6.4.35}$$

Отображения $f$ и $g$ порождены отображением $\delta$. Согласно теореме 6.3.7, отображение $f$ порождает отображение $g$. Это доказывает утверждение теоремы. □

Теорема 6.4.8. *Пусть* $A$ *- свободная конечно мерная ассоциативная алгебра над полем* $D$. *Представление алгебры* $A \otimes A$ *в алгебре* $\mathcal{L}(D; A \to A)$ *имеет конечный* **базис** $\overline{\overline{I}}$.

6.4.8.1: *Линейное отображение* $f \in \mathcal{L}(D; A \to A)$ *имеет вид*

$$f = (a_{k\cdot s_k\cdot 0} \otimes a_{k\cdot s_k\cdot 1}) \circ I_k = \sum_k a_{k\cdot s_k\cdot 0} I_k a_{k\cdot s_k\cdot 1} \tag{6.4.36}$$

6.4.8.2: *Его стандартное представление имеет вид*

$$f = a^{k\cdot ij}(e_i \otimes e_j) \circ I_k = a^{k\cdot ij} e_i I_k e_j \tag{6.4.37}$$

Доказательство. Из теоремы 6.4.7 следует, что если матрица $\mathcal{B}$ вырождена и отображение $f$ удовлетворяет равенству

$$\operatorname{rank}\begin{pmatrix} \mathcal{C}^{\cdot k}_{m\cdot ij} & f^k_m \end{pmatrix} = \operatorname{rank}\mathcal{C} \tag{6.4.38}$$

то отображение $f$ порождает то же самое множество отображений, что порождено отображением $\delta$. Следовательно, для того, чтобы построить базис представления алгебры $A \otimes A$ в модуле $\mathcal{L}(D; A \to A)$, мы должны выполнить следующее построение.

Множество решений системы уравнений (6.4.31) порождает свободный подмодуль $\mathcal{L}$ модуля $\mathcal{L}(D; A \to A)$. Мы строим базис $(\overline{h}_1, ..., \overline{h}_k)$ подмодуля $\mathcal{L}$. Затем дополняем этот базис линейно независимыми векторами $\overline{h}_{k+1}$, ..., $\overline{h}_m$, которые не принадлежат подмодулю $\mathcal{L}$, таким образом, что множество векторов $\overline{h}_1$, ..., $\overline{h}_m$ является базисом модуля $\mathcal{L}(D; A \to A)$. Множество орбит $(A \otimes A) \circ \delta$, $(A \otimes A) \circ \overline{h}_{k+1}$, ..., $(A \otimes A) \circ \overline{h}_m$ порождает модуль $\mathcal{L}(D; A \to A)$. Поскольку множество орбит конечно, мы можем выбрать орбиты так, чтобы они не пересекались. Для каждой орбиты мы можем выбрать представитель, порождающий эту орбиту. □

Пример 6.4.9. *Для поля комплексных чисел алгебра $\mathcal{L}(R; C \to C)$ имеет базис*

$$I_0 \circ z = z$$
$$I_1 \circ z = \overline{z}$$

*Для алгебры кватернионов алгебра $\mathcal{L}(R; H \to H)$ имеет базис*

$$I_0 \circ z = z$$

□

## 6.5. Линейное отображение в неассоциативную алгебру

Так как произведение неассоциативно, мы можем предположить, что действие $a$, $b \in A$ на отображение $f$ может быть представленно либо в виде $a(fb)$, либо в виде $(af)b$. Однако это предположение приводит нас к довольно сложной структуре линейного отображения. Чтобы лучше представить насколько сложна структура линейного отображения, мы начнём с рассмотрения левого и правого сдвигов в неассоциативной алгебре.

Теорема 6.5.1. *Пусть*

$$l(a) \circ x = ax \tag{6.5.1}$$

*отображение левого сдвига. Тогда*

$$l(a) \circ l(b) = l(ab) - (a, b)_1 \tag{6.5.2}$$

*где мы определили линейное отображение*

$$(a, b)_1 \circ x = (a, b, x)$$

Доказательство. Из равенств (5.1.2), (6.5.1) следует

$$\begin{aligned}(l(a)\circ l(b))\circ x &= l(a)\circ(l(b)\circ x)\\ &= a(bx)=(ab)x-(a,b,x)\\ &= l(ab)\circ x-(a,b)_1\circ x\end{aligned} \tag{6.5.3}$$

Равенство (6.5.2) следует из равенства (6.5.3). □

Теорема 6.5.2. *Пусть*

$$r(a)\circ x = xa \tag{6.5.4}$$

*отображение правого сдвига. Тогда*

$$r(a)\circ r(b) = r(ba)+(b,a)_2 \tag{6.5.5}$$

*где мы определили линейное отображение*

$$(b,a)_2\circ x=(x,b,a)$$

Доказательство. Из равенств (5.1.2), (6.5.4) следует

$$\begin{aligned}(r(a)\circ r(b))\circ x &= r(a)\circ(r(b)\circ x)\\ &= (xb)a=x(ba)+(x,b,a)\\ &= r(ba)\circ x+(x,b,a)\end{aligned} \tag{6.5.6}$$

Равенство (6.5.5) следует из равенства (6.5.6). □

Пусть

$$f:A\to A\quad f=(ax)b$$

линейное отображение алгебры $A$. Согласно теореме 6.1.8, отображение

$$g:A\to A\quad g=(cf)d$$

также линейное отображение. Однако неочевидно, можем ли мы записать отображение $g$ в виде суммы слагаемых вида $(ax)b$ и $a(xb)$.

Если $A$ - свободная конечно мерная алгебра, то мы можем предположить, что линейное отображение имеет стандартное представление в виде[6.7]

$$f\circ x = f^{ij}\ (\overline{e}_i x)\overline{e}_j \tag{6.5.8}$$

В этом случае мы можем применить теорему 6.4.8 для отображений в неассоциативную алгебру.

---

[6.7] Выбор произволен. Мы можем рассмотреть стандартное представление в виде

$$f\circ x = f^{ij}\overline{e}_i(x\overline{e}_j)$$

Тогда равенство (6.5.11) имеет вид

$$g^k_l = f^m_l g^{ij} C_{2\cdot}{}^{k}_{ip} C_{2\cdot}{}^{p}_{mj} \tag{6.5.7}$$

Я выбрал выражение (6.5.8) так как порядок сомножителей соответствует порядку, выбранному в теореме 6.4.8.

Теорема 6.5.3. *Пусть $\overline{\overline{e}}_1$ - базис свободной конечно мерной $D$-алгебры $A_1$. Пусть $\overline{\overline{e}}_2$ - базис свободной конечно мерной неассоциативной $D$-алгебры $A_2$. Пусть $C_{2\cdot}{}^{p}_{kl}$ - структурные константы алгебры $A_2$. Пусть отображение*

(6.5.9) $$g = a \circ f$$

*порождённое отображением $f \in \mathcal{L}(D; A_1 \to A_2)$ посредством тензора $a \in A_2 \otimes A_2$, имеет стандартное представление*

(6.5.10) $$g = a^{ij}(\overline{e}_i \otimes \overline{e}_j) \circ f = a^{ij}(\overline{e}_i f)\overline{e}_j$$

*Координаты отображения* (6.5.9) *и его стандартные компоненты связаны равенством*

(6.5.11) $$g^k_l = f^m_l g^{ij} C_{2\cdot}{}^{p}_{im} C_{2\cdot}{}^{k}_{pj}$$

Доказательство. Относительно базисов $\overline{\overline{e}}_1$ и $\overline{\overline{e}}_2$, линейные отображения $f$ и $g$ имеют вид

(6.5.12) $$f \circ x = f^i_j x^j \overline{e}_{2\cdot i}$$

(6.5.13) $$g \circ x = g^i_j x^j \overline{e}_{2\cdot i}$$

Из равенств (6.5.12), (6.5.13), (6.5.10) следует

(6.5.14) $$\begin{aligned} g^k_l x^l \overline{e}_{2\cdot k} &= a^{ij}(\overline{e}_{2\cdot i}(f^m_l x^l \overline{e}_{2\cdot m}))\overline{e}_{2\cdot j} \\ &= a^{ij} f^m_l x^l C_{2\cdot}{}^{p}_{im} C_{2\cdot}{}^{k}_{pj} \overline{e}_{2\cdot k} \end{aligned}$$

Так как векторы $\overline{e}_{2\cdot k}$ линейно независимы и $x^i$ произвольны, то равенство (6.5.11) следует из равенства (6.5.14). □

Теорема 6.5.4. *Пусть $A$ - свободная конечно мерная неассоциативная алгебра над кольцом $D$. Представление алгебры $A \otimes A$ в алгебре $\mathcal{L}(D; A \to A)$ имеет конечный базис $\overline{\overline{I}}$.*

(1) *Линейное отображение $f \in \mathcal{L}(D; A \to A)$ имеет вид*

(6.5.15) $$f = (a_{k\cdot s_k\cdot 0} \otimes a_{k\cdot s_k\cdot 1}) \circ I_k = (a_{k\cdot s_k\cdot 0} I_k) a_{k\cdot s_k\cdot 1}$$

(2) *Его стандартное представление имеет вид*

(6.5.16) $$f = a^{k\cdot ij}(\overline{e}_i \otimes \overline{e}_j) \circ I_k = a^{k\cdot ij}(\overline{e}_i I_k)\overline{e}_j$$

Доказательство. Рассмотрим матрицу (6.4.28). Если матрица $\mathcal{B}$ невырождена, то для заданных координат линейного преобразования $g^l_k$ и для отображения $f = \delta$, система линейных уравнений (6.5.11) относительно стандартных компонент этого преобразования $g^{kr}$ имеет единственное решение. Если матрица $\mathcal{B}$ вырождена, то согласно теореме 6.4.8 существует конечный базис $\overline{\overline{I}}$, порождающий множество линейных отображений. □

В отличие от случая ассоциативной алгебры множество генераторов $I$ в теореме 6.5.4 не является минимальным. Из равенства (6.5.2) следует, что неверно равенство (6.3.12). Следовательно, орбиты отображений $I_k$ не порождают отношения эквивалентности в алгебре $L(A; A)$. Так как мы рассматриваем только отображения вида $(aI_k)b$, то возможно, что при $k \neq l$ отображение $I_k$ порождает отображение $I_l$, если рассмотреть все возможные операции в алгебре $A$. Поэтому множество образующих $I_k$ неассоциативной алгебры $A$ не играет такой критической роли как отображение сопряжения в поле комплексных чисел. Ответ на вопрос насколько важно отображение $I_k$ в неассоциативной алгебре требует дополнительного исследования.

## 6.6. Полилинейное отображение в ассоциативную алгебру

Теорема 6.6.1. *Пусть $A_1$, ..., $A_n$, $A$ - ассоциативные $D$-алгебры. Пусть*

$$f_i \in \mathcal{L}(D; A_i \to A) \quad i=1,...,n$$

$$a_j \in A \quad j=0,...,n$$

*Для заданной перестановки $\sigma$ $n$ переменных отображение*

$$\begin{aligned}
&((a_0,...,a_n,\sigma)\circ(f_1,...,f_n))\circ(x_1,...,x_n)\\
&=(a_0\sigma(f_1)a_1...a_{n-1}\sigma(f_n)a_n)\circ(x_1,...,x_n)\\
&=a_0\sigma(f_1\circ x_1)a_1...a_{n-1}\sigma(f_n\circ x_n)a_n
\end{aligned} \tag{6.6.1}$$

*является $n$-линейным отображением в алгебру $A$.*

Доказательство. Утверждение теоремы следует из цепочек равенств

$$\begin{aligned}
&((a_0,...,a_n,\sigma)\circ(f_1,...,f_n))\circ(x_1,...,x_i+y_i,...,x_n)\\
=&a_0\sigma(f_1\circ x_1)a_1...\sigma(f_i\circ(x_i+y_i))...a_{n-1}\sigma(f_n\circ x_n)a_n\\
=&a_0\sigma(f_1\circ x_1)a_1...\sigma(f_i\circ x_i+f_i\circ y_i)...a_{n-1}\sigma(f_n\circ x_n)a_n\\
=&a_0\sigma(f_1\circ x_1)a_1...\sigma(f_i\circ x_i)...a_{n-1}\sigma(f_n\circ x_n)a_n\\
+&a_0\sigma(f_1\circ x_1)a_1...\sigma(f_i\circ y_i)...a_{n-1}\sigma(f_n\circ x_n)a_n\\
=&((a_0,...,a_n,\sigma)\circ(f_1,...,f_n))\circ(x_1,...,x_i,...,x_n)\\
+&((a_0,...,a_n,\sigma)\circ(f_1,...,f_n))\circ(x_1,...,y_i,...,x_n)
\end{aligned}$$

$$\begin{aligned}
&((a_0,...,a_n,\sigma)\circ(f_1,...,f_n))\circ(x_1,...,px_i,...,x_n)\\
=&a_0\sigma(f_1\circ x_1)a_1...\sigma(f_i\circ(px_i))...a_{n-1}\sigma(f_n\circ x_n)a_n\\
=&a_0\sigma(f_1\circ x_1)a_1...\sigma(p(f_i\circ x_i))...a_{n-1}\sigma(f_n\circ x_n)a_n\\
=&p(a_0\sigma(f_1\circ x_1)a_1...\sigma(f_i\circ x_i)...a_{n-1}\sigma(f_n\circ x_n)a_n)\\
=&p(((a_0,...,a_n,\sigma)\circ(f_1,...,f_n))\circ(x_1,...,x_i,...,x_n))
\end{aligned}$$

□

В равенстве (6.6.1), также как и в других выражениях полилинейного отображения, принято соглашение, что отображение $f_i$ имеет своим аргументом переменную $x_i$.

Теорема 6.6.2. *Пусть $A_1$, ..., $A_n$, $A$ - ассоциативные $D$-алгебры. Для заданного семейства отображений*

$$f_i \in \mathcal{L}(D; A_i \to A) \quad i=1,...,n$$

*отображение*

$$h: A^{n+1} \to \mathcal{L}(D; A_1\times...\times A_n \to A)$$

*определённое равенством*

$$(a_0,...,a_n,\sigma)\circ(f_1,...,f_n)=a_0\sigma(f_1)a_1...a_{n-1}\sigma(f_n)a_n$$

*является $n+1$-линейным отображением в $D$-модуль $\mathcal{L}(D; A_1\times...\times A_n \to A)$.*

Доказательство. Утверждение теоремы следует из цепочек равенств

$$\begin{aligned}
&((a_0,...,a_i+b_i,...a_n,\sigma)\circ(f_1,...,f_n))\circ(x_1,...,x_n)\\
=&a_0\sigma(f_1\circ x_1)a_1...(a_i+b_i)...a_{n-1}\sigma(f_n\circ x_n)a_n\\
=&a_0\sigma(f_1\circ x_1)a_1...a_i...a_{n-1}\sigma(f_n\circ x_n)a_n+a_0\sigma(f_1\circ x_1)a_1...b_i...a_{n-1}\sigma(f_n\circ x_n)a_n\\
=&((a_0,...,a_i,...,a_n,\sigma)\circ(f_1,...,f_n))\circ(x_1,...,x_n)\\
+&((a_0,...,b_i,...,a_n,\sigma)\circ(f_1,...,f_n))\circ(x_1,...,x_n)\\
=&((a_0,...,a_i,...,a_n,\sigma)\circ(f_1,...,f_n)+(a_0,...,b_i,...,a_n,\sigma)\circ(f_1,...,f_n))\circ(x_1,...,x_n)
\end{aligned}$$

$$\begin{aligned}
&((a_0,...,pa_i,...a_n,\sigma)\circ(f_1,...,f_n))\circ(x_1,...,x_n)\\
=&a_0\sigma(f_1\circ x_1)a_1...pa_i...a_{n-1}\sigma(f_n\circ x_n)a_n\\
=&p(a_0\sigma(f_1\circ x_1)a_1...a_i...a_{n-1}\sigma(f_n\circ x_n)a_n)\\
=&p(((a_0,...,a_i,...,a_n,\sigma)\circ(f_1,...,f_n))\circ(x_1,...,x_n))\\
=&(p((a_0,...,a_i,...,a_n,\sigma)\circ(f_1,...,f_n)))\circ(x_1,...,x_n)
\end{aligned}$$

□

Теорема 6.6.3. *Пусть $A_1$, ..., $A_n$, $A$ - ассоциативные $D$-алгебры. Для заданного семейства отображений*

$$f_i\in\mathcal{L}(D;A_i\to A)\quad i=1,...,n$$

*существует линейное отображение*

$$h:A^{\otimes n+1}\times S_n\to\mathcal{L}(D;A_1\times...\times A_n\to A)$$

*определённое равенством*

$$\begin{aligned}(a_0\otimes...\otimes a_n,\sigma)\circ(f_1,...,f_n)&=(a_0,...,a_n,\sigma)\circ(f_1,...,f_n)\\&=a_0\sigma(f_1)a_1...a_{n-1}\sigma(f_n)a_n\end{aligned}\tag{6.6.2}$$

Доказательство. Утверждение теоремы является следствием теорем 4.5.4, 6.6.2. □

Теорема 6.6.4. *Пусть $A_1$, ..., $A_n$, $A$ - ассоциативные $D$-алгебры. Для заданного тензора $a\in A^{\otimes n+1}$ и заданной перестановки $\sigma\in S_n$ отображение*

$$h:\prod_{i=1}^n\mathcal{L}(D;A_i\to A)\to\mathcal{L}(D;A_1\times...\times A_n\to A)$$

*определённое равенством*

$$(a_0\otimes...\otimes a_n,\sigma)\circ(f_1,...,f_n)=a_0\sigma(f_1)a_1...a_{n-1}\sigma(f_n)a_n$$

*является $n$-линейным отображением в $D$-модуль $\mathcal{L}(D;A_1\times...\times A_n\to A)$.*

Доказательство. Утверждение теоремы следует из цепочек равенств

$$\begin{aligned}
&((a_0\otimes...\otimes a_n,\sigma)\circ(f_1,...,f_i+g_i,...,f_n))\circ(x_1,...,x_n)\\
=&(a_0\sigma(f_1)a_1...\sigma(f_i+g_i)...a_{n-1}\sigma(f_n)a_n)\circ(x_1,...,x_n)\\
=&a_0\sigma(f_1\circ x_1)a_1...\sigma((f_i+g_i)\circ x_i)...a_{n-1}\sigma(f_n\circ x_n)a_n\\
=&a_0\sigma(f_1\circ x_1)a_1...\sigma(f_i\circ x_i+g_i\circ x_i)...a_{n-1}\sigma(f_n\circ x_n)a_n\\
=&a_0\sigma(f_1\circ x_1)a_1...\sigma(f_i\circ x_i)...a_{n-1}\sigma(f_n\circ x_n)a_n\\
+&a_0\sigma(f_1\circ x_1)a_1...\sigma(g_i\circ x_i)...a_{n-1}\sigma(f_n\circ x_n)a_n\\
=&(a_0\sigma(f_1)a_1...\sigma(f_i)...a_{n-1}\sigma(f_n)a_n)\circ(x_1,...,x_n)\\
+&(a_0\sigma(f_1)a_1...\sigma(g_i)...a_{n-1}\sigma(f_n)a_n)\circ(x_1,...,x_n)\\
=&((a_0\otimes...\otimes a_n,\sigma)\circ(f_1,...,f_i,...,f_n))\circ(x_1,...,x_n)\\
+&((a_0\otimes...\otimes a_n,\sigma)\circ(f_1,...,g_i,...,f_n))\circ(x_1,...,x_n)\\
=&((a_0\otimes...\otimes a_n,\sigma)\circ(f_1,...,f_i,...,f_n)\\
+&(a_0\otimes...\otimes a_n,\sigma)\circ(f_1,...,g_i,...,f_n))\circ(x_1,...,x_n)
\end{aligned}$$

$$\begin{aligned}
&((a_0\otimes...\otimes a_n,\sigma)\circ(f_1,...,pf_i,...,f_n))\circ(x_1,...,x_n)\\
=&(a_0\sigma(f_1)a_1,...\sigma(pf_i)...a_{n-1}\sigma(f_n)a_n)\circ(x_1,...,x_n)\\
=&a_0\sigma(f_1\circ x_1)a_1,...\sigma((pf_i)\circ x_i)...a_{n-1}\sigma(f_n\circ x_n)a_n\\
=&a_0\sigma(f_1\circ x_1)a_1,...\sigma(p(f_i\circ x_i))...a_{n-1}\sigma(f_n\circ x_n)a_n\\
=&p(a_0\sigma(f_1\circ x_1)a_1,...\sigma(f_i\circ x_i)...a_{n-1}\sigma(f_n\circ x_n)a_n)\\
=&p(((a_0\otimes...\otimes a_n,\sigma)\circ(f_1,...,f_i,...,f_n))\circ(x_1,...,x_n))\\
=&(p((a_0\otimes...\otimes a_n,\sigma)\circ(f_1,...,f_i,...,f_n)))\circ(x_1,...,x_n)
\end{aligned}$$

□

Теорема 6.6.5. *Пусть $A_1$, ..., $A_n$, $A$ - ассоциативные $D$-алгебры. Для заданного тензора $a\in A^{\otimes n+1}$ и заданной перестановки $\sigma\in S_n$ существует линейное отображение*

$$h:\mathcal{L}(D;A_1\to A)\otimes...\otimes\mathcal{L}(D;A_n\to A)\to\mathcal{L}(D;A_1\times...\times A_n\to A)$$

*определённое равенством*

$$(a_0\otimes...\otimes a_n,\sigma)\circ(f_1\otimes...\otimes f_n)=(a_0\otimes...\otimes a_n,\sigma)\circ(f_1,...,f_n)\tag{6.6.3}$$

Доказательство. Утверждение теоремы является следствием теорем 4.5.4, 6.6.4. □

Теорема 6.6.6. *Пусть $A$ - ассоциативная $D$-алгебра. Полилинейное отображение*

$$f:A^n\to A,a=f\circ(a_1,...,a_n)\tag{6.6.4}$$

*порождённое отображениями $I_{s\cdot1}$, ..., $I_{s\cdot n}\in\mathcal{L}(D;A\to A)$ , имеет вид*

$$a=f^n_{s\cdot0}\ \sigma_s(I_{s\cdot1}\circ a_1)\ f^n_{s\cdot1}\ ...\ \sigma_s(I_{s\cdot n}\circ a_n)\ f^n_{s\cdot n}\tag{6.6.5}$$

*где $\sigma_s$ - перестановка множества переменных $\{a_1,...,a_n\}$*

$$\sigma_s=\begin{pmatrix}a_1&...&a_n\\ \sigma_s(a_1)&...&\sigma_s(a_n)\end{pmatrix}$$

Доказательство. Мы докажем утверждение индукцией по $n$.

При $n = 1$ доказываемое утверждение является следствием утверждения 6.4.8.1. При этом мы можем отождествить[6.8]

$$f^1_{s\cdot p} = f_{s\cdot p} \quad p = 0, 1$$

Допустим, что утверждение теоремы справедливо при $n = k - 1$. Тогда отображение (6.6.4) можно представить в виде

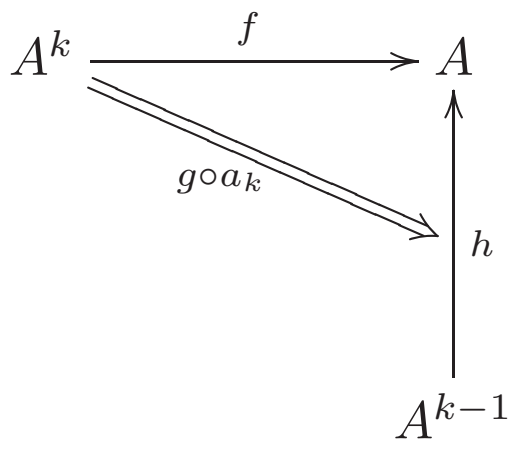

$$a = f \circ (a_1, ..., a_k) = (g \circ a_k) \circ (a_1, ..., a_{k-1})$$

Согласно предположению индукции полилинейное отображение $h$ имеет вид

$$a = h^{k-1}_{t\cdot 0}\ \sigma_t(I_{1\cdot t} \circ a_1)\ h^{k-1}_{t\cdot 1}\ ...\ \sigma_t(I_{k-1\cdot t} \circ a_{k-1})\ h^{k-1}_{t\cdot k-1}$$

Согласно построению $h = g \circ a_k$. Следовательно, выражения $h_{t\cdot p}$ являются функциями $a_k$. Поскольку $g \circ a_k$ - линейная функция $a_k$, то только одно выражение $h_{t\cdot p}$ является линейной функцией переменной $a_k$, и остальные выражения $h_{t\cdot q}$ не зависят от $a_k$.

Не нарушая общности, положим $p = 0$. Согласно равенству (6.3.10) для заданного $t$

$$h^{k-1}_{t\cdot 0} = g_{tr\cdot 0}\ I_{k\cdot r} \circ a_k\ g_{tr\cdot 1}$$

Положим $s = tr$ и определим перестановку $\sigma_s$ согласно правилу

$$\sigma_s = \sigma_{tr} = \begin{pmatrix} a_k & a_1 & ... & a_{k-1} \\ a_k & \sigma_t(a_1) & ... & \sigma_t(a_{k-1}) \end{pmatrix}$$

Положим

$$f^k_{tr\cdot q+1} = h^{k-1}_{t\cdot q} \qquad q = 1, ..., k-1$$
$$f^k_{tr\cdot q} = g_{tr\cdot q} \qquad q = 0, 1$$

Мы доказали шаг индукции. □

Определение 6.6.7. *Выражение* $f s\cdot p^n$ *в равенстве* (6.6.5) *называется* **компонентой полилинейного отображения** $f$. □

[6.8] В представлении (6.6.5) мы будем пользоваться следующими правилами.

- Если область значений какого-либо индекса - это множество, состоящее из одного элемента, мы будем опускать соответствующий индекс.
- Если $n = 1$, то $\sigma_s$ - тождественное преобразование. Это преобразование можно не указывать в выражении.

Теорема 6.6.8. *Рассмотрим $D$-алгебру $A$. Линейное отображение*

$$h : A^{\otimes n+1} \times S_n \to {}^*\mathcal{L}(D; A^n \to A)$$

*определённое равенством*

$$(a_0 \otimes ... \otimes a_n, \sigma) \circ (f_1 \otimes ... \otimes f_n) = a_0 \sigma(f_1) a_1 ... a_{n-1} \sigma(f_n) a_n \tag{6.6.6}$$
$$a_0, ..., a_n \in A \quad \sigma \in S_n \quad f_1, ..., f_n \in \mathcal{L}(D; A_n \to A)$$

*является представлением*[6.9] *алгебры* $A^{\otimes n+1} \times S^n$ *в $D$-модуле* $\mathcal{L}(D; A^n \to A)$.

Доказательство. Согласно теоремам 6.4.8, 6.6.6, $n$-линейное отображение можно представить в виде суммы слагаемых (6.6.1), где $f_i$, $i = 1, ..., n$, генераторы представления (6.3.3). Запишем слагаемое $s$ выражения (6.6.5) в виде

$$b_1 \sigma(I_{1 \cdot s} \circ x_1) c_1 b_2 ... c_{n-1} b_n \sigma(I_{n \cdot s} \circ x_n) c_n \tag{6.6.7}$$

где

$$b_1 = f^n_{s \cdot 0} \quad b_2 = ... = b_n = e \quad c_1 = f^n_{s \cdot 1} \quad ... \quad c_n = f^n_{s \cdot n}$$

Положим в равенстве (6.6.7)

$$f_i = \sigma^{-1}(b_i) I_{i \cdot s} \sigma^{-1}(c_i) \quad i = 1, ..., n$$

Следовательно, согласно теореме 6.6.3, отображение (6.6.6) является преобразованием модуля $\mathcal{L}(D; A^n \to A)$. Для данного тензора $c \in A^{\otimes n+1}$ и заданной перестановки $\sigma \in S_n$, преобразование $h(c, \sigma)$ является линейным преобразованием модуля $\mathcal{L}(D; A^n \to A)$ согласно теореме 6.6.5. Согласно теореме 6.6.3, отображение (6.6.6) является линейным отображением. Согласно определению [4]-2.1.2, отображение (6.6.6) является представлением алгебры $A^{\otimes n+1} \times S^n$ в модуле $\mathcal{L}(D; A^n \to A)$. □

Теорема 6.6.9. *Рассмотрим $D$-алгебру $A$. Представление*

$$h : A^{\otimes n+1} \times S_n \to {}^*\mathcal{L}(D; A^n \to A)$$

*алгебры $A^{\otimes n+1}$ в модуле $\mathcal{L}(D; A^n \to A)$, определённое равенством (6.6.6) позволяет отождествить тензор $d \in A^{\otimes n+1}$ и перестановку $\sigma \in S^n$ с отображением*

$$(d, \sigma) \circ (f_1, ..., f_n) \quad f_i = \delta \in \mathcal{L}(D; A \to A) \tag{6.6.8}$$

*где $\delta \in \mathcal{L}(D; A \to A)$ - тождественное отображение.*

Доказательство. Если в равенстве (6.6.2) мы положим $f_i = \delta$, $d = a_0 \otimes ... \otimes a_n$, то равенство (6.6.2) приобретает вид

$$\begin{aligned} ((a_0 \otimes ... \otimes a_n, \sigma) \circ (\delta, ..., \delta)) \circ (x_1, ..., x_n) &= a_0 \ (\delta \circ x_1) \ ... \ (\delta \circ x_n) \ a_n \\ &= a_0 \ x_1 \ ... \ x_n \ a_n \end{aligned} \tag{6.6.9}$$

[6.9] Определение представления $\Omega$-алгебры дано в определении [4]-2.1.2.

Если мы положим

$$\begin{aligned}&((a_0\otimes ...\otimes a_n,\sigma)\circ(\delta,...,\delta))\circ(x_1,...,x_n)\\ =&(a_0\otimes ...\otimes a_n,\sigma)\circ(\delta\circ x_1,...,\delta\circ x_n)\\ =&(a_0\otimes ...\otimes a_n,\sigma)\circ(x_1,...,x_n)\end{aligned}\tag{6.6.10}$$

то сравнение равенств (6.6.9) и (6.6.10) даёт основание отождествить действие тензора $d=a_0\otimes ...\otimes a_n$ и перестановки $\sigma\in S^n$ с отображением (6.6.8). □

Вместо записи $(a_0\otimes ...a_n,\sigma)$ мы также будем пользоваться записью

$$a_0\otimes_{\sigma(1)}...\otimes_{\sigma(n)}a_n$$

если мы хотим явно указать какой аргумент становится на соответствующее место. Например, следующие выражения эквивалентны

$$(a_0\otimes a_1\otimes a_2\otimes a_3,(2,1,3))\circ(x_1,x_2,x_3)=a_0x_2a_1x_1a_2x_3a_3$$

$$(a_0\otimes_2 a_1\otimes_1 a_2\otimes_3 a_3)\circ(x_1,x_2,x_3)=a_0x_2a_1x_1a_2x_3a_3$$

## 6.7. Полилинейное отображение в свободную конечно мерную ассоциативную алгебру

Теорема 6.7.1. *Пусть $A$ - свободная конечно мерная ассоциативная алгебра над кольцом $D$. Пусть $\overline{\overline{I}}$ - базис алгебры $\mathcal{L}(D;A\to A)$. Пусть $\overline{\overline{e}}$ - базис алгебры $A$ над кольцом $D$.* **Стандартное представление полилинейного отображения** *в ассоциативную алгебру имеет вид*

$$f\circ(a_1,...,a_n)=f^{i_0...i_n}_{t\cdot k_1...k_n}\ \overline{e}_{i_0}\ \sigma_t(I_{k_1}\circ a_1)\ \overline{e}_{i_1}\ ...\ \sigma_t(I_{k_n}\circ a_n)\ \overline{e}_{i_n}\tag{6.7.1}$$

*Индекс $t$ нумерует всевозможные перестановки $\sigma_t$ множества переменных $\{a_1,...,a_n\}$. Выражение $f^{i_0...i_n}_{t\cdot k_1...k_n}$ в равенстве (6.7.1) называется* **стандартной компонентой полилинейного отображения** $f$.

Доказательство. Мы изменим индекс $s$ в равенстве (6.6.5) таким образом, чтобы сгруппировать слагаемые с одинаковым набором генераторов $I_k$. Выражение (6.6.5) примет вид

$$a=f^n_{k_1...k_n\cdot s\cdot 0}\ \sigma_s(I_{k_1\cdot s}\circ a_1)\ f^n_{k_1...k_n\cdot s\cdot 1}\ ...\ \sigma_s(I_{k_n\cdot s}\circ a_n)\ f^n_{k_1...k_n\cdot s\cdot n}\tag{6.7.2}$$

Мы предполагаем, что индекс $s$ принимает значения, зависящие от $k_1$, ..., $k_n$. Компоненты полилинейного отображения $f$ имеют разложение

$$f^n_{k_1...k_n\cdot s\cdot p}=\overline{e}_i f^{ni}_{k_1...k_n\cdot s\cdot p}\tag{6.7.3}$$

относительно базиса $\overline{\overline{e}}$. Если мы подставим (6.7.3) в (6.6.5), мы получим

$$a=f^{nj_1}_{k_1...k_n\cdot s\cdot 0}\ \overline{e}_{j_1}\ \sigma_s(I_{k_1}\circ a_1)\ f^{nj_2}_{k_1...k_n\cdot s\cdot 1}\ \overline{e}_{j_2}\ ...\ \sigma_s(I_{k_n}\circ a_n)\ f^{nj_n}_{k_1...k_n\cdot s\cdot n}\ \overline{e}_{j_n}\tag{6.7.4}$$

Рассмотрим выражение

$$f^{j_0...j_n}_{t\cdot k_1...k_n}=f^{nj_1}_{k_1...k_n\cdot s\cdot 0}\ ...f^{nj_n}_{k_1...k_n\cdot s\cdot n}\tag{6.7.5}$$

В правой части подразумевается сумма тех слагаемых с индексом $s$, для которых перестановка $\sigma_s$ совпадает. Каждая такая сумма будет иметь уникальный индекс $t$. Подставив в равенство (6.7.4) выражение (6.7.5) мы получим равенство (6.7.1). □

Теорема 6.7.2. *Пусть $A$ - свободная конечно мерная ассоциативная алгебра над кольцом $D$. Пусть $\overline{\overline{e}}$ - базис алгебры $A$ над кольцом $D$. Полилинейное отображение (6.6.4) можно представить в виде $D$-значной формы степени $n$ над кольцом $D$*[6.10]

$$f(a_1, ..., a_n) = a_1^{i_1} ... a_n^{i_n} f_{i_1...i_n} \tag{6.7.6}$$

*где*

$$\begin{aligned} a_j &= \overline{e}_i a_j^i \\ f_{i_1...i_n} &= f \circ (\overline{e}_{i_1}, ..., \overline{e}_{i_n}) \end{aligned} \tag{6.7.7}$$

Доказательство. Согласно определению 6.1.2 равенство (6.7.6) следует из цепочки равенств

$$f \circ (a_1, ..., a_n) = f \circ (\overline{e}_{i_1} a_1^{i_1}, ..., \overline{e}_{i_n} a_n^{i_n}) = a_1^{i_1} ... a_n^{i_n} f \circ (\overline{e}_{i_1}, ..., \overline{e}_{i_n})$$

Пусть $\overline{\overline{e}}'$ - другой базис. Пусть

$$\overline{e}'_i = \overline{e}_j h_i^j \tag{6.7.8}$$

преобразование, отображающее базис $\overline{\overline{e}}$ в базис $\overline{\overline{e}}'$. Из равенств (6.7.8) и (6.7.7) следует

$$\begin{aligned} f'_{i_1...i_n} &= f \circ (\overline{e}'_{i_1}, ..., \overline{e}'_{i_n}) \\ &= f \circ (\overline{e}_{j_1} h_{i_1}^{j_1}, ..., \overline{e}'_{j_n} h_{i_n}^{j_n}) \\ &= h_{i_1}^{j_1} ... h_{i_n}^{j_n} f \circ (\overline{e}_{j_1}, ..., \overline{e}_{j_n}) \\ &= h_{i_1}^{j_1} ... h_{i_n}^{j_n} f_{j_1...j_n} \end{aligned} \tag{6.7.9}$$

□

Полилинейное отображение (6.6.4) **симметрично**, если

$$f \circ (a_1, ..., a_n) = f \circ (\sigma(a_1), ..., \sigma(a_n))$$

для любой перестановки $\sigma$ множества $\{a_1, ..., a_n\}$.

Теорема 6.7.3. *Если полилинейное отображение $f$ симметрично, то*

$$f_{i_1,...,i_n} = f_{\sigma(i_1),...,\sigma(i_n)} \tag{6.7.10}$$

Доказательство. Равенство (6.7.10) следует из равенства

$$\begin{aligned} a_1^{i_1} ... a_n^{i_n} f_{i_1...i_n} &= f \circ (a_1, ..., a_n) \\ &= f \circ (\sigma(a_1), ..., \sigma(a_n)) \\ &= a_1^{i_1} ... a_n^{i_n} f_{\sigma(i_1)...\sigma(i_n)} \end{aligned}$$

□

Полилинейное отображение (6.6.4) **косо симметрично**, если

$$f \circ (a_1, ..., a_n) = |\sigma| f \circ (\sigma(a_1), ..., \sigma(a_n))$$

[6.10] Теорема доказана по аналогии с теоремой в [2], с. 107, 108

для любой перестановки $\sigma$ множества $\{a_1, ..., a_n\}$. Здесь

$$|\sigma| = \begin{cases} 1 & \text{перестановка } \sigma \text{ чётная} \\ -1 & \text{перестановка } \sigma \text{ нечётная} \end{cases}$$

Теорема 6.7.4. *Если полилинейное отображение $f$ косо симметрично, то*

$$f_{i_1,...,i_n} = |\sigma| f_{\sigma(i_1),...,\sigma(i_n)} \tag{6.7.11}$$

Доказательство. Равенство (6.7.11) следует из равенства

$$\begin{aligned} a_1^{i_1} ... a_n^{i_n} f_{i_1...i_n} &= f \circ (a_1, ..., a_n) \\ &= |\sigma| f \circ (\sigma(a_1), ..., \sigma(a_n)) \\ &= a_1^{i_1} ... a_n^{i_n} |\sigma| f_{\sigma(i_1)...\sigma(i_n)} \end{aligned}$$

□

Теорема 6.7.5. *Пусть $A$ - свободная конечно мерная ассоциативная алгебра над кольцом $D$. Пусть $\overline{\overline{I}}$ - базис алгебры $\mathcal{L}(D; A \to A)$. Пусть $\overline{\overline{e}}$ - базис алгебры $A$ над кольцом $D$. Пусть полилинейное над кольцом $D$ отображение $f$ порождено набором отображений $(I_{k_1}, ..., I_{k_n})$. Координаты отображения $f$ и его компоненты относительно базиса $\overline{\overline{e}}$ удовлетворяют равенству*

$$f_{l_1...l_n} = f^{i_0...i_n}_{t \cdot k_1...k_n} I_{k_1 \cdot}{}^{j_1}_{l_1} ... I_{k_n \cdot}{}^{j_n}_{l_n} C^{k_1}_{i_0 \sigma_t(j_1)} C^{l_1}_{k_1 i_1} ... B^{k_n}_{l_{n-1} \sigma_t(j_n)} C^{l_n}_{k_n i_n} \overline{e}_{l_n} \tag{6.7.12}$$

$$f^{p}_{l_1...l_n} = f^{i_0...i_n}_{t \cdot k_1...k_n} I_{k_1 \cdot}{}^{j_1}_{l_1} ... I_{k_n \cdot}{}^{j_n}_{l_n} C^{k_1}_{i_0 \sigma_t(j_1)} C^{l_1}_{k_1 i_1} ... C^{k_n}_{l_{n-1} \sigma_t(j_n)} C^{p}_{k_n i_n} \tag{6.7.13}$$

Доказательство. В равенстве (6.7.1) положим

$$I_{k_i} \circ a_i = \overline{e}_{j_i} I_{k_i \cdot}{}^{j_i}_{l_i} a_i^{l_i}$$

Тогда равенство (6.7.1) примет вид

(6.7.14)

$$\begin{aligned} f \circ (a_1, ..., a_n) &= f^{i_0...i_n}_{t \cdot k_1...k_n} \overline{e}_{i_0} \sigma_t(a_1^{l_1} I_{k_1 \cdot}{}^{j_1}_{l_1} \overline{e}_{j_1}) \overline{e}_{i_1} ... \sigma_t(a_n^{l_n} I_{k_n \cdot}{}^{j_n}_{l_n} \overline{e}_{j_n}) \overline{e}_{i_n} \\ &= a_1^{l_1} ... a_n^{l_n} f^{i_0...i_n}_{t \cdot k_1...k_n} I_{k_1 \cdot}{}^{j_1}_{l_1} ... I_{k_n \cdot}{}^{j_n}_{l_n} \overline{e}_{i_0} \sigma_t(\overline{e}_{j_1}) \overline{e}_{i_1} ... \sigma_t(\overline{e}_{j_n}) \overline{e}_{i_n} \\ &= a_1^{l_1} ... a_n^{l_n} f^{i_0...i_n}_{t \cdot k_1...k_n} I_{k_1 \cdot}{}^{j_1}_{l_1} ... I_{k_n \cdot}{}^{j_n}_{l_n} C^{k_1}_{i_0 \sigma_t(j_1)} C^{l_1}_{k_1 i_1} \\ &\quad ... C^{k_n}_{l_{n-1} \sigma_t(j_n)} C^{l_n}_{k_n i_n} \overline{e}_{l_n} \end{aligned}$$

Из равенства (6.7.6) следует

$$f \circ (a_1, ..., a_n) = \overline{e}_p f^{p}_{i_1...i_n} a_1^{i_1} ... a_n^{i_n} \tag{6.7.15}$$

Равенство (6.7.12) следует из сравнения равенств (6.7.14) и (6.7.6). Равенство (6.7.13) следует из сравнения равенств (6.7.14) и (6.7.15). □

# Список литературы

# Предметный указатель

# Специальные символы и обозначения